\documentclass[12pt, a4,oneside]{amsart}
\let\latexput\put
\usepackage{amsmath,amsfonts,amssymb,amsthm,bm,wasysym}%mathtext

\usepackage{subfigure}
\usepackage{verbatim}
\usepackage{enumerate}
\usepackage{enumitem} 
\usepackage{version,etex}
\usepackage[top=1in, bottom=1in,left=1.5in, right=0.5in]{geometry}
\usepackage[pdftex]{graphicx,color}
\usepackage{graphicx}
\usepackage[percent]{overpic}
\usepackage{epsfig}
\usepackage{hyperref}
\usepackage{caption}
\usepackage[normalem]{ulem}

\DeclareGraphicsExtensions{.pdf_t}
\reversemarginpar
\catcode`@=11 \catcode`!=11

\expandafter\ifx\csname fiverm\endcsname\relax
  \let\fiverm\fivrm
\fi
  
\let\!latexendpicture=\endpicture 
\let\!latexframe=\frame
\let\!latexlinethickness=\linethickness
\let\!latexmultiput=\multiput
\let\!latexput=\put
 
\def\@picture(#1,#2)(#3,#4){%
  \@picht #2\unitlength
  \setbox\@picbox\hbox to #1\unitlength\bgroup 
  \let\endpicture=\!latexendpicture
  \let\frame=\!latexframe
  \let\linethickness=\!latexlinethickness
  \let\multiput=\!latexmultiput
  \let\put=\!latexput
  \hskip -#3\unitlength \lower #4\unitlength \hbox\bgroup}

\catcode`@=12 \catcode`!=12
\input pictex.tex
\catcode`@=11 \catcode`!=11
  
\let\!pictexendpicture=\endpicture 
\let\!pictexframe=\frame
\let\!pictexlinethickness=\linethickness
\let\!pictexmultiput=\multiput
\let\!pictexput=\put

\def\beginpicture{%
  \setbox\!picbox=\hbox\bgroup%
  \let\endpicture=\!pictexendpicture
  \let\frame=\!pictexframe
  \let\linethickness=\!pictexlinethickness
  \let\multiput=\!pictexmultiput
  \let\put=\!pictexput
  \let\input=\@@input   % \@@input is LaTeX's saved version of TeX's primitive
  \!xleft=\maxdimen  
  \!xright=-\maxdimen
  \!ybot=\maxdimen
  \!ytop=-\maxdimen}

\let\frame=\!latexframe

\let\pictexframe=\!pictexframe

\let\linethickness=\!latexlinethickness
\let\pictexlinethickness=\!pictexlinethickness

\let\\=\@normalcr
\catcode`@=12 \catcode`!=12
\usepackage{tikz}
\let\put\latexput
\usetikzlibrary{matrix}
\def\D{\mathbb D} 
 
\def\Re{\operatorname{Re}}
\def\Im{\operatorname{Im}}
\newenvironment{proofof}[1]{\medskip
\noindent{\bf Proof of #1.}}{ \hfill\qed\\}
\def\E{E} 
\def\V{\hat E} 
\newtheorem{Extension Lemma}{Extension Lemma}[section] 
\def\Comp{\mbox{cc}}

\newtheorem{question}{Question}

\begin{document}
\setcounter{secnumdepth}{4}

\newcommand{\proclaim}[2]{\medbreak {\bf #1}{\sl #2} \medbreak}

\newcommand{\ntop}[2]{\genfrac{}{}{0pt}{1}{#1}{#2}}

\let\newpf\proof \let\proof\relax \let\endproof\relax
\newenvironment{pf}{\newpf[\proofname]}{\qed\endtrivlist}

\def\PC{{\mathrm{PC}}}
\def\R{{\mathbb{R}}}
\def\C{{\mathbb{C}}}
\def\N{{\mathbb{N}}}
\def\H{{\mathbb{H}}}
\def\P{{\mathcal{P}}}
\def\Crit{{\mathrm{Cr}}}
\def\Per{{\mathrm{Per}}}
\def\Y{{\mathcal{Y}}}
\def\dbar{{\bar{\partial}}}
\def\Im{{\mathrm{Im}}}
\def\diam{{\mathrm{diam}}}
\def\Dom{{\mathrm{Dom}}}
\def\dist{{\mathrm{dist}}}
\def\gap{{\mathrm{Gap}}}
\def\Space{{\mathrm{Space}}}
\def\spa{{\mathrm{Space}}}
\def\cen{{\mathrm{Cen}}}
\def\Land{{\hat{\mathcal{L}}}}

\def\comp{{\mathrm{Comp}}}
\def\LL{{\mathcal{L}}}
\def\clos{{\mathrm{Cl}}}
\def\crit{{\mathrm{Crit}}}
\def\mod{{\mathrm{mod}}}
\newcommand{\I}{\mathcal{I}}
\newcommand{\II}{\bm{I}}
\newcommand{\dombf}{\operatorname{\bf D}}
\newcommand{\DD}{\operatorname{\Dom}}
\newcommand{\U}{\bm{U}}
\newcommand{\J}{\bm{J}}
\newcommand{\W}{\bm{W}}
\newcommand{\PP}{\bm{P}}
\newcommand{\Z}{\bm{Z}}
\newcommand{\G}{\bm{G}}
\newcommand{\A}{\mathcal{A}}
\newcommand{\HH}{\mathbb{H}}
\newcommand{\frm}{R}
\newcommand{\compl}{\operatorname{compl}}
\newcommand{\interior}{\operatorname{int}}
\newcommand{\UU}{\bm{\mathcal{U}}}
\newcommand{\VV}{\bm{\mathcal{V}}}

\newtheorem{thm}{Theorem}[section]
\newtheorem{thmx}{Theorem}
\renewcommand{\thethmx}{\Alph{thmx}}
\newtheorem{cor}[thm]{Corollary}
\newtheorem{lem}[thm]{Lemma}
\newtheorem{lemma}[thm]{Lemma}
\newtheorem{claim}[thm]{Claim}
\newtheorem{prop}[thm]{Proposition}
\newtheorem{conj}[thm]{Conjecture}
\newtheorem{schw}{Schwarz Lemma} 
\newtheorem{sectl}{Sector Lemma}
\theoremstyle{remark}
\newtheorem{rem}{Remark}[section]
\newtheorem{notation}{Notation}
\newtheorem{example}{Example}[section]

\numberwithin{equation}{section}
\newcommand{\thmref}[1]{Theorem~\ref{#1}}
\newcommand{\propref}[1]{Proposition~\ref{#1}}
\newcommand{\secref}[1]{\S\ref{#1}}
\newcommand{\lemref}[1]{Lemma~\ref{#1}}
\newcommand{\corref}[1]{Corollary~\ref{#1}} 
\newcommand{\figref}[1]{Fig.~\ref{#1}}

\theoremstyle{definition}
\newtheorem{defn}{Definition}[section]
\newtheorem{definition}{Definition}[section]
\def\proof{\bn {\bf Proof.} }

\def\trevor{\textcolor{red}{.....QQQ....}}

\thanks{
The authors were supported by ERC AdG RGDD No 339523.  We are grateful to K. Drach for many 
helpful comments, and D. Preiss and M. Lyubich for several useful discussions. Part of this  paper was written  
during the spring 2022 programme at MSRI on holomorphic dynamics.}

\title[Conjugacy classes of real analytic one-dimensional maps] {Conjugacy classes of real analytic one-dimensional maps are analytic connected manifolds}
\author{Trevor Clark and Sebastian van Strien   }
\address{Mathematics Department, Imperial College London} 

\begin{abstract} An important question is to describe  topological conjugacy classes of 
dynamical systems. Here we show that within the space of real analytic one-dimensional maps with critical points of prescribed order,
the conjugacy class of a map is a real analytic manifold.  This extends results of Avila-Lyubich-de Melo \cite{ALM} 
for the quasi-quadratic unimodal case and of Clark \cite{C} for the more general unimodal case. Their methods fail in the case where there are several critical points, and 
for this reason we introduce the new notions of {\em pruned Julia set} of a real analytic map, 
and associate to a real analytic map an {\em external map} of the circle {\em with discontinuities} and   {\em a pruned polynomial-like} complex extension
of the real analytic map. Using this we are also able to show that topological conjugacy classes are connected (something 
which was not even known in the general unimodal setting). 
Even more,  this space is contractible. 

In a companion paper, further applications of this paper will be given. It will be shown that  
%(1) even 
% the space of maps which are merely partially conjugate (in the sense that the itineraries of certain critical points agree)
% to a given map forms a real analytic manifold;
%(2) 
 within any real analytic family of real analytic one-dimensional maps, hyperbolic parameters form an open and dense subset.

\end{abstract} 
\date{\today} 
\maketitle
 \setcounter{tocdepth}{1} 
\tableofcontents

\def\Crit{\mathrm{Cr}} 
\def\f{\hat f} 
\def\g{\hat g} 
\def\F{\hat F} 
\def\G{\hat G} 

\part*{} %Part A: Pruned polynomial-like maps and their external maps}
\section{Introduction and statement of results}
One of the main contributions of Sullivan into the field of dynamical systems is
his introduction of quasiconformal mappings, and in particular the Measurable
Riemann Mapping Theorem into the subject. For example, in his proof of the
absence of wandering domains for rational maps, he shows that if a rational map
has a wandering domain, then it is possible to construct an infinite dimensional
space of distinct deformations of this map, contradicting the fact that the
space of rational maps is finite dimensional, see \cite{Su1}. Similarly, he
proposed a strategy for proving density of hyperbolicity for real polynomial
maps by establishing quasisymmetric rigidity. This was successfully implemented
in the real quadratic case in \cite{Lyubich-quadratic dynamics, GS2} and in the
general case in \cite{KSS-rigidity, KSS-density}. For a survey on the
techniques developed in the latter papers, see \cite{CDKvS}. 
One of the limitations of quasisymmetric rigidity is that there is no analogue of the
Measurable Riemann Mapping Theorem in the real setting.

In this paper we will show how to overcome this in a number of settings, by
showing that two interval maps which are topologically conjugate can be
connected by a real-analytic path of maps in the same topological conjugacy class. Moreover 
we will show that topological conjugacy classes form real analytic manifolds, generalising 
a result in the unimodal setting due to Avila-Lyubich-de Melo \cite{ALM} and Clark \cite{C}.  
It turns out that these results imply density of hyperbolicity within full families. 

In a sequel to this paper, we will use the results
obtained here to show  that in the space of real analytic maps with {\em small basins}, 
 topological conjugacy classes laminate
the entire space. Here we say that a real analytic map has small basins, 
if the complex extension of the real basins of its periodic attractors are compactly 
contained in the domain of analyticity of the map; for a formal definition see 
Definition~\ref{def:smallbasin}.
This is quite a surprising result, as topological conjugacy classes 
{\em do not} laminate the space of all unimodal real analytic mappings. 
One application of this laminar structure is that the topological entropy of real analytic families of unimodal maps which are close to the 
family of quadratic maps depends monotonically on the parameter.  

In another sequel of this paper we will use the results in this paper to show that 
hyperbolic maps form a dense subset within real analytic families of real analytic one-dimensional maps.

\medskip
\label{page1} Let us state our results more precisely. Let $I$ be a compact interval; to be definite, let us
take $I=[-1,1]$. Let $\nu\in\mathbb N$, and fix a vector $\underline
\nu=(\ell_1,\ell_2,\dots,\ell_{\nu}) \in\mathbb N^{\nu}$ with each $\ell_i\geq
2$, and let $\mathcal A^{\underline \nu}$ denote the space of real analytic
mappings $f:I\to I$ with $f(\partial I)\subset\partial I$, with precisely $\nu$ 
critical points 
$-1<c_1<c_2<\dots<c_{\nu}<1$ such that $Df(c_i)=\dots = Df^{\ell_i-1}(c_i)=0$ and $Df^{\ell_i}(c_i)\ne 0$. 
For convenience, let us also assume that $\partial I$ is hyperbolic repelling and that 
$f$ is in the class $\mathcal A^{\underline \nu}_a$ of maps in $\mathcal A^{\underline \nu}$
which have a holomorphic extension to $\Omega_a=\{z\in \C: \dist(z,I)<a\}$, with precisely 
$\nu$ critical points in $\Omega_a$ of order $\ell_1,\dots,\ell_\nu$ and  which extends 
continuously to $\overline \Omega_a$.
The space $\mathcal A^{\underline \nu}_a$ is endowed with the supremum metric $d(f,g)=\sup_{z\in \overline{\Omega}_a}|f(z)-g(z)|$. 
Following Appendix 2 of \cite{Lyubich}, we will also introduce a topology on the space $\mathcal A^{\underline \nu}$ and the 
 concept  of {\em real analytic manifold modelled on Banach spaces} in Section~\ref{manifold-structure}.

If $f\in\mathcal A^{\underline \nu}$ has only hyperbolic periodic points (i.e. with multipliers not $0,\pm 1$), then the {\em topological conjugacy class} $\mathcal T^{\underline \nu}_{f}$ 
\label{Tf} denotes the
set of mappings $g\in \mathcal A^{\underline \nu}$ with {\em only} hyperbolic periodic points 
that are topologically conjugate to $f$ on $I$ by an order preserving topological conjugacy which maps
the critical points of $f$  to those of $g$. Since $f,g\in  \mathcal A^{\underline \nu}$,  this conjugacy necessarily preserves 
the order  $\ell_i$ of the critical points. 
Let $\zeta(f)$ be the maximal number of critical points {\em in the basins} of periodic attractors
of $f$ with pairwise disjoint infinite orbits.

\begin{thmx}[Manifold structure] \label{thm:manifold}
  Assume that all periodic orbits of $f\in \mathcal A^{\underline \nu}$ are hyperbolic. Then
  \begin{enumerate}
  \item   $\mathcal T_{f}^{\underline \nu}$   is  an embedded real analytic  submanifold of $\mathcal A^{\underline \nu}$ modelled on a family of Banach spaces of codimension $\nu-\zeta(f)$;
  \item for each $a>0$, $\mathcal T_{f}^{\underline \nu}\cap \mathcal  A^{\underline \nu}_a$   is  an embedded real Banach  submanifold of $\mathcal A^{\underline \nu}_a$ of codimension $\nu-\zeta(f)$.
  \end{enumerate} 
   \end{thmx}

In Section~\ref{manifold-structure}  we will give a formal definition of the term {\em embedded real analytic  submanifold  modelled on a family of Banach spaces}. 
%In fact, in Theorem~\ref{thm:hybridbanach}
%we improve this statement to obtain that $\mathcal T_{f}\cap \mathcal A^{\underline \nu}_a$ is a real Banach submanifold of $\mathcal A^{\underline \nu}_a$.

\begin{thmx}[Topological conjugacy classes are contractible] \label{thm:connected}
   Assume that all periodic orbits of $f\in \mathcal A^{\underline \nu}$ are hyperbolic. Then 
   \begin{enumerate}
     \item   for each $g\in \mathcal T_{f}^{\underline \nu}$  there exists a one-parameter family of real analytic maps
  $f_t\colon I\to I$ in $\mathcal A^{\underline \nu}_a$ with $f_0=f$, $f_1=g$ so
  that $f_t$ depends analytically on $t$ and $f_t$ is topologically conjugate on $I$ to
  $f$ for each $t\in [0,1]$. 
   \item Moreover, if $f\in \mathcal A^{\underline \nu}_a$ then the  manifold $\mathcal T_{f}^{\underline \nu}\cap  \mathcal A^{\underline \nu}_a$ is contractible. 
     \end{enumerate} 
\end{thmx}

In Sections~\ref{sec:manifoldstructure}-\ref{sec:contractible} we will state  corresponding theorems 
with the {\em real-hybrid class}  $\mathcal H_{f}^{\R}$ of $f$, where (in the case that all periodic orbits of $f\colon I\to I$
are hyperbolic)  $\mathcal H_{f}^{\R}$ is defined as the subset of $\mathcal T_{f}$ where the  topological conjugacy extends as a holomorphic conjugacy in a small complex neighbourhood  
of the periodic attractors, or equivalently in a complex neighbourhood of the real 
basin of the periodic attractors.  In particular,   the multipliers at periodic attractors for each map $g\in  \mathcal H_f$ are the same as those for the corresponding periodic attractors of $f$. 
The assumption that all periodic attractors are hyperbolic is not required in the {\lq}hybrid{\rq} analogues of Theorems \ref{thm:manifold} and \ref{thm:connected}, 
but the presence of parabolic periodic points  requires additional assumptions 
(namely that  the parabolic periodic points are {\lq}simple{\rq}) and arguments in the proof.

\begin{thmx}[Partial lamination] \label{thm:lamination}
   Assume that all periodic orbits of $f\in \mathcal A_a^{\underline \nu}$ are hyperbolic. 
  Then $f$ has a neighbourhood which is laminated by hybrid conjugacy classes. More precisely, for
  each neighbourhood $\mathcal V_2$ of $f$  in $\mathcal{A}_a^{\underline{\nu}}$ there exists a   
%	neighbourhood $\mathcal V$ of $f$ in $\mathcal{A}_a^{\underline{\nu}}$ so that for any $g\in \mathcal V$
%	the intersection of $\mathcal H_g$ with $\mathcal V$ is pathconnected. 
neighbourhood 
	$\mathcal V_1\subset \mathcal V_2$ of $f$ in $\mathcal{A}_a^{\underline{\nu}}$ so that for each 
	$g\in \mathcal V_1$ and each $g_0,g_1\in \mathcal V_1\cap \mathcal H_g^\R$ there exists a path 
	$g_t\in \mathcal{A} _a^{\underline{\nu}}$, $t\in [0,1]$ inside 
	$\mathcal V_2\cap \mathcal H_g^{\R}$  connecting $g_0,g_1$.  
\end{thmx}

As mentioned,  Theorems~\ref{thm:manifold}  and  \ref{thm:lamination} were obtained in \cite{ALM} for  
hybrid classes in the context of unimodal maps
with quadratic critical points.  In \cite{GSm}, using completely different methods, an analogous result to Theorem~\ref{thm:manifold} is proved in the setting  of piecewise expanding mappings. 

%\textcolor{red}{QQQ Later in the paper we also will state a corresponding theorem which allows the presence of 
%periodic attractors or neutral periodic orbits. QQQ }
%

\medskip 
One of the main steps in the proof of Theorems~\ref{thm:manifold} and \ref{thm:connected} is 

\begin{thmx}[External maps and pruned polynomial-like structure]\label{thm:prunedpl}
  Associated to each map  $f\in \mathcal A^{\underline \nu}$ with only hyperbolic periodic points 
  there exist 
  \begin{enumerate}
    \item a circle map with discontinuities (called an {\em external map}) associated to $f$;
  \item a  pruned-polynomial-like  mapping $F\colon U\to U'$ which is an extension
  of the real analytic $f$ and so that its domain $U$ is a neighbourhood of $I$. 
  \item For each $f \in \mathcal A_a^{\underline \nu}$ there exists a neighbourhood
  $\mathcal U\subset  \mathcal A_a^{\underline \nu}$ so that each map $g\in \mathcal U$ has 
  a pruned polynomial-like extension $G\colon U_g\to U_g'$ which is obtained from $F\colon U\to U'$
  by holomorphic motion, see Theorem~\ref{thm: pruned-pol-like-map} for a more precise statement. 
  (Note that $g$ is not required to be conjugate to $f$.) 
  \end{enumerate} 
  \end{thmx}
  
The  definition of the notion of  pruned polynomial-like map will be given in Section~\ref{sec:ppl}.
The main point is that they share essentially all useful properties of polynomial-like maps; 
an example is shown in Figure~\ref{fig1}. 
For a more formal statement, see  Theorem~\ref{thm: pruned-pol-like-map}. 
The main ingredients for proving these results is   quasisymmetric rigidity which was proved in 
\cite{CvS} and the complex bounds from \cite{CvST}. These results build on the enhanced nest construction from 
\cite{KSS-rigidity}, see also \cite{KvS,CDKvS}.

%\begin{thmx}

%Theorem A-E also hold for circle mappings, but this will be shown in a separate paper. 
%(The main new aspect in that setting is to develop the notion of {\lq full family} for circle maps.)

%
% Indeed, let $S^1$ be the circle, let $\underline \nu =
%(\ell_1,\dots,\ell_\nu)\in \N^\nu$ with $\ell_i\ge 2$, and $\mathcal A^{\underline \nu}(S^1)$ be the 
%set of real analytic mappings $f\colon S^1\to S^1$ with precisely $\nu$ critical points
%$c_1<c_2<\dots<c_\nu$ such that $Df(c_i)=\dots=Df^{\ell_i-1}(c_i)=0$ and $Df^{\ell_i}(c_i)\ne 0$. 
%Then if $f\in \mathcal B^{\underline \nu}$ has only hyperbolic periodic points,  then $\mathcal T^{\underline \nu}_f$ i
%denotes the set of $g\in \mathcal A^{\underline \nu}(S^1)$ with only hyperbolic periodic points,
%which are  topologically conjugate to $f$ on $S^1$ by an order preserving topological conjugacy which maps the critical points  $c_{i,f}$ of $f$ to the critical point $c_{i,g}$ of $g$. This means that  the critical points are marked.
%
%\begin{thmx}\label{thm:circledense} 
%Theorem A-E also hold when replacing $\mathcal A^{\underline \nu}$ and $I$ 
% by  $\mathcal A^{\underline \nu}(S^1)$ and $S^1$.
%\end{thmx}
%
%The proof of Theorem F will be given elsewhere. 
%Theorem~\ref{{thm:circledense} is a generalisation of Theorem 1.2 
%in \cite{Rempe-Strien}. 

%\begin{thmx}
%\end{thmx} 

% Corollary \ref{cor:polynomials} does not follow from \cite{Koz} or \cite{KSS-density},
% because there polynomials are approximated by hyperbolic polynomials of possibly much larger degree. 

The theorems above are stated for maps $f\in \mathcal A^{\underline \nu}_a$
without parabolic periodic points. Most theorems also go through when $f$ does have
parabolic periodic points, see for example Section~\ref{subsec:parabolic}. 
The lamination structure from  Theorem~\ref{thm:lamination} does not 
hold near maps with parabolic periodic points, and this is the topic of a companion paper
which is currently in preparation.  

\subsection{Further applications} 
In a companion paper it will be shown that the results in this paper imply that one
has density of hyperbolicity within real analytic families of real analytic one-dimensional maps. 

\subsection{Notation and terminology} \label{notation} 
We let $\mathbb C$ denote the complex plane and we say that a subset $A\subset \C$ is {\em real symmetric} if $z\in A$ iff $\bar z\in A$. A function $g\colon \partial \D\to \partial \D$ is called real symmetric if
$g(\bar z)=\overline{g(z)}$. 
If $U\subset\mathbb C$ is an open set, we let
$\mathcal B_U$ denote the set of holomorphic mappings on $U$ 
which are
continuous on $\overline U$. For $F\in \mathcal B_U$ we will define 
$$B_U(F,\epsilon)= \{G\in \mathcal B_U; |G(z)-F(z)|<\epsilon \mbox{ for all }z\in \overline U\}.$$
As before, for $a>0$, we let
$\Omega_a=\{ z\in\mathbb C: \dist(z,I)<a\}$ where $\dist$ is the Euclidean distance on $\C$.

Let $\mathcal A^{\underline\nu}, \mathcal A^{\underline \nu}_a$ are defined as in the introduction.
%and let $\mathcal A^{\underline\nu}_a$
%denote the space of all real analytic mappings $f\in\mathcal A^{\underline \nu}$ 
%which extend holomorphically to $\Omega_a$, continuously to  $\overline {\Omega_a}$ and 
%%with $\nu$ distinct real
%%critical points   and 
%so that for each $x\in (\R\cap \Omega_a) \setminus I$ there exists $n>0$ so that $f^n(x)\notin \Omega_a$. 
%Since, $\partial I$ is repelling and $f$ is real analytic, we can always choose such $a>0$.
We let $\Crit(f)$ (or simply $\Crit$) denote the set of 
critical points of a differentiable mapping $f$. We should emphasise that $f\in\mathcal A^{\underline \nu}$ has at most a finite number of attracting 
periodic points, see \cite{MMvS,dMvS}, but attracting periodic points do not necessarily contain 
a critical point in their basin. We will denote by $B_0(f)\subset I$ the
union of immediate basins of the periodic attractors
of $f\colon I\to I$. Note that some periodic attractors may not contain critical points in their immediate basin.

Given a set $K\subset \C$ we let $\Comp_{x}(K)$  (or $\mbox{cc}_xK$)
denote the connected component of $K$ containing $x$
and if $A\subset \C$ then $\Comp_A(K)$ is the union of the connected components of $K$ containing
points of $A$. 

\label{comp}

\section{Organisation and outline of the ideas used in this paper} 

In Part A  of this paper we introduce the notion of pruned polynomial-like mapping. 
As we will show,  having a pruned polynomial-like mapping is almost as good 
as having a polynomial-like mapping. Indeed, it allows us to extend the  results
of \cite{ALM} for non-renormalizable unimodal maps to the general setting. In particular, 
we do not need to assume {\lq}big bounds{\rq} as in that paper. 
In Section~\ref{sec:ppl} we formulate the notion of pruned polynomial-like mapping
 and 
state Theorem~\ref{thm: pruned-pol-like-map} which asserts 
that {\em each} real analytic interval map has a complex
analytic extension which is a pruned polynomial-like mapping.
Thus we obtain a good Markov structure in a complex neighbourhood 
of the entire dynamical  interval $I$. 
The proof of  the existence of such a pruned polynomial-like extension 
of a real analytic map, goes in several steps:
\begin{enumerate}
\item First we associate to a real analytic map a pruned Julia set.
If the real analytic map is in fact a polynomial, then this pruned Julia set
$K_X$  is simply a subset of the usual Julia set, but pruned so that it is in a small 
neighbourhood of the the interval $I$. Where the Julia set is pruned depends
on a set $X$, where $X$ consists of the boundary points of intervals around the 
critical values of $f$. 
\item Using the uniformisation of $\C\setminus K_X$ and the map $f$ near $K_X$
we obtain a real analytic ordering preserving circle map (whose {\lq}degree{\rq} depends on the degree of the critical points of $f$) with discontinuities (each critical point of $f$ corresponds to a discontinuity of 
$\hat f_X$). This circle map is called the external map $\hat f_X\colon \partial \D \to \partial \D$ associated to $f$ and the pruning data $X$.  Later on we will encode the pruning data in a more combinatorial way by a subset $Q\subset \partial \D$ which is forward invariant under the map $z\mapsto z^{|\nu|}$ where $|\nu|=\sum_{i=1,\dots,\nu} \ell_i$. 
\item The expanding properties of the map $\hat f_X$ can be used to obtain a pruned polynomial-like
structure near the interval $I$.
\end{enumerate} 
This construction is given in  Sections~\ref{sec:prunedJulia}-\ref{sec:existence-pruned-plmap}
and in Section~\ref{sec:exter-qc-conj} it is then shown that the pruned polynomial-like extensions of
two  topologically conjugate interval maps are qc conjugate and that this conjugacy 
preserves the structure of these pruned polynomial-like mappings.

In Part B of the paper we will then show that pruned polynomial-like mappings can be 
treated more or less like polynomial-like mappings. This means that with this structure in hand, 
we no longer need to use a necklace neighbourhood around the $I$ as in \cite{ALM},  but obtain a quasiconformal conjugacy between  two conjugate interval maps 
on a full complex neighbourhood of $I$. This allows us to prove a mating result as in  
\cite{Lyubich-quadratic dynamics}:
associated to two pruned polynomial-like mappings with the same (combinatorial)
pruning data there exists a new pruned polynomial-like mapping which is 
qc conjugate to the first mapping and has the same external mapping
as the second one, see Section~\ref{sec:mating}.  This gives that one 
conjugacy class inherits the manifold structure from another conjugacy class. 
To determine the codimension of these manifolds, we extend the infinitesimal pullback argument 
and key lemma of \cite{ALM} to our setting. Once these techniques are in place, the theorems readily follow.  

The main differences with \cite{ALM} are:
\begin{enumerate}
\item  We do not (and cannot) assume that  there are big bounds.  
In \cite{C} a  polynomial-like argument is used for the case when the geometry is bounded (that argument here would not work when the post-critical set is non-minimal) and arguments
similar to \cite{ALM}  when the geometry is unbounded. Here we do not need to make this distinction.
We overcome the lack of big bounds  by introducing the notion of pruned polynomial-like map
and showing that each real analytic interval map has a complex extension with such a structure.
\item In \cite{ALM} a manifold structure $\mathcal H_f^{\R}$ is obtained 
for each unimodal $f\in \mathcal A^{\underline \nu}_a$ by first 
obtaining the manifold structure for hyperbolic unimodal maps (using the implicit function theorem) 
and then using the density of hyperbolic maps and a lambda-lemma argument (to complex 
perturbations of these maps) to extend this manifold structure for non-hyperbolic maps. 
Here we cannot use this lambda-lemma as our manifolds are of higher codimension (as our maps 
have $\nu$ critical points).
Instead we use our pruned polynomial-like structure to obtain a mating result: 
given two maps one can find a third map which is hybrid-conjugate to the first one
and which has the same external map as the 2nd one. Then we proceed essentially
as in \cite{Lyubich} to inherit a manifold structure on the hybrid class of non-hyperbolic
maps from the manifold structure of a given nearby hyperbolic map. 
\item In \cite{ALM} vertical vectors are obtained, as in \cite{Koz},
by constructing first smooth vertical vector fields and then using a polynomial approximation. 
This is one of the most subtle arguments in \cite{ALM}, for which 
the notion of puzzle maps is introduced. These are 
maps whose domains form a necklace neighbourhood of the interval $I$
(rather than an actual neighbourhood). 
 Here   we can avoid 
this discussion and argue as in the polynomial-like case more or less in the spirit of \cite{Lyubich}. 
\end{enumerate}
A difference with the polynomial-like case dealt with in \cite{Lyubich} is that
\begin{enumerate} 
\item[(4)] in our case (of pruned polynomial-like mappings) the space of these external maps is more complicated than in the case of polynomial-like mappings, due to the existence of 
discontinuities. For that reason we do not attempt  to prove nor use that the space of external
maps (in our setting) has a manifold structure. This means that, for example, our proof of the 
contractibility of conjugacy classes is quite different from  \cite{AL}. Note that such a result is not claimed in \cite{ALM}. 

An important step in our proof of Theorem B2 (that the hybrid class is contractible), is to assign a pruned polynomial-like map to each 
map in a real hybrid class, using a quasiconformal motion (which is obtained using a 
partition of unity of unity argument). Here we use the notion of quasiconformal motion 
(the continuous analogue of a holomorphic motion) introduced in \cite{ST}. 
\end{enumerate}

%\section{A pruned polynomial-like extension of a real analytic map}\label{sec:ppl} 
\part*{Part A: Pruned polynomial-like maps and their external maps}
\label{partA} 
\section{Pruned polynomial-like mappings}\label{sec:ppl} 
% Let $I\subset \mathbb R,$ be a compact interval and $f\colon I\to I$ be a real analytic map with critical points $c_1,\dots,c_\nu\in I$, and assume that $f(\partial I)\subset \partial I$.

\begin{figure}
  \begin{center}
    %\includegraphics[width=10.3cm]{truncated-pol-like-map.pdf}
%    \begin{overpic}[width=12cm]{Aeroplane-julia-set.png}  
%\put (60,17) {\tiny $X$} 
% \arrow <3pt> [.2,.67] from 200  61 to 175  61
%\put (60,5) {\tiny $X$} 
% \arrow <3pt> [.2,.67] from 200  17 to 175  17
%\end{overpic} 
    \begin{overpic}[width=12cm]{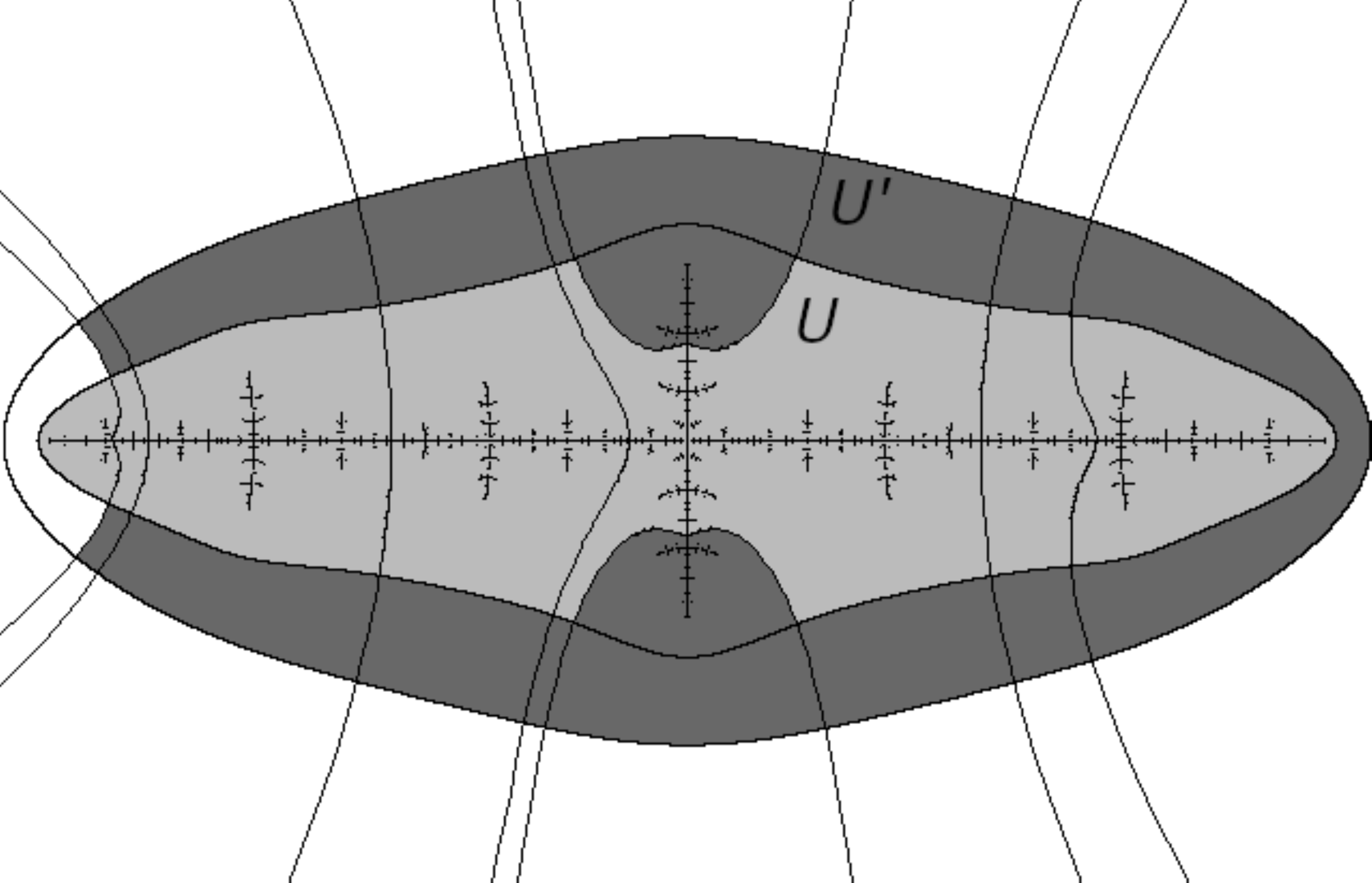}  
\arrow <3pt> [.2,.67] from 190  127 to 174  135
\put (56,36) {\tiny $X$} 
 \arrow <3pt> [.2,.67] from 190  92 to 175  84
  \put (56,26) {\tiny $X$} 
   \end{overpic} 
  \end{center}
\caption{\label{fig1} A pruned polynomial-like map $f\colon U\to U'$ associated to a quadratic interval map $f\colon I\to I$  with only repelling periodic points and whose pruned Julia set  $K_X(f)$ consists of the closure of pre-images of $J$, i.e. pruned at pre-images of the points $X=\partial J$.  The set $\Gamma$ consists of 
the union of the drawn curves (intersected with $\overline U$), consisting of preimages of a
ray landing at a periodic point.  Note that 
$U\setminus U'$ is equal to  the part of $U$ to the left of
$\gamma$. In this figure, the curves drawn are rays and exponentials obtained
from the B\"ottcher coordinates associated  to the quadratic map, 
but in general these curves we will need to be constructed using local methods.}
\end{figure}

One of the main technical innovations in this paper is the introduction
of the notion of pruned polynomial-like map, and the theorem 
that shows that to each real analytic map one can associate such a 
pruned polynomial-like map.

\begin{defn}[Pruned polynomial-like mapping]  \label{prunedpollike}
 We say that a holomorphic map $F\colon U\to U'$ between 
two open sets $U, U'$ consisting of finitely many components 
together with a finite union of arcs or closed curves $\Gamma,\Gamma_a^*,\Gamma_a$  is 
a {\em pruned polynomial-like mapping} if 
\begin{itemize}
\item %$F$ extends continuously to the boundary of $U$; 
$F$ extends holomorphically to a neighbourhood of $\overline U$.
\item 
$U'=F(U)$ and $F(\partial U)\subset \partial U'\cup \Gamma$;  
\item  $U\cap U'\ne \emptyset$, $\partial U\cap \partial U' \subset \Gamma$ and 
$\partial U'\cap   U \subset \Gamma \cup \Gamma^*_a\cup \Gamma_a$;
\item each component of $\Gamma$ is either a piecewise smooth arc in $U'$ connecting boundary points of $U'$
or is contained  in $\partial U'\cup \partial U$;
\item $F(\Gamma\cap \overline{U}) \supset \Gamma$; 
\item each component of $\Gamma^*_a$ and $\Gamma_a$
is contained in the basin of a periodic attractor $F$; 
\item each component of $U'\setminus \Gamma'$, 
$U'\setminus (U\cup \Gamma')$,  $U\setminus \Gamma'$ and $U\setminus (U'\setminus \Gamma'$)
is a quasidisk. Here $\Gamma'=\Gamma\cup\Gamma^*_a\cup \Gamma_a$. 
\end{itemize}
We say that a pruned polynomial-like mapping is {\em real} if $U,U'$ are real-symmetric and $F$ commutes with complex conjugation. 
\end{defn}

\begin{rem} When $U\subset U'$ and $\Gamma=\emptyset$ then $U$ is compactly contained in $U'$
and this definition reduces to that of a polynomial-like map. Moreover, if $F$ has no periodic attractors,
 then $\Gamma_a=\Gamma^*_a=\emptyset$. 
\end{rem} 

\begin{rem} An example of such a pruned polynomial-like map is shown in Figure~\ref{fig1}. This example shows that in general 
  $U\setminus U'\ne \emptyset$ and $U'\setminus U\ne \emptyset$.  
%If $B_0(f)$ is equal to the interior of $I$, as is the case for $f(x)=2x(1-x)$, then $U=\emptyset$. 
In this paper the curves of $\Gamma$ %in  the proof of Theorem~\ref{thm: pruned-pol-like-map}
will consist of {\lq}rays{\rq} landing on real periodic or pre-periodic points of $F$. Note that 
the pruned polynomial-like map from Figure~\ref{fig1} is not proper. 
\end{rem}

\begin{rem} If $f$ has a periodic attractor with a critical point in its basin, 
then there exists a bit of $\partial U'$ which is mapped to a curve which 
iterates to the attracting fixed point, see Figure~\ref{fig:basin}.   
\end{rem} 

If $F:U\to U'$ is a pruned polynomial-like mapping, we call
$$K_F:=\{z\in \C; F^n(z)\in \overline{U} \mbox{ for all }n\ge 0\}$$
the {\em pruned filled Julia set} of $F:U\to U'$, and  $J_F:=\partial K_F$ its 
{\em pruned Julia set}. 
 Observe that while $K_F\subset\overline U$, $K_F$ need not be contained in $U$
and also does not need to be backward invariant, see for example Figure~\ref{fig1}.
We will see that for the pruned polynomial-like mappings we construct  in 
the theorem below,  one has $K_F\cap\partial U=K_F\cap\partial U'$.

\medskip 

Let us now show that one can associate a pruned polynomial-like map to a real analytic map 
$f\in\mathcal A^{\underline \nu}_a$. 
 Since $\Omega_a$ might be a very small neighbourhood of $I$, one cannot expect $f^{-1}(\R)$ (or even $f^{-1}(I)$) to be contained in the interior of $\Omega_a$, and so one should not expect $f\colon I\to I$ to have a genuine polynomial-like extension. In this section we will construct something almost as good as a polynomial-like extension, namely the extension given by the following theorem (which contains part of  Theorem~\ref{thm:prunedpl}).

\begin{thm}[Pruned polynomial-like mapping associated to $f$] \label{thm: pruned-pol-like-map}
  Suppose that $f\in\mathcal A^{\underline \nu}_a$ has only hyperbolic periodic points. 
  Associated to $f\colon I\to I$ there exist 
open sets $U, U'$ in the complex plane and a finite union of curves $\Gamma,\Gamma_a^*,\Gamma_a$ so that
  \begin{enumerate}
  \item  $f$ extends to a map 
  $F\colon U\to U'$ which together with $\Gamma,,\Gamma_a^*,\Gamma_a$ forms  a pruned polynomial-like mapping.
%  \item 
%$U$ is a  tubular neighbourhood of $I\setminus B_0(f)$, i.e. 
%for each point in $\partial B_0(f)$ there is a curve in $\Gamma$  landing on it and
%near such a point $\partial U$ is equal to a union of two such curves. 
%\marginpar{Do we want to $U\subset I$. Can we skip 2?} 
\item $U$ is  a neighbourhood of $I$ and for each $\delta>0$ the sets $U,U'$ can be chosen so that they are 
contained in a $\delta$-neighbourhood of $I$.
%\item  For each $f \in \mathcal A_a^{\underline \nu}$ there exists a pruned polynomial-like map  
%$F\colon U\to U'$ which is a complex extension of $f$.
\item   there exists a neighbourhood  $\mathcal U\subset  \mathcal A_a^{\underline \nu}$ of $f$
so that each map $g\in \mathcal U$ has 
  a pruned polynomial-like extension $G$ which is obtained from $F\colon U\to U'$
  by holomorphic motion. More precisely, there exists 
  a holomorphic motion $H_\lambda$, $\lambda \in \mathcal U$ of $U,U'$ where $H_f=id$ so that  each $g\in \mathcal U$ has a complex extension
  $G\colon H_g(U)\to H_g(U')$ which is a pruned polynomial-like mapping. 
  Moreover, $G(z)= H_g^{-1}\circ F\circ H_g(z)$ for $z\in \partial U$. 
 \end{enumerate}   
 \end{thm}

To obtain the pruned polynomial-like extension $F\colon U\to U'$ of the real analytic map 
$f\colon I\to I$ we will associate to $f$:
\begin{enumerate}
\item[-] A {\lq}pruned Julia set{\rq} $K_X\subset \C$ associated to $f$.
\item[-]  An {\lq}external{\rq} map $\hat f_X\colon \partial \D \to \partial \D$ which is real analytic, 
except in a finite number of points where it is discontinuous. 
\item[-]  Exploit the expanding structure  of $\hat f_X$
together with the attracting structure near the immediate basins of the periodic attractors of $\hat f_X$. 
\item[-] The persistence of the pruned polynomial-like structure is proved in Proposition~\ref{prop:persistenceprunedstructure}.  
\end{enumerate}

\begin{rem} If $f$ has no periodic attractors, then $\Gamma_a^*=\Gamma_a=\emptyset$ and 
 $\Gamma$ consists of {\lq}rays{\lq} landing at points which are eventually 
mapped onto repelling periodic points. If $f$ the periodic attractors are {\em small} in the sense
of Definition~\ref{def:smallbasin} (roughly speaking this means that the basins are contained
in the domain of analyticity of $f$) then we also can ensure that 
$\Gamma_a^*=\Gamma_a=\emptyset$, see Theorem~\ref{thm:smallbasins}. 
\end{rem}
\begin{rem} If $f$ has periodic attractors, then 
$\Gamma$ may also consist of pieces of curves which land on periodic attractors, see Section~\ref{polstructure-attractors}. However, if the basins of periodic attractor of $f$ are compactly contained in the domain of analyticity of  $f$, then one can ensure that these basins are compactly contained in the domain of $F\colon U\to U'$ and one can take again 
$\Gamma^*_a=\Gamma_a=\emptyset$, see Subsection~\ref{subsec:KXattractors}. 
\end{rem} 
\begin{rem}
Even if $f$ is allowed to have simple
parabolic periodic points (see  Definition~\ref{def:simplepara}),  the previous theorem still holds, see
Subsection~\ref{subsec:parabolic} except that the statement of part (4) of Theorem~\ref{thm: pruned-pol-like-map} needs to be amended, see Subsection~\ref{subsec:holoparabolic}. 
\end{rem}

A somewhat similar result to Theorem~\ref{thm: pruned-pol-like-map} is obtained in \cite[Appendix B]{ALM} for unimodal maps, with non-degenerate critical points, which are at most finitely many times renormalizable. Their proof requires the mappings have a large scaling factor at some level of the principal nest;
it does not  go through for multimodal mappings or unicritical mappings of higher degree. Consequently, our construction is necessarily quite different and more abstract, since we need to deal with multimodal maps and infinitely renormalizable maps.

To construct the associated pruned-polynomial-like map, we will first define  the notion of a
pruned Julia set, show that this is near the real line 
and then associate to this object (and the map $f$) some circle map with discontinuities. 

\subsection{Notation used for interval maps and their extensions} Interval maps will be denoted by   $f\colon I\to I$ \label{notationfandfX} 
as well as its complex extension $f\colon U\to \C$.  An extension which will have 
the structure of a pruned polynomial-like mapping will be denoted by $F\colon U\to U'$. 
The map $f\colon U\to \C$ will induce (via some Riemann mapping $\phi_X$) 
an {\em external map}  $\hat f_X\colon \partial \D \to \partial \D$. The  complex extension 
of this map will again be denoted by $\hat f_X$ but for the part where this 
extension to a neighbourhood of $\partial \D$
has an expanding Markov structure  we will write $\hat F_X\colon \V\to \V'$. 
The sets $\V,\V'$ will be used to construct the {\lq}expanding{\rq} part
of the pruned polynomial-like extension of $f$, denoted by  $F\colon \E\to \E'$. 
We also have a suitable structure near basins of $f$, which gives an extension $F\colon B\to B'$
of $f$ near basins. Taking $U:=\E\cup B$  and $U'=\E'\cup B'$ we will obtain
the required pruned polynomial-like map $F\colon U\to U'$. 
%of a pruned polynomial-like mapping it will be denoted by $\hat F_X\colon \hat U \to 
%\hat U'$. In order to construct the sets $U,U',\hat U,\hat U'$ we will 
%use subsets $\hat E,\hat E'$ corresponding to the expanding part
%of $\hat f_X$ and sets $\hat B,\hat B'$ corresponding
%to the basins of $f$ and $\hat f_X$. So for example, $\hat F_X\colon \hat E \to \hat E'$ will be an 
%expanding Markov map. 

\section{A pruned Julia set associated to $f\colon \Omega_a\to \C$} \label{sec:prunedJulia} 
As mentioned, to prove Theorem~\ref{thm: pruned-pol-like-map} we will first
define the notion of a pruned Julia set. Let $1\leq \nu'\leq\nu$ denote the
number of distinct critical values of $f\colon I\to I$. Take (real) disjoint interval
neighbourhoods $J_1,\dots,J_{\nu'}\subset I$ of the critical values $f(c_1),\dots,f(c_\nu) \in I\subset\R$
so that the closure of the component of $f^{-1}(J_i)$ containing $c_i$ is contained in $\Omega_a=\{z\in \C:
\dist(z,I)<a\}$. 
Next set $J=\cup J_i$ and define
\label{XJi} 
$$J^{-1}=\overline{f^{-1}(J)\setminus \R}\mbox{ and } X=\partial J^{-1}.$$
For later use, also define 
$\Crit'(f)$ to be the subset of critical points of $f$ which are periodic.
%and let 
%\textcolor{red}{$$X=X_0\cup \{ c_i; c_i \mbox{ is periodic}\}.$$} %in the basin of a periodic attractor}\}.$$
We shall call $J_i$ {\em pruning intervals} and the points $X$ {\em pruning points}. 
By choosing the intervals $J_i$ sufficiently small, we can ensure that 
\begin{enumerate}
\item each $J_i$ contains exactly one critical value of $f$, 
\item  
%   By the real bounds, see \cite[Theorem A']{vSV}, we can and will assume that each $J_i$ contains exactly one critical value of $f$.
 $f^{-1}(J)$ consists of arcs whose closures are contained in $\Omega_a$;
 \item if the critical value $f(c_i)$ in $J_i$ is contained in the basin of a periodic attractor, 
 then we assume that $J_i$ is contained compactly in this basin, and moreover  does not any point 
 in the backward orbit of a non-periodic  critical point.
  \end{enumerate}
The intervals $J$ (and the finite set $X$)  will be used to define the pruned Julia set associated to $f\colon \Omega_a \to \C$.

To define the pruned Julia set $K_X(f)$, 
inductively define $K_n$ by taking $K_0=I$, 
\begin{equation}
K_1=K_0\cup    \mbox{cc}_{K_0} J^{-1} 
\end{equation} 
and for $n\ge 1$, 
\begin{equation}
K_{n+1}:= K_n\cup    \mbox{cc}_{K_n'} f^{-{n}}( J^{-1}) . \label{def:Kn} \end{equation} 
where $\mbox{cc}_{K_n'}$ stands for the connected components which intersect
$K_n':=K_n\setminus \Crit'(f)$. 
So 
in the definition of $K_{n+1}$ we exclude any preimage of  $J^{-1}$ which once again intersects a critical point. 
Implicit in the definition (\ref{def:Kn}) is that $J$ can be chosen 
so that $\mbox{cc}_{K_n'} f^{-{n}}( J^{-1})$ is again in $\Omega_a$ for each $n\ge 0$. 
In Theorem~\ref{thm:KX}  we will show that one can choose $J$
so that this property indeed holds. 

For example,  the set $K_1$ consists of $I$ together with 
$2(\ell_i-1)$ arcs  in $\mathbb C\setminus\mathbb R$ attached to $c_i$ 
for each $i$. 
Note that $K_{n+1}\setminus
K_n\subset \mathbb C\setminus\mathbb R,$ and is contained in
$f^{-n}(K_1\setminus I)$. 
It is immediately clear from the definition that for all $n\ge 0$, 
$$f^{-1}(K_n)\supset K_n\mbox{ and }K_{n+1}\supset K_n.$$ 

Next define 
\begin{equation} K_X(f)=\mbox{ closure of } \bigcup_{n=0}^\infty K_n .\label{KXf} \end{equation}
Note that $K_0\subset K_1\subset K_2\subset \dots $ are finite trees with finite degree 
whose {\lq}complex endpoints{\rq} or {\lq}pruning points{\rq} 
are preimages of $X$.   We will call $K_X(f)$ the {\em pruned Julia set} associated to $f\colon \Omega_a \to \C$ associated to the set $X$.
Note that for the pruned polynomial-like extension $F\colon U\to U'$ that we will construct
we may not have  $K_X(f)\subset K_F$. However, 
if $f$ has only repelling periodic orbits then we will obtain  $K_F=K_X(f)$, see Lemma~\ref{lem:equalityKXKFX}.

The reason we use the notation $K_X(f)$ rather than $K_{J_i}(f)$
is that for real polynomials without periodic critical points, $K_X(f)$ is obtained 
from the filled Julia set of $f$ by {\lq}pruning{\rq} the filled Julia set at
all preimages of $X$, see Example~\ref{example:K_X}(ii). 
%This point of view is also
%made precise in the alternative definition of $K_X(f)$ given in Lemma~\ref{def:alternative}. 
If $f$ has periodic critical points, 
then we prune more because otherwise we would not be able to ensure that $K_X(f)$
is close to the real line.

\begin{figure}
\centering 
 \subfigure[]{\includegraphics[width=.4\linewidth]{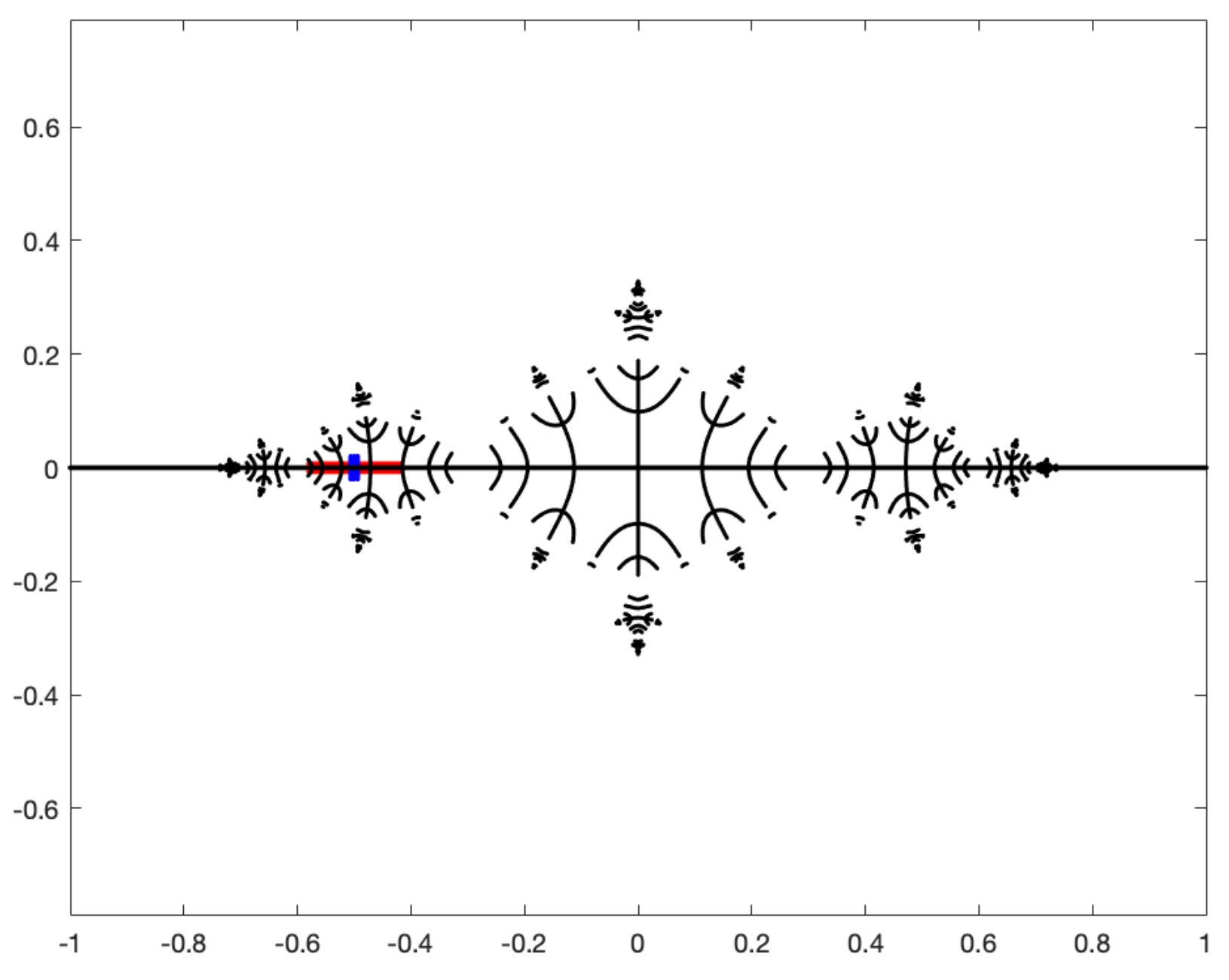}}
\subfigure[]{\includegraphics[width=.4\linewidth]{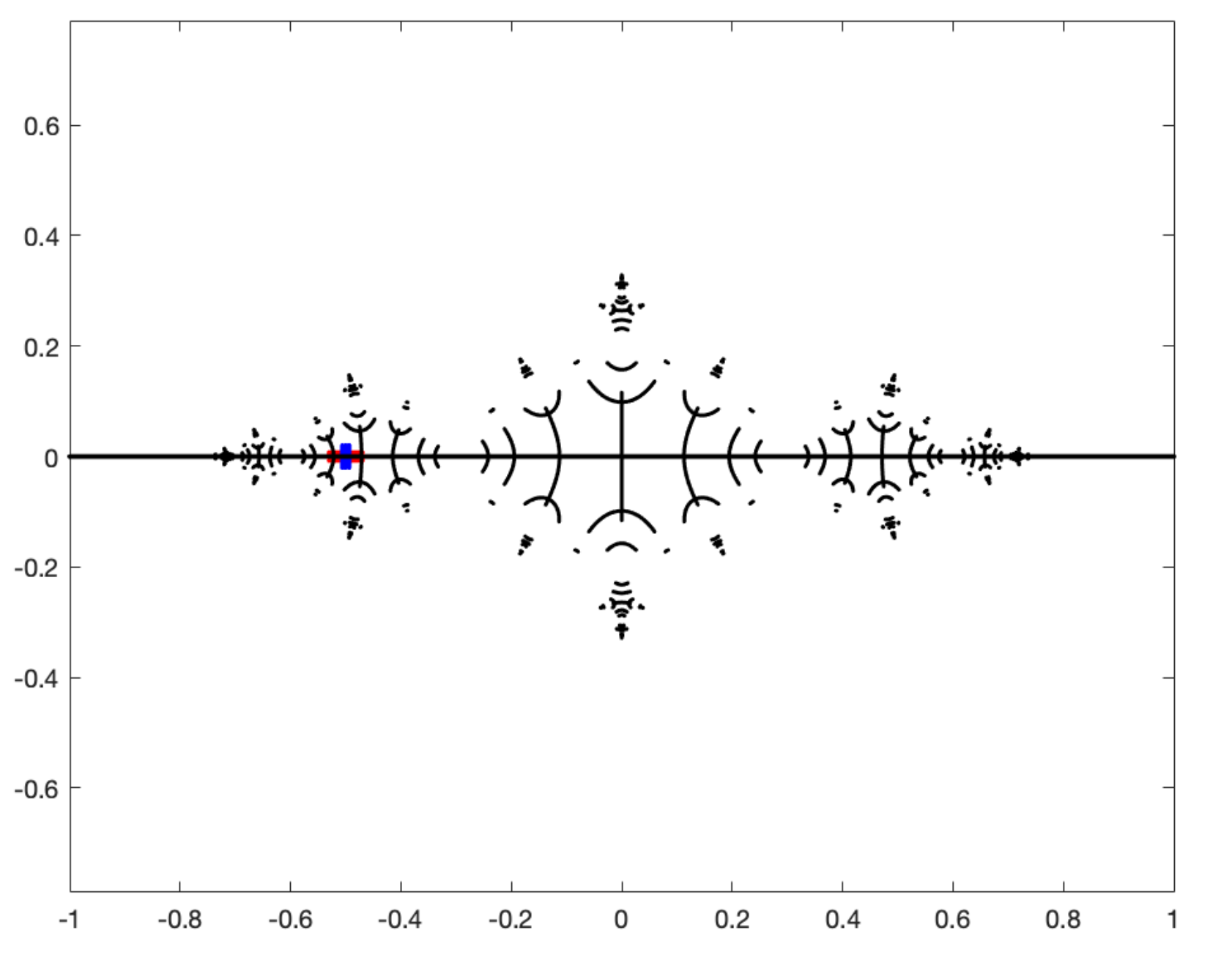}}
\caption{The $\cup_{i=1,\dots,8}  f^{-i}(J^{-1})$ corresponding to two slightly different small {\lq}pruning{\rq} 
intervals $J$ (drawn in red) around   the critical value of a real  quadratic map $f$. Note the different places where the tree structure is {\em pruned}. 
 If we denote by $X=f^{-1}(\partial J)\setminus \R$, then the pruned Julia set $K_X$  corresponds to the complex pre-images of $J$ connected to $I$. By taking the interval $J$ small enough, we can ensure that $K_X$ lies near $I$ (and so in the domain of definition of an extension of $f$ to the complex plane).}
\label{fig:truncated}
\end{figure}

\begin{figure}
\centering 
 \subfigure[]{\includegraphics[width=.4\linewidth]{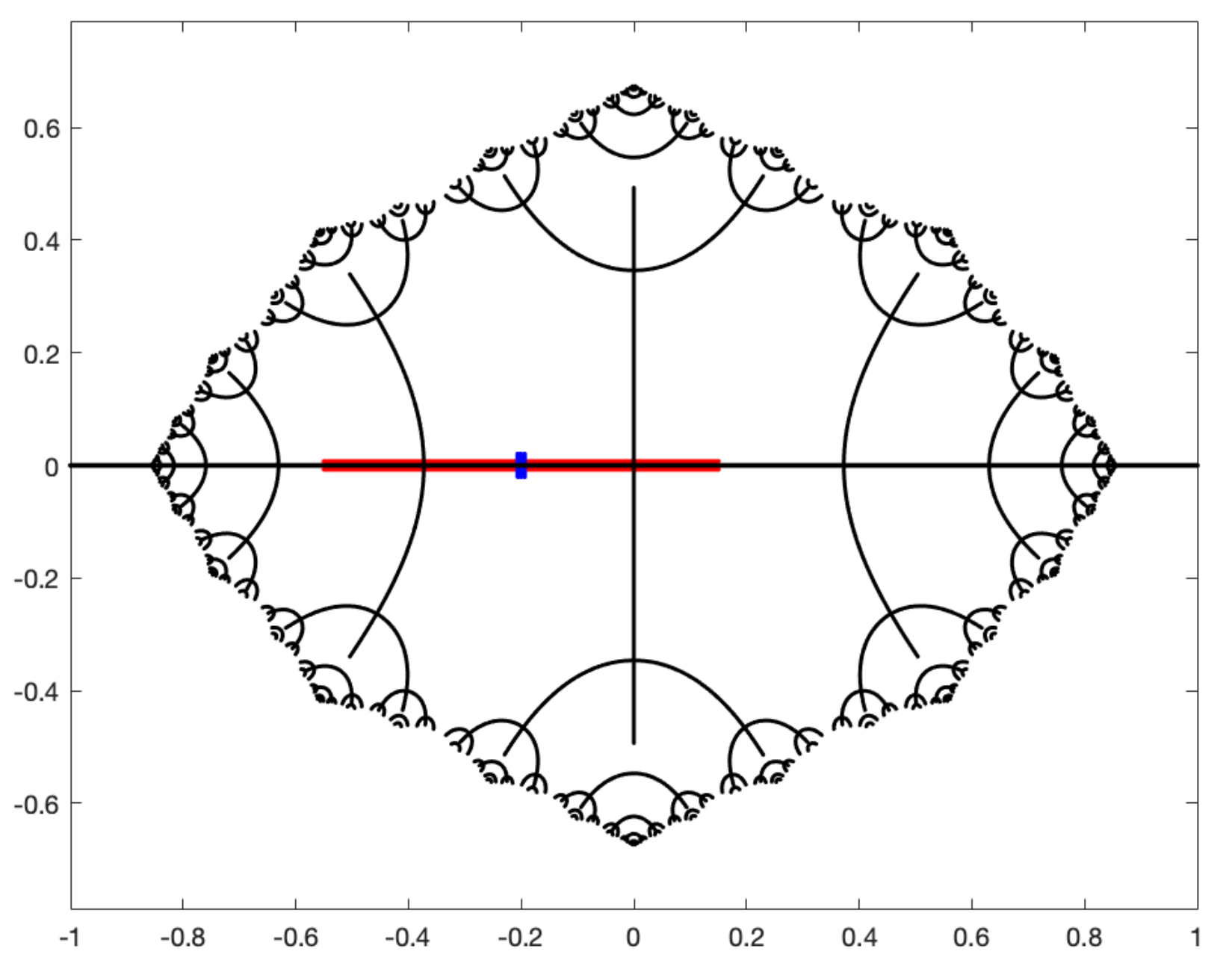}}
\subfigure[]{\includegraphics[width=.4\linewidth]{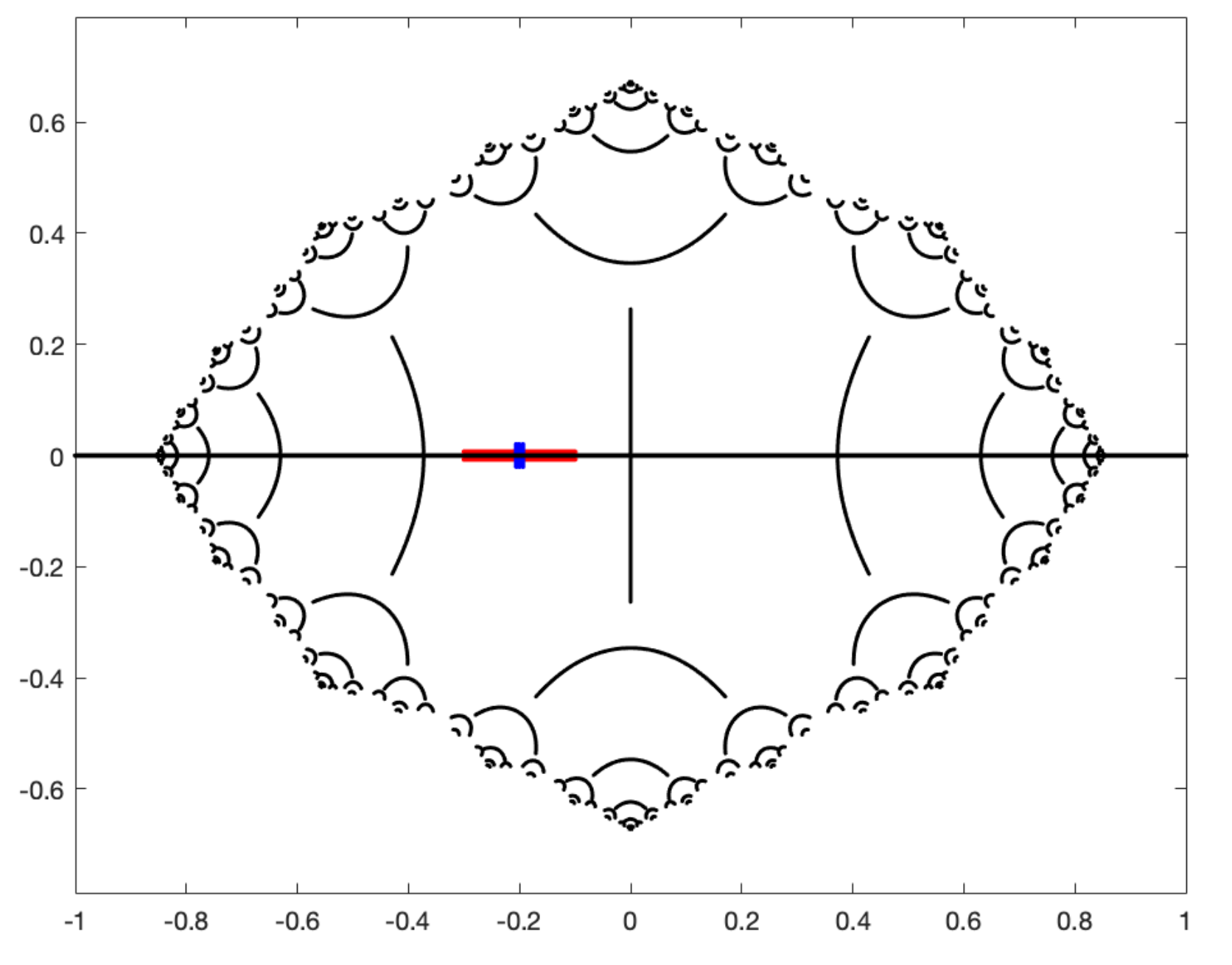}}
\caption{The quadratic map $F_c(z)=(-c+1)z^2+c$, normalised so that $F_c(\pm 1)=1$ with $c=-0.2$. 
Here we take an interval $J$  around the critical value $c$ in the basin of an  attracting fixed point for two choices of $J$: The set 
$K_X(F)$ shrinks as $J$ gets smaller, and in particular when $J$ does not contain a backward iterate
of the critical point. Here, as everywhere, $X=\partial J^{-1}$. }
\label{fig:truncated2}
\end{figure}

\begin{figure}
\centering 
\begin{overpic}[width=12cm]{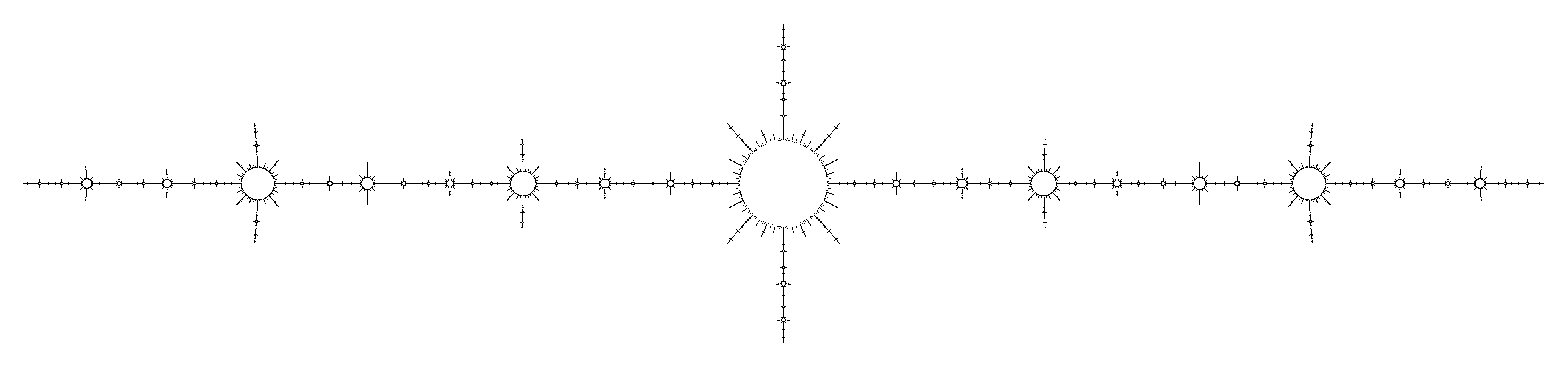}  
\put (60,17) {\tiny $X$} 
 \arrow <3pt> [.2,.67] from 200  61 to 175  61
\put (60,5) {\tiny $X$} 
 \arrow <3pt> [.2,.67] from 200  17 to 175  17
\end{overpic} 
% \subfigure[]{\includegraphics[width=.4\linewidth]{Aeroplane-julia-set.png}}
\caption{If the basins are near the interval $I$, then one can prune along all preimages of the finite set $X=\partial J^{-1}$
suggested in this figure. 
If $f$ is real analytic, 
then $K_X(f)$ is the analogue of the Julia set {\lq}near{\rq} the interval $I$.} 
\label{fig:truncated3}
\end{figure}

\begin{example}[The set $K_X$] \label{example:K_X} 
Let us illustrate the notion of a pruned Julia set in some cases: \\
(i) If $I=[-1,1]$, $J=(-\epsilon,\epsilon)$ with $\epsilon<1$ and $f(x)=x^\ell$ then $K_0=I$ and for $n\ge 1$, 
$$K_n= \bigcup_{k=1,\dots,\ell-1} e^{\pi i k/\ell} \cdot  J'   \,\, \bigcup \,\, I 
\mbox{ where }J'=[-\epsilon^{1/\ell},\epsilon^{1/\ell}].$$ 
Note that $K_n=K_1$ for $n\ge 1$ because the only pre-image of the critical point
$c_1=0$ is itself and because of the 
definition of $K_1'=K_1\setminus \{0\}$. So in this case $K_X(f)$ consists of $\ell$ curves through $0$. 
If we replaced in the definition (\ref{def:Kn}) $\mbox{cc}_{K_n'} f^{-n}(J)$ 
by $\mbox{cc}_{K_n} f^{-n}(J)$ then we would have obtained that $K_X(f)$
is equal to the closed unit disc, which is undesirable for our purposes. 
Note that the relevance of $K_n'$ in the inductive definition of $K_{n+1}$
only occurs when $f$ has a periodic critical point. 
\\ 
(ii) 
Let $I=[-1,1]$, $f(z)=(-c+1)z^2+c$, normalised so that $f(\pm 1)=1$. If $c<0$ then 
$f$ has an attracting orientation reversing fixed point $a\in (c,0)$.
Let $x\in (-1,a)$ be so that $f(x)=0$. 
If the pruning interval $J\ni c$ is chosen so large that it contains $[x,c]$ then 
$\C\setminus K_X$ has bounded components, for example the region $A$ which lies NW of $0$, see 
Figure~\ref{fig:truncated2}(a). 
Note that $A$ is mapped by $f$ to its complex-conjugate $\overline A$, and that $A$ and $\overline A$
each contain a periodic point of period two in its boundary. Moreover, $f$ maps 
$B:=\interior(A\cup \overline A)$ into $B\setminus [x,c]$ and so is contracting w.r.t. the hyperbolic metric. 
Note that we require that $J$ does not contain backward iterates of non-periodic critical points, 
so in actual fact, $\C\setminus K_X$ will have no bounded components, see  Figure~\ref{fig:truncated2}(b). 
\\
(iii)  Assume that $I=[-1,1]$ and $f\colon I\to I$ is a polynomial so that all its critical points are real 
and  non-periodic, 
then for each $n\ge 0$ we have $K_n=\{z\in \C; f^n(z)\in I\}$  when $J=(-1,1)$ and so $K_X$
is the usual Julia set of $f$.  On the other hand,
if $J\Subset (-1,1)$ then  $K_n$ is a subset 
of the Julia set of $f$ pruned  at preimages of $X=\partial f^{-1}(J)\setminus \R$. Again $K_n$ 
consists of a finite union of smooth curves and $K_X$ is the usual Julia set
but pruned at infinitely many preimages of $X$. 
The sets $K_n$ are illustrated in Figure~\ref{fig:truncated}, \ref{fig:truncated2} and \ref{fig:truncated3}.  \\
(iv)  if one critical points
$c_1$ of degree $\ell_1$ is mapped after $k$ steps onto another critical point $c_2$ of degree $\ell_2$, then $K_1$ has $\ell_i$ curves emanating from $c_i$, $i=1,2$ but then $K_{k+1}$ has $\ell_1\cdot\ell_2$
curves emanating from $c_1$. \\
(v) if $f$ has no periodic attractor,  then the backward orbits of critical points are dense 
in $I$ (this follows from the absence of wandering intervals).  If $z\in I$ is such a point, 
and $f^n(z)$ is a critical point then $f^{-(n+1)}(J)$ will contain a curve through
$z$ transversal to $I$.   On the other hand, 
if  $f$ has a periodic attractor containing a critical point $c$ and $J\subset I$ is a very small neighbourhood of $f(c)$, then there may not be any points in the backward orbit of any critical point in $J$. In that case $K_n$ has a particularly simple tree structure. 
\end{example}

\begin{figure}[htp]
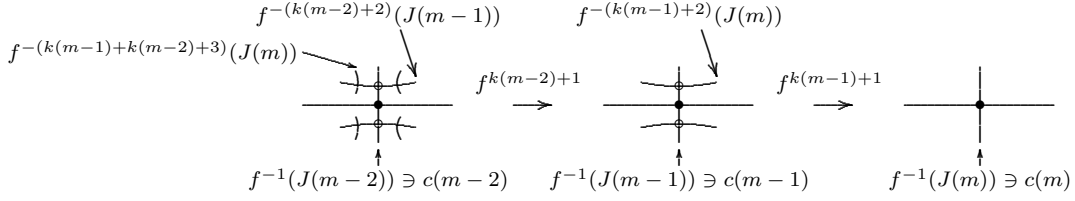
  \hfil \beginpicture \dimen0=1cm 
  \setcoordinatesystem units <\dimen0,\dimen0> point at 50 0 
  \setplotarea x from 1 to 20, y from -1.3 to 1.3  
  \setlinear 
  \plot 10 0 12 0 / 
  \plot 11 -0.5 11 0.5 / 
\put {\tiny $f^{-1}(J(m))\ni c(m)$} at 11 -1 
 \put {\tiny $\bullet$} at 11 0 
    \arrow <3pt> [.2,.67] from 11  -0.8 to 11  -0.6 
\plot 6 0 8 0 / 
\plot 7 -0.5 7 0.5 / 
  \setquadratic
  \plot 6.5 0.3 7 0.25 7.5 0.3 / 
    \plot 6.5 -0.3 7 -0.25 7.5 -0.3 / 
    \put {\tiny $\circ$} at 7 0.25
      \put {\tiny $\circ$} at 7 -0.25
    \put {\tiny $f^{-(k(m-1)+2)}(J(m))$} at 7  1.2 
    \arrow <5pt> [.2,.67] from 7.3  1 to 7.5  0.38 
        \arrow <3pt> [.2,.67] from 7  -0.8 to 7  -0.6 
    \put {\tiny $f^{-1}(J(m-1))\ni c(m-1)$} at 7 -1 
     \put {\tiny $\bullet$} at 7 0 
          \arrow <5pt> [.2,.67] from 4.8  0 to 5.3  0
  \put {\tiny $f^{k(m-2)+1}$} at 5 0.3 
                 \arrow <5pt> [.2,.67] from 8.8  0 to 9.3  0
                             \put {\tiny $f^{k(m-1)+1}$} at 9 0.3 
    \setlinear 
    \plot 2 0 4 0 / 
     \plot 3 -0.5 3 0.5 / 
 \setquadratic
  \plot 2.5 0.3 3 0.25 3.5 0.3 / 
   \plot 2.5 -0.3 3 -0.25 3.5 -0.3 / 
       \put {\tiny $\circ$} at 3 0.25
      \put {\tiny $\circ$} at 3 -0.25
   \put {\tiny $f^{-(k(m-2)+2)}(J(m-1))$} at 3  1.2 
    \arrow <5pt> [.2,.67] from 3.3  1 to 3.5  0.38 
    \arrow <3pt> [.2,.67] from 3  -0.8 to 3  -0.6 
    \put {\tiny $f^{-1}(J(m-2))\ni c(m-2)$} at 3 -1 
       \put {\tiny $\bullet$} at 3 0 
     \setquadratic
  \plot 2.7 0.45 2.73 0.3 2.7 0.15 / 
  \plot 2.7 -0.45 2.73 -0.3 2.7 -0.15 / 
  \plot 3.3 0.45 3.25 0.3 3.3 0.15 / 
  \plot 3.3 -0.45 3.25 -0.3 3.3 -0.15 /   
\put {\tiny $f^{-(k(m-1)+k(m-2)+3)}(J(m))$} at 0 0.7 
    \arrow <3pt> [.2,.67] from 2  0.7  to 2.7  0.5 

%        \put {\tiny $f^{-((m-1)+2)}(J(m))$} at 7  1.2 
  %  \arrow <10pt> [.2,.67] from 7.3  1 to 7.5  0.38 
  % -4 4 -4 4 4 -4 4 -4 -4 / \plot -4 -4 4 4 / \setquadratic
%%   \plot -4 -4 -2 -0.3 0 4 /
%\plot 0 -4 1 -2.6 2 -1.3 / \plot 2 1.3 3 2.3 4 4 / \plot 0 4 -1 2.6 -2 1.3 / \plot -2 -1.3 -3 -2.3 -4 -4 / 
%\setdots \setlinear \plot 2 1 2 -4 / 
%\put {$\f_X$} at 1 -1 
%\put {$\hat x_1^f$} at -4.5 -1.4 
%\put {$\hat x_2^f$} at -4.5 1.4 
%\put {$\hat x$} at 2 -4.6
%\setsolid \plot 1.6 -3.95 2.4 -3.95 / \plot 1.6 -3.9 2.4 -3.9 / \put {$\hat Y$} at 2.4 -3.3 \setdots 
%\setlinear \plot -2 1 -2 -4 / \put {$\hat x'$} at -2 -4.6 \setsolid \plot -1.6 -3.95 -2.4 -3.95 / \plot -1.6 -3.9 -2.4 -3.9 /
% \put {${\bar \C} \setminus  K_X$} at -12 -4.4
%\arrow <10pt> [.2,.67] from -10.5 -4.4 to -7.1 -4.4 \put {${\bar \C} \setminus {\bar {\mathbb D}}$} at -6 -4.4
% \put {$\phi_X$} at -9 -3.2
\endpicture
\caption{\label{fig:tree} The tree structure of the sets $K_n$. The critical points
are marked in the figure by the symbol {\tiny $\bullet$}. The symbol {\small $\circ$}
indicates the points in $f^{-1}(z(m-1))$ 
in the components of $f^{-1}(J(m-1))$ containing $c(m-1)$ respectively the points
in $f^{-1}(z(m-2))$ 
in the components of  $f^{-1}(J(m-2))$ containing $c(m-2)$.}
\end{figure}

\subsection{Maps whose periodic attractors have small 
basins} \label{subsec:KXattractors} 
In general, the basin of a periodic attractor containing a critical point $c$
in its basin could intersect the boundary of the domain of analyticity of $f$. 
In this case, we will choose the pruning interval $J\ni f(c)$ so it is contained
in the basin of a periodic attractor. However, if $f$ has periodic attractors whose basins are small enough to fit 
compactly  inside the domain of analyticity of $f$, then we can define a larger version of $K_X(f)$
as follows.  More precisely, 

\begin{defn}[Small basins] \label{def:smallbasin}  %Let $O$ be a union of periodic attractors of $f$. 
We say that  $f$ has {\em small basins} 
or that the {\em  basins of $f$ fit inside the domain of analyticity} 
if for each periodic attractor which contains a real critical point in its basin the following holds:
\begin{enumerate}
\item[$\bullet$]  For each periodic attractor 
$a$, let  $B_0(a)$ be the  immediate (complex) basin  of a periodic attracting point $a$ (possibly parabolic) 
and let $B_\R(a)$ be the union of the preimages of $B_0(a)$ that intersect the real line. 
Then $$\overline{B_\R(a)}\subset \Omega,$$
where $\Omega$ is a domain of analyticity of $f$ on which $f|\Omega$  has only critical points on the real line.
\end{enumerate} 
%Note that periodic attractors without critical points in their basin do not contribute to $K_{X,O}$.  
\end{defn} 

In this definition  it is allowed that $f$ has parabolic periodic points.
If $f$ has small basins,  we take the pruning intervals $J_i$ so that for 
each critical point $c_i$ of $f$  one has 
that $\comp_{f(c_i)}B_0(a)\cap \R$ is compactly contained in $J_i\ni f(c_i)$ 
and so that the set $K_{X,O}(f)$ is compactly contained in a domain of analyticity of $f$, where 
$$
K_{X,O}(f):=K_X \cup \bigcup_i \mbox{cc}_{K_X} \overline{B_{i,O}} 
$$
and where $B_{i,O}$ are the connected components of $B_O$. 
Implicit in this definition 
is that $J_i$ can be chosen 
so that $\mbox{cc}_{K_X} \overline{B_{i,O}}$ is again in $\Omega$. 
In the proof of Theorem~\ref{thm:KX}  we will show that if $f$ has small basins, 
then we can indeed choose $J_i$ so that this is the case. 

\bigskip

If $f$ has small basins, then we have the following analogue of Theorem~\ref{thm: pruned-pol-like-map}:

\begin{thm}\label{thm:smallbasins} 
 If $f\in \mathcal A^{\underline \nu}$ has small basins, 
then  the properties of  Theorem~\ref{thm: pruned-pol-like-map}
hold with $\Gamma^*=\Gamma_a=\emptyset$. 
\end{thm}

\begin{rem} The above enlargement $K_{X,O}(f)$  ensures that the attracting periodic orbits 
 are not visible for the external map $\hat f_X$ that we define in the next section. 
Indeed, if  $f$ has  small basins, 
%and 
%can be included in $K_{X,O}$, 
then the circle map $\hat f_X$ that we will construct below will have 
no periodic attractors, see Lemma~\ref{lem:fXproperties} below. 
\end{rem} 

\begin{rem} Assume that $f$ has only repelling periodic points.
Then not only we will obtain a pruned polynomial-like extension $F$ of $f$, 
but we will even obtain (via holomorphic motion) a  pruned polynomial-like extension $\tilde F$
 for all maps $\tilde f$ near $f$. 
 Of course $\tilde f$ may have attracting periodic orbits (which will necessarily have small basins)
and the basins of these attractors will be compactly contained in the domain 
of the pruned polynomial-like extension $\tilde F\colon \tilde U\to \tilde U'$, as we will show
in Section~\ref{sec:expandingpruned}.
\end{rem}

%
%
%define $\hat B_i=B_i$ 
%and otherwise define $\hat B_i=\emptyset$.  Then as before let
% $$J^{-1}=\overline{f^{-1}(J)\setminus \R}$$
%\textcolor{red}{Check QQQ} 
%\begin{equation} \J = J \cup B_{\small} \mbox{ where }B_{\small}=\cup_i \hat B_i. \end{equation} 
%\begin{equation}
%K_1=K_0\cup \mbox{closure of }   \mbox{cc}_{K_0} f^{-1}(\J)
%\end{equation} 
%and for $n\ge 1$, 
%\begin{equation}
%K_{n+1}:= K_n\cup \mbox{closure of }   \mbox{cc}_{K_n'} f^{-n}(\J) . \label{def:Kn} \end{equation} 
%where $\mbox{cc}_{K_n'}$ stands for the connected components which intersect
%$K_n':=K_n\setminus \Crit'(f)$.

\subsection{How the tree structure of $K_n$ is created}  Although not needed in our discussion, let
us discuss the tree structure of $K_n$. 
Take  $z\in K_n\setminus \R$. Then the component $L_n(z)$ of $K_n\setminus \R$  containing $z$ contains a path connecting $z$ to $\R$
consisting of a finite number of smooth arcs which cross each other transversally  which are generated as follows:  
there exist $0\le m\le n$ and $k\ge 0$, and for each $1\le i\le m$ there exist a critical point $c(i)$  of $f$, an interval 
$J(i)$  from the collection $\{J_1,\dots,J_{\nu'}\}$  containing the critical value
$f(c(i))$, a point $z(i)\in J(i) \setminus \{c(i)\}$  and an integer $k(i)\ge 1$  with the following properties.

Take the smooth curve $f^{-1}(J(m))$ through $c(m)$,  consider $z(m-1)\in J(m-1)$ and $k(m-1)\ge 0$ so that $f^{k(m-1)}(z(m-1))=c(m)$.  Then $f^{-(k(m-1)+1)}(J(m))$ contains a smooth arc through $z(m-1) \in J(m-1)$ transversal to $\R$ and 
$f^{-(k(m-1)+2)}(J(m))$ contains a smooth arc transversally crossing the curve $f^{-1}(J(m-1))$ through $c(m-1)$.

Next consider $z(m-2)\in J(m-2)$ and 
$k(m-2)\ge 0$ so that $f^{k(m-2)}(z(m-2))=c(m-1)$. Then $f^{-(k(m-1)+k(m-2) + 2)}(J(m))$ contains a smooth arc, which crosses the smooth arc $f^{-(k(m-1)+1)}(J(m-1))$ transversally, which  in turn intersects $J(m-2)$ transversally. 
Hence   
$f^{-(k(m-1)+k(m-2)+3)}(J(m))$ crosses a component of $f^{-(k(m-2)+2)}(J(m-1))$ transversally, which in turn intersects $f^{-1}(J(m-2))\ni c(m-2)$ transversally, see Figure~\ref{fig:tree}. 

Continuing in this way, 
$f^{-(k(m-1) +\dots + k(1)  + m) }(J(m))$ contains a smooth arc, which crosses a smooth arc in $f^{-(k(m-2)+\dots + k(1) + (m-1))}(J(m-1))$
transversally, which in turns crosses $f^{-(k(m-3)+\dots + k(1) +(m-2))}(J(m-2))$, etc, until $f^{-(k(1)+2)}(J(2))$ which intersects $f^{-1}(J(1))\ni c(i)$  transversally. If the component  $L_n(z)$ contains a critical point, then the shortest path  in $L_n(z)$ connecting $z$ to $\R$
is of this form. If $L_n(z)$ does not contain a critical point then it is any path connecting $z$ to $\R$
is  mapped onto such such a path by some map $f^k$.

\subsection{The set $K_X(f)$ can be chosen arbitrarily near the real line and is locally connected}

%Note that the above description considers backward orbits, rather than forward orbits. 
%In particular, we need to  take into account complex backward orbits of the intervals $J_i$. To emphasise this point, 
%consider for example a unimodal map with a non-recurrent critical point. Then the set $J^{-1}$ 
%through $c$ is not contained in any real domain of a complex box mapping
%with range (containing $c$) so small that it does intersect the orbit of the critical value.

To obtain control on the geometry of $K_X$ we have to  develop a construction  which incorporates 
certain complex pullbacks in the complex box construction from 
\cite{CvST}. This will be done in Appendix A, where also the proof of the following 
theorem will be given. 

\begin{thm}\label{thm:KX} 
For each $a>0$ and each $f\in \mathcal A^{\underline \nu}_a$ 
there  exists $\delta>0$ so that if the pruning intervals $J_i$ are disjoint and have length $<\delta$
then 
\begin{enumerate} 
\item $K_X(f)$ is well-defined:  in  the inductive definition (\ref{def:Kn}) we have that 
 $ \mbox{cc}_{K_n'} f^{-{n}}( J^{-1})\subset \Omega_a$. Moreover, 
$K_X(f)\Subset \Omega_a$ and $K_X(f)$ is real symmetric. 
\item  $f^{-1}(K_X(f))\supset K_X(f)$ and $K_X(f)$ is connected;
\item $K_X(f)$ has no interior, is  full (i.e. $\C\setminus K_X(f)$ is connected) and locally connected;
\item   each point of $K_X(f)$ has finitely many accesses;
\item every periodic point of $f\colon K_X(f)\to K_X(f)$ is either repelling or in the real line
(where $f$ also denotes the complex extension of the real analytic interval map). 
\end{enumerate}
Moreover,   let $O$ be a union of periodic attractors and assume that $O$
has a small basin   in the sense of Definition~\ref{def:smallbasin},
then 
\begin{enumerate}
\item[(6)] properties (1)-(5) hold for $K_{X,O}(f)$ and except that in this case
$K_{X,O}(f)$ will have interior; 
\item[(7)] all periodic points of $\hat f_X$ in  $\partial K_{X,O}(f)$ are repelling. 
\end{enumerate}
\end{thm}

As mentioned, the proof of this theorem will be given at the end of Appendix A.

\begin{rem} Note that periodic points are in general not dense in $K_X(f)$.  For example,
consider a unimodal map with a maximum at $c$. Then points to the right of $J\ni f(c)$
do not contain pre-images  in $K_X(f)$, and therefore any component of $K_X(f)\setminus I$ which is 
{\lq}to the right{\rq} of $J$ does not contain periodic points.  In the next subsection we will associate to 
$f$ and $K_X(f)$ a circle map $\f_X$. This circle map will have discontinuities 
and periodic points of this circle map are again not dense  in the circle. 
\end{rem}

\begin{rem} 
We should emphasise that $f\in\mathcal A^{\underline \nu}$ has at most a finite number of attracting  or parabolic  periodic points, see \cite{MMvS,dMvS}. 
%Hence one can always ensure that 
%properties (5)-(6) hold when taking $\delta>0$ sufficiently small. Property (7) requires that 
%these properties (5)-(6) also hold for nearby maps, making this no longer a trivial statement. 
\end{rem}

%\textcolor{red}{\begin{rem}  One  reason for including  in $X$ critical points which are in the basin of 
%periodic attractors in $X$ is that $K_X(f)$ could contain 
%a neighbourhood of such a critical point which is not compactly contained in $\Omega_a$. 
%\end{rem}} 

\section{The associated external map: a circle mapping with discontinuities}\label{subsec:ppl}
The external mapping associated to a polynomial-like mapping $f$ is an analytic
expanding map of the circle $\f:\partial\mathbb D\to \partial\mathbb D$
that has an extension to a holomorphic mapping $\f:A\to A'$ where $A$ and $A'$ are annuli so that $\partial\mathbb D\Subset A\Subset A'$. In this section, we will construct an external mapping which is associated to a real analytic map $f\colon I\to I$ together with its pruned Julia set $K_X$, and use it to prove 
Theorem~\ref{thm: pruned-pol-like-map}. We will go on to discuss this class of external mappings in general.

To do this we will  make the following assumption from now on: 
\begin{equation*}
\left\{ \begin{array}{cl}
&\mbox{ the pruning intervals $J$ and the corresponding pruned Julia set }\\
&\mbox{$K_X$ satisfy the conclusions of Theorem~\ref{thm:KX}.}\end{array}
\right. 
\end{equation*} 
%\subsubsection{A pruned external map associated to $f|_{K_X(f)}$}

Let $\psi_X\colon \overline \C \setminus \overline{ \mathbb D} \to \overline \C
\setminus K_X$ be the uniformising map that fixes $\infty$ and with real
derivative at $\infty,$ and write $\phi_X=\psi_X^{-1}\colon \overline \C \setminus
K_X\to \overline \C \setminus \overline{ \mathbb D}$. Set 
\begin{equation} \f_X=\phi_X\circ f\circ\psi_X, \end{equation} 
defined on $\psi_X^{-1}(\Omega_a\setminus K_X)\subset \overline \C\setminus 
\overline{ \mathbb D}$ which is an open annulus with inner boundary $\partial \mathbb D$. Since $K_X$ is locally connected and each point  of $K_X(f)$ has finitely many accesses, by
Carath\'eodory's Theorem $\psi_X$ extends continuously to the boundary of $\mathbb
D$ and $\phi_X$ extends as a multi-valued function to the boundary of $K_X.$
So we can write
$$\psi_X\colon \partial \D \to K_X \mbox{ and } \phi_X\colon K_X\to \partial \D$$
where the latter map is multi-valued. 
Note that the normalization
of $\psi_X$, the assumption that $I=[-1,1]$ and that $K_X$ is real-symmetric, 
imply that $\phi_X(-1)=-1$ and $\phi_X(1)=1.$  
Since there are finitely many accesses at each point in $K_X$, 
\begin{equation} \hat X= \phi_X(X)  \end{equation}  
is a finite subset
of $\partial \mathbb D$.
%, where as before $\Crit'(f)$ is the subset of critical points of $f$ which are periodic. 
Note that $\pm 1 \notin X$ and therefore 
$$\pm 1 \notin \hat X.$$

\begin{figure}
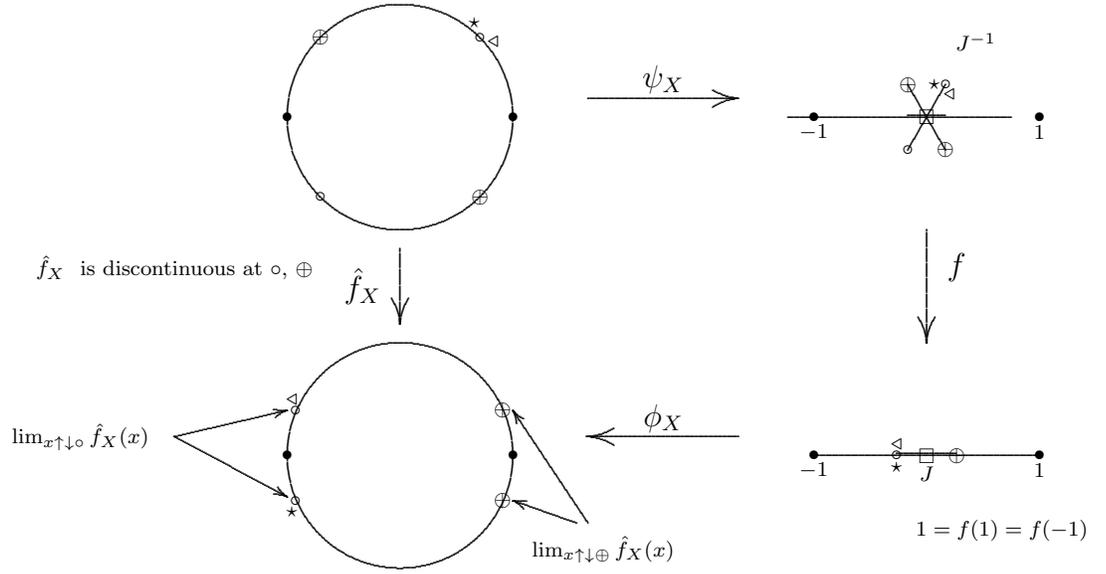
 \hfil \beginpicture \dimen0=0.5cm 
\setcoordinatesystem units <\dimen0,\dimen0> point at 19 0 
\setplotarea x from -12 to 5, y from -12 to 4 
\setlinear \plot -3 0 3 0 / 
%\setdots <0.5mm> 
%\plot -1 0.1 -12 0.1 / 
%\plot -5 0.05 -7 0.05 / 
\put {\tiny $J^{-1}$} at 2 2 / 
\plot 0 0 0.5  0.866 / 
\plot 0 0 0.5  -0.866 /
\plot 0 0.05 -0.5 0.05 / %
 \put {\tiny $\circ$} at 0.5  0.866 
%   \put {\tiny $\circ$} at 0.5  -0.866 
  \put {\tiny $\oplus$} at 0.5  -0.866 
%  \put {\tiny $\star$} at 0.6 0.6 
    \put {\tiny $\triangleleft$} at 0.6 0.6 
    \put {\tiny $\star$} at 0.2 0.866 
   \put {\tiny $\oplus$} at -0.5  0.866 
%   \put {$\$} at -0.5  0.866    
  \put {\tiny $\circ$} at -0.5  -0.866 
   \put {\tiny $\bullet$} at 3 0 
      \put {\tiny $\bullet$} at -3 0 
      \put{\tiny $1$} at 3 -0.4
      \put {\tiny $-1$} at -3 -0.4   
%\setdashes  <0.5mm> 
\plot 0 0 -0.5  0.866 / 
\plot 0 0 -0.5  -0.866 / 
\plot 0 0.05 0.5 0.05 / 
 \circulararc 360 degrees from -17 0  center at -14  0
 \put {\tiny $\bullet$} at -11 0 
  \put {\tiny $\bullet$} at -17 0 
     \put {\tiny $\star$} at  -12.03  2.51
          \put {\tiny $\triangleleft$} at  -11.53  2.01
  %  \put {\tiny $\star$} at  -11.53  2.01
    \put {\tiny $\circ$} at  -11.8786  2.121
        \put {\tiny $\oplus$} at  -11.8786  -2.121
%                \put {$\circ$} at  -11.8786  -2.121
            \put {\tiny $\oplus$} at   -16.1213  2.121
%                     \put {$\circ$} at   -16.1213  2.121
        \put {\tiny $\circ$} at   -16.1213  -2.121
         \arrow <10pt> [.2,.67] from -9  0.5 to -5 0.5 
         \put {$\psi_X$} at -7 1 
             \arrow <10pt> [.2,.67] from -14  -3.5 to -14 -5.5 
                     \arrow <10pt> [.2,.67] from 0  -3 to 0 -6 
                     \put {$f$} at 0.8 -4 
          \circulararc 360 degrees from -17 -9  center at -14  -9
%        \put {\tiny $\star$} at  -16.88 -7.5 
                \put {\tiny $\triangleleft$} at  -16.88 -7.5 
         \put {\tiny $\star$} at  -16.88 -10.5 
            \put {\tiny $\circ$} at  -16.7786  -7.8        
           \put {\tiny $\circ$} at  -16.7786  -10.2
%                    \put {$\oplus$} at  -11.2786  -7.8
            \put {\tiny $\oplus$} at   -11.2786  -7.8
%                             \put {$\circ$} at  -11.2786  -10.2
            \put {\tiny $\oplus$} at   -11.2786  -10.2

         \put {$\f_X$} at -15 -4.5      
                 \arrow <10pt> [.2,.67] from -5  -8.5 to -9 -8.5 
                         \put {$\phi_X$} at -7 -8 
  \plot -3 -9 3 -9 /  
     \put {\tiny $\bullet$} at 3 -9  %%%
      \put {\tiny $\bullet$} at -3 -9 %%%
            \put{\tiny $1$} at 3 -9.4
            \put {\tiny $1=f(1)=f(-1)$} at 2 -11
                      \put {\tiny $J$} at 0 -9.5
      \put {\tiny $-1$} at -3 -9.4 
  \plot -0.8 -8.95 0.8  -8.95 / 
%      \put {\tiny $\star$} at  -0.8  -8.7       
       \put {\tiny $\triangleleft$} at  -0.8  -8.7   
      \put {\tiny $\star$} at  -0.8  -9.3
          \put {\tiny $\circ$} at  -0.8  -9 
             \put {\tiny $\oplus$} at  0.8  -9 
%                    \put {$\circ$} at  0.8  -9 
              \put {\tiny $\bullet$} at -11 -9 
  \put {\tiny $\bullet$} at -17 -9 
  \put {\tiny $\square$} at 0 0 
   \put {\tiny $\square$} at 0 -9 
  \put {\tiny $\hat f_X$ \mbox{ is discontinuous at $\circ$, $\oplus$}} at -20 -4 
  \put {\tiny $\lim_{x\uparrow\downarrow \circ} \hat f_X(x)$} at -22.5 -8.5   
%  \put {$\oplus$} at -16.55 -4.09   
                   \arrow <5pt> [.2,.67] from -20  -8.5 to -17 -7.8 
 \arrow <5pt> [.2,.67] from -20  -8.5 to -17 -10.1 
   \put {\tiny $\lim_{x\uparrow \downarrow \oplus} \hat f_X(x)$} at -8.6 -11.5   
%  \put {\tiny $\circ$} at -8.64 -10.9   
                   \arrow <5pt> [.2,.67] from -9  -10.8 to -11 -7.8 
 \arrow <5pt> [.2,.67] from -9.3  -10.8 to -11 -10.2 
\endpicture
\caption{The pruned Julia sets $K_X$ associated to the map $f$ defined by $z\mapsto z^3$ 
when $J=(-\epsilon,\epsilon)$. This diagram shows how the
 external map $\f_X=\phi_X\circ f\circ\psi_X\colon  A\to \C\setminus \overline \D$ is constructed,  where 
$A\subset \C \setminus \overline \D$ is an annulus with inner boundary $\partial \mathbb D$. 
Points in the set $X=f^{-1}(\partial J)\setminus \R$ are marked with the symbol $\circ$
in the top right. The corresponding points $\hat X=\phi_X(X)\subset \partial \D$ are also marked with $\circ$ 
and bold $\circ$  in the figures on the left. Two points near a point in $X$ are mapped by $f$ to different sides
of $I$, and consequently the map $\f_X$ has discontinuities at $\hat X$. 
For the same reason $\f_X$ also  has discontinuities at all six points $\phi_X(0)\subset \partial \D$. 
} 
\end{figure}

\begin{example}\label{example:hatf} Let us describe the external maps corresponding to 
 the examples from above: \\
\t(i) $f$ has a super-attractor. Let $I=[-1,1]$, $J=(-\epsilon,\epsilon)$ with $\epsilon<1$ and $f(x)=x^3$ then
as we saw in Example~\ref{example:K_X} 
$$K_X= \bigcup_{k=1,\dots,\ell-1} e^{\pi i k/\ell} \cdot  J'   \,\, \bigcup \,\, I 
\mbox{ where }J'=[-\epsilon^{1/\ell},\epsilon^{1/\ell}].$$ 
 Since the pruned filled Julia is not backward invariant, the associated external map $\f_X$ is not
 going to extend continuously to $\partial \D$.
  In the next lemma we will show that nevertheless  $\f_X$ extends to a piecewise real analytic map 
 map of $\partial \D$ with discontinuities at pre-images under $\psi_X$ of $X=f^{-1}(\partial J)\setminus \R$ and of the critical point.
\\
(ii)  Assume that $I=[-1,1]$, $J_n=(-1,1)$ and $f\colon I\to I$ is a polynomial so that all its critical points are real 
and  non-periodic. 
Then $K_n=\{z\in \C; f^n(z)\in I\}$, $\forall n\ge 0$ and therefore $K_X$ is 
equal to the filled Julia set of $f$.  On the other hand,
if $J\Subset (-1,1)$ then  $K_n$ is a strict subset 
of the Julia set of $f$ pruned  at preimages of $X=\partial J^{-1}$. 
%$K_n$ again  consists of a finite union of smooth curves.  
\end{example} 

\begin{figure}
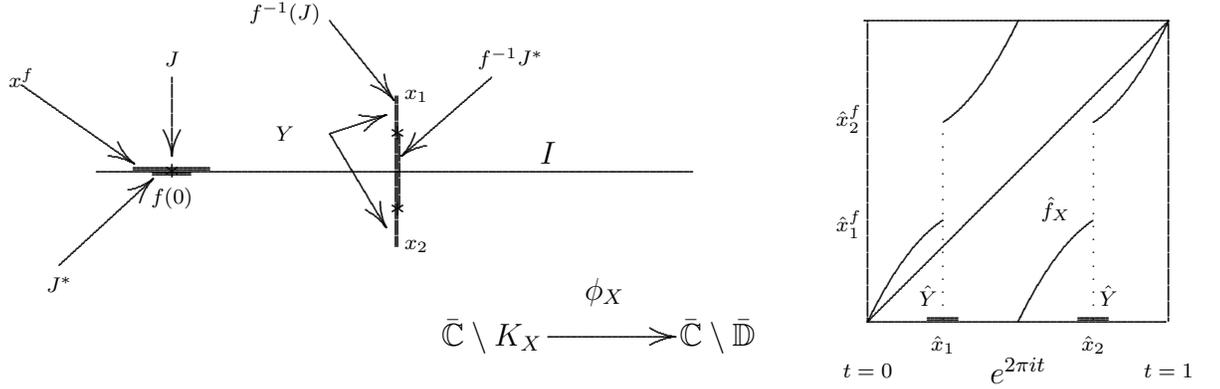
 \hfil \beginpicture \dimen0=0.5cm 
\setcoordinatesystem units <\dimen0,\dimen0> point at 19 0 
\setplotarea x from -4 to 4, y from -6 to 6 \setlinear 
\plot -8 0 8 0 / \plot -5 0.1 -7 0.1 / \plot -5 0.05 -7 0.05 / 
\plot 0 2 0 -2 /
\plot -0.05 2 -0.05 -2 /
\plot 0.05 1 0.05 -1 /
 \plot -5.5 -0.05 -6.5 -0.05 / 
\plot -5.5 -0.08 -6.5 -0.08 / 
\put {$*$} at 0 1 
\put {$*$} at 0 -1
\setdots <0.7mm> %\plot 0 -3 0 3 / 
\put {\tiny $x_1$} at 0.5 2 %\put {$X$} at -0.7 1.75 
%\put {\tiny $Y$} at -1 1.5 
\put {\tiny $x_2$} at 0.5 -2 
%\put {\tiny $*$} at 0 1.3 
%\put {\tiny $*$} at 0 -1.3 
\put {\tiny $f(0)$} at -6 -0.7 
\put {\tiny $J$} at -6 3  \put {$*$} at -6 0 
\setsolid \arrow <10pt> [.2,.67] from -6 2.5 to -6 0.4 
\put {\tiny $J^*$} at -9 -3  
\setsolid \arrow <10pt> [.2,.67] from -9 -2.5 to -6.5 -0.2 
\put {\tiny $x^f$} at -10 2.5  
\setsolid \arrow <10pt> [.2,.67] from -10 2.3 to -7 0.2 
\arrow <10pt> [.2,.67] from 2.5 2.5 to 0.1 0.4 \put {\tiny $f^{-1}J^*$} at 3 3
\arrow <10pt> [.2,.67] from -1.8 1 to -0.3 1.5 \put {\tiny $Y$} at -3 1
\arrow <10pt> [.2,.67] from -1.8 4 to -0.1 1.9 \put {\tiny $f^{-1}(J)$} at -3 4.2
\arrow <10pt> [.2,.67] from -1.8 1 to -0.3 -1.5
\put {$I$} at 4 0.5 
  %     \setquadratic
  \setcoordinatesystem units <\dimen0,\dimen0> point at 2.5 0 
  \setplotarea x from -4 to 4, y from -6 to 6 \setlinear \plot -4 -4 4 -4 4 4 -4 4 -4 -4 / \plot -4 -4 4 4 / 
  \put {$e^{2\pi i t}$} at 0 -5.3
  \setquadratic
%   \plot -4 -4 -2 -0.3 0 4 /
%\plot 0 -4 1 -2.6 2 -1.3 
 \plot 2 1.3 3 2.3 4 4 / \plot 0 4 -1 2.3 -2 1.3 / \plot -2 -1.3 -3 -2.3 -4 -4 / 
\plot 0 -4 1 -2.3  2 -1.3 / 
\setdots \setlinear \plot 2 1 2 -4 / 
\put {\tiny $\f_X$} at 1 -1 
\put {\tiny $t=0$} at -4 -5.3 
\put {\tiny $t=1$} at 4 -5.3 
\put {\tiny $\hat x_1^f$} at -4.5 -1.4 
\put {\tiny $\hat x_2^f$} at -4.5 1.4 
\put {\tiny $\hat x_2$} at 2 -4.6
\setsolid \plot 1.6 -3.95 2.4 -3.95 / \plot 1.6 -3.9 2.4 -3.9 / \put {\tiny $\hat Y$} at 2.4 -3.3
\put {\tiny $\hat Y$} at -2.4 -3.3
 \setdots \setlinear \plot -2 1 -2 -4 / \put {\tiny $\hat x_1$} at -2 -4.6 \setsolid \plot -1.6 -3.95 -2.4 -3.95 / \plot -1.6 -3.9 -2.4 -3.9 /
 \put {${\bar \C} \setminus  K_X$} at -14 -4.4
\arrow <10pt> [.2,.67] from -12.5 -4.4 to -9.1 -4.4 \put {${\bar \C} \setminus {\bar {\mathbb D}}$} at -8 -4.4
 \put {$\phi_X$} at -11 -3.2
%   \put {$\hat Y$} at 2.4 -3.3
%%      \put {$a_k'$} at -5.3 -4.5 \put {$b_{k+1}$} at 1.5 -4.5 \put {\tiny $\bullet$} at 4.5 -4 \put {\tiny $\bullet$} at 4 -4 \put {\tiny $\bullet$} at -4 -4 \put {\tiny $\bullet$} at -4.5 -4 \put {\small $0$} at -2.2 -4.5 \put {\small $d_k$} at -1.3 -4.5 \put {\small $e_k$} at -0 -4.5 \put {\small $b_k$} at -5 4 \put {\small $B_k$} at -5 8.6 \put {\small $A_k$} at -5 -6
%        \put {\small $0$} at -2 -4.5 \setlinear \plot -4.5 -5.8 -4 -4 -1.5 5 -0.5 8.6
 %           / % 2.5 9 i.e. 1 3.6
%             \setquadratic
%%                \plot -2 3 0 1.5 4 -4 /
%                  \plot -2 3 1.3 1.3 4 -4 / \setlinear \plot 4 -4 4.5 -6 / \setdots <0.8mm> \setlinear
%%                     \plot -0.5 -4 -0.5 8.6 / \plot -1.7 -4 -1.7 4 / \plot -4 8.6 -0.2 8.6 /
%                       \plot -2.5 -2.5 -2.5 -5 / \plot 1.3 -2.5 1.3 -5 / \put {$a_1$} at -2.5 -5.5 \put {$b_1$} at 1.7 -5.5 \plot -2 -4 -2 3 /
%%                          \plot -4 -2 -2 -2 /
\endpicture
\caption{\label{fig2} On the left, the set $K_1$ is drawn in the case of a quadratic map $f(x)=(1-a)x^2+a$ 
without periodic attractors for $f(c)=a\in (-1,0)$. 
  In this case, $J$ is an interval containing $f(c)$,  $X=\{x_1,x_2\}=f^{-1}(\partial J)\setminus \R$ and $x^f=f(x_1)=f(x_2)\in \partial J$. 
  %   The set $K_X(f)$ is constructed by taking successive
%     preimages of $f^{-1}(J)$.
On the right, the mapping $\f_X \colon \partial \D \to \partial \D$ is drawn  - it is discontinuous at $\hat X=\{\hat x_1,\hat x_2\}$.
Note that $x^f$ corresponds to two distinct points $\hat x^f_1, \hat x^{f}_2$ in $\D$ and that these
are the limits of  $\f_X(x)$ as $x$ tends from the left or right to $\hat x_1$ or to $\hat x_2$. 
We also draw the curves $f^{-1}(J^*)\subset f^{-1}(J)$ (the non-real preimages of $J^*\subset J$) 
where $J^*$ is a curve in $I$ which compactly contains $J$. The set $X^*$ consists of the two boundary points of $f^{-1}(J^*)\setminus \R$ and is marked with the symbol $*$.  The set $Y=K_X\setminus X^*$ contains the two intervals
which are drawn.  The two intervals $\hat Y=\phi_X(Y)$ are also marked. 
%$\hat Y=\phi_X(Y)$ where $Y\subset K_X$ is a set containing $X$, corresponding to a further pruning at points $X^*$. 
%$\lim_{x\uparrow \hat x} \f_X(x)= \hat x_1^f$ and  $\lim_{x\downarrow \hat x} \f_X(x)= \hat x_2^f$.
}
\end{figure}

\subsection{A semi-conjugacy with an expanding circle map} 

\begin{defn} \label{def:signf} 
 Define the {\em sign} $\epsilon(f)\in \{-1,1\}$ 
of $f\colon I\to I$ by $$\epsilon(f)=\begin{cases}  1 &\mbox{ if }f(1)=1 \\
  -1 &\mbox{ if }f(1)=-1. \end{cases} $$ 
\end{defn} 

In this section we will show that associated to $\hat f_X$ is is a map 
in the class of maps $\mathcal E^d_\epsilon$  defined  below. 

\begin{defn}[$\mathcal E^d_\epsilon$] \label{def:Ed}   Assume that $\epsilon\in \{-1,1\}$. 
We say that $g \in \mathcal E^d_\epsilon$ if 
\begin{enumerate}
\item there exists a subset $\hat X\subset \partial \D$ of finite  cardinality so that $g$ is a real-symmetric analytic map  $g\colon \partial \D\setminus \hat X \to \partial \D$, so that $g$ has discontinuities 
at each point of $\hat X$ and  so that $\pm 1 \notin \hat X$; 
\item $g(\bar z)=\overline{g(z)}$ for all $z\in \partial \D\setminus \hat X$; 
\item $g(1)=\epsilon$. 
\item there exists a  map $g^*\colon \R\setminus \pi^{-1}(\hat X) \to \R$ so that  $\pi \circ  g^* (x) = g \circ \pi(x)$ for all $x\in \R \setminus \pi^{-1}(\hat X) $
 where  $\pi\colon \R\to \partial \D$ is the covering map $t\mapsto e^{2\pi  i t}$ so that
$g^*\colon \R\setminus \pi^{-1}(\hat X) \to \R$ is strictly monotone increasing on each component 
of $\pi^{-1}(\hat X)$  and so that  the jumps at discontinuities of $g^*$ are of size $<1$;
\item $g^*(x+1)-g^*(x)=d$ for all $x\in \R \setminus \pi^{-1}(\hat X)$;  
\end{enumerate} 
 \end{defn}

\begin{lemma}\label{lem:fXproperties}
  Assume that $f\in  \mathcal A^{\underline \nu}_a$ has only hyperbolic periodic points and assume that 
  the intervals $J_i\ni f(c_i),$ $c_i\in\Crit(f)$, are disjoint and are compactly contained in $I$. 
  Then   $\f_X$ extends to a map on $\partial \D$ (minus a finite number of points):
  \begin{enumerate}
  \item $\f_X \colon \partial \mathbb D \setminus \hat X\to \partial \mathbb
    D$ which is real analytic, order preserving and has no critical point;
  \item $\f_X$ has a jump discontinuity at {\em each} point $\hat X$ and $\phi_X(\Crit(f))\cap \hat X=\emptyset$;
  \item each periodic point of $\f_X\colon \partial \mathbb D \setminus
    \hat X\to  \partial \mathbb D$ is either hyperbolically repelling or corresponds to a (real) periodic attractor or parabolic point of $f\colon I\to I$. 
      \end{enumerate}
      Moreover, $\f_X$ is in the class $ \mathcal E^d_\epsilon$ defined in the definition above,  where $d=\ell_1+\dots+\ell_\nu$, where
      $\ell_i$ is the $i$-th critical point of $f$,  and where $\epsilon=\epsilon(f)$ (defined above). 

If $f$ has no periodic attractors then each periodic orbit of $\f_{X}\colon \partial \mathbb D \setminus     \hat X\to  \partial \mathbb D$ 
is repelling. Moreover, if $O$ is the set of all (real) periodic attractors of $f$ and assume that $O$ is 
 as in the 2nd part of Theorem~\ref{thm:KX}, 
  then  each periodic point of 
$\f_{X,O}\colon \partial \mathbb D \setminus     \hat X\to  \partial \mathbb D$ is repelling. 
    Here  $\hat f_{X,O}$ is the  external map corresponding to $K_{X,O}$. 
\end{lemma} 
\proofof{Lemma~\ref{lem:fXproperties}}
 Let us for simplicity  assume in the proof below (except at the end of the proof) 
 that no critical point of $f$ is eventually 
mapped to itself or to another critical point. 

To prove the first statement of the lemma  let us for simplicity write $K=K_X$ and show 
that $\f_X$ extends to a real analytic map outside $\hat X=\psi_X(X)$ and has discontinuities
at each point in $\hat X$. Indeed,  for each $x\in K\setminus X$, there exists a neighbourhood $W\subset\C$ of $x$ so that $f$ maps each component $W'$ of $W\setminus K$ homeomorphically  onto some component $W''$ of 
$f(W)\setminus K$.  
It follows that $\f_X$ maps $\hat W'=\psi^{-1}_X(W')$ homeomorphically onto $\hat W''=\psi^{-1}_X(W'')$.
Taking the reflection of the sets $\hat W'$ and $\hat W''$  in $\partial \D$,
and using Schwarz reflection principle, it follows that 
$\f_X$ is a well-defined analytic map in a
neighbourhood of each point of $\partial \mathbb D\setminus X,$ and
$\f_X$ maps $ \partial \mathbb D\setminus X$  into $\partial \mathbb D$. 
Another way of seeing the previous sentence is by considering curves in $\C\setminus \overline \D$
which land at a point in $\overline \D$. For example, consider a critical 
point $c_i\in K$ which is not eventually mapped into a critical point. 
Then, locally,  each of the $2\ell_i$ sectors near $c_i$ of $\C\setminus K_X$  is mapped by $f$ to a half-sector near $f(c_i)$. Hence the above argument shows that $\f_X$ extends continuously to $\phi_X(c_i)$, and 
therefore as a real analytic map near $\phi_X(c_i)$. 
If $c_i$ is mapped by $f$ to another non-periodic critical point, then by construction the number of sectors of 
$\C\setminus K_X$ at 
$c_i$ will be $\ell_i$ times the number of sectors of $\C\setminus K_X$ at $f(c_i)=c_j$ and so the previous argument still goes through. In particular, even in this case,  $\f_X$ is real analytic at $\phi_X(c_i)$.
% (here we use that $c_i$ is not eventually mapped into a critical point).
%From this Property (4) of Definition (\ref{def:Ed}) follows.}} 

To prove the second statement of the lemma, we need to analyse what happens at 
points $\hat X$ and show they are discontinuities of $\f_X$. This is described in the following  
 {\bf claim:} For each point $x\in X$ there exists a unique point $\hat x$ (which is in $\hat X$) so that $\phi_X(\hat x)=x$
but there are (precisely) two points $\hat x_1^f,\hat x_2^f$ so that $\phi_X(\hat x_1^f)=\phi_X(\hat x_2^f)=x^f=f(x)$.

To prove this claim, note that for {\em each} small  neighbourhood $W$
of $x$, we have that $W\setminus K$ consists of a single component, but
the set $f(W)\setminus K$  consists of two components, one in the upper
half-plane and one in the lower half plane, see Figure~\ref{fig2}. This 
proves the claim.  

It follows that  $\f_X\colon \partial \mathbb D \setminus \hat X\to  \partial \mathbb D$ has a discontinuity 
at each point $\hat x \in \hat X$ and its left and right limit values correspond to 
$\hat x_1^f$ and $\hat x_2^f$. Since $\phi_X(\hat x_1^f)=\phi_X(\hat x_1^f)$ it follows that 
$\hat x_1^f$ and $\hat x_2^f$ are each others complex conjugates.
Again we can also see that $\f_X$ has a discontinuity at such a point by considering two curves $\gamma_1,\gamma_2$ in $\C\setminus K_X$ landing  on $x_1$ from the right respectively from the left of the curve $f^{-1}J'$ as in Figure~\ref{fig2} . Then $f(\gamma_1)$ lands on $f(x_1)$ from above
and $f(\gamma_2)$ lands on $f(x_2)$ from below. The same argument also 
shows that one has discontinuities at $\phi_X(c_i)$ when $c_i$ is a periodic 
critical point of $f$. 
%If $f$ has periodic critical point $c_i$ whose orbit does not
%contain other critical points, then as $f(c_i)$ is a pre-image of $c_i$, by  definition 
%the set $K_X$ near $c_i$ and near $f(c_i)$ are homeomorphic.  Since $f$ has local degree $\ell_i>1$, it follows  
%that $\f_X$ has also discontinuities at each point in $\phi_X^{-1}(c_i)$,  see also example~\ref{example:hatf}.}

%Because of
%this, $\f_X$ is discontinuous at each $\hat x\in \psi_X^{-1}(X)$. Indeed,
%arbitrarily close to $x\in X$ there are points $z_1,z_2\in \mathbb C\setminus
%K,$ so that $f(z_1)$ and $f(z_2)$ are in the upper half-plane and the lower half-plane,
%respectively.   Thus, $\f_X$ has a jump discontinuity at $\hat x$. 

To prove the third  and fourth statement of the lemma, 
observe that each periodic point
$\hat p$ of period $s$ for $\f_X$ corresponds to a periodic point $p$ of period
$s$ for $f\colon \Omega_a\to \C$. By part (4) of Theorem~\ref{thm:KX},   $p$ is 
repelling unless it lies on the real line.
Let us first assume that $p$ is  a repelling periodic 
 point of $f\colon \Omega_a\to \C$.  Then we can 
choose a neighbourhood $W$ of $p$, so that $f^s(W)=W'\Supset W.$
We have that $\psi_X^{-1}(W)$ is a union of disjoint topological disks each with the
property that its boundary can be decomposed into two components one in
$\phi_X(\partial W)$ and one in $\partial\mathbb D.$ Reflect each of these disks
about $\partial \mathbb D$ to obtain a union of topological disks $W_\phi$.
One of these components contains $\hat p$; label it by $\hat W.$ Repeat this
procedure with $W'$ and to obtain a topological disk $\hat W'\Supset \hat W\owns
\hat p,$ so that $\f_X^s(\hat W)=\hat W'.$ Thus we have that $\hat p$ is a
repelling periodic point of $\f_X$.
If $p$ is an attracting periodic point of $f$, then we can find a neighbourhood
$W$ of $p$ so that $f^s(W)=W'\Subset W$ and repeat a similar argument. The parabolic case
also goes similarly. 
If $O$ is as in  (6) of Theorem~\ref{thm:KX}, then all periodic points of $f$ in $\partial K_{X,O}$
are repelling, and so  $\hat f_{X,O}$ has only repelling periodic points. 
 
 Now it is easy to see that the properties from  Definition (\ref{def:Ed}) hold. 
 For example, take a point $p\in f^{-1}(J_i)\setminus \R$. There are two curves landing at $p$ approaching 
 $f^{-1}(J_i)\setminus \R$ from opposite sides.  The image under a complex extension of $f$ of one of these curves  lies in the upper half plane and the other one  in the lower half plane, but land 
 at the same point. The number of curves $f^{-1}(J)\setminus \R$ is determined by $\ell_1+\dots+\ell_\nu$. From all this we obtain that $g$ has degree $d=\ell_1+\dots+\ell_\nu$
 as claimed in property (5) of Definition (\ref{def:Ed}).
\hfill\qed

 \bigskip 

Later on it will be useful to apply a further pruning in a combinatorial 
well-defined way. For this we will use the following 

\begin{lemma}[Semi-conjugacy with expanding covering map of $\partial \D$] \label{lem:comb-semiconj}
Let  $g\in \mathcal{E}^d_\epsilon$ with $\epsilon\in \{-1,1\}$ and let
$Q_g$ be a finite forward invariant set disjoint from 
the set $\hat X$ of discontinuities of $g$.
 Then 
\begin{enumerate}
\item there exists an orientation preserving degree $d$ continuous covering map $g_*\colon \partial \D \to \partial \D$ which agrees with 
$g$ outside a neighbourhood $U$ of $\hat X$ with $U\cap Q_f=\emptyset$
and so that the image of $g_*$ of each component of $U$ has length $<1$;
\item there exists a unique semi-conjugacy $h_g$ of $g_*$  with  the map $ \partial \D \to \partial \D$
defined by  
$$\begin{cases} z\mapsto z^d &\mbox{ if }\epsilon=+1  \\
z\mapsto -z^d &\mbox{ if }\epsilon=-1.\end{cases}$$ 
Moreover,   $h_g(1)=1$ and 
$h_g$ is real symmetric, i.e.,  
$h_g(\overline z)=\overline{h_g(z)}$; 
\item $Q=h_g(Q_g)$ only depends on $g$ and $Q_g$ and not on the extension $g_*$ of $g$ (and so 
not on $h_g$).
\end{enumerate} 
\end{lemma} 
\begin{pf} That $g_*$ exists immediately follows from the properties of $\mathcal{E}^d$. Parts (2) and (3) immediately follow from \cite{Sh}. 
\end{pf}

Using this lemma, we can give the following definition. 

\begin{defn}  $\mathcal E_{\pm,Q}^d$. Take $\epsilon\in \{-1,1\}$ and assume that 
$Q\subset \partial \D$ is a finite forward invariant set under $z\mapsto z^d$ if $\epsilon=1$ 
and invariant under $z\mapsto -z^d$ if $\epsilon=-1$.
Then   $\mathcal E_{\pm,Q}^d$ is defined  
 to be the class of maps $g\in \mathcal E^d$ so that $g$ has an invariant subset $Q_g$ which does not intersect 
the set of discontinuities of $g$, so that $h_g(Q_g)=Q$.
% and so that between any two  points of discontinuities of $g$ lie two points of $Q_g$.  
\end{defn} 

Sometimes we tacitly use $\mathcal E^d$ for the space of circle maps
which do actually do arise from interval maps. 

\subsection{Further pruning and the forward invariant set $\Lambda_N$} % and combinatorial pruning data}
\label{subsec:furtherpruning} 
Since the analytic behaviour of the map $\f_X$ near its discontinuities is difficult to describe, we will take 
for each critical value, $f(c_i),$ a closed, nice interval $J_i^*\Subset J_i$, containing $f(c_i),$ so that its boundary points are periodic, preperiodic or in the basin of a periodic attractor.
Here $J_i$ are the pruning intervals used in the construction of the pruned Julia sets.  
Such intervals $J_i^*$ exist  by the real bounds, see \cite[Theorem A$'$]{vSV}. Let $J^*=\cup_i J_i^*$
and $$X^*=f^{-1}(\partial J^*)\setminus \R$$ and let 
$Y$  be the union of the components of 
$K_X\setminus X^*$ which do not intersect $I$.
Associated to each critical point
$c_i$  there are $2(\ell_i-1)$ such components.
Note that 
$$\hat Y=\phi_X(Y),$$ 
consists of open intervals around the discontinuity points of $\f_X$, 
i.e.,  around  the finite set $\hat X=\phi_X(X)\subset \partial \mathbb D$.

Let $B_{0,i}\subset \C$ be the immediate basins of the (real) periodic attractors of  
$f\colon \Omega_a \to \C$, let $B_0=\cup B_{0,i}$ and let 
$B_{0,i}^*$ be the connected components of $K_X\setminus \partial B_0$ containing points of $B_{0,i}$.
Then define 
$$\hat B_{0,i}=\phi_X(B_{0,i}^*)$$ and let $B,B^*,\hat B_0$ be the union of these sets. 
Note that $f$ has at most a finite number of real periodic attractors, see \cite{MMvS}. 

Next we define the $\f_X$-forward invariant set 
\label{defLambdaN} 
$$\Lambda_N:=\{z\in \partial \mathbb D; \f_X^n(z)\notin 
\hat Y \mbox{ for all }0\le n\le N\} \mbox{ and } \Lambda_\infty= \cap_{N\geq 0} \Lambda_N.\label{lambda}  $$

\begin{lemma}\label{lemQ}  The map $\phi_X$, $\hat f_X$ and $\hat Y$ have the following properties. 
\begin{enumerate}
\item $\phi_X(I) \subset \Lambda_\infty$ and therefore $\Lambda_\infty$ corresponds to the real points; 
\item $\partial \hat Y= \phi_X(X^*)$, $\partial \hat Y\subset \Lambda_\infty$ and $\partial \hat Y$
is contained in a finite forward invariant 
set which avoids the set of discontinuities of $\f_X$;
\item $\partial \hat Y\subset \partial \D$ consists
of periodic or eventually periodic points or is in the basin of a periodic attractor
of $\f_X$; 
\item $\hat B_0\subset \partial \D$ consists of a finite number of intervals, each of which
containing  a point of  $\phi_X(a)$ for some 
attracting periodic point $a$ of $f$. Moreover,  $\partial \hat B_0$ is a finite forward invariant set. 
\end{enumerate} 
\end{lemma}

\begin{rem}\label{rem:smallattractorsB}
If each component of the basin of a (real) periodic attractor
of $f \colon \Omega_a \to \C$ has small diameter,  as in  (6) of Theorem~\ref{thm:KX}, then 
by Lemma~\ref{lem:fXproperties} all periodic points of $f$ in $\partial K_{X,O}$
are hyperbolically repelling, and so  $\hat f_{X,O}$ has only repelling periodic points. In that case, 
it will not be necessary to use the set $\hat B_0$ we just defined in the next section
(provided we use $\hat f_{X,O}$ instead of $\f_X$). 
\end{rem} 

\begin{pf} By construction the closure of $Y$ is disjoint from $I$. In particular, 
$f$-forward iterates of $x\in I$ never intersect the closure of $Y$. It follows that $\f_X$-forward iterates points
of $z\in \phi_X(I)$ do not intersect $\hat Y$, proving assertion (1). $\partial \hat Y= \phi_X(X^*)$ is obvious. 
Since $f(\hat X^*)\subset I$ it follows that $\partial \hat Y\subset \Lambda_\infty$. 
For each periodic point $p$, $K_X\setminus \{p\}$ consists of at most finitely many components. 
So at most finitely many accesses (or {\lq}rays{\rq} land) at each periodic point of $f\colon K_X\to K_X$. Hence
periodic (preperiodic) points of $f\colon K_X\to K_X$ correspond to periodic (resp. preperiodic) points of $\f_X \colon \partial \mathbb D\setminus \hat X\to \partial \mathbb D$ and in particular $\partial \hat Y$ consists of periodic or preperiodic points of $\f_X$. This implies the final part of (2). 
That $\partial \hat Y$ consists of periodic or eventually periodic orbit, follows from assumption that 
$\partial J^*$ consists of periodic, or eventually periodic points (or the basin of a periodic attractor), proving 
(3).  Let us prove (4). Let $K_1$ be a connected component of $K_X\setminus \partial B_0$ that intersects $B_0.$ Then $\phi_X(K_1)$ consists either of
  \begin{enumerate}
    \item one component, and $\partial K_1\cap\mathbb R$ contains an endpoint of $I,$ or
    \item two disjoint, connected sets one corresponding to curves through the upper half-plane landing at $K_1$, and  the other to curves passing through the lower half-plane.
    \end{enumerate}
    In either case, $ \psi_X (\partial \hat B_0)\subset\partial K_1\cap I\subset\partial B_0\cap I,$ which is a finite forward invariant set, and so $\partial \hat B_0$ is a finite forward invariant set too.
\end{pf}

%(and so correspond to the boundary of the  pruning intervals $Y$ and the basins $B_0$ of $f$).
%We will indicate this by saying that 
%$$\f_X\in \mathcal{E}_Q^d.$$

%   Since $\f$ is analytic and defined on a neighbourhood in $\bar \C \setminus \overline{\mathbb D}$ of $\partial \mathbb D\setminus \hat Y$, it follows by Schwarz reflection, $\f$ is analytic on a neighbourhood of $\overline{\partial \mathbb D\setminus \hat Y}$ mapping $\partial \mathbb D$ (where it is defined) into $\mathbb D$. In particular, $\f \colon \overline{\partial \mathbb D\setminus \hat Y} \to \partial \mathbb D$ is analytic circle map without critical points and with only repelling periodic points. Since $p\in \partial \hat Y$ is a repelling periodic point of $f$, there exists a curve $\gamma$ through $p$ (transversal to the $\partial \mathbb D$) so that $\gamma,\dots,\f^{n-1}(\gamma)$ are disjoint and so that $\f^n(\gamma)\supset \hat \gamma$. Let $\gamma'$ be the corresponding curve through $y$. Since $\f_X$ is analytic and defined on a neighbourhood in $\overline \C \setminus \overline{\mathbb D}$ of $\partial \mathbb D\setminus X_{\phi}^*$, and preserves $\partial \mathbb D$ (wherever it is defined), it follows by Schwarz reflection that $\f_X$ has an analytic extension to a neighbourhood of $\partial \mathbb D\setminus X_{\phi}^*$.

\section{An expanding Markov structure of the external  map}  \label{sec:expandingpruned}
%and its complex extension}
%Let $\hat B_0$ be the basin of the periodic attractors of  $\f_X\colon \partial \D\setminus \hat X \to \partial \D$.
%By  Lemma~\ref{lem:fXproperties}  and Properties (1),(2) of Lemma~\ref{lemQ} the immediate basin 
%of each periodic attractor of $\f_X$ is disjoint from $\hat Y$ if this periodic attractor does not 
%contain a critical point in its basin, and contains a component of the interior of $\hat Y$ if the periodic attractor
%contains a critical point in its basin.
%For each finite forward invariant subset $Q$ of $z\mapsto z^d$ define 
%$\mathcal E^d_Q$ to be the set of maps $\f_X$ so that $Q=Q(\f_X,X,X^*)$. 
Assume that $f$ has only hyperbolic periodic points
(or that the basins of the parabolic points are sufficiently small so that 
they be included in a set $K_{X,O}$ as in Theorem~\ref{thm:KX}).
\label{def:Lambda'} 
Then define  
$$\Lambda'_N=\{z\in  \partial D ; \f^n_X(z)\notin \hat Y \cup \hat B_0 \mbox{ for all }0\le n\le N\} \mbox{ and }\Lambda'_\infty=\cap_{N\ge 0}  \Lambda_N$$
where we should note that $\hat Y$ contains a neighbourhood of the discontinuity 
of $\f_X$.  By Lemma~\ref{lemQ}(3)  each point in $\partial (\hat Y \cup \hat B_0)$ is periodic 
or eventually periodic. Therefore we define 
$$Q_{\hat f_X}=\partial (\hat Y \cup \hat B_0)$$
and set 
$$Q_f :=  h_{\hat f_X}(Q_{\hat f_X} )$$
where $h_{\hat f_X}$ is the semi-conjugacy from Lemma~\ref{lem:comb-semiconj}.

\begin{rem} In fact, if the attractors of $B$ are have sufficiently small diameter, 
as described in the 2nd part of Theorem~\ref{thm:KX} then we can replace 
$\f_X$ by $\f_{X,O}$. In this case,  as pointed out in Remark~\ref{rem:smallattractorsB}, all periodic points
of $\f_{X,O}$ will be repelling, we  do {\em not}  
need to include $\hat B_0$ in the above definition of $\Lambda'_N$.
In this case in the next lemma one can replace $\Lambda_\infty'$ 
by $\Lambda_\infty$.  It follows that in this case the basins of the periodic attractors  of $f\colon \Omega_a\to \C$
will be well-inside the domain of the pruned polynomial-like mapping
$F\colon \E \to \E'$ that we will construct.
\end{rem} 

\medskip

\begin{lemma}[$\Lambda'_\infty$ is expanding]
  \label{lem:expmetric}
%   Let $\Lambda=\{z\in \partial \mathbb D; \f^n(z)\notin \interior \hat Y\mbox{ for all }n\ge 0\}$. 
The set $\Lambda'_\infty$ is a forward invariant hyperbolic repelling set and there exists a Riemannian metric $|\cdot |_x$ on $\Lambda'_\infty$ and $\lambda>1$ so that
$|D\f_X(x)v|_{\f_X(x)}\ge\lambda |v|_x$ for $v\in T_x\partial \mathbb D$, $x\in \Lambda'_\infty$. Moreover, 
\begin{enumerate}
\item 
%It follows that there exists $N<\infty$ so that $f$ is expanding on $\Lambda_N'$. In other words,
 there exist $N<\infty$ and  $\lambda>1$ so that 
 $|D\f_X(x)v|_{\f_X(x)}\ge\lambda |v|_x$ for each 
$v\in T_x\partial \mathbb D$, $x\in \Lambda_N'$. 
\item Let $I_1',\dots,I_l'\subset \partial \mathbb D$ 
be the components of $\Lambda_N'$ and let $I_1,\dots,I_k\subset \partial \mathbb D$
be  so that for each $i$ there exists $j_i$ so that  $\f_X(I_i)=I_{j_i}'$ and so that 
$\Lambda_{N+1}'= \cup I_i\subset \cup I'_i$.
\item Each boundary points of $I_i$ and $I_j'$ is periodic or eventually periodic. 
\end{enumerate} 
\end{lemma}
\begin{pf} That $\Lambda'$ is hyperbolic follows from Ma\~n\'e, see \cite[Theorem III.5.1]{dMvS}. For example, one can modify  $\f_X$ near $\hat X$ to obtain a new $C^2$ covering map of the circle 
  $\partial \mathbb D$. Since the statement only concerns points which stay away outside 
  a neighbourhood of $\hat X$ one can indeed apply the above theorem of Ma\~n\'e. 
Part (1) of this lemma therefore follows by taking an adapted metric, see \cite[Lemma III.1.3]{dMvS}. 
Parts (2) and (3) follow from the fact that each point in $\partial (\hat Y\cup \hat B_0)$ 
is eventually periodic. 
\end{pf}

%a finite union of intervals
%$I_1,\dots,I_k\subset \partial \mathbb D$
%\label{page:I} 
%and $I_1',\dots,I_l'\subset \partial \mathbb D$ covering $\Lambda$ and $\lambda>1$ so that for each $i$
%there exists $j_i$ so that 
%$\f_X(I_i)=I_{j_i}',$ and so that
%$|D\f_X(x)v|_{\f_X(x)}\ge\lambda |v|_x$ for each 
%$v\in T_x\partial \mathbb D$, $x\in I_i$. 

\subsection{The complex extension of the external map and its Markov structure} 

Since $\f_X\colon \partial  \mathbb D \setminus \hat X \to \partial \mathbb D $ is real analytic, it is defined  on a complex neighbourhood of its domain. 
%A complex extensions of this map
%will be denoted by $\F_X$. We also sometimes simply write $\f$ and $\F$
%instead of $\f_X$ and $\F_X$. 
%%Abusing notation, we let
%$\f_X$ also denote this extension.

\begin{lemma}\label{lem:raysrepelling}  For each $p\in Q_{\hat f_X}$ there exists a smooth curve  $\gamma_p$ through $p$ which is transversal to $\partial \D$, so that 
$\f_X(\gamma_p)\supset \gamma_{\f(p)}$ and so that the curves $\gamma_p$, $p\in Q_{\hat f_X}$
are pairwise disjoint.  
\end{lemma} 
\begin{pf} To construct a ray $\gamma_p$ through
$p\in Q_{\hat f_X}$, observe that $p$ is a finite forward or backward
iterate of a repelling periodic point $q\in\partial \mathbb D$ of
period $s$ under $\f_X$.
Using a local linearizing coordinate at $q$ we have that 
$\f_X^s$ is conjugate in a neighbourhood of $q$ to
$z\mapsto\lambda z$ in a neighbourhood of 0,
for some $\lambda>1$. The mapping $z\mapsto \lambda z$ preserves the
line landing at $0$ that corresponds under the inverse of the
linearization to a ray $\gamma_q$ transverse to the circle at $q$. 
We transfer this ray to $p$ through the appropriate iterate of $\f_X$.
Moreover, since $\f_X$ is conformal, we have that
each $\f_X^i(\gamma_p)$ is transverse to $\partial\mathbb D$, so
provided that $\gamma_q$ was chosen short enough, we have that
the $\f_X^i(\gamma_p)$ are pairwise disjoint.
\end{pf}

Now define $$\Gamma_{\f_X}=\cup_{p\in Q_{\hat f_X}} \gamma_p.$$ 

\begin{defn}\label{def:expmarkov} 
 We say that  $\F_X\colon \V\to \V'$  
has  an {\em expanding Markov structure}  if the following holds: 
\begin{itemize}[leftmargin=3mm,labelsep=2mm,itemsep=1mm]
%\item aa 
%\end{itemize}
%\begin{itemize} 
\item[-]   $\F_X\colon \V\to \V'$  is a locally univalent covering map; 
\item[-] $\V\cap\partial \mathbb D=I_1\cup\dots\cup I_k$, 
  $\V'\cap\partial \mathbb D=I'_1\cup\dots\cup I'_l$
  where $I_i,I_j'$ are intervals; 
 \item[-]  $\F_X$ maps $I_i$ onto some $I_j'$ and each boundary point of $I_i$ and $I_j'$ 
 is periodic or eventually periodic.  If such a point $p$
 is periodic, then  $\psi_X(x)\in I$. 
\item[-] $\partial \V\cap\partial \V'\subset\Gamma_{\f_X}$;
\item[-] $\F_X(\partial \V\setminus \Gamma_{\f_X}) \subset
  \partial \V'\setminus \Gamma_{\f_X}$; 
\item[-] the set $\Gamma_{\f_X}$ is forward invariant
in the sense that $\F_X(\Gamma_{\f_X}\cap \overline{\V} )=\Gamma_{\f_X}\cap \V'$;
\item[-]  $\partial \V\cap \partial \V'\subset \Gamma_{\f_X}$;
\item[-] the diameters of puzzle pieces, i.e. components of $\F_X^{n}(\V')$,  tend to zero as $n\to \infty$;
\item[-] every point in $\V \setminus \partial \D$ eventually escapes $\V$; 
\item[-] $\V$ contains a tubular neighbourhood of $\partial \D \setminus \hat B_0$ where $\hat B_0\subset \partial \D$ is the set corresponding to the immediate basins of periodic attractors of $f$.
\end{itemize}
\end{defn} 

%Let $\F_X$ be the complex extension of $\f_X$.

\begin{prop}[Existence of expanding Markov structure] \label{prop:expandingcircle}
  Let $\f_X$ and $\Gamma_{\f_X}$ be as above. 
There exist open sets $\V, \V'$ near $\partial \D\subset \C$ so that  $\f_X$
extends to a map $F\colon \V\to \V'$ which has an expanding Markov structure, in the sense
defined above. 
 \end{prop}

%\begin{defn}  The extension  $\F_X\colon V\to V'$ of $\f_X$ between two open subsets $V\subset V'$ of 
%$\C$  together with a finite collection of 
%rays $\Gamma_{\f_X}$ is said to have an  {\em  expanding Markov structure} 
%%{\em pruned external map} 
%if  the following hold:
%\begin{itemize}
%\item  $\F_X$  is a locally univalent covering map; 
%\item $V\cap\partial \mathbb D=I_1\cup\dots\cup I_k$, 
%  $V'\cap\partial \mathbb D=I'_1\cup\dots\cup I'_l$
%  \item  $\F_X$ maps $I_i$ onto some $I_j'$ 
%\item $\partial V\cap\partial V'\subset\Gamma_{\f_X}$;
%\item $\F_X(\partial V\setminus \Gamma_{\f_X}) \subset
%  \partial V'\setminus \Gamma_{\f_X},$ and
%\item the set $\Gamma_{\f_X}$ is forward invariant
%in the sense that $\F_X(\Gamma_{\f_X}\cap \overline{V} )=\Gamma_{\f_X}\cap V'$.
%\end{itemize}
%%If $f$ has only repelling periodic points, then $V\supset \phi_X(I)$.
%\end{defn}

\begin{pf} For each $i$ take a {\lq}rectangular{\rq} set 
  $\V_i'\supset I_i'$ bounded by two arcs from 
$\Gamma_{\f_X}$ and two curves
from $\{z\in \C; d(z,\partial \mathbb D)=\tau\}$ and let $\V'=\cup \V_i'$.
Here  $d$ is the metric coming induced from the norm from Lemma~\ref{lem:expmetric}. 
Let $\V=\f_X^{-1}(\V')$. 
By the previous lemma, provided $\tau>0$ is chosen small enough, $\f_X$
maps
$\V$ onto $\V'$ (locally univalently). 
\end{pf}

\begin{defn}\label{def:ray-equipotential}   The parts of  $\partial \V ,\partial \V'$ consisting of curves transversal
to $\partial \D$ will be called {\em rays} and the other curves will be called {\em roofs} or the 
 {\em equipotentials} of  $\V$. Each ray is eventually mapped to a ray through 
a periodic point. We shall also use the corresponding terminology for the boundary curves
of $\E:=\psi_X(\V)$ and $\E':=\psi_X(\V')$. 
\end{defn} 

\begin{rem}  Note that we use the notation $\hat F_X$ to emphasise
this complex extension of $\hat f_X$ has additional structure, see 
Remark~\ref{notationfandfX}.  Also note that $\V$ contains $\Lambda'_\infty$. 
\end{rem}

\begin{figure} \hfil 
\beginpicture 
\dimen0=0.32cm 
  \setcoordinatesystem units <\dimen0,\dimen0> point at 50 13 
  \setplotarea x from 0 to 35, y from -1.3 to -40
%      \put {\tiny expanding structure} at -2 -4 
%      \put {\tiny together with the attracting structure} at -2 -5
%        \put {\tiny $\hat F \colon V \to V'$} at -2 -6
  \setlinear
  \plot 0.5 0 14 0 /  
  \plot 0.5 -2 0.5 2 2.5 2 2.5 -2 0.5 -2 /
  \put {\tiny $V_1$} at 1 1  
  \plot 4 -2 4 2 6 2 6 -2 4 -2 / 
  \put {\tiny $V_2$} at 5 1    
  \plot 8 -2 8 2 10 2 10 -2 8 -2 / 
  \put {\tiny $\star$} at 7 0 
    \put {\tiny $V_3$} at 9 1    
 %  \plot 6 -2 6 2 8 2 8 -2 6 -2 / 
   \plot 12 -2 12 2 14 2 14 -2 12 -2 / 
\setdashes <0.5mm> 
   \plot 14.1 -2 14.1 2 /
   \plot 0.4 -2 0.4 2 / 
   \setsolid  
       \put {\tiny $V_4$} at 13 1  
       \setlinear 
   \setdots <0.5mm>  
       \plot 4.05 -3 4.05  3  14.65  3 14.65 -3 4.05 -3 / 
       \plot 5.95 -3 5.95 3 / 
              \plot 7.95 -3 7.95 3 / 
                         \plot 9.95 -3 9.95 3 / 
                                                \plot 11.95 -3 11.95 3 / 
                                                \put {\tiny $\F_X(\V_3)$} at 5 4.3 
                                                \setsolid 
                                                \arrow <5pt> [.2,.4] from 4 3.5 to 6  3.5
                                                 \arrow <5pt> [.2,.4] from 6 3.5 to 4  3.5
                                                 \put {\tiny $\F_X(\V_1)$} at 13 -4.4 
                                                       \arrow <5pt> [.2,.4] from 12 -3.5 to 14  -3.5
                                                 \arrow <5pt> [.2,.4] from 14 -3.5 to 12  -3.5
%                                                                   \put {\tiny $B_0$} at 7 -1.5 
%                  \put {\tiny $B_1$} at 3.3 -1.5 
%                             \put {\tiny $B_2$} at 11 -1.5 
                 \betweenarrows {\tiny $\F_X(\V_2)$}  from 4 -4 to 10 -4  %{\tiny $\gamma$} 
                  \betweenarrows {\tiny $\F_X(\V_4)$}  from 8 4 to 13.5 4 
                     \circulararc 360 degrees from 28 0  center at 24  0
                       \put {\tiny $\star$} at 24 4
                       \put {\tiny $\star$} at 24 -4    
                       \put {\tiny $\bullet$} at 28 0  
                          \put {\tiny $\oplus$} at 20 0  
                          \setdashes <0.5mm> 
                          \plot 28 0 29 0 / 
                          \plot 19 0 20 0 / 
                          \setsolid 
                \circulararc 10 degrees from 29 1  center at 24  0
%                     \plot 28 0 29 1 / 
                       \plot 27.8 1.3 28.8 1.9 / 
                  \circulararc -10 degrees from 25.2 5 center at 24 0 
                             \plot  24.8 4 25.2 5 / 
                             \plot  25.6 3.6 26 4.6 / 
                    \circulararc 10 degrees from 22.8 5 center at 24 0 
			\plot 23.2 4 22.8 5 / 
			\plot 22.4 3.6 22 4.6 / 
  \circulararc -10 degrees from 19 1  center at 24  0
    %                 \plot 20 0 19 1 / 
                       \plot 20.2 1.3 19.2 1.9 /  
               %        \setdots <0.3mm> 
                               \circulararc -30 degrees from 29 1  center at 24  0
                         \circulararc  30 degrees from 19 1  center at 24  0
                         
                       \setsolid 
                            
                                       \circulararc -10 degrees from 29 -1  center at 24  0
       %              \plot 28 0 29 -1 / 
                       \plot 27.8 -1.3 28.8 -1.9 / 
                  \circulararc 10 degrees from 25.2 -5 center at 24 0 
                             \plot  24.8 -4 25.2 -5 / 
                             \plot  25.6 -3.6 26 -4.6 / 
                    \circulararc -10 degrees from 22.8 -5 center at 24 0 
			\plot 23.2 -4 22.8 -5 / 
			\plot 22.4 -3.6 22 -4.6 / 
  \circulararc 10 degrees from 19 -1  center at 24  0
  %                   \plot 20 0 19 -1 / 
                       \plot 20.2 -1.3 19.2 -1.9 / 
                       \put {\tiny $\V_4^+$} at     29.4 2.3     
                                              \put {\tiny $\V_4^-$} at     29.4 -2.3  
                                                                  \put {\tiny $\V_3^+$} at     25.6 6     
                                              \put {\tiny $\V_3^-$} at     25.6 -6  
                                                                  \put {\tiny $\V_2^+$} at     22.4 6     
                                              \put {\tiny $\V_2^-$} at     22.4 -6                               
                                      \put {\tiny $\V_1^+$} at     18.2 1     
                                              \put {\tiny $\V_1^-$} at     18.2 -1    
                                                             \setsolid 
                  \plot -10  -2  -10 2  -6 2 -6 -2 -10 -2    / 
%                  \plot -10 -2 -6 2 / 
                                    \setquadratic 
                  \plot -10 2 -8 -1.3 -6 2 / 
                  \put {\tiny graph of $f\colon I\to I$} at -8 -3 
                           \endpicture
                            \caption{A expanding structure $\hat F_X\colon \V\to \V'$ 
                            for a unimodal map $f\colon I\to I$ with a period three attractor. The upper  parts of 
                            $\V_1^+,\V_2^+$  are mapped to the lower half plane by $\hat F_X$, while $\V_3^+,\V_4^+$ are mapped to the upper half plane. In the middle panel the  action of the corresponding regions $\V_i$ 
                            is drawn schematically, by identifying $z$ and $\bar z$ near $\partial \D$ (and identifying
                            the top right and bottom sides of $\V_4$ and also the top left and bottom side of $\V_1$ in the middle figure).  
                            If we can consider the external map $\hat F_{X,O}$ for the attracting periodic orbit $O$,  then the circle map $\hat F_{X,O}$ is expanding (outside a neighbourhood
                            of the jump points), and then the regions $\V_3$ and $\V_4$ will have a ray in common, 
                            and the  same holds for $V_1$ and $V_2$.  In that case, $\V,\V'$ and $U:=\E,U':=\E'$ will look like 
                            Figure~\ref{fig-prunedplmap2}, where $\V_3\cup \V_4$ will form the right half of $\V$
                            and $\V_1\cup \V_2$ the left half of $\V$.
%                             The interval map 
%has a corresponding global polynomial-like extension $F\colon U \to U'\cup B'$. 
} 
\label{fig;ExtendPLE}
                           \end{figure} 
                           
\begin{center}
\begin{figure}[htbp]
\begin{overpic}[width=4cm]{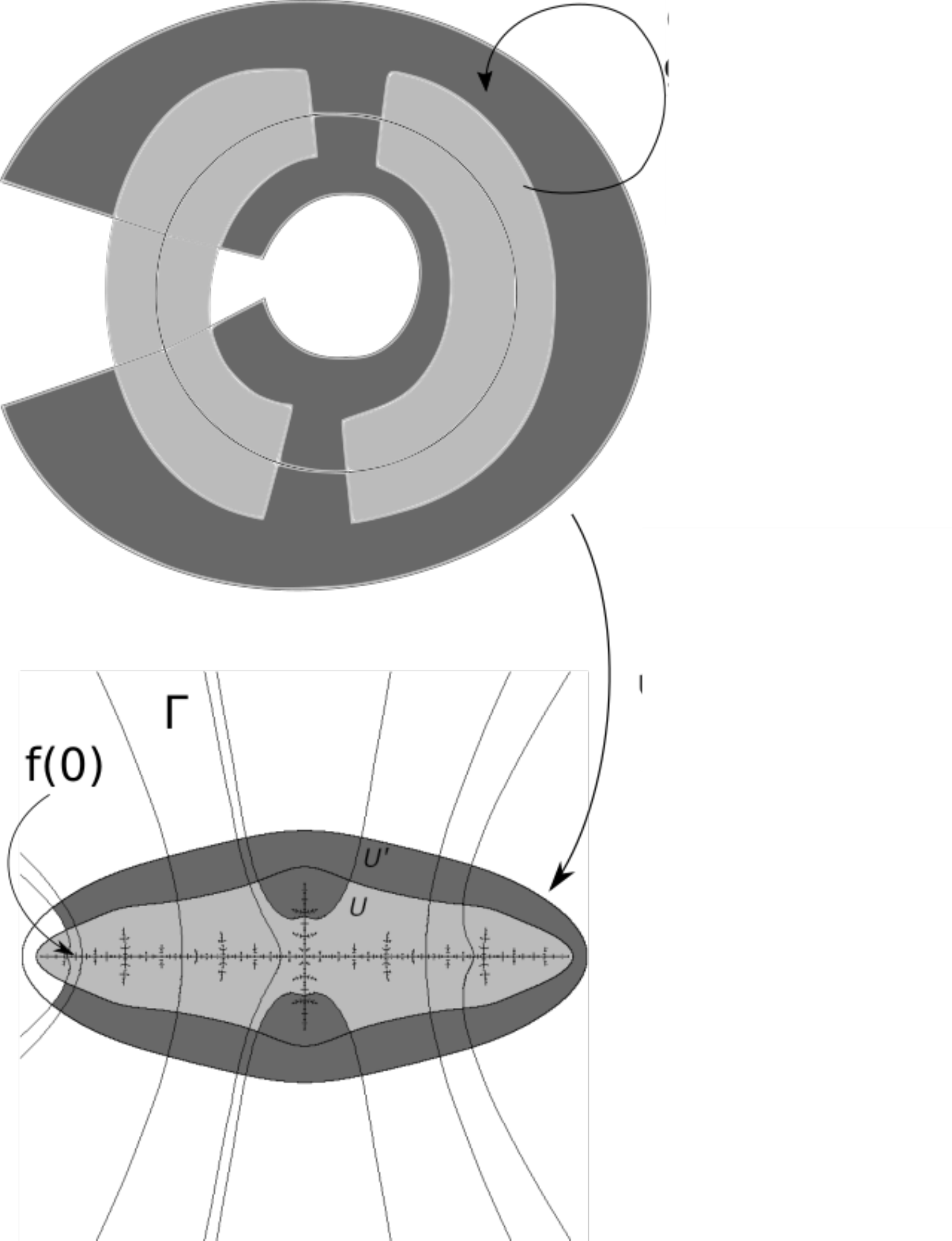}  \put (60,100) {\tiny $\hat f_X\colon \partial \D\to \partial \D$} \put (60,90) {\tiny $\hat f_X= \phi_X \circ f\circ \psi_X$}
\put (60,60) {\tiny $\psi_X\colon \overline \C \setminus \overline{ \mathbb D} \to \overline \C
\setminus K_X$}
\put (60,50) {\tiny $\phi_X\colon  \overline \C
\setminus K_X \to \overline \C \setminus \overline{ \mathbb D} $}
\put (40,73) {\tiny $\V$} 
\put (10,73) {\tiny $\V$} 
\put (54,69) {\tiny $\V'$} 
\put (2,10) {\tiny $U'$} 
\put (25,20) {\tiny $U$} 
\end{overpic} 
\caption{\label{fig-prunedplmap2} The pruned polynomial-like map $F\colon U\to U'$
is constructed from the expanding structure of the external map $\hat F_X\colon \V\to \V'$.
The closure of the $\phi_X$-image of the two components of $\V\setminus \D$ form a connected set because 
the boundary points of $\V\cap \partial \D$ are identified under this map. The curves $\Gamma$ are images
under $\psi_X$ of curves which are transversal to $\partial \D$ 
and which go through periodic and eventually periodic points in $\hat Y\subset \partial \D$. 
This figure corresponds to a simplified version of Figure~\ref{fig;ExtendPLE} where the 
regions $V_3$ and $V_4$ are merged.}
\end{figure}
\end{center}

%\textcolor{blue}{DELETE???
%We will call mapping $\f_X:V\to V'$ 
%%associated to $f,$ 
%%$K_X(f)$, $X^*$ and $Q\supset \phi_X(X^*)$ 
%%by 
%from Proposition~\ref{prop:expandingcircle} 
%an expanding Markov structure associated to the external 
%map $\hat f_X$ together with $f,$  $K_X(f)$, $X^*$ and $Q\supset \phi_X(X^*)$.
%%{\em pruned external mapping associated to} $f$. 
%%
%%
%It will be useful to also  call {\em any} such  holomorphic mapping $F:V\to V'$, together with
%a set of rays $\Gamma_F$, which satisfies the conclusion of 
%Proposition~\ref{prop:expandingcircle} a {\em pruned external
%  mapping}
%  and denote this
%  space by $\mathcal E_Q^d$.
%\textcolor{red}{QQQ SKIP  defining this space?}   
%   In particular, 
%$$\f_X\colon V\to V'\in \mathcal E_Q^d.$$  
%%
%Let us emphasize that a pruned external mapping depends on both the map $f$, 
%as well as on the sets $V,V',\Gamma_{\f_X}$ and to emphasise this dependence 
%we sometimes write  $(f,V,V',\Gamma_{\f_X})$. 
%END DELETE???}

%\textcolor{red}{Perhaps it is not necessary to introduce a
%particular notation for the expanding part of a pruned polynomial-like map? } 

\subsection{Associating a pruned polynomial-like map to a real 
analytic map: the proof of Theorem~\ref{thm: pruned-pol-like-map} when $f$ has only
repelling periodic points.}
\label{subsec:pruned-proof}
Let $\E'=\psi_X(\V')$, $\E=\psi_X(\V)$ 
and $\Gamma_F=\psi_X(\Gamma_{\f_X})$.
Let us show that 
$$F\colon \E \to \E'$$
is a pruned polynomial-like map if $f$ has no periodic attractors. 
By Proposition~\ref{prop:expandingcircle},
$(K_{X}(f)\setminus X)\subset \E$. 
Moreover, we have that if
$z\in\partial E\cap\partial E'$, then
$z\in \Gamma_F$. By Lemma~\ref{lemQ}(1) and since $\V$ contains 
$\Lambda_\infty$ and since we assume 
that there are no periodic attractor (or parabolic periodic points), 
we have that $\E$ contains a neighbourhood of $I$. 
So in that case we define $U=\E$ and $U'=\E'$ and 
 obtain a pruned polynomial-like extension $F\colon U\to U'$
of $f\colon I\to I$. 
\qed

\begin{rem} If $f$ does have periodic attractors, then $\E$ does not necessarily contain $I$. 
To address
this issue   we will define sets  $B,B'$ corresponding 
to the basins of the periodic attractors 
so that $F\colon \E \cup B \to \E'\cup B'$ is a pruned polynomial-like map. 
\end{rem} 

\section{The pruning data $Q$} 
Since the boundary points of $\hat Y\cup \hat B_0$ are (eventually) periodic
 there exists a smallest finite $\hat f_X$-forward invariant set
$Q(\hat f_X,\hat Y, \hat B_0) \subset \partial \D$ so that  $$\partial I_i,\partial I_j', \partial \hat B_0\subset Q(\hat f_X,\hat Y,\hat B_0)\quad \forall i,$$
where $I_i$ are the intervals from the previous section.  

We call $Q(\hat f_X,\hat Y,\hat B_0) $ the {\em pruning data} of the external map 
$\f_X$. The set $Q(\hat f_X,\hat Y,\hat B_0)$ gives some finite combinatorial data
about  $f$ together with the choices made for $J$ and $J^*$. 

\begin{defn}\label{def:admispruning} 
A set $Q\subset \partial \D$ is an {\em admissible pruning set  for }$f\in \mathcal A^{\underline \nu}$ if 
it is  of the above form. 
\end{defn}

We  have already shown that if $f$ has only repelling periodic orbits, 
then it has a pruned polynomial-like extension $F\colon U\to U'$ with respect an admissible 
pruning set $Q$. 

\begin{rem}
If  the basins of all real periodic attractors of $f$ have small diameter (as discussed before), then 
we can replace $\hat f_X$ by $\hat f_{X,O}$. As $\hat f_{X,O}$ will 
have no periodic attractor, we have that $\hat B_0=\emptyset$ in this case.  
On the other hand, if $f$ has a periodic attractor whose diameter is not small, then 
it may be necessary to choose the pruning interval 
$J_i$ and $J_i^*$ to be contained in the basin of this attractor. 
In that case,  the boundary points of $\hat Y$ would be in the basin of $\hat B_0$.
Nevertheless the set $Q$ defined above would be finite.  
\end{rem}

From Lemma~\ref{lem:comb-semiconj} the map  is semi-conjugate via 
a map $h_{\f_X}$ 
to either $z\mapsto z^d$ or to $z\mapsto -z^d$, depending on the sign $\epsilon(f)$ of $f$. 
This means that $h_{\f_X}(Q(\hat f_X,\hat Y,\hat B_0))$ is a finite forward set which is invariant 
under $z\mapsto z^d$ or under $z\mapsto -z^d$ and so this set is {\lq}combinatorially{\rq} defined. 

\begin{defn} 
We say  that the external 
maps $\f_X$ and $\g_X$ associated to two interval maps are {\em pruning equivalent}
\label{pruning-equivalence} 
if the degrees of $\f_X$ and $\g_X$ are the same, the signs of $f,g$ are the same and if 
moreover 
 $$h_{\f_X}(Q(\f_X,\hat Y,\hat B_0))=h_{\g_X}(Q(\g_X,\hat Y_g,\hat B_{0,g})).$$ 
We will write  $Q(\f_X)\sim Q(\g_X)$
if one can make choices for $\hat Y$ and $\hat Y_g$ to that the above equality holds. 
 \end{defn} 

\begin{rem}
$Q(\f_X)\sim Q(\g_X)$ 
does not give any information on the (orbits of the) discontinuities of $\f_X$ and $\g_X$, 
and in particular  does not imply that $\f_X,\g_X$ are topologically conjugate, nor that 
$f,g$ are topologically conjugate. 
\end{rem}

\begin{defn}[The set $Q(F)$ associated to a pruned polynomial-like map] 
Given a pruned polynomial-like map $F\colon U\to U'$ with rays $\Gamma$,
we also have a set $K_F$, a conformal map $\phi_X\colon \C\setminus K_F 
\to \C \setminus \D$, and an   external map $\hat F$ with discontinuities. 
The set $\Gamma\cap K_X$ is forward $F$-invariant and eventually periodic. Thus 
the set $\phi_X(\Gamma\cap K_X)$ corresponds an $\hat F$ invariant subset
of $\D$ which is eventually periodic.  Then $Q(F)$ is defined as the finite subset $\partial \D$
 the image of $\phi_X(\Gamma\cap K_X)$ under the semi-conjugacy with the linear map 
 on $\D$, as in Lemma~\ref{lem:comb-semiconj}.
\end{defn}

\section{An attracting structure near hyperbolic attracting basins} \label{polstructure-attractors} 

In this subsection, we will add some additional structure to the pruned polynomial-like map
which takes care of basins of %periodic attractors of $\f_X$, and correspondingly 
the {\lq}large{\rq} basins of  the interval map $f\colon I\to I$
(namely the ones that cannot be treated by considering the sets $K_{X,O}$).

\begin{defn}  We say that $f \colon B\to B'$ has an {\em attracting structure near attracting basins} if 
there exists a finite union  $\Gamma_a$ of curves so that 
\begin{itemize}[leftmargin=6mm,labelsep=2mm,itemsep=1mm]
\item[-] $B\cap \partial \R$ agrees with the immediate basin of  the periodic attractors of $f$
%(outside a neighbourhood of the points of discontinuities) 
and each component of $B$ is contained in the basin of a periodic attractor of $f$;
\item[-]  $B$ has finitely many components and $B\subset \Omega_a$; 
\item[-] $\partial B\cup \partial B' \subset \Gamma_a$, 
$B'=f(B)\subset B$ and $f(\partial B)\subset \partial B'\cup \Gamma_a$; 
\item[-] $\Gamma_a=\Gamma_a^* \cup \Gamma_a^v\cup \Gamma_a^a$, where 
\item[-] $\Gamma_{a}^*$ is so that $f(\Gamma_a^*\cap B) = \Gamma_a^*$  and each component of $\Gamma_a^*$ is a smooth curve connecting an attracting periodic point $a$
to a repelling periodic points in the boundary of the basin of $a$,  or pre-images of these point; 
\item[-]  each component $\gamma$ of $\Gamma_a^v$ is a piecewise smooth arc in $B$ connecting boundary points of $B$ and iterates of $\gamma$ are disjoint from $\Gamma_a^v$;
\item[-] each component $\gamma$ of   $\Gamma_a^a$ bounds a disk $D_0(a)\ni a$, where $a$
is an attracting periodic point of  $\F_X$ so that $D_0(a)$ is in the basin of $f$ 
and  so that   $\C\setminus (f^n(\gamma)\cup \gamma)$ bounds an annulus  where $n$ is the period of $a$; 
\item[-] each component of $B'\setminus (\Gamma_a \cup \R)$ and of $B'\setminus \Gamma_a$ is a quasidisk.
\end{itemize}
\end{defn}

\begin{figure}
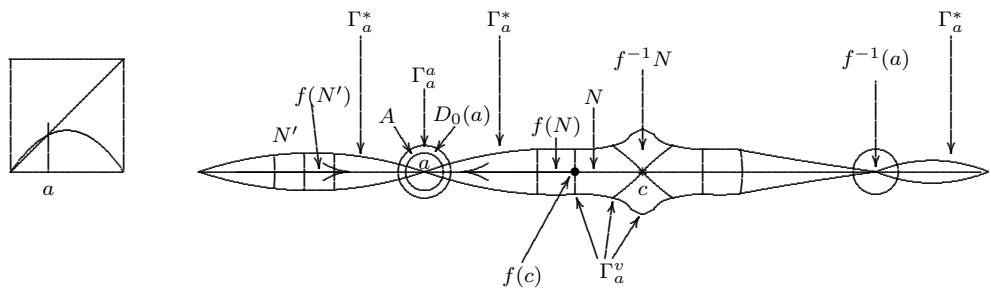
 \hfil \beginpicture \dimen0=0.5cm 
  \setcoordinatesystem units <\dimen0,\dimen0> point at -10 0 
  \setplotarea x from 5 to 25, y from -1.3 to 5
  \setlinear 
  \plot 2 0 23 0 /  
 \plot -3 0 0 3 / 
  \plot -3 0 0 0 0 3 -3 3 -3 0 / 
 \put {\tiny $a$}  at -2   -0.5 
% \setdots <1mm> 
 \plot -2 0 -2 1.3 / 
 \setsolid 
   \setquadratic 
  \plot -3 0 -1.5 1.1 0 0 / 
  \put {\tiny $N'$} at 4.3 1  
    \put {\tiny $f(N')$} at 5.3 2  
         \arrow <4pt> [.2,.67] from 5  1.8  to 5.2  0.1
           \put {\tiny $N$} at 12.5 2
                  \arrow <4pt> [.2,.67] from 12.5  1.7  to 12.5  0.1
                             \put {\tiny $f^{-1}N$} at 13.8 3
                  \arrow <4pt> [.2,.67] from 13.8  2.7  to 13.8  0.6
                            \put {\tiny $f(N)$} at 11.5 1.2
                  \arrow <4pt> [.2,.67] from 11.5  0.9  to 11.5  0.1
   \setquadratic 
  \plot 2 0 5 0.5 8 0 / 
\plot  2 0 5 -0.5 8 0 / 
  \put {\tiny $a$} at 8 0.2 
  \plot 4 -0.42 4.02 0 4 0.42  /
    \plot 4.8 -0.48 4.81 0 4.8 0.48  /
       \plot 5.6 -0.46 5.6 0 5.6 0.46  /
    \plot 8 0 10 0.5 12 0.6 / 
        \plot 8 0 10 -0.5 12 -0.6 / 
               \plot 11 -0.6 11 0 11 0.6 / 
        \plot 12 -0.6 12 0 12 0.6 / 
        \plot 12 0.6 12.5 0.61 13 0.7 / 
         \plot 12 -0.6 12.5 -0.61 13 -0.7 / 
        \plot 13 0.7 13.5 0.3 13.8 0 /                
        \plot 13 -0.7 13.5 -0.3 13.8 0 / 
        \plot 14.6 0.7 14.1 0.3 13.8 0 /                
        \plot 14.6 -0.7 14.1 -0.3 13.8 0 / 
        %%% above $c$
        \plot 13 0.7 13.3 0.8 13.5 1 / 
        \setlinear 
        \plot 13.5 1     13.8 1.14 / 
            \plot 14.1 1     13.8 1.14 / 
            \setquadratic 
            \plot 14.1 1 14.3 0.8 14.6 0.7 / 
               \plot 13 -0.7 13.3 -0.8 13.5 -1 / 
        \setlinear 
        \plot 13.5 -1     13.8 -1.14 / 
            \plot 14.1 -1     13.8 -1.14 / 
            \setquadratic 
            \plot 14.1 -1 14.3 -0.8 14.6 -0.7 /        
        \plot 15.4 0.6 15.1 0.61 14.6 0.7 / 
         \plot 15.4 -0.6 15.1 -0.61 14.6 -0.7 /          
         \plot 15.4 -0.6 15.4 0 15.4 0.6 / 
         \plot 15.4 0.6 16  0.61 16.4 0.58 / 
           \plot 15.4 -0.6 16  -0.61 16.4 -0.58 / 
           \plot 16.4 0.6 16.45 0 16.4 -0.58 / 
           \plot 16.4 0.58 18 0.3 20 0 / 
                   \plot 16.4 -0.58 18 -0.3 20 0 / 
                   \plot 20 0 21.5 0.3 23 0 / 
   \plot 20 0 21.5 -0.3 23 0 / 
%    \put {$p'$} at 2 0.2  
   \put {\tiny $*$} at 13.8  0 
      \put {\tiny $\bullet$} at 12  0 
  \arrow <4pt> [.2,.67] from 10.6  -2.4  to 11.85  -0.15
                            \put {\tiny $f(c)$} at 10.6 -2.8
      \put {\tiny $c$} at 13.8 -0.5 
      \put {\tiny $\Gamma^v_a$} at 13 -2.7 
                \arrow <4pt> [.2,.67] from 13  -2.3  to 13.7   -1.2
          \arrow <4pt> [.2,.67] from 12.8  -2.3  to 13   -0.8
                  \arrow <4pt> [.2,.67] from 12.6  -2.3  to 12.05   -0.7
%    \setdots <2pt>
%    \plot   4.2 -0.43 4.22 0 4.2 0.43  /
       \setsolid 
     \arrow <10pt> [.2,.67] from 3  0  to 6   0
          \arrow <10pt> [.2,.67] from 12  0  to 9   0
%          \put {\tiny $\Gamma$} at 9.2 1.5 
          \put {\tiny $A$} at 7 1.5 
           \put {\tiny $\Gamma_a^a$} at 8 2.5 
            \arrow <4pt> [.2,.67] from 8  2.2  to 8   0.69
          \put {\tiny $D_0(a)$} at 9 1.5 
            \arrow <4pt> [.2,.67] from 7.2  1.2  to 7.65   0.49
                    \arrow <4pt> [.2,.67] from 8.8  1.2  to 8.3   0.56
 %    \arrow <3pt> [.2,.67] from 9  1.5  to 8.3   0.475 %% arrow near a
% \setdots <0.5mm> 
     \circulararc 360 degrees from 8.7 0  center at 8  0
     \circulararc 360 degrees from 8.5 0  center at 8  0
          \circulararc 360 degrees from 20.6 0  center at 20  0
       \setsolid 
            \arrow <4pt> [.2,.67] from 20  2.4  to 20  0.15
                            \put {\tiny $f^{-1}(a)$} at 20 3
                            \put {\tiny $\Gamma^*_a$} at 10 4 
 \arrow <4pt> [.2,.67] from 10  3.6  to 10  0.7
                         \put {\tiny $\Gamma^*_a$} at 22 4 
 \arrow <4pt> [.2,.67] from 22  3.6  to 22  0.5                            
                          \put {\tiny $\Gamma^*_a$} at 6.3 4 
 \arrow <4pt> [.2,.67] from 6.3  3.6  to 6.3  0.5   
\endpicture 
\caption{\label{fig:basin} The attracting structure $f\colon B\to B'$ near an attracting fixed point $a$ associated to a map $f$ whose graph is
as shown on the left, with the rectangular fundamental domain $N$
and its image $f(N)$ and with a curve through the critical point $c$ marked with the symbol {\tiny $*$} and its image $f(c)$ is marked with the symbol 
{\tiny $\bullet$}. The set $B_0^*$ is the set bounded by the curves $\Gamma_a^*$. 
By adding a disk $D_0(a)$ round $a$ and its preimage under $f$, one obtains an attracting structure.  }
\end{figure}

An example of such a structure near a periodic attractor for $f\colon I\to I$ is shown in \ref{fig:basin}.

\begin{prop}\label{prop:expanding-attractor-circle}
Near each hyperbolic periodic attractor of   $f\colon I\to I$, 
the % $a$ whose immediate basin  is $B_0(a)\subset \partial \D$. 
 complex extension $f\colon \Omega_a\to \C$
has an attracting structure $f\colon B\to B'$  in the sense of the previous definition.  
There is a corresponding attracting structure $\F_X\colon \hat B \to \hat B'$ near 
$\partial \D \cap \hat B$. % where as before $\hat B$ is the set corresponding to $B$. 
\end{prop}

\begin{rem}
Note that  by \cite{MMvS,dMvS} the period of periodic attractors of 
$f\colon I\to I$ are bounded. Since $f$ is real analytic,  it follows that $f$
 can have at most a finite number of periodic attractors. \end{rem}

\begin{pf}  Assume that $f$  has only hyperbolic periodic points, 
and assume that it has one or more periodic attractors. Let $a$ be one of these periodic attractors 
and assume it has period $n$.  Let  $B_0(a)$ be the  component of the basin 
which contains $a$.  Choose  fundamental domains  $N,N'$ in $B_0(p)\subset I$, 
so that all critical points that are in the basin of $p$ have forward iterates in $N\cup N'$
and so that $c\in \partial N$.  Here we take $N=N'$ if $Df^n(a)<0$  or if $Df^n(a)=0$
and $x\mapsto f^n(x)$ has a local extremum at $x=a$, 
 while if $Df^n(a)>0$ or $Df^n(a)=0$ and $x\mapsto f^n(x)$ has an inflection point at 
 $x=a$ then 
$N$ and $N'$ are on opposite sides of $a$, see Figure~\ref{fig:basin}.
Abusing notation, we also denote by $N,N'$  rectangular sets in $\C$ 
whose real traces agree with the fundamental domains in $\R$,  whose boundaries consist of  pieces of 
smooth curves and so that $N,f^n(N)$ (respectively $N,f^{2n}(N)$) 
have a smooth common boundary in the orientation  preserving (resp. reversing) 
case. Choose $N,N'$ so that their {\lq}vertical{\rq} boundaries are smooth curves, orthogonal 
to $\partial \D$
and add an additional smooth vertical curve $\gamma\in \Gamma^v$ in $N,N'$ through each iterate of 
a critical point in $N,N'$, orthogonal to $\partial \D$ and connecting the top and bottom boundaries
of $N,N'$.  

Adding forward and backward iterates of $N$ (intersecting $\R$) we obtain a
set $B_0^*\subset \C$ whose boundary is a union of smooth curves landing on the boundary 
points of $\partial B_0\cap I$ and on $a$ (and possibly some of its preimages in the immediate basin
of $a$ as is shown in
Figure~\ref{fig:basin}). Note that preimages of the curves $\gamma\in \Gamma^v$ which go
through a critical value will no longer be orthogonal, and thus we 
get that some preimages of the fundamental domain $N$ will be {\lq}complex{\rq},  see Figure~\ref{fig:basin}.

Now add a topological disk $D_0(a)$  around $a$  so that $f^n(D_0(a))\subset D_0(a)$ and so that $B_0^*\cap \partial D_0(a)$ coincides with an iterate of the {\lq}vertical curve{\rq} 
in $\partial N$. Next add preimages $D_0'$ of $D_0(a)$ centered at preimages of $a$ in the 
immediate basin of $a$. Thus we obtain a set 
$$B=B_0^*\cup D_0(a)\cup D_0'$$ 
 so that $f^n(B)\subset B$ and so that $f^n(\partial B)\subset \partial B \cup \Gamma$  
 where $$\Gamma= f^n(\partial D_0(a))\cup \partial B_0^*.$$
%  is the union of a circle around the orbit of $a$ and a curve going through the 
% \textcolor{red}{critical value???.} 
Taking forward iterates of $f$,  we have shown that around the entire immediate basin of the periodic attractor $\{a,\dots,f^{n-1}(a)\}$ we have an {\em attracting structure} completing the construction.
Using the map $\phi_X$ the attracting structure $f\colon B\to B'$
induces an attracting structure for $\hat f_X\colon \hat B \to \hat B'$. 
\end{pf}

\section{A global pruned polynomial-like structure associated to interval maps} 
\label{sec:existence-pruned-plmap}

Let us now combine the expanding Markov structure and the attracting structure
near basins of periodic attractors, to a global structure on a neighbourhood
of the dynamical interval:

\begin{defn}\label{def:globalpruned} 
 We say that $\F_X\colon \V\cup \hat B \to \V'\cup \hat B'$
is a {\em global pruned polynomial-like structure}, if   $\F_X\colon \V\to \V'$  is an 
expanding Markov structure and $\hat F_X\colon \hat B\to \hat B'$ is 
attracting structure and if these are {\em compatible} in the sense that 
\begin{enumerate}
\item $\hat B\cup \V$ contains a tubular neighbourhood of $\phi_X(I)\subset \partial \D$;
\item each component of  $\hat B$ is either strictly contained in a component of $\V'\setminus \Gamma$ or 
disjoint from $\V'$ and in the former case $\V'\setminus \hat B$ forms a quasidisk.
\end{enumerate} 
\end{defn}

\begin{figure}
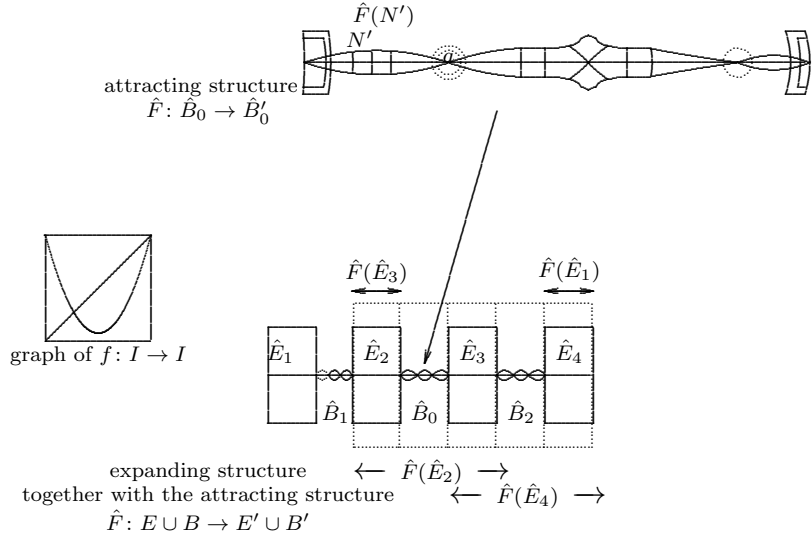
 \hfil 
\beginpicture 
\dimen0=0.32cm 
  \setcoordinatesystem units <\dimen0,\dimen0> point at 10 0 
  \setplotarea x from -5 to 25, y from -10 to -10
    \put {\tiny attracting structure} at -2.3 -1 
           \put {\tiny $\hat F\colon \hat B_0\to \hat B_0'$} at -2 -2
  \setlinear
  \plot 2 0 23 0 /  
% \plot -3 0 0 3 / 
%  \plot -3 0 0 0 0 3 -3 3 -3 0 / 
% \put {\tiny $a$}  at -2   -0.5 
% \setdots <1mm> 
% \plot -2 0 -2 1.3 / 
% \setsolid 
%   \setquadratic 
%  \plot -3 0 -1.5 1.1 0 0 / 
 \setdots <0.5mm> 
     \circulararc 360 degrees from 8.7 0  center at 8  0
     \circulararc 360 degrees from 8.5 0  center at 8  0
          \circulararc 360 degrees from 20.6 0  center at 20  0
          \setsolid 
  \put {\tiny $N'$} at 4.3 1  
    \put {\tiny $\hat F(N')$} at 5.3 2  
  %       \arrow <4pt> [.2,.67] from 5  1.8  to 5.2  0.1
%           \put {\tiny $N$} at 12.5 2
%                  \arrow <4pt> [.2,.67] from 12.5  1.7  to 12.5  0.1
%                             \put {\tiny $f^{-1}N$} at 13.8 3
%                  \arrow <4pt> [.2,.67] from 13.8  2.7  to 13.8  0.6
%                            \put {\tiny $f(N)$} at 11.5 1.2
%                  \arrow <4pt> [.2,.67] from 11.5  0.9  to 11.5  0.1
   \setquadratic 
  \plot 2 0 5 0.5 8 0 / 
\plot  2 0 5 -0.5 8 0 / 
  \put {\tiny $a$} at 8 0.2 
  \plot 4 -0.42 4.02 0 4 0.42  /
          \plot 2.8 1 2.9 0 2.8 -1  /
      \plot 3 1.3 3.1 0 3 -1.3  /
      \plot 1.95 -1.3 2 0 1.95 1.3 / 
          \plot 22.8 1.3 23 0 22.8 -1.3  /
      \plot 22.5 1 22.6 0 22.5 -1  /
      \plot 22 -1.3 22.2 0 22 1.3 /      
      \setlinear 
      \plot 1.95 -1 2.8 -1 / 
          \plot 1.95 1 2.8 1 / 
           \plot 1.95 -1.3 3 -1.3 / 
          \plot 1.95 1.3 3 1.3 / 
          \plot 22 -1.3 22.8 -1.3 / 
     \plot 22 1.3 22.8 1.3 / 
              \plot 22.5 -1 22.8 -1 / 
     \plot 22.5 1 22.8 1 / 
      \setquadratic 
    \plot 4.8 -0.48 4.81 0 4.8 0.48  /
       \plot 5.6 -0.46 5.6 0 5.6 0.46  /
    \plot 8 0 10 0.5 12 0.6 / 
        \plot 8 0 10 -0.5 12 -0.6 / 
               \plot 11 -0.6 11 0 11 0.6 / 
        \plot 12 -0.6 12 0 12 0.6 / 
        \plot 12 0.6 12.5 0.61 13 0.7 / 
         \plot 12 -0.6 12.5 -0.61 13 -0.7 / 
        \plot 13 0.7 13.5 0.3 13.8 0 /                
        \plot 13 -0.7 13.5 -0.3 13.8 0 / 
        \plot 14.6 0.7 14.1 0.3 13.8 0 /                
        \plot 14.6 -0.7 14.1 -0.3 13.8 0 / 
        %%% above $c$
        \plot 13 0.7 13.3 0.8 13.5 1 / 
        \setlinear 
        \plot 13.5 1     13.8 1.14 / 
            \plot 14.1 1     13.8 1.14 / 
            \setquadratic 
            \plot 14.1 1 14.3 0.8 14.6 0.7 / 
               \plot 13 -0.7 13.3 -0.8 13.5 -1 / 
        \setlinear 
        \plot 13.5 -1     13.8 -1.14 / 
            \plot 14.1 -1     13.8 -1.14 / 
            \setquadratic 
            \plot 14.1 -1 14.3 -0.8 14.6 -0.7 /                
        \plot 15.4 0.6 15.1 0.61 14.6 0.7 / 
         \plot 15.4 -0.6 15.1 -0.61 14.6 -0.7 /          
         \plot 15.4 -0.6 15.4 0 15.4 0.6 / 

         \plot 15.4 0.6 16  0.61 16.4 0.58 / 
           \plot 15.4 -0.6 16  -0.61 16.4 -0.58 / 
           \plot 16.4 0.6 16.45 0 16.4 -0.58 / 
           \plot 16.4 0.58 18 0.3 20 0 / 
                   \plot 16.4 -0.58 18 -0.3 20 0 / 
                   \plot 20 0 21.5 0.3 23 0 / 
   \plot 20 0 21.5 -0.3 23 0 /  
   \arrow  <5pt> [.2,.4] from 10 -2 to 7  -12.5
%\endpicture
%
%\beginpicture
%\dimen0=0.3cm 
  \setcoordinatesystem units <\dimen0,\dimen0> point at 10 13 
  \setplotarea x from 0 to 25, y from -1.3 to -40
      \put {\tiny expanding structure} at -2 -4 
      \put {\tiny together with the attracting structure} at -2 -5
        \put {\tiny $\hat F \colon \E\cup B \to \E'\cup B'$} at -2 -6
  \setlinear
  \plot 0.5 0 14 0 /  
  \plot 0.5 -2 0.5 2 2.5 2 2.5 -2 0.5 -2 /
  \put {\tiny $\V_1$} at 1 1  
  \plot 4 -2 4 2 6 2 6 -2 4 -2 / 
  \put {\tiny $\V_2$} at 5 1    
  \plot 8 -2 8 2 10 2 10 -2 8 -2 / 
    \put {\tiny $\V_3$} at 9 1    
 %  \plot 6 -2 6 2 8 2 8 -2 6 -2 / 
   \plot 12 -2 12 2 14 2 14 -2 12 -2 / 
       \put {\tiny $\V_4$} at 13 1 
       \setquadratic 
       \setdots < 0.35mm> 
       \plot 2.5 0 2.75 0.2 3 0 / 
        \plot 2.5 0 2.75 -0.2 3 0 /  
        \setsolid   
            \plot 3 0 3.25 0.2 3.5 0 / 
            \plot 3 0 3.25 -0.2 3.5 0 / 
     \plot 3.5 0 3.75 0.2 4 0 / 
            \plot 3.5 0 3.75 -0.2 4 0 / 
        \plot 6 0 6.3 0.2 6.7 0 /  
             \plot 6.7 0 7 0.2 7.3 0 / 
                     \plot 7.3 0 7.7 0.2 8 0 / 
                             \plot 6 0 6.3 -0.2 6.7 0 /  
             \plot 6.7 0 7 -0.2 7.3 0 / 
                     \plot 7.3 0 7.7 -0.2 8 0 / 
                      \plot 10 0 10.3 0.2 10.7 0 /  
             \plot 10.7 0 11 0.2 11.3 0 / 
                     \plot 11.3 0 11.7 0.2 12 0 / 
                             \plot 10 0 10.3 -0.2 10.7 0 /  
             \plot 10.7 0 11 -0.2 11.3 0 / 
                     \plot 11.3 0 11.7 -0.2 12 0 / 
       \setlinear 
   \setdots <0.5mm>  
       \plot 4.05 -3 4.05  3  13.95  3 13.95 -3 4.05 -3 / 
       \plot 5.95 -3 5.95 3 / 
              \plot 7.95 -3 7.95 3 / 
                         \plot 9.95 -3 9.95 3 / 
                                                \plot 11.95 -3 11.95 3 / 
                                                \put {\tiny $\hat F(\V_3)$} at 5 4.3 
                                                \setsolid 
                                                \arrow <5pt> [.2,.4] from 4 3.5 to 6  3.5
                                                 \arrow <5pt> [.2,.4] from 6 3.5 to 4  3.5
                                                 \put {\tiny $\hat F(\V_1)$} at 13 4.4 
                                                       \arrow <5pt> [.2,.4] from 12 3.5 to 14  3.5
                                                 \arrow <5pt> [.2,.4] from 14 3.5 to 12  3.5
                                                                   \put {\tiny $\hat B_0$} at 7 -1.5 
                  \put {\tiny $\hat B_1$} at 3.3 -1.5 
                             \put {\tiny $\hat B_2$} at 11 -1.5 
                 \betweenarrows {\tiny $\hat F(\V_2)$}  from 4 -4 to 10 -4  %{\tiny $\gamma$} 
                  \betweenarrows {\tiny $\hat F(\V_4)$}  from 8 -5 to 14 -5 
               
               \dimen1=0.1cm 
                 \setcoordinatesystem units <\dimen1,\dimen1> point at 80 30 
  \setplotarea x from 20 to 34, y from -1.3 to -40
                  \setsolid 
                  \plot 20 -7 34  -7 34 7 20 7 20 -7  / 
                  \plot 20 -7 34 7 / 
                  \setquadratic 
                  \plot 20 7 27 -6 34 7 / 
                  \put {\tiny graph of $f\colon I\to I$} at 27 -9 
%                  \setlinear 
%                  \plot 22.5 -7.1 24 -7.1 / 
%                   \plot 26 -7.1 28 -7.1 / 
%                     \plot 30 -7.1 32 -7.1 /
%                \put {\tiny $B_1$} at 23.2 -8 
%                  \put {\tiny $B_0$} at 27 -8 
%                             \put {\tiny $B_2$} at 31 -8
                         %    \put {\tiny the expanding and attracting structure in the setting of a period three attractor } at  13 -9 
                           \endpicture
                            \caption{A global structure $\hat F\colon \V\cup \hat B \to \V'\cup \hat B'$ 
                            which is compatible with the expanding structure 
and the  attracting-structure in the setting of a period three attractor. The interval map 
has a corresponding global polynomial-like extension $F\colon \E\cup B \to \E'\cup B'$. 
} 
\label{fig;ExtendPL}
                           \end{figure}

Similarly, an extension  $F\colon \E\cup B\to \E'\cup B'$ of an interval map $f\colon [-1,1]\to [-1,1]$ is said to have a {\em global pruned polynomial-like  structure} if $\E \cup B\supset [-1,1]$ and if 
 $F\colon \E\to \E'$ (coming from the expanding map $\hat F_X\colon V\to V'$)
matches the {\lq}attracting structure{\rq} $F\colon B\to B'$ we constructed 
in the previous section near basins in the sense above. 

\begin{thm} \label{thm:prunedexpanding} 
There exists a global  pruned polynomial-like structure $F\colon \E \cup B \to \E'\cup B'$.
\end{thm} 
\begin{pf}
To do this, 
extend $\F_X\colon \V \to \V'$  near the boundary of $\partial  B_0$ to a map $\F_X \colon \V_1\to \V_1'$, 
so that $\V_1\setminus \V$, $\V_1'\setminus \V'$ are {\lq}rectangular{\rq}
regions, so that  $(\V_1'\setminus \V_1)\cap \hat B$ agrees with one of the 
sets $\F_X^{-i}(N)$ where $N$ is from the construction on the previous subsection.
This is illustrated in   Figure~\ref{fig;ExtendPL}.
We will add these sets to $\V$ and $\V'$ and denote the resulting map again
by $\F_X\colon \V \to \V'$. Taking the corresponding
sets $\E, \E', B,B'$ we thus obtain a global pruned polynomial-like extension
$F\colon \E\cup B\to \E'\cup B'$ of $f$. 
%The construction of the global pruned polynomial-like map $F\colon \E\cup B\to \E'\cup B'$ goes exactly the same as before. 
\end{pf}

\subsection{ Associating a pruned polynomial-like map to a real 
analytic map: the proof of Theorem~\ref{thm: pruned-pol-like-map} when there are only hyperbolic periodic points.}
\label{subsec:pruned-proof-bis}
As in Subsection~\ref{subsec:pruned-proof} we obtain a global pruned polynomial-like mapping
even in the presence of hyperbolic periodic attractors.

\subsection{Summarising the proof of  Theorem~\ref{thm: pruned-pol-like-map}(1-3)  in a diagram} 

In the previous subsection we showed how to assign to a real analytic map $f$ and 
intervals $J_i$ around the critical values $f(c_i)$ of $f$,  a pruned polynomial-like mapping. 
This was done by considering the pruned Julia set $K_X$, where $X=\partial J^{-1}$. 
Using the corresponding Riemann mapping of $\bar \C\setminus K_X\to \bar \C\setminus \D$
we obtained a circle map and, if $f$ has no periodic attractors, 
an expanding structure near $\partial \D$ (away from the 
discontinuities of $\hat f_X\colon \partial \D \to \partial \D$). 
Thus we obtain a pruned polynomial-like map $F\colon \E\to \E'$ which is a complex extension of 
$f$. In this case we write $U=\E$ and $U'=\E'$. 

If $f$ has periodic attractors whose basins are sufficient small (in the sense of 
Theorem~\ref{thm:KX}) then we can replace $K_X$ by a larger set $K_{X,O}$
and then we again get an expanding structure near $\partial \D$ and a pruned polynomial-like 
map $F\colon \E\to \E'$. In this case the basin of the periodic attractors of $F$ are compactly 
contained inside  of $\E$, and we again obtain a pruned polynomial-like 
map $F\colon U\to U'$ when setting $U=\E, U'=\E'$. In this case we have
that $U$ contains $I$. 

If $f$ has attractors with larger basins, then we 
need to treat the attracting structure near those basins 
as in  Section~\ref{polstructure-attractors} and Figure~\ref{fig:basin}
and ensure that it is compatible with the expanding structure as in Figure~\ref{fig;ExtendPL}.
We thus obtain a global pruned polynomial-like map $F\colon \E \cup B \to \E'\cup B'$. 
In this case we set $U=\E\cup B$ and $U'=\E'\cup B'$.

We should also remark that if $f\in \mathcal{A}^{\underline \nu}$, and given sufficiently small 
interval neighbourhoods $J_i$, the following diagram commutes, 

%Observe that a pruned polynomial-like mapping depends on $(f,U,U',\Gamma_f).$ 
\begin{center}
\begin{tikzpicture}
  \matrix (m) [matrix of math nodes,row sep=3em,column sep=6em,minimum width=2em]
  {
     \f_X \colon \partial \D \to \partial \D\in \mathcal{E}^d  & \F_X \colon \V\to \V'  \\
     f\colon I\to I \in \mathcal{A}^{\underline \nu}  & F_X\colon \E \to \E' .  \\ };
  \path[-stealth]
    (m-2-1) edge node [left] {$\Phi_X$}  (m-1-1) 
    (m-1-1)  edge  node [above] {extension} (m-1-2)
    (m-2-2) edge node [above] {restriction} (m-2-1) 
%   (m-2-1.east|-m-2-2) edge node [below] {$\psi_X$}
 %           node [above] {$\exists$} (m-2-2)
    (m-1-2) edge node [right] {$\Psi_X$} (m-2-2) ; 
  %          edge [dashed,-] (m-2-1); 
  \end{tikzpicture}\\
\end{center}
Here, $\f_X:=\Phi_X(f)=\phi_X\circ f\circ \psi_X$,  $F_X:=\Psi_X(\F_X)=\psi_X\circ \hat F_X\circ \phi_X|_{\psi_X(\V)}$. Here the restriction is to $I\subset \C$ or $\partial \D\subset \C$
and extension refers to a complex extension. Note that $\V,\V'$ and therefore $\E,\E'$
depends on the choice of the set $Y$ and the intervals $I_i\supset \Lambda'$ in $\partial D$. 

%Here $Q$ depends on $(f,X,X^*)$.

\medskip 

\begin{rem} Let us motivate the terminology:  the sets $\V,\V'$ come from the {\em expanding} set $\hat F_X \colon \V\to \V'$, whereas $\hat B,\hat B'$ are associated to the {\em basins} of 
attractors. The set $\E,\E',B,B'$ are the corresponding sets for $F$, 
and $U=\E\cup B$ is a neighbourhood  of $I$. 
\end{rem} 

It may be useful to add the following: 

\begin{lemma}\label{lem:equalityKXKFX}  
Let $f$ be a real analytic map with only repelling periodic point. 
Given sufficiently small  interval neighbourhoods $J_i$ of  
around the critical values $f(c_i)$ of $f$, and setting
$X=\partial  J^{-1}$, we obtain 
$$K(F_X)=K_X(f)$$ 
where we let $F_X\colon U\to U'$ be the pruned polynomial-like mapping constructed above.
\end{lemma}
\begin{pf} This follows from part 3 of Proposition~\ref{prop:expandingcircle}. 
\end{pf}

%Define  
%$$Q_{\partial \Gamma} (\f_X) :=\cup_i   \partial I_i  =\partial \Gamma_{\F}\cap \partial \D.$$
%and 
%$$Q(\f_X) :=Q_{\partial \hat Y}\cup Q_{\partial \hat B_0}\cup Q_{I_i}.$$
%We will call $Q(\f_X)$ an \textcolor{red}{ {\em combinatorial structure} associated to $\f$}. 
%By definition,  $$\f_X \in \mathcal E_Q^d \mbox{ for }Q=Q(\f_X).$$ %and we let $X_{\f_X}$ denote the landing points on $\partial \mathbb
%   D$ of rays in $\Gamma_{\f_X}$.

%\textcolor{red}{Notation $V\to V'$ replace by $V_e\to V_e'$?} 

\section{Holomorphic motions of pruned polynomial-like mappings and line fields}  
\label{sec:holomotion} 

\begin{defn} A {\em holomorphic motion} of $X\subset \C$ over a complex Banach manifold $T\ni 0$
is a family of maps 
$$h_\lambda\colon X\to \C$$
so that 
\begin{enumerate}
\item  $\lambda \mapsto h_\lambda(x)$ depends holomorphically on $\lambda \in T$;
\item $x\mapsto h_\lambda(x)$ is injective; 
\item $h_0=id$. 
\end{enumerate} 
\end{defn} 

Holomorphic motions have very useful properties:

\begin{thm}\cite{BR,ST} \label{thm:extension}
A holomorphic motion $h_\lambda$  of $X\subset \C$ over a Banach ball $T=B_r$
 admits an extension to a holomorphic motion $h_\lambda\colon \C\to \C$ over $T'=B_{r/3}$. 
\end{thm} 

\begin{thm}\cite{BR,MMS} \label{thm:holomotion} 
A holomorphic motion $h_\lambda$  of $X\subset \C$ over $T=\D$ admits an extension to a holomorphic  of
$\C$ over  $T=\D$. Moreover,  the quasiconformal dilatation $\varkappa (h_\lambda)$ of $h_\lambda$ 
is at most $\dfrac{1+|\lambda|}{1-|\lambda|}$ for $\lambda\in \D$. 
\end{thm} 
\begin{pf}  See page 209 of \cite{AIM}. 
\end{pf}

%\trevor{\textcolor{red}{can you obtain this also in a Banach ball, by taking slices????. More precise estimate??.
%Page 209 of Astala $\kappa_\lambda \le \frac{1+|\lambda|}{1-|\lambda|} $ (Extended $\lambda$-Lemma). }}

At the end of this paper,  we will also need to consider deformations over some infinite dimensional manifold,
Therefore  we will also consider quasiconformal motions in Section~\ref{sec:contractible} and in that
setting we will no longer use the above theorems. 

%\textcolor{red}{See also page 11 \cite{ALM}} 

\medskip

%\trevor{\textcolor{red}{Check for $U_f$ vs $U_F$}} 

\begin{defn} \label{defBnua} 
Let  $\mathcal B^{\underline \nu}_a$  be the set of maps which 
are holomorphic on $\Omega_a$, with precisely $\nu$ critical points in $\Omega_a$ 
of order $\ell_1,\dots,\ell_\nu$ and which extend continuously to $\overline{\Omega_a}$. 
\end{defn} 

 Obviously $\mathcal A^{\underline \nu}_a\subset \mathcal B^{\underline \nu}_a$.
 It is easy to see that there exists an open set $\mathcal U\ni f$  in  $\mathcal B^{\underline \nu}_a$
which is complex analytic manifold, see Lemma~\ref{lem:anumanifold} below.
The box mappings we constructed above move holomorphically 
over open subsets of  $\mathcal B^{\underline \nu}_a$.

\begin{prop}\label{prop:persistencepruned}
[Persistence of pruned polynomial-like maps via holomorphic motion] \label{prop:persistenceprunedstructure} 
Let $f\in \mathcal A_a^{\underline \nu}$ and assume that all periodic points
points of $f$  hyperbolic.
%Let $\mathcal{B}_{\Omega_a}^{\underline \nu}$ be the subset of 
%$\mathcal{B}_{\Omega_a}$ of maps $\nu$ critical points $c_1(g),\dots,c_\nu(g)$
%of order $\ell_i$.  
Then there exists an open complex neighbourhood  $\mathcal U$ in $\mathcal B^{\underline \nu}_a$ 
of  $f$ %a neighbourhood
%  $\mathcal U\subset \mathcal A^{\underline \nu}_a$ of $f$ 
%  so that $\mathcal U$ in \textcolor{red}{$\mathcal{B}_{\Omega_a}^{\underline \nu}$ 
% is a complex manifold, QQQQQ and 
 so that there exists a  pruned polynomial-like extension $F\colon U\to U'$ 
   of $f$ with rays $\Gamma$
   so that each map $g\in \mathcal U\cap \mathcal A_a^{\underline \nu}$ %\cap \mathcal{B}_{\Omega_a}^{\underline \nu}$ 
   has 
  a pruned polynomial-like extension $G\colon U_G \to U_G'$
   which is obtained from   $F\colon U\to U'$ by holomorphic motion over $\mathcal U$.
   %\cap \mathcal{B}_{\Omega_a}^{\underline \nu}$.     
  More precisely, 
  \begin{itemize}
  \item there exists 
  a holomorphic motion $h_G$ of $\partial U\cup \partial U'\cup \Gamma$ over 
  $G\in \mathcal U$ %\cap \mathcal{B}_{\Omega_a}^{\underline \nu}$ 
  where 
  $$h_F=id,$$
 so that   if we take $U_G,U_G'$ as the regions bounded by the $h_G$-images
of $\partial U,\partial U'$
%  $$U_G=h_G(U)\mbox{ , }  U'_G=h_G(U') \mbox{ , } \Gamma_G=h_G(\Gamma)$$
  then $G\colon U_G\to U'_G$ is a pruned polynomial-like map;
\item we have 
  $$h_G\circ G(z)=  F\circ h_G(z)\mbox{ for all } z\in \partial U;$$ 
%\item  When $$G\colon U_G \to U'_G$$ for each  
%  $g\in \mathcal U\cap \mathcal A^{\underline \nu}_a$ 
%so that 
%
%
%Let $f_\lambda\in \mathcal A_a^{\underline \nu}$ 
%be a family of real analytic maps depending real analytically 
%on $\lambda\in (-1,1)$ so that all periodic points of $f_0$ are hyperbolic.  
%Let $F_0:U_0\to U'_0$  (or $F_0\colon U_0\cup B_0\to U_0'\cup B_0'$ if $f_0$
%has attracting periodic orbits)  be the pruned polynomial-like mapping $F_0$. 
%
%Then there exists a pruned polynomial-like extension 
%$F_\lambda\colon  U_\lambda   \to U'_\lambda$ of $F_0$
%and  a holomorphic motion $h_\lambda$ of 
%$\partial U'\cup\partial U$ over a neighbourhood of 
%$0\in \D$ such that 
%$$U_\lambda=h_\lambda(U_0)\mbox{ , }  U'_\lambda=h_\lambda(U'_0) \mbox{ , } \Gamma_\lambda=h_\lambda(\Gamma_0)$$ and 
%$$h_{\lambda}\circ F_0(z)=F_{\lambda}\circ
%h_{\lambda}(z)\mbox{ for }z\in \partial U$$ for all $\lambda$ in this neighbourhood of $0\in \D$. 
\item if $f$ has periodic attractors, then in the previous statement we can take
$U_G=\E_G\cup B_G$ and $U'_G=\E_G'\cup B_G'$. 
\item We can choose $\mathcal U$  so that $h_G$ extends to a holomorphic motion of $\C$ over 
  $G\in \mathcal U$. 
\end{itemize} 
\end{prop}
\begin{pf} 
Let us first consider the case that all periodic points of $f$ are repelling.
Choose a pruned polynomial-like extension $F\colon U\to U'$ so that 
$U,U',\Gamma\Subset \Omega_a$.  Note that the boundary of $U$ and $U'$
consists of rays and equipotentials (using the terminology of Definition~\ref{def:ray-equipotential}).  
The finitely many rays (which are in the set $\Gamma$) are  eventually mapped to an  invariant curve 
through hyperbolic repelling periodic points. So choose the neighbourhood $\mathcal U$ of $f$ in $\mathcal A^{\underline \nu}_a$
so that each of these (finitely many) periodic points remains repelling. 

As mentioned, each curve in $\Gamma$ is eventually mapped to a ray $\gamma$ 
through some repelling periodic point $p\in I$. 
For each such $\gamma$ pick an arc $\alpha_t \subset \gamma$, $t\in [0,1]$ which is a fundamental domain (so each orbit hits 
this arc at most once) so that $\alpha_0\in \partial U'$
and so that $F^n(\alpha_1)=\alpha(0)$ or $F^{2n}(\alpha_1)=\alpha_0$ 
where $n$ is the period of the periodic orbit $p$. Here the case $F^{2n}$ corresponds to the situation
that the multiplier of $p$ is negative. Now choose a family of arcs $\alpha_{g,t}\in  \C$, $t\in [0,1]$ 
depending 
analytically on $g$, so that this arc has no self-intersections 
and so that $\alpha_{g,0}=\alpha_0$ and $G^n(\alpha_{g,1})=\alpha_{g,0}$. 
So this defines a holomorphic motion of $\alpha$ over a small neighbourhood $\mathcal U$.

Now extend this
holomorphic motion over each $\gamma\subset \Gamma$ by  considering $F^{-n}(\alpha)$ for any $n$
and thus over all of $\Gamma$. Also define $h_G|\partial U'=id$
and define $h_G$ restricted to the equipotential arcs in $U$ so that there 
$h_G\circ G(z)=  F\circ h_G(z)$. 
Since the equipotential arcs in $\partial U, \partial U'$ have a positive distance from each other,
one can choose a neighbourhood $\mathcal U'\subset \mathcal U$ of $f$ in $\mathcal B^{\underline \nu}_a$ 
so that this defines a  holomorphic motion of $\partial U,\partial U',\Gamma$ over $\mathcal U'$.
If we define $U_G,U_G'$ as the regions bounded by the $h_G$-images
of $\partial U,\partial U'$ then we immediately get that the complex
extension of $g\in \mathcal U' \cap \mathcal A^{\underline \nu}_a$ forms a pruned polynomial-like map $G\colon U_G\to U_G'$
with rays $\Gamma_{G}=h_G(\Gamma)$, 
with the above properties. 

If $f$ has also hyperbolic periodic attractors, 
then one argues similarly considering the arcs in $\Gamma_a$ which are in the basin 
of the periodic attractors. 

Theorem~\ref{thm:extension} implies that there exists a neighbourhood $\mathcal U^*\subset \mathcal U'$ of $f$ in 
$\mathcal B^{\underline \nu}_a$  so that the holomorphic motion $h_G$ of $\partial U,\partial U',\Gamma$ over 
$\mathcal U'$ extends to a holomorphic motion of $\C$ over $\mathcal U^*$. Relabeling $\mathcal U^*$  by $\mathcal U$ 
gives  the last assertion. 
\end{pf}

Combining the results of the previous sections gives the proof of Theorem~\ref{thm: pruned-pol-like-map}. 

\begin{rem} Of course $f$ may only have repelling periodic points, 
whereas a nearby map $g$ might have periodic attractors or parabolic periodic points. 
In that case the immediate 
basin of these periodic attractors of $G\colon U_G\to U'_G$
will be compactly contained in $U_G$ and will not be relevant for the above discussion. 
\end{rem} 
In fact, later on we will also need that the pruned polynomial-like structure persists 
over the entire topologically conjugacy class (or hybrid class) of $f_0$. 
For this we cannot use holomorphic motions, but will use quasiconformal motions, see
Section~\ref{sec:contractible}.

\section{A global pruning structure in the presence of simple parabolic periodic points} \label{subsec:parabolic} 
Lemma~\ref{lem:expmetric} no longer holds in the presence of parabolic periodic points.
However, we can still treat in a more or less similar way  provided all 
parabolic points are simple, in the following sense:

\begin{defn}\label{def:simplepara}
We say that a periodic point $a$ of $f\colon I\to I$ with (minimal) period $n$ is a 
{\em simple} parabolic periodic point if $|Df^n(a)|=1$ and it is of saddle-node,
period-doubling or of pitchfork type. 
Here we say that $a$ is of
\begin{itemize}
\item  {\em saddle-node} type if 
one can write 
$$ Df^n(a)=1  \mbox{ and }  f^n(x)\,\, =a + (x-a) + \tau (x-a)^2 + O((x-a)^3) \mbox{ for }x\approx a$$
with $\tau\ne 0$.
\item {\em period-doubling} type if
$$Df^n(a)=-1 \mbox{ and } f^{2n}(x)=a + (x-a)  - \tau (x-a)^3 + O((x-a)^4) \mbox{ for }x\approx a,$$
with $\tau\ne 0$.
\item  {\em pitchfork type} if  
$$ Df^n(a)=1  \mbox{ and }  f^n(x)\,\, =a + (x-a) + \tau_- (x-a)^3 + O((x-a)^4)  \mbox{ for } x\approx a  $$
 with $\tau_-<0$. 
 \end{itemize} 
\end{defn} 

\begin{lemma} If the Schwarzian derivative of $f$ is negative then each parabolic periodic 
point of $f$ is simple. Moreover, each parabolic periodic point is attracting from (at least) one side. 
\end{lemma}
\begin{pf} This follows from a simple computation. (One can show that if $f$ has only 
one critical point then only the first case can hold when $Df^n(a)=1$.)
\end{pf}  

\begin{rem}\label{rem:perioddoubling}  If $Df^n(a)=-1$ then necessarily 
$D^2f^{2n}=0$. More precisely, if $f^n(x)=-x+ax^2+bx^3+O(x^4)$
then $f^{2n}(x)=x-2(b+a^2)x^3+O(x^4)$. 
\end{rem}

So let us assume that all parabolic periodic points of $f$ are simple. 
Let $B_{0,hyp}$ be the immediate basins of hyperbolic periodic points and $B_{0,par}$
be the immediate basins of parabolic periodic points of $f$. Let $\hat B_{0,hyp}$
and $\hat B_{0,par}$ be the corresponding immediate basins for $\hat f_X$. 
Choose a set $\hat Z_{0,par}$ consisting
of intervals whose union forms a neighbourhood of $\hat B_{0,par}$ so that each boundary point
of $\hat B_{0,par}$ is eventually periodic. Next for each 
 $N\in \N$, define 
$$\Lambda'_{N,\delta}=\{z\in  \partial D ; \f^n_X(z)\notin \hat Y \cup \hat B_{0,hyp} \cup  \hat Z_{0,par} \mbox{ for all }0\le n\le N\}$$ 
and 
$$\Lambda'_{\infty,\delta}=\cap_{N\ge 0}  \Lambda_{N,\delta}.$$

\begin{lemma}[$\Lambda'_\infty$ is semi-expanding]
  \label{lem:semi-expmetric} 
  %   Let $\Lambda=\{z\in \partial \mathbb D; \f^n(z)\notin \interior \hat Y\mbox{ for all }n\ge 0\}$. 
The set $\Lambda'_\infty$ is a forward invariant semi-expanding repelling set:  there exists a Riemannian metric $|\cdot |_x$ on $\Lambda'_\infty$ and $\lambda>1$ so that
$$|D\f_X(x)v|_{\f_X(x)}\ge \lambda'_\delta (x) |v|_x\mbox{ for all }v\in T_x\partial \mathbb D, x\in \Lambda'_\infty$$
where 
$$\lambda_\delta'(x)=\begin{cases} \lambda & \mbox{ if $x$ has distance $\ge \delta$ from any parabolic periodic point} \\
1& \mbox{ otherwise.}\end{cases}
$$ 
Moreover, 
\begin{enumerate}
\item 
%It follows that there exists $N<\infty$ so that $f$ is expanding on $\Lambda_N'$. In other words,
 there exist $N<\infty$ and  $\lambda>1$ so that 
 $|D\f_X(x)v|_{\f_X(x)}\ge\lambda'(x) |v|_x$ for each 
$v\in T_x\partial \mathbb D$, $x\in \Lambda_{N,\delta}'$. 
\item Let $I_1',\dots,I_l'\subset \partial \mathbb D$ 
be the components of $\Lambda_{N,\delta}'$ and let $I_1,\dots,I_k\subset \partial \mathbb D$
be  so that for each $i$ there exists $j_i$ so that  $\f_X(I_i)=I_{j_i}'$ and so that 
$\Lambda_{N+1,\delta}'= \cup I_i\subset \cup I'_i$.
\item Each boundary points of $I_i$ and $I_j'$ is periodic or eventually periodic. 
\end{enumerate} 
\end{lemma}
%
%
%\begin{lemma}
%  \label{lem:expmetric-parabolic}
%%   Let $\Lambda=\{z\in \partial \mathbb D; \f^n(z)\notin \interior \hat Y\mbox{ for all }n\ge 0\}$. 
%The set $\Lambda'_\infty$ is a forward invariant semi-hyperbolic  repelling set: there exist for each $\delta>0$
%constants  $C>0$ and $\lambda>1$ so that 
%$|D\f_X^n(x)|\ge C\lambda^{n'}$ for each $x\in \Lambda$ for which $n$ and $x$ so that
%the distance of $f^n(x)$ to any parabolic periodic point is at least $\delta$. 
%Moreover, there exists a Riemannian metric $|\cdot |_x$ on $\Lambda$ and $\lambda>1$ so that
%$|D\f_X(x)v|_{\f_X(x)}\ge\lambda |v|_x$ whenever  $v\in T_x\partial \mathbb D$, $x\in \Lambda'$,
%and the distance of $x$ to any parabolic periodic  is at least $ \delta$. 
%\end{lemma}
\begin{pf} This lemma is obtained by a minor modification of Ma\~n\'e's Lemma
and the proof of Lemma~\ref{lem:expmetric}.
\end{pf}

Using this lemma we then obtain the analogues of Lemma~\ref{lem:raysrepelling} 
and Proposition~\ref{prop:expandingcircle} in the setting when simple parabolic periodic 
points are allowed. The only difference is that now
we do not obtain an expanding structure on $\partial \D$ outside $\hat B_0$ but outside 
$\hat B_{0,hyp}\cup \hat Z_{0,par}$. Here  $\hat Z_{0,par}\subset \partial \D$
is defined as $\phi_X(Z_{0,par})$.

\begin{figure}[htp]
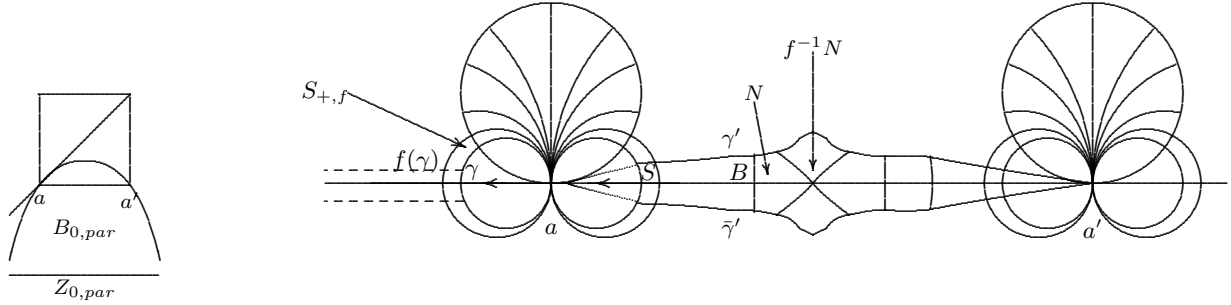
 \hfil
  \beginpicture
\dimen0=0.4cm
\dimen1=0.6cm
\setcoordinatesystem units <\dimen0,\dimen0> point at 38 5
\setplotarea x from -2 to 2, y from 0 to 0
\setlinear
\setsolid 
 \plot -4 -1 -3 0 0 3 / 
  \plot -3 0 0 0 0 3 -3 3 -3 0 / 
 \put {\tiny $a$}  at -3   -0.5 
  \put {\tiny $a'$}  at 0   -0.5 
 \put {\tiny $B_{0,par}$} at -1.5 -1.5 
 \plot -4 -3 1 -3 / 
  \put {\tiny $Z_{0,par}$} at -1.5 -3.5 
% \setdots <1mm> 
% \plot -2 0 -2 1.3 / 
 \setsolid 
   \setquadratic 
  \plot  -4 -2.5 -3 0 -1.5 0.8 0 0 1 -2.5 / 
\setcoordinatesystem units <\dimen1,\dimen1> point at 16 -0.7
\setplotarea x from -2 to 10, y from -2 to 0
\setlinear 
\plot -4 -4 12 -4 /
 \arrow <5pt> [.2,.67] from 3 -4 to 1 -4
  \arrow <5pt> [.2,.67] from -1 -4 to -1.5 -4
   %      \setdashes <0.5mm>
\circulararc 360 degrees from 0 -4 center at 0 -2
   %      \setshadegrid span <4pt>
   %      \hshade -6 -7 2 4 -7 2 / 
   %      \setshadegrid span <2pt>
   %      \hshade 2 -10 5 6 -10 5 / 
%%         \hshade -11 11 27 0 11 20 /
%      \arrow <10pt> [.2,.67] from 4.5 -2 to 2.2 -3.2
\put {\scriptsize $S_{+,f}$} at -5 -2 
      \arrow <5pt> [.2,.67] from -4.5 -2 to -1.9 -3.2
%\put {\scriptsize $S_{-,f}$} at -5 -2 
   %      \setdots <0.5mm>
      \circulararc 43 degrees from 0 -4 center at  -5 -4   
   \circulararc 70 degrees from 0 -4 center at  -2.8 -4   
\circulararc 102 degrees from 0 -4 center at  -1.6 -4   
   \plot 0 -4 0 0 / 
   \circulararc -43 degrees from 0 -4 center at  5 -4   
   \circulararc -70 degrees from 0 -4 center at  2.8 -4   
\circulararc -102 degrees from 0 -4 center at  1.6 -4   
\circulararc 360 degrees from 0 -4 center at -1 -4
\circulararc -360 degrees from 0 -4 center at 1 -4
\circulararc 360 degrees from 0 -4 center at -1.2 -4
\circulararc -360 degrees from 0 -4 center at 1.2 -4
%\plot 1.82 -3.4 2.3 -3.4 / 
%\plot 1.82 -4.6 2.3 -4.6 / 
\put {\scriptsize $\gamma$}  at -1.75 -3.75
\put {\scriptsize $f(\gamma)$}  at -3 -3.5
\plot -5 -4 0 -4 / 
\setdashes <1mm> 
\plot -1.95 -3.7 -5 -3.7 / 
\plot -1.95 -4.4 -5 -4.4 / 
\setsolid 
\put {\scriptsize $S$}  at 2.15 -3.75
\put {\scriptsize $B$}  at 4.15 -3.75
\put {\scriptsize $a$} at 0 -5 
%\plot -4 -3.4 -2.3 -3.4 / 
%\plot 4 -3.4 2.3 -3.4 / 
   %      \put {\scriptsize $-2\epsilon$} at -2.5 -3.5 
   %      \put {\scriptsize $2\epsilon$} at 2.5 -3.5 
   %      \put {\tiny $\bullet$} at -2 -4
   %      \put {\tiny $\bullet$} at 2  -4
\setsolid
   %      \arrow <10pt> [.2,.67] from -6 -2 to -2.3 -3.8
   %      \put {\scriptsize $\tau_f(\Re(w_k))$} at -7 -1.5
   %      \put {\tiny $\bullet$} at -1.8 -3.35
   %      \setsolid
   %      \arrow <10pt> [.2,.67] from -5 0 to -2 -3.22
   %      \put {\scriptsize $\tau_f(w_k)$} at -5 0.7
   %      \put {\scriptsize $C_1$} at -1 -5.5 
   %      \put {\scriptsize $C_2$} at 0 0.2
   \setcoordinatesystem units <\dimen1,\dimen1> point at 24 3.3
\setplotarea x from 14 to 16, y from -2 to 0
%             \setquadratic 
%        \plot 15.6 0.6         15.1 0.61 14.77 0.65 / 
%         \plot 15.4 -0.6 15.1 -0.61 14.77 -0.65 /          
%     %    \plot 15.4 -0.6 15.4 0 15.4 0.6 / 
%         \plot 15.4 0.6 16  0.61 16.4 0.58 / 
%           \plot 15.4 -0.6 16  -0.61 16.4 -0.58 / 
%%           \plot 16.4 0.6 16.45 0 16.4 -0.58 / 
% \plot 15.7 0.6 15.8 0 15.7 -0.58 / 
%           \plot 16.4 0.58 18 0.3 20 0 / 
%                   \plot 16.4 -0.58 18 -0.3 20 0 / 
%                   \plot 20 0 21.5 0.3 23 0 / 
%   \plot 20 0 21.5 -0.3 23 0 / 
      \setquadratic 
%  \plot 2 0 5 0.5 8 0 / 
%\plot  2 0 5 -0.5 8 0 / 
%  \put {\tiny $a$} at 8 0.2 
%  \plot 4 -0.42 4.02 0 4 0.42  /
%    \plot 4.8 -0.48 4.81 0 4.8 0.48  /
%       \plot 5.6 -0.46 5.6 0 5.6 0.46  /
    \plot 9.9 0.45 11 0.5 12 0.6 / 
        \plot 9.9 -0.45 11 -0.5 12 -0.6 / 
   %            \plot 11 -0.6 11 0 11 0.6 / 
  %      \plot 12 -0.6 12 0 12 0.6 / 
        \plot 12 0.6 12.5 0.61 13 0.7 / 
         \plot 12 -0.6 12.5 -0.61 13 -0.7 / 
        \plot 13 0.7 13.5 0.3 13.8 0 /                
        \plot 13 -0.7 13.5 -0.3 13.8 0 / 
        \plot 14.6 0.7 14.1 0.3 13.8 0 /                
        \plot 14.6 -0.7 14.1 -0.3 13.8 0 / 
        %%% above $c$
        \plot 13 0.7 13.3 0.8 13.5 1 / 
        \setlinear 
        \plot 13.5 1     13.8 1.14 / 
            \plot 14.1 1     13.8 1.14 / 
            \setquadratic 
            \plot 14.1 1 14.3 0.8 14.6 0.7 / 
               \plot 13 -0.7 13.3 -0.8 13.5 -1 / 
        \setlinear 
        \plot 13.5 -1     13.8 -1.14 / 
            \plot 14.1 -1     13.8 -1.14 / 
            \plot 12.5 -0.65 12.5 0.65 / 
            \setquadratic 
            \plot 14.1 -1 14.3 -0.8 14.6 -0.7 /        
        \plot 15.4 0.6 15.1 0.61 14.6 0.7 / 
         \plot 15.4 -0.6 15.1 -0.61 14.6 -0.7 /          
         \plot 15.4 -0.6 15.4 0 15.4 0.6 / 
         \plot 15.4 0.6 16  0.61 16.4 0.58 / 
           \plot 15.4 -0.6 16  -0.61 16.4 -0.58 / 
           \plot 16.4 0.6 16.45 0 16.4 -0.58 / 
           \plot 16.4 0.58 18 0.3 20 0 / 
                   \plot 16.4 -0.58 18 -0.3 20 0 / 
                     \arrow <4pt> [.2,.67] from 12.5  1.8  to 12.8  0.2
           \put {\tiny $N$} at 12.5 2
                  \arrow <4pt> [.2,.67] from 13.8  2.9  to 13.8  0.3
                             \put {\tiny $f^{-1}N$} at 13.8 3
%                   \plot 20 0 21.5 0.3 23 0 / 
%   \plot 20 0 21.5 -0.3 23 0 / 
   \put {\scriptsize $a'$} at 20 -1 
   \put {\tiny $\gamma'$} at 12 1 
      \put {\tiny $\bar \gamma'$} at 12 -1 
 %        \circulararc 43 degrees from 20 0 center at  25 0   
     \circulararc 360 degrees from 20 0 center at  20 2     
     \circulararc 43 degrees from 20 0 center at  15 0   
  \circulararc 70 degrees from 20 0 center at  17.2 0   
\circulararc 102 degrees from 20 0 center at  18.4 0   
\setlinear 
 \plot 20 0 20 4 /
 \plot 20 0 23 0 / 
  \circulararc -43 degrees from 20 0 center at  25 0   
  \circulararc -70 degrees from 20 0 center at  22.8 0   
\circulararc -102 degrees from 20 0 center at  21.6 0   
\circulararc 360 degrees from 20 0 center at 19 0
\circulararc -360 degrees from 20 0 center at 21 0
\circulararc 360 degrees from 20 0 center at 18.8 0
\circulararc -360 degrees from 20 0 center at 21.2 0
   \setdots <0.2mm> 
   \setquadratic 
   \plot 8 0 8.3 0.005 8.4 0.013 8.8 0.13 9.9 0.4 / 
      \plot 8 0 8.3 -0.005 8.4 -0.013 8.8 -0.13 9.9 -0.4 / 
\endpicture
\caption{\label{fig:parabolic} The attracting structure near a parabolic one-sided attractor with quadratic order
of contact. The region bounded by the dashed curves come from the expanding structure.
The set $B$ (which is bounded between the curves $\gamma'$ and $\bar \gamma'$ is not a quasidisk.}
\end{figure}

To construct an attracting structure in the the basin of a parabolic periodic point, we
start with a  ray $\gamma$ through $\partial Z_{0,par}$ (on the repelling side
of parabolic point).  We then extend the ray $\gamma$
 to a curve connecting to the parabolic periodic point $a$, as in  Figure~\ref{fig:parabolic}.
More detail on how to construct $\gamma$ is given in the proof of Theorem~\ref{thm:qc-rigiditynearsaddle}. 
Next consider the iterates of $\gamma$ inside a disc $D$ tangent to $\R$ at $a$.
The region $S_{+,f}$ between $\gamma$ and $f(\gamma)$ forms a fundamental crescent:
each orbit near $a$ passes precisely once through $S_{+,f}$. 
One can get a precise description on the shape of these regions using Fatou coordinates. 
Next complete the construction as in the proof of Proposition~\ref{prop:expanding-attractor-circle}.
Let $\gamma'$ be the invariant curve  as shown in Figure~\ref{fig:parabolic}
(analogous to the curves $\Gamma_a^*$  in Figure~\ref{fig:basin}).
Then the region $B$ bounded between $\gamma$ and $\bar \gamma'$ is no longer a quasidisk. 
This is because the curves $\gamma$ and $\gamma'$ will be tangent at $a$. 
Nevertheless if we add the disc discussed above (with the crescent-like regions inside)
we obtain an attracting structure near a simple parabolic periodic point. 
Thus we get a global pruning structure for maps $\mathcal A^{\underline \nu}$
for which all parabolic periodic points are simple,  analogously as in 
 Theorem~\ref{thm:prunedexpanding}.

%To make the structure compatible with the expanding structure from  Lemma~\ref{lem:semi-expmetric} above, 
%we will choose the crescent so that one the crescent coincides with the curve $\gamma$ 
%which comes from the expanding structure outside $\hat B_0$ (going through the eventually 
%periodic point in $\partial \hat Z_{0,par}$). 

\subsection{Holomorphic motion (in a restricted sense) in the presence of parabolic periodic points} 
\label{subsec:holoparabolic} 

If $f$ has parabolic periodic points then the pruned polynomial-like structure that we 
constructed does not  persist under a holomorphic deformation. However, some 
version of this holds nevertheless. Assume that $f$ has a simple parabolic periodic point
of period $n$. Let $\mathcal{B}_a^{\underline \nu}$ be defined as in 
Definition~\ref{defBnua}. 
%Let $\mathcal{B}_{\Omega_a}^{\underline \nu}$ be the subset of 
%$\mathcal{B}_{\Omega_a}$ of maps with $\nu$ critical points $c_1(g),\dots,c_\nu(g)$
%of order $\ell_i$.  \trevor{ \textcolor{red}{Already define?? $\mathcal{B}_a$????} }
Let  $\mathcal{PAR}_f$ the set of  maps $g\in \mathcal{B}_a^{\underline \nu}$
so that $g$ has a simple parabolic periodic point of period $n$. It is not hard to 
show that there exists a neighbourhood $\mathcal U$ of $f$ in $\mathcal{B}_a$
so that $\mathcal U \cap \mathcal{B}_a^{\underline \nu}$ is a complex analytic 
manifold,  see Subsection~\ref{subsec:para-transv}.   

It follows that the analogue of Proposition~\ref{prop:persistenceprunedstructure} holds:
there exists a neighbourhood  $\mathcal U$ of $f$ in   $\mathcal{B}_a$
  so that   there exists a holomorphic motion  $h_G$ of $\partial \E, \partial \E',\partial B,\partial B'$ over  $g\ni \mathcal U \cap \mathcal{PAR}_f$ the conclusion of Proposition~\ref{prop:persistenceprunedstructure} 
  holds. In particular  each $g\in \mathcal A^{\underline \nu}_a\cap  \mathcal{PAR}_f$ 
near $f$ has a pruned polynomial-like extension 
   $$G\colon U_G := \E_G\cup B_G \to \E_G'\cup B_G'=: U'_G$$ 
   where $B_G,B_G'$ has the parabolic structure as described in the previous subsection. 
% \trevor{\textcolor{red}{MORE DETAILS?}}

\section{Absence of line fields}  
\label{sec:linefields} 
%\subsection{Some additional comments}

It is useful to observe that using the arguments of \cite{McM1, Shen}, see also 
\cite{CDKS}, the complex bounds of \cite{CvST, CvS} imply:

\begin{prop}
  Suppose that $F:U\to U'$ is a pruned polynomial-like mapping with a filled pruned Julia set $K_F$. Then the boundary of  $K_F$ does not support a measurable invariant line field.
\end{prop}

\section{Two external maps with the same pruning data 
are qc-conjugate} \label{sec:exter-qc-conj} 

\begin{prop} \label{prop:qc-conjucy-external} 
Let $\f_X,\g_X$ be two external maps with 
global pruned polynomial-like extensions  $\F_X\colon \V_{\F}\cup B_{\F} \to  \V'_{\F}\cup B'_{\F}$
 and $\G_X\colon \V_{\G}\cup B_{\G} \to  \V'_{\G}\cup B'_{\G}$
 so that so that  $Q(\f_X)=Q(\g_X)$. 

Then there exists a $\partial \D$-symmetric quasiconformal map
$H\colon \V_{\F}\cup \V_{\F}'\cup B_{\F} \to \V_{\G}\cup \V_{\G}'\cup B_{\G}$
which sends   $(\V'_{\F},\V_{\F},B_{\F},B'_{\F})$ to the corresponding sets
$(\V'_{\G},\V_{\G},B_{\G},B'_{\G})$ and which conjugates $\F_X$ and $\G_X$. 
%
% $\F_X\colon V_{\F}\cup B_{\F} \to  V'_{\F}\cup B'_{\F}$
% and $\G_X\colon V_{\G}\cup B_{\G} \to  V'_{\G}\cup B'_{\G}$.
% \textcolor{red}{Delete: In particular, $\f_X$ and $\g_X$ are qc conjugate on a neighbourhood 
% of their invariant Cantor sets $\Lambda_\infty(\f_X)$ and $\Lambda_\infty(\g_X)$
% and between $\F_X|\Gamma_{\F}$ to $\G_X|\Gamma_{\G}$?} 
 %If the multipliers of corresponding periodic attractors are the same, 
 %then $f_X$ and $\g_X$ are hybrid conjugate. 
\end{prop} 
\begin{pf} Since $Q(\f_X)=Q(\g_X)$, and since all the relevant sets are quasidisks, 
there exists a  $\partial \D$-symmetric map  which is quasiconformal and maps the following 
sets to corresponding sets
$$H:(\V_{\F},\V'_{\F},B_{\F},B'_{\F},\Gamma_{\F})\rightarrow (\V_{\G},\V'_{\G},B_{\G},B'_{\G},\Gamma_{\G})$$ 
and so that  $H$ is a conjugacy between $\F_X$ and $\G_X$ restricted to $(\partial \V_{\F},\partial \V'_{\F},B_{\F},B'_{\F},\Gamma_{\F})$.
% \to (\partial V_{\G},\partial V'_{\G},B_{\G},B'_{\G},\Gamma_{\G})$ is a conjugacy between $\F$ and $\G$ on these sets. 
Now set $H_0=H$ and define $H_{n+1}$ by $\F_X\circ H_{n+1} = H_n \circ \G_X$.
Since $\F_X$  and $\G_X$ do not have critical points, 
we obtain a qc conjugacy 
%$H\colon (V_{\F},V'_{\F},\Gamma_{\F})\to (V_{\G},V'_{\G},\Gamma_{\G})$
between the external mappings 
$\F:\V_{\F_X}\cup B_{\F_X} \to \V'_{\F_X}\cup B'_{\F_X}$ and 
$\G:\V_{\G_X}\cup B_{\G_X} \to \V'_{\G_X}\cup B'_{\G_X}$.   
\end{pf} 

The previous result shows that two maps $f,g\in \mathcal A^{\underline \nu}_a$ 
with the same pruning data   are {\lq}externally{\rq} conjugate:

\begin{prop}[{\lq}External{\rq} qc-conjugacy of the pruned polynonial-like maps] 
\label{prop:qcconjpartial}
 Assume that  all periodic points of $f,g\in \mathcal A^{\underline \nu}_a$ are hyperbolic and 
 that $Q(\f_X)=Q(\g_X)$. Then the following holds. 
 \begin{enumerate}
\item If $f,g$ have either no periodic attractors,
or all periodic attractors have small basins, 
these maps  have pruned polynomial-like
complex extensions
$F\colon U_F\to U_F'$ and $G\colon U_G \to U'_G$ and  
there exists a qc map $h\colon U_F\cup U_F'\to U_G\cup U_G'$
so that 
$$U_G=h(U_F)\mbox{ , }  U'_G=h(U'_F) \mbox{ , } \Gamma_g=h(\Gamma_F)$$
and $$h\circ F(z)= G\circ h(z)\mbox{ for }z\in \partial U_F\cup \Gamma_F.$$
Any such qc map extends to a conjugacy  $H\colon (U_F \cup U'_F)\setminus K_F 
\to (U_G\cup U_G')\setminus K_G$ between 
$$F\colon U_F\setminus K_F\to U_F' \mbox{ and } G\colon U_G\setminus K_G\to U_G' $$  
with the same dilatation as $h$. 
\item  If $f,g$ do have periodic attractors, then they have a global pruned polynomial-like complex 
extensions
$F\colon \E_F\cup B_F \to \E_G'\cup B_G'$ and $G\colon \E_G\cup B_G \to \E'_G\cup B'_G$
so that if we set $U_F:=\E_F\cup B_F$ and $U_G:=\E_G\cup B_G$ then 
$$F\colon U_F\setminus K_F\to U'_F \mbox{ and }
G\colon U_G\setminus K_G\to U_G' $$  
are qc-conjugate
and so that this qc-conjugacy maps  $\Gamma_F$ to $\Gamma_{G}$.
\end{enumerate} 
%\begin{enumerate}
%\item 
%\item if $f,g$ are topologically conjugate then 
%$F\colon U_F\to U_F'$ and $G\colon U_G\to U_G'$ are qc conjugate.
%\end{enumerate}  
\end{prop} 
\begin{pf} Let  $\E_F,\E'_F,B_F,B_F',\Gamma_F$ and $\E_G,\E'_G,B_G,B_G',\Gamma_G$ to be the sets
corresponding to $\V_{\F},\V'_{\F}$, $\hat B_{\F},\hat B_{\F}',\Gamma_{\F}$ and $
\V_{\G},\V'_{\G},\hat B_{\G},\hat B_{\G}',\Gamma_{\G}$ then the proposition follows immediately from  Proposition~\ref{prop:qc-conjucy-external}.
% then 
%$F\colon U_F\setminus K_F\to U'_F\setminus K_F$ and 
%$G\colon U_G\setminus K_G\to U_G'\setminus K_G$  are qc-conjugate
%such that  $\Gamma_F$ is mapped to $\Gamma_{\tilde F}$.
\end{pf} 

\begin{rem} 
Note that $H$ does not necessarily map $\phi_{X_F}(c_i)$ to $\phi_{X_{\tilde F}}(c_i)$
(where $\phi_{X_F}\colon K_F\to \partial \D$ and $\phi_{X_{\tilde F}}\colon K_{\tilde F} \to \partial \D$ are the boundary maps 
associated to the Riemann mappings  from above)
and so the above proposition definitely does NOT imply that $F\colon U_F\to U'_F$ and $G\colon U\to U'_G$ are conjugate.  
\end{rem} 

\section{Hybrid conjugacy}

\begin{defn}\label{Hf}  Assume that all periodic points of $f,g\in \mathcal A^{\underline \nu}$ are 
either hyperbolic or simple parabolic.  We say that $f,g\colon I\to I$ are {\em real-hybrid conjugate}  if they are 
qc conjugate on a neighbourhood of  $I$, and  there exist a choice for this qc-conjugacy and 
a neighbourhood  $W$ of the attracting periodic points of $f$, 
so that $h$ is holomorphic restricted to the intersection of $W$ with the basin of the attracting 
periodic points of $f$.  We denote  by $\mathcal{H}_{f}^{\R}$
the real hybrid class of $f$.
\end{defn} 

The above definition does not impose anything about the domain 
on which the qc conjugacy is defined. By contrast, the next definition requires
the qc conjugacy to be defined in the domain of pruned polynomial-like mappings.

\begin{defn}\label{HF} 
 We say pruned polynomial-like maps $F\colon U_F\to U_F'$ and $G\colon U_G\to U_G'$ 
 are {\em hybrid conjugate} if there exists 
 a qc topological conjugacy $H$ between $F$ and $G$ such that $\bar{\partial} H=0$ almost everywhere on $K_F$.
 We denote by $\mathcal{H}_{F}$  the hybrid class of $F$, i.e. the set of pruned polynomial-like mappings which are hybrid conjugate to $F$. 
 If $F$ is real then $\mathcal{H}_{F}^{\R}$ denotes the set of real pruned polynomial-like mappings which are qc conjugate to $F$.
\end{defn}

\begin{rem}  If all periodic points of $f,g$  are repelling then \cite{CvS}, \cite{CvST} and absence of line fields 
imply that these maps have pruned polynomial-like extensions $F$ and $G$ which are hybrid conjugate, 
see Theorem~\ref{thm:qcconj}.
\end{rem}

\begin{rem}
If two maps $F,G$ are hybrid conjugate and they have periodic attractors, 
then the conjugacy is conformal on the basins of these attractors. 
\end{rem}

The domain and range $U_F,U_F'$ of a pruned polynomial-like extension of a real analytic
are not uniquely defined.  Hence: 

\begin{defn}\label{HF-equiv} 
 We say pruned polynomial-like maps $F\colon U_F\to U_F'$ and $G\colon U_G\to U_G'$ 
 are {\em equivalent}  if $F=G$ on $U_F\cap U_G$ and 
 $Q(F)=Q(G)$.
 \end{defn}
 
 \begin{rem}
 Note that this definition does {\em not} consider the germ equivalence class.  Indeed 
one way to obtain an equivalent pruned polynomial-like map from 
 $F\colon U_F\to U_F'$ is by keeping, but shortening, the rays $\Gamma_F$
 and lowering the roofs (i.e. equipotentials) in  $\partial \V ,\partial \V'$
 see Figure~\ref{fig;ExtendPLE} corresponding to part of the boundary 
 $\partial U_F,\partial U_F'$ of the domains
 which is not contained in $\Gamma_F$. Since the pruning data of two 
 equivalent pruned polynomial-like maps is the same, their Julia sets also are the same. 
 \end{rem} 
 
\begin{lemma}\label{lem:equiv-hybrid} 
 If $F\colon U_F\to U_F'$ and $G\colon U_G\to U_G'$ are equivalent, then there exists domains 
 $U,U'$ so that $U,U'\subset U_F\cap U_G$ and so that 
 $F=G\colon U\to U'$ is a pruned polynomial-like mapping.
\end{lemma}
\begin{pf} Let $U'$ be the connected component of $U_F'\cap U_G'$
which intersects $I$. Then let $U$ be the connected component of $F^{-1}(U')$
containing $I$. 
\end{pf}

\section{Topologically conjugate maps have   qc-conjugate pruned extensions}\label{sec:qc-conj} 

One of the main consequences of having pruned polynomial-like structure
is that topologically conjugate maps  in $\mathcal A^{\underline \nu}$ (without parabolic periodic points)
are qc conjugate on a neighbourhood of the interval $I$. The next theorem is based on the usual 
pullback argument.

The most important part of the next theorem is point 3, where it is shown 
that the dilation of the qc conjugacy between $F$ and $G$ is {\em only} depends 
on dilatation of the qc conjugacy  between
$U'_F\setminus (U_F\cup \Gamma_F')$,  $U_F\setminus \Gamma_F'$ and $U_F\setminus (U_F'\setminus \Gamma_F'$)
and $U'_G\setminus (U_G\cup \Gamma_G')$,  $U_G\setminus \Gamma_G'$ and $U_G\setminus (U_G'\setminus \Gamma_G'$).
A consequence of this is that nearby hyrbid conjugate maps, are $\kappa$-qc conjugate
with $\kappa>1$ close to one, see Corollary~\ref{cor:keytypelemma}.

\begin{thm}[Pullback argument] 
\label{thm:qcconj} \label{thm:pullback} 
Assume that $f,g\in \mathcal A^{\underline \nu}_a$ are topologically conjugate on $I$
and do not have parabolic periodic points. Then the following holds. 
\begin{enumerate}
\item If $f,g$ have either no periodic attractors,
or all periodic attractors have small basins, 
these maps  have pruned polynomial-like
complex extensions
$F\colon U_F\to U_F'$ and $G\colon U_G \to U'_G$
which are also qc conjugate.
\item  If $f,g$ do have periodic attractors, then they have a global pruned polynomial-like complex 
extensions
$F\colon U_F \to U_F'$ and $G\colon U_G\to U_G'$
which are qc conjugate.
In this case $U_F:=\E_F\cup B_F$, $U_F':=\E_F'\cup B_F'$, $U_G:=\E_G\cup B_G$ and $U_G':=\E'_G\cup B'_G$. 
\item 
%In particular, $f,g$ are qc conjugate on a neighbourhood of the real interval $I$.
Moreover, if $f,g$ are real-hybrid conjugate, and if the external qc conjugacy 
between their pruned polynomial-like extensions $F,G$  from Proposition~\ref{prop:qcconjpartial}
 has dilatation at most $\varkappa$ outside $K_F$,
 or equivalently on
 $U'_F\setminus (U_F\cup \Gamma_F')$,  $U_F\setminus \Gamma_F'$ and $U_F\setminus (U_F'\setminus \Gamma_F'$)
  %$U'\setminus U,U \setminus U',\Gamma\cap (U'\setminus U)$, 
  then $F,G$ are hybrid-conjugate and 
 $\varkappa$-qc conjugate on their entire domain. 
\end{enumerate} 
\end{thm} 
\begin{pf}  Since $f,g$ are topologically conjugate, one can choose pruning intervals
$J_{i,f}, J^*_{i,f} \ni f(c_i)$ and $J_{i,g},J_{i,g}^* \ni g(c_i)$ so that $Q(\f_X)=Q(\g_X)$. 
Provided we take these pruning intervals sufficiently small, 
there exist pruned polynomial-like extensions  $F,G$ which are complex extensions of $f,g$
and so that $U_F,U_F',U_G,U_G'$ are compactly contained in $\Omega_a$. By  Proposition~\ref{prop:qcconjpartial}, 
we also have that there exists a qc conjugacy $h_0$ between $F,G$  outside $K_F$ and $K_G$.

%In fact, if $f,g$ are close to each other then $G$ be
%can be obtained from $F$  by a holomorphic motion $h_G$ over a neighbourhood $\mathcal U$ of $f$ in $\mathcal B^{\underline\nu}_a$, 
%see Proposition~\ref{prop:persistenceprunedstructure}.
%Here $h_G$ is a qc map which sends $(U_F,U'_F,\Gamma_F,B_F)$ to $(U_G,U_G',\Gamma_G,B_G)$. 

By the Main Theorem in \cite{CvS} the maps $f,g$ are qs conjugate. 
Extend this qs conjugacy on the interval $I$ to $U\cup U'$ so that 
it agrees with $h_0$ on $U'\setminus U$, $U\setminus U'$, $\Gamma\cap (U'\setminus U), B$.
%and $B$ with the qc map from the previous proposition. 
Here we do not 
claim any bound on the qc bound of this extension.  Denote this new 
qc map by $h_{0,g}$ and let $K$ be its quasiconformal dilatation. 
Since $h_{0,g}$ is a conjugacy between the postcritical sets of $f,g$, one can 
define a sequence of qc maps $h_{n,g}$ so that $F\circ H_{n+1,g}=H_{n,g}\circ G$
and so that the qc dilatation of $h_{n+1,g}$ is the same as that of $H_{n,g}$. 
Since the space of $\kappa$-qc maps is compact,  $H_{n,g}$ has a convergent subsequence. 
Moreover, for each $z\in U\setminus K_F$ there exists $n$ so that $F^n(z)\notin U$
and therefore $H_{n+i,g}(z)=H_{n,g}(z)$ for all $i\ge 0$, 
Since $h_0\colon B_F\to B_G$ is already a holomorphic conjugacy, 
and since $K_F\setminus B_F$ has no interior, 
in fact $h_n$ converges to a $\kappa$-qc map $H_\infty$.  
Hence $h_\infty$ conjugates $F,G$. 
The last assertion follows immediately since $F$ has no invariant line fields on $\partial K_F$
and so the quasiconformal dilatation of $h_\infty$ is equal to $\varkappa(h_g)$. 

If $f,g$ have no periodic attractors then the second assertion holds: since 
then the filled Julia set  of $F\colon U\to U'$ has no interior, and since
there are no invariant line fields on the boundary  of the filled Julia set. 
So assume $f,g$ do have periodic attractors and that $f,g$ are real-hybrid conjugate. 
 Then by definition there exists a qc conjugacy $\varphi$ 
which is holomorphic on a neighbourhood $W$ of the set of periodic attractors. 
Now choose the sets $B,B'$ of the pruned polynomial-like maps $F\colon \E\cup B\to \E'\cup B'$ 
so that $W$ contains the closed curves $\Gamma_a^a$ surrounding the periodic attractors, see Figure~\ref{fig:basin}.
Now extend the conjugacy $\varphi$ to the basin $B_0$ of the set of periodic attractors. 
Use this extension of $\varphi$ to construct $B_G, B'_G$ (so define $B_G:=\varphi(B)$, $B'_G:=\varphi(B')$). 
Thus we get a pruned polynomial-like map $G\colon \E_G\cup B_G \to \E'_G\cup B_G'$.
Now modify the qc-conjugacy $H_0$ near the periodic attractors  $O$ 
so that it agrees with $\varphi$ on $W$. Then using the pullback argument 
from above, we obtain a $\kappa$-qc conjugacy $H_\infty$ which is holomorphic on the basins 
of $O$. This completes the proof of the second and third assertion. 
\end{pf} 

%\begin{rem} 
%The last statement is  Theorem 6.5 in \cite{ALM}.
%\end{rem} 

\begin{rem} Using the pruned polynomial-like structure, one can also immediately extend the proof 
in \cite{KSS-rigidity} to obtain a qc conjugacy between $F,G$. 
\end{rem} 

The previous theorem gives: 

%%%%KEY ESTIMATE 
\begin{cor}[qc bound] 
\label{cor:keytypelemma}  For each $f\in \mathcal A^{\underline \nu}_a$ so that all its periodic 
points are hyperbolic %and each $\epsilon>0$ 
there 
exist $\delta>0$  and $L>0$ so that the following holds.
Assume that  $g_0,g_1\in \mathcal A^{\underline \nu}_a$ are real-hybrid conjugate to each
other and that $||g_i-f||_\infty<\delta$. Then there exist
pruned polynomial-like extensions $G_i\colon U_{G_i}\to U_{G_i}'$ 
of $g_i$, $i=0,1$ and 
%$L$ %depending on $\delta$ and  
a qc conjugacy $h_{G_0,G_1}$ between   them 
%$G_0\colon U_{G_0}\to U_{G_0}'$ and $G_1\colon U_{G_1}\to U'_{G_1}$ 
whose qc dilatation 
$$\varkappa(h_{G_0,G_1}) \le  1 + L ||g_0-g_1||_\infty.$$  
Here $||\cdot ||_\infty$ is the supremum norm on $\overline{\Omega_a}$. 
\end{cor}
%\trevor{\textcolor{red}{Does $L$ depend on $g_1,g_2$. If so, then no use. So question is whether we can apply BR-MMS result in the high dimensional case?}}
\begin{pf} 
Let  $G_i\colon U_{G_i}\to U'_{G_i}$ the pruned polynomial extensions of $g_i$  obtained
by holomorphic motion from $F\colon U_F\to U_F'$ on some neighbourhood $\mathcal U$ of $f$ 
in $\mathcal B^{\underline \nu}_a$, see Proposition~\ref{prop:persistenceprunedstructure}.  The previous theorem,
Theorem~\ref{thm:pullback},  implies that the dilatation of $h_{G_1,G_2}$
is bounded by the dilatation of the qc conjugacy between the sets
 $U'_{G_i}\setminus U_{G_i},U_{G_i} \setminus U_{G_i}',\Gamma\cap (U_{G_i}'\setminus U_{G_i})$.
Let us now obtain an upper bound for the qc-dilation for (some) such conjugacy. 

%
%
%
% It follows from 
%Proposition~\ref{prop:persistenceprunedstructure}  and Theorem~\ref{thm:holomotion} 
%that the external conjugacy 
%$$h\colon (U_{G_1}\cup U'_{G_1})  \setminus  K(G_1) \to (U_{G_2}\cup U'_{G_2}) \setminus K(G_2)$$
%has dilatation $\le 1+O(||g_1-g_2||)$ provided $g_1,g_2$ are sufficiently close to $f$ on $\Omega_a$.
Assume that $\mathcal U$ is a $\delta_0$ ball in $\mathcal B^{\underline\nu}_a$
and assume $||g_i-f||_\infty<\delta_0/2$. It follows that $\delta_0/2$ ball around $g_0$
is contained in $\mathcal U$.   Let $A=(\delta_0/2)/||g_0-g_1||_\infty$ and 
consider the complex one-dimensional slice
$T=\{g^t; g^t = g_0+ At(g_1-g_0) \mbox{ with }t\in \D\}$ parametrised by $t\in \D$.
It follows that $T\subset \mathcal U$.
Denote by $G_t$ the pruned polynomial-like extension of $g_t$ coming from the holomorphic motion. 
By Theorem~\ref{thm:holomotion}, there exists a holomorphic motion $h^t$ of $\C$ over $\D$ so that 
its quasiconformal dilatation is at most 
$$\varkappa(h^t) \le \dfrac{1+|t|}{1-|t|} \quad \forall t\in \D.$$
Note that  $g_1=g^{t_1}$ for $t_1=1/A$. 
Hence
\begin{equation}
\varkappa(h_{G_0,G_1}) \le  \varkappa(h^{t_1}) \le  \dfrac{1+t_1}{1-t_1} \le 
\dfrac{1+2||g_0-g_1||_\infty/{\delta_0}}{1-2||g_0-g_1||_\infty/{\delta_0}} \le 1+L||g_0-g_1||_\infty \end{equation} 
for some $L$ which depends on $\delta_0$. 
\end{pf}

\begin{rem} In the previous statement it is allowed for $g_i \in \mathcal A^{\underline \nu}_a$ to have attracting or 
parabolic periodic points that are not present for $f$, as their basins would be sufficiently 
small and be compactly contained in $U_{g_i}$ and be contained in $K(G_i)$.
\end{rem} 

\section{Hybrid-conjugacy  in the parabolic case} 

For maps with parabolic periodic orbits we have the following analogue of 
Corollary~\ref{cor:keytypelemma}. 

%%%%KEY Estimate 

\begin{thm}[QC-rigidity of maps with simple parabolic periodic points] \label{thm:qc-rigiditynearsaddle} 
Consider a  maps $f$ with simple parabolic periodic orbits $O_1,\dots,O_p$. Then there exist $\delta>0$ and $L>0$
so that for any two topologically conjugate maps $g_1,g_2$ with precisely one simple parabolic orbit in a $\delta$-neighbourhood 
of $O_i$ for each $i=1,\dots,p$ and with $||g_1-g_2||_\infty<\delta$ the following holds. 
Let $G_1,G_2$ be pruned polynomial-like 
 extensions of $g_1,g_2$.  Then there exists a $\varkappa$-qc $h$ between $G_1,G_2$ so that  
 \begin{equation} 
 \varkappa(h)\le \varkappa \le 1+ L ||g_1-g_2||_\infty.\label{eq:Khineq} 
 \end{equation} 
 If $g_1,g_2$ are hybrid conjugate, then this $\bar \partial h=0$ on the basin the periodic attractors
 and parabolic attractors  of $G_1,G_2$. 
\end{thm}

\begin{pf} Note that $g_1,g,g_2$ are all contained in an infinite dimensional  manifold  $\mathcal P_f$
consisting of maps with parabolic periodic points near $O_i$, $i=1,\dots,p$. 

Consider a simple parabolic periodic point of $f$. For simplicity assume it is attracting from
one side and repelling from the other side, and that this parabolic point is a fixed point at $x=a$.
 This means that $f(x)=(x-a)-(x-a)^2 +h.o.t.$ for $x$ near $a$. As usual one can take the coordinate change
 $w=1/(x-a)$ near $x\in \C$ near $x=a$ and in these coordinates one gets that 
$f$  takes the form $T(w)=w+1 + O(1/w)$ for $w\approx \infty$. 
Take a curve $\gamma$ through the preperiodic  point $\partial Z_{0,par}$ (on the repelling side of parabolic point) which forms part of the boundary of the expanding set $E$. Note that 
$\gamma$ is  a preimage of an invariant curve through a periodic point, 
and we can assume that $\gamma$ (and all its forward iterates) curves are orthogonal to $I$. 
Let $\hat \gamma$ be the corresponding curve in $w$ coordinates. 
Extend $\gamma$ to a curve $\tau$ containing  the fixed point $a$ 
in such a manner that the corresponding 
curve $\hat \tau \supset \hat \gamma$  in $w$ coordinates is disjoint from $T(\hat \tau)$. 
This implies that the region $\hat S_-$ between  $\hat \tau$ and  $T(\hat \tau)$ forms a fundamental domain. The corresponding region $S_-$ in the $z$-plane
has a crescent shape. Now fix some large $R>0$ and 
consider the forward iterates of the region $\hat S'=\hat S\cap \{w;\Im(w)\ge R\}$.
This region is contained in $H=\{w; \Re(w)\ge  R   - C  \}$ where $C>0$ is a constant which is independent of $R$.  
The  corresponding region in the $z$ plane is contained in a disc $D$ whose boundary 
goes through $a$.  Choose $R>0$ so large so that $D$ is disjoint from $\gamma$ and $E'\supset E$. 
Note that the basin of $a$ is contained in $D$ and forward iterates of a crescent shaped region $S_+$.

 Let us show that there exists  $\kappa$-qc conjugacy  $h$ between 
$g_1,g_2$ with the required property.  Denote the regions corresponding to these maps
by $E_i',E_i,S_{i,0}$.  To construct $h$ first define $h \colon E_1'\setminus E_1 \to 
E_2'\setminus E_2$ and choose a hybrid conjugacy $h$ on the immediate basin of the parabolic fixed point $a$. 
Next extend this to a homeomorphism $S_{1,-} \to S_{2,-}$ 
as a conjugacy on the boundary of $S_{1,-}$, by extending  the conjugacy 
in the region between $S_{1,-}\cap D$ and $S_{1,-}\cap E'$.  Next extend $h$ to a neighbourhood of $a$
by pulling back $h|S_{i,-}$. Note that this can be done using the (restricted) holomorphic motion, 
and that is why, as before, we obtain the desired inequality (\ref{eq:Khineq}).  
%\trevor{ \textcolor{red}{ QQQQ } }
\end{pf}

\begin{figure}[htp]
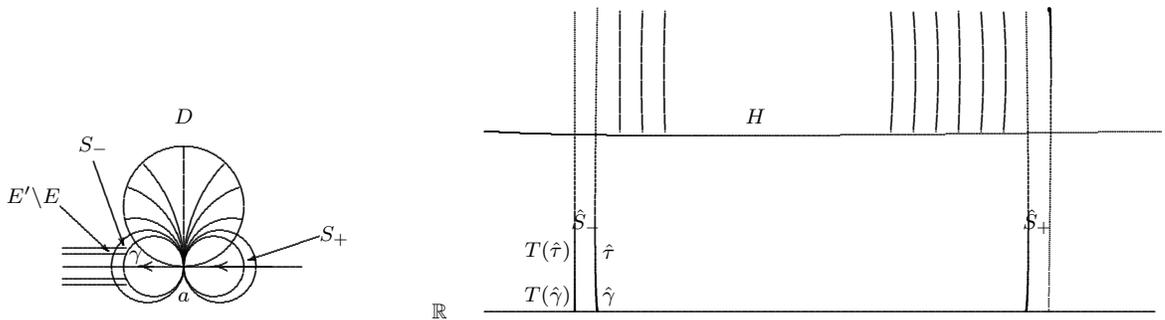
 \hfil
  \beginpicture
\dimen0=0.3cm
\dimen1=0.2cm
\setcoordinatesystem units <\dimen0,\dimen1> point at 10 1
\setplotarea x from -10 to 28, y from -6 to 14
\setlinear
\plot -10 -6 20 -6 / 
%\plot -10 6 20 6 / 
\put {$\scriptstyle  \hat S_-$} at -5.5 0
\put {$\scriptstyle  \hat S_+$} at 14.5 0
%\put {$\scriptstyle R_2$} at 14.5 2.2
%\plot 14 -4 15 -4 / 
\setsolid 
\setlinear 
   %      \plot -2 6 -2 -6 / 
  \put {$\scriptstyle \hat \gamma$} at -4.5 -5 
    \put {$\scriptstyle \hat \tau$} at -4.5 -2 
      \put {$\scriptstyle T(\hat \gamma)$} at -7.2 -5 
    \put {$\scriptstyle T(\hat \tau)$} at -7.2 -2 
    \put {$\scriptstyle H$} at 2  7     
 % \put {$\scriptstyle \bullet$} at -6 -6 
%    \put {$\scriptstyle \tau$} at 14.6 -6.7 
%  \put {$\scriptstyle \bullet$} at 14.6 -6 
  \setquadratic 
\plot  -10 6  0 5.7 20 6 / 
\plot 14 -6 14.1 0 14 14 / 
\plot 15 -6 15.1 10 15 14 / 
%\plot 16 -6 16.03 10 16 14 / 
\setquadratic 
\plot 13 6 13.1 10 13 14 / 
\plot 12 6 12.1 10 12 14 / 
\plot 11 6 11.06 10 11 14 / 
\plot 10 6 10.1 10 10 14 / 
\plot 9 6 9.1 10 9 14 / 
\plot 8 6 8.1 10 8 14 / 
\plot -6 -6 -6 0 -6 14 / 
\plot -5 -6 -5.1 0 -5 14.2 / 
\plot -4 6 -4.0 10 -4 14 / 
\plot -3 6 -3.05 10 -3 14 / 
\plot -2 6 -2.05 10 -2 14 / 
\setlinear
%\plot 13 14 14 14 /
%\plot 12 14 13 14 /
%\plot -6 9.4 -5 9.4 /
%\plot -5 9 -4 9 /
%\plot -4 8.6 -3 8.6 /
%\plot 15 5 20 5 / 
%\plot -10 5 -6 5 / 
\setquadratic 
%\plot -6 5 4  5.3 15 5  / 
\setlinear 
   %      \setdashes <2pt>
   %      \plot -5.6 -6 -5.6 6 / 
%\setshadegrid span <2pt>
%\hshade -6 -5 14 6 -5 14 / 
%\setshadegrid span <1pt>
%\hshade -6 14 15 6 14 15 / 
%\hshade 5 15 20 6 15 20 / 
%\hshade 5 -10 -6 6 -10 -6 / 
%\hshade -6 15 20 -4 15 20 / 
%\hshade -6 -10 -6 -4 -10 -6 / 
\setsolid 
\setlinear 
%\plot 15 -6 15 6 /
   %      \put {$w_0$} at -11 2 
%\put {$\scriptstyle R''$} at 18 5.4 
%\put {$\scriptstyle R''$} at -8 5.4 
%\arrow <10pt> [.2,.67] from 21.5 4 to 20 5
%\put {$\scriptstyle F(\partial E_{\epsilon})$} at 23 3.7 
\put {$\scriptstyle \R$} at -12 -6
%\put {$\scriptstyle B$} at 4 2
%\put {$\scriptstyle \Im(w)=1/\epsilon$} at 17 7
%\put {$\scriptstyle V$} at 21 1
%\put {$\scriptstyle F^{-1}(V)$} at 22.3 2
%\arrow <10pt> [.2,.67] from 20 1 to 15 2
%\arrow <10pt> [.2,.67] from 20.3 2 to 14 3.3
%\arrow <10pt> [.2,.67] from 18 0 to 15.2 1
%\put {$\scriptstyle F^{n+1}(\partial E_\epsilon)\cap R$} at 20 -1
\setquadratic 
%\plot 14.2 1 14.6 1.024 15 1 / 
\setlinear 
%\put {$\scriptstyle F^{-n-1}(V)$} at -11 -0.9
%\put {$\scriptstyle F^{-n}(R)$} at -8.8 1 
%\arrow <10pt> [.2,.67] from -7 1 to -5.7 1
%\arrow <10pt> [.2,.67] from -8.5 -1 to -6.3 -1
   %      \arrow <10pt> [.2,.67] from 21 -2 to 14 -2
   %      \put {$F^{-1}(L)$} at 23.2 -2
%\put {\tiny $\bullet$} at -10 -6 
%\put {$\scriptstyle -\pi/|\sigma|$} at -10 -8 
%\put {\tiny $\bullet$} at 20 -6 
%\put {$\scriptstyle \pi/|\sigma|$} at 20 -8 
\dimen1=0.4cm
\setcoordinatesystem units <\dimen1,\dimen1> point at 25 -2
\setplotarea x from -4 to 4, y from -5 to 0
\setlinear
\setsolid 
\plot -4 -4 4 -4 /
 \arrow <5pt> [.2,.67] from 3 -4 to 1 -4
  \arrow <5pt> [.2,.67] from -1 -4 to -1.5 -4
   %      \setdashes <0.5mm>
\circulararc 360 degrees from 0 -4 center at 0 -2
%\circulararc 360 degrees from 0 -4 center at 0 -1.9
%\put {$\scriptstyle \tau_f(E_{\epsilon})$} at 0 -1
%\put {\tiny $\bullet$} at 0  -3.7 
%\put {$\scriptstyle \sigma$} at 0.4 -3.5 
% \put {\tiny $\bullet$} at 0  -4.3 
%\put {$\scriptstyle \bar \sigma$} at 0.5 -4.3 
\put {$\scriptstyle D$} at 0 1 
\put {$\scriptstyle a$} at 0 -5
      \arrow <3pt> [.2,.67] from -3 -0.5 to -2 -3.3
\put {$\scriptstyle S_-$} at -3 0 
      \arrow <3pt> [.2,.67] from -4.1 -2 to -2.5 -3.5
\put {$\scriptstyle E'\setminus E$} at -5 -1.7 
\put {$\scriptstyle \gamma$} at -1.6 -3.7
   %      \setdots <0.5mm>
         \circulararc 43 degrees from 0 -4 center at  -5 -4   
   \circulararc 70 degrees from 0 -4 center at  -2.8 -4   
\circulararc 102 degrees from 0 -4 center at  -1.6 -4   
   \plot 0 -4 0 0 / 
   \circulararc -43 degrees from 0 -4 center at  5 -4   
   \circulararc -70 degrees from 0 -4 center at  2.8 -4   
\circulararc -102 degrees from 0 -4 center at  1.6 -4   
\circulararc 360 degrees from 0 -4 center at -1 -4
\circulararc -360 degrees from 0 -4 center at 1 -4
\circulararc 360 degrees from 0 -4 center at -1.2 -4
\circulararc -360 degrees from 0 -4 center at 1.2 -4
\plot -1.92 -3.6 -4 -3.6 / 
\plot -1.92 -3.4 -4 -3.4 / 
\plot -1.92 -4.4 -4 -4.4 / 
\plot -1.92 -4.6 -4 -4.6 / 
%\plot 1.92 -3.6 4 -3.6 / 
%\plot 1.92 -4.4 4 -4.4 / 
   \arrow <3pt> [.2,.67] from 4.5 -3 to 2.15 -3.79
\put {$\scriptstyle S_+$}  at 5 -3
\endpicture
\caption{\label{fig:winoutstructure} The construction of the qc conjugacy between $g_1$ and $g_2$ 
  on a {\lq}large{\rq} neighbourhood of $a$.  
  %The strip between and $\partial E_\epsilon$, $F(\partial E_\epsilon)$
  %corresponds to a fundamental annulus around the disc  $\tau_f(E_\epsilon)$. 
  } 
\end{figure}

\newpage 

\part*{Part B: Applications}

\section{Topological and analytic structure on the space of real analytic functions and pruned polynomial-like mappings} \label{manifold-structure}

As before we define a metric on the space  $\mathcal A_a$
by the supremum metric 
$$d(f,g):=||f-g||_{\overline{\Omega}_a}:=\sup_{z\in \overline{\Omega}_a}|f(z)-g(z)|.$$ 
This makes $\mathcal A^{\underline \nu}_a$ a real analytic Banach submanifold 
of the linear real Banach space $\mathcal A_a$ (of real analytic maps on $\Omega_a$ without any conditions on critical points), see Section~\ref{sec:anumanifold}. Using the Implicit Function Theorem, we will show in Section~\ref{sec:manifoldstructure} that 
real-hybrid classes for so-called semi-hyperbolic  maps are  real analytic Banach submanifolds of $\mathcal A^{\underline \nu}_a$. To show that real-hybrid classes are also real analytic 
 Banach submanifolds of $\mathcal A^{\underline \nu}_a$ is one of the main aims of the rest of this paper. 

The {\em tangent space to} $\mathcal A_a^{\underline \nu} $ {\em at} $f$, denoted by
$T_f\mathcal{A}_a^{\underline \nu} $,  consists of 
holomorphic vector fields $v$ defined on $\Omega_a$ 
which vanish at $\partial I$ (since for any $f\in \mathcal A_a$ we have
$f(\partial I)\subset \partial I$)  and so that at each critical point
$c_i$ of $f$, $1\leq i\leq \nu,$ we have
$v^{(j)}(c_i)=0$ for $1\leq j< \ell_i-2$. %, where $\ell_i$ is the order of $c_i$.
The condition on the derivatives $v^{(j)}(c_i)$ ensures that
the critical points do not change order. %bifurcate infinitesimally. 
Note that critical points are allowed to vary within the space $\mathcal A_a^{\underline \nu} $, 
and also that $\mathcal A_a^{\underline \nu} $ is not a linear space. Indeed, if $f\in \mathcal A_a^{\underline \nu} $
and $v\in T_f\mathcal{A}_a^{\underline \nu} $ then in general $f+tv\notin \mathcal A_a^{\underline \nu} $ for $t\ne 0$.

\medskip

\label{realanalytictopology} 
\subsection{The real analytic topology} 
By definition the space $\mathcal A^{\underline \nu}$ is the union of $\mathcal A^{\underline \nu}_a$  over all  $a>0$. 
Of course, it is also equal to the union of $\mathcal A^{\underline \nu}_U$  for all open sets
$U\supset I$, where  $\mathcal A^{\underline \nu}_U$ denotes the set of real analytic
maps in  $\mathcal A^{\underline \nu}$ which have an analytic extension to $U$ and extend continuously 
to $\overline U$. 
The {\em real analytic $C^\omega$} topology on $\mathcal A^{\underline \nu}$ 
(also called the {\em inverse limit topology})  is defined by saying that 
a set $\mathcal O$ is open if $\mathcal O\cap \mathcal A^{\underline \nu}_a$ is open for all $a>0$ (in the topology 
induced by the supremum norm on $\mathcal A^{\underline \nu}_a$). This defines a Hausdorff topology on $ \mathcal A^{\underline \nu}$.
In Appendix~\ref{Append:toponA} we will describe more properties of this topology, in particular:
\begin{itemize}
\item  for  $f_n,f\in \mathcal A^{\underline \nu}$ we have
$f_n\to f$  (in the 
real analytic topology) 
 if and only if there exists $a>0$ so that $f_n,f\in \mathcal A^{\underline \nu}_a$ and $f_n\to f$ uniformly
 on $\overline \Omega_a$.
%We say that $\mathcal O\subset  \mathcal A^{\underline \nu}$ is open 
%if for each $f\in \mathcal O$ and each $f_n\in \mathcal A^{\underline \nu}$ with $f_n\to f$
%we have that $f_n\in \mathcal O$ for all but finitely many $n$. 
\item  for any $f\in \mathcal A^{\underline \nu}_a$,  any $a'>0$ and any open set $\mathcal O\ni f$ 
in the real analytic topology, there exists  $g\in  \mathcal O \setminus \mathcal A^{\underline \nu}_{a'}$;
in particular $\mathcal A^{\underline \nu}_a$ is not an open subset of $ \mathcal A^{\underline \nu}$.
\item  $f\in \mathcal A^{\underline \nu}_a$ can be approximated
(in the $C^k$ topology on $[-1,1]$) by real analytic maps
$f_i\in \mathcal A^{\underline \nu}_{a_i}$ which are not 
contained in $f_i\in \mathcal A^{\underline \nu}_{a'_i}$ where $0<a_i<a_i'$
is so that $a_i'\to 0$. 
\end{itemize} 
The space $\mathcal A^{\underline \nu}$ can be viewed 
as {\em germs} of real analytic interval maps, in the usual sense,  and $\mathcal A^{\underline \nu}=\lim \mathcal A^{\underline \nu}_a$. 
In Appendix~\ref{Append:toponA}  we will elaborate on all this and discuss some further properties of the real analytic topology on 
$\mathcal A$.  This topology on $\mathcal A^{\underline \nu}$ is useful 
as the spaces $\mathcal T_f$ and $\mathcal H_f$ do not  specify the domain of these maps.

\subsection{Real-analytic manifolds} \label{subsec:real-analytic} 
$\mathcal M\subset \mathcal A^{\underline \nu}$ is called a {\em real analytic manifold modelled on a family
of Banach spaces} (or a {\em real-analytic manifold} in short),
if $\mathcal M$ is the union of sets of the form $j_U( \mathcal O_U)$, 
where  $j_U\colon \mathcal O_U\to \mathcal M$ is a family of (the canonical) injections, 
 $U\in \mathcal U$ and  $\mathcal O_U$ is an open subset of  the Banach space $\mathcal A_U$. 
The set  $j_U( \mathcal O_U)$ is called a {\em Banach slice} of $\mathcal M$. 
%if $\mathcal M\cap \mathcal A^{\underline \nu}_a$ is a real Banach manifold for each small  $a>0$. 
A set $\mathcal X\subset \mathcal A^{\underline \nu}$ is called {\em an immersed submanifold} if there exists an analytic manifold
$\mathcal M$ and an analytic  map $i\colon \mathcal M\to  \mathcal A^{\underline \nu}$  so that $Di(m)$
is a linear homeomorphism onto its range, and if $\mathcal X=i(\mathcal M)$. We say that 
$\mathcal X$ is an {\em embedded manifold} if $i\colon \mathcal M \to \mathcal X$ is a homeomorphism with the topology on $\mathcal X$ coming from the one on $\mathcal A^{\underline \nu}$.

\subsection{Topology and analytic structure on the space of germs of  pruned polynomial-like mappings} 
Following \cite[\S4]{McM1} we will also endow a topology on the space of germs of pruned polynomial-like mappings:

\begin{defn} \label{CaratheodoryTopology}  $F_n\colon U_n\to U_n'$ converges to $F\colon U\to U'$ if 
\begin{enumerate}
\item  $(U_n,u_n)\to (U,u)$ in the  {\em sense of  Carath\'eodory}, i.e.  
(i) $u_n\to u$,  (ii) for each compact $K\subset U$, $K\subset U_n$ holds
for $n$ large and (iii) for any open connected set $K'$ containing $u$, 
if $K'\subset U_n$ for infinitely many $n$ then $K'\subset U$. Here we will take $u_n=u=0$. 
\item $F_n\to F$ on compact subsets of $U$. 
\end{enumerate} 
\end{defn} 

Given a pruned polynomial-like mapping $F\colon U\to U'$ so that all 
its periodic points are hyperbolic, and a neighbourhood $\tilde U$ of $\overline U$
on which $F$ is holomorphic and extends continuously to its closure, 
there exists a neighbourhood $B_{\tilde U}(F,\epsilon)$ consisting of 
all holomorphic $G\colon \tilde U\to \C$ 
with $||G-F||_{\tilde U} < \epsilon$ so that each $G\in B_{\tilde U}(F,\epsilon)$
has a pruned polynomial-like map
$G\colon U_G\to U_G'$. This is obtained via holomorphic motion,  see Proposition~\ref{prop:persistencepruned}. Let 
\begin{equation} 
j\colon B_{\tilde U}(F,\epsilon) \to \mathcal{PPL} :=\{\mbox{space of pruned polynomial-like mappings}\}/\!\! \sim 
\label{injbanach} 
\end{equation} 
be the corresponding injection. Here $\sim$ is the germ-equivalence
relation from Definition~\ref{HF-equiv}. 
From the properties
of the holomorphic motion $j$ is continuous. 
 Because of Lemma~\ref{lem:equiv-hybrid}  it follows that if 
 $j(B_{\tilde U} (F,\epsilon)) \cap j(B_{\tilde U'} (F',\epsilon'))\ne \emptyset$
the composition map $(j|B_{\tilde U'}(F',\epsilon'))^{-1}\circ ( j|B_{\tilde U}(F,\epsilon) )$ is the identity map.
 It follows that 
the map $j$ from (\ref{injbanach}) defines an analytic structure on the space of pruned polynomial-like mappings.

 In fact, if we denote by  $\mathcal{PPL}_U$ the space of pruned polynomial-like
 whose domain contains $U$, then for 
 each $F\in \mathcal{PPL}_U$  there exists $a>0$ so that $F\in \mathcal A^{\underline \nu}_a$.
Thus we have the map 
$$\mathcal{PPL}_U  \to  \mathcal A^{\underline \nu}_a \to \mathcal A^{\underline \nu}.$$
The image of 
$B_{U}(F,\epsilon)$ under this map 
$$B_{U}(F,\epsilon) \to \mathcal{PPL}  \to  \mathcal A^{\underline \nu}_a \to \mathcal A^{\underline \nu}.$$
is called a {\em slice} in $\mathcal A^{\underline \nu}$.

\medskip

\subsection{Organisation of this part of the paper} 
In Sections~\ref{sec:anumanifold}-\ref{sec:manifoldstructure}  we will show that the real-hybrid class of any special (hyperbolic) map $f\in  \mathcal A^{\underline \nu}_a$
forms a real analytic  Banach manifold in $\mathcal A^{\underline \nu}_a$. Using this and 
the so-called Mating Theorem, we will then obtain in Sections~\ref{sec:mating}-\ref{sec:topmanifolds} 
an immersed manifold structure 
on the hybrid class of an arbitrary function $f\in  \mathcal A^{\underline \nu}$. 
In Sections~\ref{sec:infinitesimal}-\ref{sec:codimension}  we will determine
the codimension of this manifold. 
In Section~\ref{sec:embedded} we will show that this manifold is even an embedded manifold. 
In  Section~\ref{sec:hybridbanach} we then show that $\mathcal H^\R\cap \mathcal A^{\underline \nu}_a$
is a real Banach space, which means that it is not necessary to work within the real analytic 
topology or within a germ setting. 
%The remaining sections will show that
In Sections~\ref{sec:pathconnected}-\ref{sec:contractible} we then 
show that it is contractible and that these manifolds form a partial lamination.

\section{The space $\mathcal A_a^{\nu} $  is a Banach manifold}\label{sec:anumanifold} 
\begin{lemma}\label{lem:anumanifold} 
 The space $\mathcal A_a^{\underline \nu} $  is a real analytic submanifold
of the linear space $\mathcal A_a$ (of real analytic maps on $\Omega_a$ without any conditions on critical points). 
\end{lemma} 

\begin{pf} If $\ell_i=2$ for all $i=1,\dots,\nu$ then $\mathcal A_a^{\underline \nu}$
is an open subset of $\mathcal A_a$ and so there is nothing to prove. Otherwise
 one can prove this lemma  as was done in Theorem 2.1 in \cite{LSvS3}.
Indeed, since $f\in \mathcal A_a^{\underline \nu}$, there exist $\nu$, $\ell_1,\dots,\ell_\nu$
so that for each $i = 1,2,\dots,\nu$, 
$$f'(c_i) = f''(c_i) = \dots = f^{(\ell_i-1)}(c_i) = 0, f^{(\ell_i)}(c_i)\ne  0.$$ 
Applying the Implicit Function Theorem to the maps $\mathcal A_a\times \C \ni (g,\zeta_i)\mapsto g^{(\ell_i-1)}(\zeta_i)$ for $(g,\zeta_i)$ near $(f,c_i)$,
gives that there exists a neighborhood $W$ of $f$ in $\mathcal A_a$ and uniquely defined functions
$\zeta_i \colon W \to \C$ which are holomorphic  such that $\zeta_i(f) = c_i$ and $g^{(\ell_i-1)}(\zeta_i(g)) = 0, g^{(\ell_i)}(\zeta_i(g))\ne 0$ for each $g \in W$.
Replacing $W$ by a smaller neighborhood, for each $g\in W$ the equation $g'(\zeta)=0$
has $\ell_i-1$ solutions $\zeta$ (counting multiplicity) near $c_i$. It follows that
 for any $g\in W\cap \mathcal A_a^{\underline \nu} $, \,\,  $g'(\zeta)=0$
has a unique solution near $c_i$ (with multiplicity $\ell_i$); hence
  $\zeta_i(g)$ is the only critical point of $g\in W\cap \mathcal A_a^{\underline \nu}$ near $c_i$ and it has multiplicity $\ell_i-1$.

For $g\in W$, write
$$\zeta_i^0(g) = g(\zeta_i(g)), \zeta_i^1(g) = g'(\zeta_i(g)), \zeta_i^2(g) = g''(\zeta_i(g)), \dots.$$
Thus $\zeta_i(g)$ is a critical point of $g$ with multiplicity $\ell_i-1$ if and only if
$\zeta_i^j(g)=0$ for all $1\le j\le \ell_i-2$ (note that
$g^{(\ell_i-1)}(\zeta_i(g)) = 0, g^{(\ell_i)}(\zeta_i(g))\ne 0$ holds automatically for all $g\in W$).
Define $G\colon W\to \C^{(\ell_1-2)+\dots+(\ell_\nu-2)}$ by
$$g  \to (\zeta_1^1(g), \dots, \zeta_1^{(\ell_1-2)}(g), \dots, \zeta_\nu^1(g), \dots , \zeta_\nu^{(\ell_\nu-2)}(g)).$$
The map $(DG)_f\colon T_fW\to \C^{(\ell_1-2)+\dots+(\ell_\nu-2)}$ has maximal rank. Indeed,  if we 
take $1\le i_0\le \nu$,
$1\le j\le \ell_{i_0}-2$ and  the  family $g_t=f+tv$ where $v(z)= \prod_{i\ne i_0}(z-c_i)^{\ell_i} (z-c_{i_0})^j$
then $(DG)_f(v)=se_i\ne 0$ where $i=(\ell_1-2)+\dots+(\ell_{i_0-1}-2)+j$, $e_i$ is the standard unit vector
 in the $i$-th coordinate and $s=j! \prod_{i\ne i_0}(c_{i_0}-c_{i})^{\ell_i}\ne 0$. 
Hence the lemma follows
from the implicit function theorem. 
\end{pf}

\section{Conjugacy classes of semi-hyperbolic  maps form Banach manifolds} 
\label{sec:manifoldstructure}

\begin{defn}\label{def:semihyperbolic} 
 We say that  $f\in \mathcal A^{\underline \nu}_a$ is {\em hyperbolic} 
if all its periodic points are hyperbolic  and all critical points are in basins of periodic attractors.
In this paper, we say that $f$ is {\em semi-hyperbolic } if all its periodic points are hyperbolic and
  so that each critical point of  $f$ is 
  either 
  \begin{enumerate}
   \item[-] eventually mapped onto  a hyperbolic repelling periodic point, or 
  \item[-] periodic, or 
   \item[-] is in the basin of 
  a hyperbolic periodic attractor.  
  \end{enumerate}  
  \end{defn} 
  
   In this section 
  we will show that the hybrid classes of such, and similar, maps form Banach manifolds. 
  
\subsection{A necessary and sufficient condition for a map
to be topologically or hybrid conjugate to a semi-hyperbolic  map} \label{subsec:condhybrid} 

%  In this subsection we will define a function $\Psi$ which determines
%  whether a given map is inside a conjugacy class of $f_0$.  £
%  In this subsection we will associate to each hyperbolic interval map $f$
%  a collection of vector fields $v$ which is transversal to the hybrid class of $f$.
%  , 
%in the sense   that $f_t=f+tv\in\mathcal A^{\underline \nu}_a$ is not contained in the hybrid class of $f$
%for $|t|$ small. 
   
The theorem below
  shows, amongst other things, when a map $g\in \mathcal A^{\underline \nu}_a$
  is real-hybrid conjugate to a semi-hyperbolic  map $f$.

\begin{thm}\label{thm:definitonPsi} 
  Let $f\in\mathcal A^{\underline \nu}_a$ 
  be a  real analytic interval map which is semi-hyperbolic. 
%  so that 
%  all its periodic points are hyperbolic and so that each critical point of  $f$ is 
%  either 
%  \begin{enumerate}
%   \item[-] eventually mapped onto  a hyperbolic repelling periodic point, or 
%  \item[-] periodic, or 
%   \item[-] is in the basin of 
%  a hyperbolic periodic attractor.  
%  \end{enumerate}   
Then 
  \begin{enumerate}
  \item There exists a neighbourhood $\mathcal O$ of $f$ in $\mathcal A^{\underline \nu}_a$ 
and a smooth  function $\Psi_H\colon \mathcal O\to \R^{\nu'}$  so that 
$$\Psi_H(g)=\Psi_H(f) \iff g \mbox{ is real-hybrid conjugate to }f.$$ 
Here $\nu'=\nu+\nu_{noness-att}$ where $\nu$ is the number of critical points of $f$
and $\nu_{noness-att}$ is the number of periodic attractors of $f$ without 
critical points in their basin. 
\item There there exists a neighbourhood $\mathcal O$ of $f$ in $\mathcal A^{\underline \nu}_a$ 
and a smooth  function $\Psi_T\colon \mathcal O\to \R^{\nu_T}$   so that 
$$\Psi_T(g)=\Psi_T(f) \iff g \mbox{ is topologically conjugate to }f.$$ 
Here  $\nu_T=\nu-\zeta(f)$, where $\nu$ is the number of critical points of $f$
and $\zeta(f)$ is the maximal number of critical points in the basins of periodic 
attractors with pairwise disjoint infinite orbits. 
\end{enumerate}
\end{thm}
\begin{pf}  
%By assumption, each critical point is either eventually mapped
%onto a repelling periodic point, or 
Let us first define $\Psi_H$ as in assertion (1). Choose a set  $Atr$ of periodic attracting points, so 
that each periodic attracting orbit intersects $Atr$  in precisely one point. 
For each $p\in Atr$ which contains critical points in its basin, choose
one such critical point $c(p)$ and let $\Crit'_{at}(p)$ be the remaining critical 
points in the basin of $p$. Note that since $f,g$ are real analytic, there may not be 
any critical point in the basin of $p$. Also note that any real analytic map $f$
has at most a finite number of periodic attractors, see \cite{MMvS}. 

By assumption, for  each $c\in\Crit(f)$ exactly one of the following holds 
 \begin{enumerate}
 \item[(i)] there  exists $q_c>0$ so that $c'(c):=f^{q_c}(c)\in\Crit(f)$
 and so that $f^{k}(c)\notin\Crit(f)$ for $1\le k<q_c$;  

 \item[(ii)] there exists $1\le l_c < q_c$ so that $f^{q_c}(c)=
 f^{l_c}(c)$ and $f^k(c)\notin \Crit(f)$ for $1\le k \le q_c$;

 \item[(iii)]  $c$ has an infinite orbit and is in the basin of a  hyperbolic periodic attractor $p\in Atr$
so that $f^k(c)\notin \Crit(f)$ for all $k> 0$.   \end{enumerate}

For each $c$ as in  case (ii), we have by assumption that $Df^{q_c-l_c}(f^{l_c}(c))\ne 1$
where $f^{l_c}(c)$ is a periodic point of period $q_c-l_c$. 
Let us denote the set of critical points of $f$ for which (i), (ii) and (iii) holds by 
$\Crit_{ec}$, $\Crit_{ep}$ and $\Crit_{at}$. Note
that $\Crit_{at}=\cup_{p\in Atr} [\Crit'_{at}(p)\cup \{c(p)\}]$.

%Choose some subset $\Crit'_{at}\subset \Crit_{at}$
%so that each periodic attractor contains at most one critical point from $\Crit'_{at}$
%and for $c\in \Crit'_{at}$ let $\Crit'_{at}(c)$ be subset of $\Crit_{at}\setminus \Crit'_{at}$
%which are attracted the same periodic attractor as $c$. 

Let $\mathcal O$ be a neighbourhood of $f$ in $\mathcal A^{\underline \nu}_a$
so that  critical points $c_g$ of $g$ depends smoothly on $g\in \mathcal O$. In particular
there exists a critical point $c'_g(c)$ of $g$ near $f^{q_c}(c)$ for each $c\in \Crit_{ec}$.

 For $p\in Atr$,  let $r$ be the period of $p$ and let $\varphi_{p}\colon (\R,p) \to (\R,0)$ be the 
 holomorphic map so that $\varphi_{p}\circ f^{r}(z) = \lambda \varphi_{p}(z)$ near $p$,  see \cite{Mil}. 
(The map $\varphi_{p}$ is unique up to a multiplicative constant.)
For each critical point $c$ in the basin of $p$, let $n_c$ be 
so that $f^{n_c}(c)$ is in the immediate basin of $p$. 
Normalise $\varphi_p$ so that 
$\varphi_{p}(f^{n_c}(c(p)))=\pm 1$  where $c(p)$ is a {\lq}preferred{\rq} critical point 
in the basin of $p$ and  where the sign depends on whether $f^{n_c}(c(p))$ is to 
the left or to the right of $p$. 
Note that
$\varphi_p$ is defined by 
$$\varphi_p(z)=\lim_{n\to \infty}  (f^{nr}(z)-p_c)/\lambda^n$$
where $\lambda=Df^r(p)$ and that  $\varphi_p\circ f^r = \lambda \varphi_p$. 
Normalise $\varphi_{g,p}$ as above by $\varphi_{p,g}(g^{n_c}(c(p)))=\pm 1$.  
Now define 
$$\Psi_H\colon \mathcal O\to \C^{\nu'}$$ by 
\begin{equation} \label{eq:Psi-trans} 
\begin{array}{rl}  
&\Psi_H(g) =
 \left( ( g^{q_c}(c_g)-c'_g(c))_{c\in \Crit_{ec}} , (g^{q_c}(c_g)-g^{l_c}(c_g))_{c\in \Crit_{ep}}, \right.  \\
& \\ 
 &  \left. \qquad (Dg^{r_c}(p_c),(\varphi_{p,g}(g^{n(c(p))}(c(p)))-\varphi_{p,g}(g^{n(c')}(c')))_{c'\in \Crit_{at}(p)\setminus \{c(p)\}}) _{p\in Atr}  \right)\in \R^{\nu'}.
 \end{array}
 \end{equation} 
 Here if $p\in Atr$ has only at most one critical point in its basin,
%$c'\in \Crit_{at}(p) \setminus \{c(p)\}=\emptyset$ 
the term %corresponding term 
$$(\varphi_{p,g}(g^{n(c(p))}(c(p)))-\varphi_{p,g}(g^{n(c')}(c')))_{c'\in \Crit_{at}(p)\setminus \{c(p)\}}$$
should be ignored. 
That $\Psi_H(f)\in \R^{\nu'}$  follows from $\Crit=\Crit_{ec}\cup \Crit_{ep}\cup \Crit_{at}$ and since
 $\Crit_{at}=\cup_{p\in Atr} [\Crit'_{at}(p)\cup \{c(p)\}]$.

Let us now show that there exists a neighbourhood $\mathcal O$ of $f$ so that for any $g\in \mathcal O$ one
has that $\Psi_H(f)=\Psi_H(g)$ if and only if $f,g$ are hybrid equivalent. 
To see that $\Psi_H$ is constant on hybrid classes, first note that by definition the coordinates corresponding to 
$\Crit_{ec}$ and $\Crit_{ep}$ of $\Psi_H(f)$ are zero. 
%For $g$ to be conjugate to 
%$f$ it is necessary that those coordinates of $\Psi_H(g)$ are also zero. 
If $f$ has a (hyperbolic) periodic attractor $p$ then denote by $p_g$ the corresponding
periodic attractor for $g$.  Then for $f,g$ to be hybrid conjugate
it is necessary that the conjugacy $h$ is  holomorphic  map from the immediate basin 
of $f$ at $p$ to the immediate basin of $g$ at $p_g$. This implies that the 
multiplier of $f$ at $p$ is equal to the multiplier of $g$ at $p_g$. Let $H$
be the map $h$  in terms of the $\varphi_{p,f}$ and $\varphi_{p,g}$ linearising coordinates.
It follows that $H\colon (\C,0) \to (\C,0)$  is also univalent and since the only univalent map from $\C$ to $\C$
is linear, it follows that this (real) linear map is completely determined by 
the condition that $h(f^{n(c(p))}(c(p)))=g^{n(c(p))}(c(p))$ which implies that for the lift 
$H(\varphi_{g,c}(g^{n(c(p))}(c(p))))= \varphi_{g,c}(g^{n(c(p))}(c(p)))$. 
Here $c(p)$ is  the preferred critical point in the basin of $p$. 
Because  $\varphi$ is unique up to a multiplicative constant, 
but is normalised  so that $\varphi_{p,g}(g^{n_c}(c(p)))=\pm 1$ 
it follows that $H=id$. Hence  
$$\mathcal O\ni g\mapsto \left(\varphi_{p,g}(g^{n(c(p))}(c(p)))- \varphi_{p,g}(g^{n(c')}(c'))_{c'\in \Crit_{at}(p)\setminus \{c(p)\}}) _{p\in Atr}  \right)\in \C^{\nu'}
$$
is constant for all maps which are hybrid conjugate to $f$. 

Similarly, if $\Psi_H(f)=\Psi_H(g)$ and $g\in \mathcal O$ then all the critical relations are preserved. 
Moreover, there exists a qc-conjugacy $h$ between pruned-polynomial-like extensions
 $F\colon U\to U'$  and $G\colon U_g\to U_g'$ of $f$ and $g$,
whose dilatation vanishes on the basins of the periodic attractors and so that $\bar \partial h=0$ on $K_F$.
Notice that we may need to 
adjust the choice of boundaries of the basins of attraction of $G\colon U_g:=\E_g\cup B_g\to \E'_g\cup B_g'=:U_g'$ 
to obtain a conformal conjugacy on basins of attraction.   
Thus,   $F\colon U\to U'$  and $G\colon U_g\to U_g'$  are hybrid conjugate.

Let us now construct the analogous map $\Psi_T$ for the topological conjugacy classes. In that case, the conjugacy is more flexible inside the basin of periodic attractors.
However, $\Psi_T$ now needs to take account of critical relations within the basin of periodic attractors. 
To do this, choose a maximal subset of $\Crit^*_{at}$ of $\Crit_{at}$ so that the forward orbits of 
the critical orbit in $\hat \Crit_{at}$ are pairwise disjoint. Moreover, for each $c\in \Crit^*_{at}$
choose the subset of $\Crit_{at}(c)$ of  $\Crit_{at}$ of critical points $c'\ne c$ so that 
there exist $l_c,l_{c'}\ge 1$ so that  $f^{l_c}(c)=
 f^{l_{c'}}(c')$. 

Now define 
$$\Psi_T\colon \mathcal O\to \C^{\nu-\zeta(t)}$$ by 
\begin{equation} \label{eq:Psi-transT} 
\begin{array}{rl}  
&\Psi_T(g) =
 \left( ( g^{q_c}(c_g)-c'_g(c))_{c\in \Crit_{ec}} , (g^{q_c}(c_g)-g^{l_c}(c_g))_{c\in \Crit_{ep}}, \right.  \\
& \\ 
 &  \left. \qquad (g^{l_c}(c_g)-g^{l_{c'}}(c'_g))_{c\in \Crit^*_{at}, c\in \Crit_{at}(c)}   \right) \in \R^{\nu-\zeta(f)}.
 \end{array}
 \end{equation} 
 Here, if for some  $c\in \Crit^*_{at}$ the set $\Crit_{at}(c)$ is empty then that term should be ignored. 
% 
%\begin{enumerate}
% \item[(iii)] there exists $1\le l_c < q_c$ so that $f^{q_c}(c)=
% f^{l_c}(c)$ and $f^k(c)\notin \Crit(f)$ for $1\le k \le q_c$;
% \end{enumerate} semi
\end{pf}

\subsection{Conjugacy 
  classes of semi-hyperbolic  mappings are real analytic Banach manifolds}
  
  Let $\nu,\nu',\zeta(f)$ be as 
in Theorem~\ref{thm:definitonPsi}.  

\begin{thm}\label{thm:manifold-specialmaps} 
  Let $f\in\mathcal A^{\underline \nu}_a$ be semi-hyperbolic. 
%  be a  real analytic interval map so that 
%  all its periodic points are hyperbolic and so that each critical point of  $f$ is 
%  either 
%  \begin{enumerate}
%  \item[-] periodic, 
%  \item[-] eventually mapped onto  a hyperbolic periodic point or 
%  \item[-] is in the basin of 
%  a hyperbolic periodic attractor.  
%  \end{enumerate} 
  Then 
%   Let $F\colon U\to U'$
%  be a pruned polynomial-like extension of $f$.  Then
  \begin{enumerate}
%  \item  the hybrid conjugacy class of $F$ is 
%a complex-analytic   Banach manifold with complex codimension 
%$\nu'$.
\item    the real-hybrid conjugacy class of $f$ is 
a real analytic  Banach manifold in $\mathcal A^{\underline \nu}_a$ with real codimension  $\nu'$.
\item the topological conjugacy class of $f$ is 
a real analytic Banach manifold in $\mathcal A^{\underline \nu}_a$ with real codimension 
$\nu-\zeta(f)$.
\end{enumerate} 
\end{thm}
\begin{pf} To prove (1) it is enough to show that the map $\Psi_H\colon \mathcal O \to \R^{\nu'}$ defined
in equation (\ref{eq:Psi-trans}) (in the proof of  the previous theorem) has full rank. 
So take a unit vector $e_j \in \C^{\nu'}$ and show how to 
pick a family  $f_t$ through $f$ so that $\dfrac{d}{dt}\Psi_H(f_t)=e_j$. 
For convenience  we will choose $f_t(x)=f(x)+tv(x)$ where
$v$ is a polynomial function which vanishes at $c\in \Crit(f)$ of 
order $\ell_c$ and which, in order to ensure that $f_t\in \mathcal \A_a$, also vanishes at $\pm 1$.   
 This implies that the critical points $c_t$ of $f_t$ do not 
depend on $t$ and remain of the same order,  and so $v\in T_{f}\mathcal A^{\underline \nu}_a$. 
Note that  $\dfrac{d}{dt}f_t(x)=v(x)$.  Let us consider each coordinate of $\Psi_H$ in turn and show
that one can choose $v$ so that coordinate of $D\Psi_H(v)$ is non-zero while the 
other ones are zero. Define $\Crit_{ec},\Crit_{ep},\Crit_{att}$ as in the proof of Theorem~\ref{thm:definitonPsi}. 

\noindent 
{\bf Case 1: $c\in \Crit_{ec}$.} This corresponds to the first component of $\Psi_H$, see  
equation (\ref{eq:Psi-trans}), i.e. to  critical 
points which are eventually mapped onto other critical points. 
%For each $c\in\Crit(F_0)$, let $j_c>0$ be {\em minimal} so that 
%$F_0^{j_c}(c)\in\Crit(F_0)$. 
%For any $v\in T_{F_0}\mathcal{A}_a^{\nu}$,
%let $v_{c,k}=v(F_0^k(c_i))$. 
%Now, that $v\in E^h_{F_0}$ implies that
%there exists a qc vector field $\alpha$ on the interval so that 
%for $\alpha_{c,k}=\alpha(F_0^k(c))$ the following equations hold: 
%\textcolor{red}{CHECJK QQQ}
%$$\alpha_{c,k+1}  = DF_0(F_0^k(c))\alpha_{c,k} + v_{c,k}  \mbox{ for } 
%0\leq k <q_c, c \in \Crit(F_0). 
%$$
% Now $v\in E^h_{F_0}$ if and only if $\alpha_{c,q_c}=0$ for $c\in \Crit_{ec}$
%and $\alpha_{c,q_c}=\alpha_{c,l_c}$ for $c\in \Crit_{ep}$.  
% The previous recurrence equations imply that 
By the choice of $f_t$ we have $\frac{d}{dt}c_t\vert_{t=0}=0$
and a simple calculation shows 
$$\dfrac{f_t^{q_c}(c)}{dt}\vert_{t=0}   
=v(f^{q_c-1}(c))+ Df(f^{q_c-1}(c)) v(f^{q_c-2}(c)) + \dots + 
 Df^{q_c-1}(f(c))v(c) .$$
 It follows that if we choose the function $v$ so that 
 $v(c)\ne 0$ and $v(f^i(c))=0$ for $0<i\le q_c-1$ then 
 \begin{equation} \left. \dfrac{d}{dt} (f_t^{q_c}(c_t) - c_t) \right|_{t=0} \ne 0. \label{eq:trans132} \end{equation} 
 Thus, with this choice of $v$,  the coordinate of $D\Psi_H(v)$  corresponding
 to $c\in \Crit_{ec}$  will be non-zero, 
 and as we will see if we choose $v(f^i(\tilde c))=0$ for all $\tilde c\ne c$  in $\Crit_{ec}$ and finitely many $i\ge 0$,  
 then the other  coordinates of $D\Psi_H(v)$ will be zero. 
The inequality   (\ref{eq:trans132}) implies that $v$ is transversal to the hybrid class.  
 
 \noindent 
{\bf Case 2: $c\in \Crit_{ep}$.} 
Next we will consider the 2nd component of  $\Psi_H$ corresponding
to  critical points which are eventually mapped to a hyperbolic periodic point: 
$$\left. \dfrac{d}{dt} (f_t^{q_c}(c)-f_t^{l_c}(c_t))\right|_{t=0}
 = \left(v(f^{q_c-1}(c))+ \dots +   Df^{q_l-l_c-1} (  f^{l_c+1}(c) ) v(f^{l_c}(c_t)) \right) +  $$
 $$ \left. +  \left( Df^{q_c-l_c}(f^{l_c}(c)) - 1 \right) \dfrac{f_t^{l_c}(c_t)}{dt}\right|_{t=0} $$
where  
 $$\dfrac{f_t^{l_c}(c_t)}{dt}\vert_{t=0}   =v(f^{l_c-1}(c))+ Df(f^{l_c-1}(c)) v(f^{l_c-2}(c)) + \dots + 
 Df^{l_c-1}(f(c))v(c) .$$ 
From this  it follows that, as  $Df^{q_c-l_c}(f^{l_c}(c)) \ne 1$,  if we take $v(c)\ne 0$, $v(f^i(c))=0$ for all $i>0$ 
and $v(f^i(\tilde c))=0$ for all critical points $\tilde c\ne c$ and all $i\ge 0$, then the coordinate corresponding to 
$c\in  \Crit_{ep}$ in  $\dfrac{d}{dt}\Psi_H(f_t)$ will be non-zero whereas the others will be zero
(as we will see). 
Note that if $c\in \Crit_{ep}$  then by assumption $Df^{q_c-l_c}(f^{l_c}(c)) \ne 1$.
Hence there exists $p_{c,t}$ depending smoothly  on $t$ so that  $p_{c,0}=f^{q_c-l_c}(f^{l_c}(c))$ and so that $f_t^{q_l-l_c}(p_{c,t})=p_{c,t}$ for $t$ near $0$. Case 2 therefore follows from 
 
 \noindent 
{\bf Claim:}  Let $c\in \Crit_{ep}$. Then under the assumption that  $Df^{q_c-l_c}(f^{l_c}(c)) \ne 1$
we have that   $\dfrac{d}{dt} (f_t^{l_c}(c)- p_{c,t})\vert_{t=0}\ne 0$
 is equivalent to $\dfrac{d}{dt} (f_t^{q_c}(c)- f_t^{l_c}(c) ) \vert_{t=0} \ne 0$.
 
 \noindent 
 {\bf Proof of Claim:} 
$$\dfrac{d}{dt} (f_t^{q_c}(c)-f_t^{l_c}(c))\vert_{t=0}  = \dfrac{d}{dt} (f_t^{q_c-l_c}(c))\vert_{t=0} 
+  \left( Df^{q_c-l_c}(f^{l_c}(c)) - 1 \right) \dfrac{d}{dt} f_t^{l_c}(c_f)\vert_{t=0} $$
Differentiating the equality $f_t^{q_l-l_c}(p_{c,t})=p_{c,t}$  gives 
$$\dfrac{d}{dt} f_t^{q_c-l_c}(p_{c,0})\vert_{t=0} + Df^{q_c-l_c}(f^{l_c}(c)) \dfrac{dp_{c,t}}{dt}\vert_{t=0} = \dfrac{dp_{c,t}}{dt}\vert_{t=0} $$ 
and so
$$\begin{array}{rl} 
&(1-Df^{q_c-l_c}(f^{l_c}(c))) \dfrac{dp_{c,t}}{dt} = \frac{d}{dt} f_t^{q_c-l_c}(p_{c,0})\vert_{t=0}  =\\
&\qquad\qquad \dfrac{d}{dt} (f_t^{q_c}(c)-f_t^{l_c}(c))\vert_{t=0} -\left( Df^{q_c-l_c}(f^{l_c}(c)) - 1 \right) \dfrac{d}{dt} f_t^{l_c}(c_f)\vert_{t=0}.  \end{array} $$
Hence 
$$\left( Df^{q_c-l_c}(f^{l_c}(c)) - 1 \right)  
\left( \dfrac{d}{dt} \left(f_t^{l_c}(c)- p_{c,t}\right)\big\vert_{t=0} \right) 
= \dfrac{d}{dt} (f_t^{q_c}(c)-f_t^{l_c}(c))\vert_{t=0}$$
This implies the claim. 
%So if $Df^{q_c-l_c}(f^{l_c}(c)) \ne 1$ then 
%  $\dfrac{dx}{dt}\vert_{t=0}  \ne \dfrac{df_t^{l_c}(c_f)}{dt}\vert_{t=0}$
% is equivalent to $\dfrac{d}{dt} (f_t^{q_c}(c)-f_t^{l_c}(c))\vert_{t=0}\ne 0$. 
 \checkmark

 \noindent 
{\bf Case 3: $p\in Atr$.}  Let us now consider $p\in Atr$ and assume that $Df^r(p)\ne 1$.
 Note that 
 $$\dfrac{d}{dt} f_t^r(p_t)= v(f^{r-1}(p)) + Df(f^{r-2}(p))v(f^{r-2}(p)) + 
\dots +  Df^{r-1}(f(p))v(f(p)) + Df^r(p) \dfrac{d}{dt} p_t$$
is equal to $\dfrac{d}{dt} p_t$.
It follows that $v(p)=v(f(p))=\dots=v(f^{r-1}(p))=0$ implies $ \dfrac{d}{dt} p_t=0$
and similarly $\dfrac{d}{dt} f^i(p_t)=0$ for $i\ge 0$. 
Hence
$$\dfrac{d}{dt} Df^r_t(p_t)= Df^r(p) \left( \frac{v'(f^{r-1}(p))}{Df(f^{r-1}(p))} +  \frac{v'(f^{r-2}(p))}{Df(f^{r-2}(p))} +  \dots + \frac{v'(p)}{Df(p)}\right).$$
Taking $v'(p)\ne 0$ and $v'(f^i(p))=0$ for $0<i<r$ gives $\dfrac{d}{dt} Df^r_t(p_t)\ne 0$. This is therefore
equal to third term component of $D\Psi_H(v)$. 

 \noindent 
{\bf Case 4: $p\in Atr$, $c'\in \Crit_{at}(p)\setminus \{c(p)\}$.}  
Let us finally consider the term 
\begin{equation}  A:= \left. \dfrac{d}{dt} [ \varphi_{p_t,f_t}(f_t^{n(c(p))}(c(p))) -  \varphi_{p_t,f_t}(f_t^{n(c'(p))}(c'(p)))] \right|_{t=0}
\label{eq:varphitrans} 
\end{equation} 
which is the final component of $D\Phi_H(v)$,  see definition (\ref{eq:Psi-trans}). This term  is
related to the position of orbits of the critical points in the basin 
of the periodic attractor $p$ (which contains $c$ in its basin) and here  $\varphi_{p_t,f_t}$ is the linearisation 
at $p_t$ and so 
$\varphi_{p_t,f_t}(z)=\lim_{k\to \infty} f_t^{rk}(z)/a_t^k$ where $a_t=Df^r_t(p_t)$,
and $|a|<1$. Let us choose a real analytic function $v$ so that $v(x)=v'(x)=0$ for each $x$ in the forward orbit  of $p$. This then implies that $\frac{d}{dt}p_t \vert_{t=0}=p$ and $\frac{d}{dt}a_t\vert_{t=0}=a$. 
Note that  
\begin{equation} \dfrac{d}{dt} f^k_t(z)= v(f^{k-1}(z)) + Df(f^{k-1}(z)) v(f^{k-2}(z)) + \dots + Df^k(z)v(z). \label{eq:speediterate} \end{equation}
Let us take $|v|<\epsilon$ along forward iterates of $f_t^{n(c(p))}(c(p))$ and $f_t^{n(c'(p))}(c'(p))$.
This is possible as $v(x)=v'(x)=0$ for each $x$ in the forward orbit  of $p$, by taking 
$v$ zero in a finite number of forward iterates of  $f_t^{n(c(p))}(c(p))$ and $f_t^{n(c'(p))}(c'(p))$.
Since $\varphi(z)=\lim_{k\to \infty} f^{rk}(z)/a^k$ and using equation (\ref{eq:speediterate}) it follows that 
\begin{equation} - \delta <  \left.  \dfrac{d}{dt} [ \varphi_{p_t,f_t}(f^{n(c(p))}(c(p))) -  \varphi_{p_t,f_t}(f^{n(c')}(c'))] \right|_{t=0} < \delta
\label{eq:varphitrans2} 
\end{equation} 
where $\delta>0$ is determined by $\epsilon$ and the multiplier at the periodic point $p$; in fact, $\delta \approx \epsilon/(1-|a|)$. 
Hence (\ref{eq:varphitrans}) is, by the chain rule, equal to some $s\in [-\delta,\delta]$ plus 
\begin{equation*}  B:=  \left[ D\varphi_{p,f}(f^{n(c(p))}(c(p)))   \left. \dfrac{d}{dt} f_t^{n(c(p))}(c(p))\right|_{t=0} 
- D\varphi_{p,f}(f^{n(c')}(c'))  \left. \dfrac{d}{dt} f_t^{n(c')}(c')\right|_{t=0} \right] 
\label{eq:varphitrans3} 
\end{equation*} 
Since $D\varphi_{p,f} \ne 0$, one can choose $A\ne 0$ and $|B|>C\delta$  by choosing  
$v(f(c))=\dots=v(f^{n(c(p))}(c(p)))=0$ 
and $v(f(c'))=\dots=v(f^{n(c'(p))}(c'(p)))=0$ and $v(c),v(c')\ne 0$ appropriately. 
Indeed, then 
$$B= D\varphi_{p,f}(f^{n(c(p))}(c(p))) Df^{n(c(p))} v(c) -  D\varphi_{p,f}(f^{n(c')}(c'(p))) Df^{n(c')} v(c') .$$
Since the factors in front of $v(c)$ and $v(c')$ are non-zero, 
we can choose $v(c),v(c')$ so that 
$$|A| = |B|- \delta > 0.$$ 
%
% which means that $\alpha_{c,q_c}=\alpha_{c,l_c}$ is equivalent to o
%$$\begin{array}{cc}  &v_{c,l_c-1}+ DF_0(F_0^{l_c-1}(c)) v_{c,l_c-2} + \dots + 
% DF_0^{l_c-1}(F_0(c))v_{c,0}=\alpha_{c,l_c} = 
% \\ 
%&  
%\dfrac{ \left(v_{c,q_c-1}+ \dots +   DF_0^{q_l-l_c-1} (  F_0^{l_c+1}(c) ) v_{c,l_c}\right)}{ 1-DF_0^{q_l-l_c-1}}.
%\end{array} $$
%
%
%
%
% 
%For $c\in $\Crit_{ec}$
%the homogeneous system of linear
%equations:
%$$
%\begin{array}{rll} \alpha_{c,k+1}-DF_0(F_0^k(c))\alpha_{c,k} &= v_{c,k}  & \mbox{for}\;
%1\leq k <q_c, c \in \Crit(F_0) \\
%& \\
% \alpha_{c,1} &= v_{c,0} &\mbox{for } c\in \Crit(F_0) \end{array} $$
%%$$\mbox{and}\; v_{c,0}=\alpha_{c,1},\;\mbox{where}\; 
%%c\in\Crit(F_0)$$
%has a unique solution, $(\alpha_{c,k})_{c\in\Crit(F_0), 1\leq
%  k\leq q_c}$ because $DF_0(F_0^k(c))\ne 0$ for each $1\le k<q_c$. 
%Observe that 
%$$L_1: (v_{c,k})_{c\in\Crit(F_0), 0\leq
%  k<q_c}\mapsto (\alpha_{c,k})_{c\in\Crit(F_0), 1\leq
% k\leq q_c}, \;\mbox{and}$$
%$$L_2:v\mapsto \left( (\alpha_{c,q_c})_{c\in\Crit(F_0)_{ec}}, (\alpha_{c,q_c}-\alpha_{c,l_c})_{c\in\Crit(F_0)_{ep}} \right) \in \C^\nu
%$$
%are continuous linear mappings and that 
%$v\in E^h_{F_0}$ if and only if $v\in\mathrm{ker}(L_2)$. 
In particular, it follows that if $v$ so that the values of $v$ along the forward orbits of 
$c$ and $c'$  as above, then  the fourth component of $\Psi_H(v)$ is non-zero. 

\noindent 
Combining the previous cases, we see that for each coordinate of 
$\Psi_H$ one can choose $v$ 
so that coordinate of $D\Psi_H(v)$ is non-zero while 
the others are zero.  Hence $D\Psi$  has maximal rank and 
therefore Assertion (1)  follows from the  Implicit Function Theorem. 

The proof of Assertion (2) is entirely analogous.  
\end{pf}

\begin{rem} Of course one could  
denote $$\dfrac{f_t^{n}(c)}{dt}\vert_{t=0}   
=v(f^{n-1}(c))+ Df(f^{n-1}(c)) v(f^{n-2}(c)) + \dots + 
 Df^{n-1}(f(c))v(c) $$
by $v^{n}(c)$ and obtain shorter expressions. The full expressions
show how to choose $v$ so that the required assumptions are satisfied.  
\end{rem} 
\subsection{The simple parabolic case} \label{subsec:para-transv} 

\begin{thm} \label{thm:parabolicmanifold} 
Let $a>0$ and $f\in \mathcal A^{\underline \nu}_a$. Assume that 
$f$ has a simple parabolic periodic point $a_0$ of period $n$. 
Then  there exists a real analytic Banach manifold $\mathcal P_{f,a}^{\underline \nu}$ in $\mathcal A^{\underline \nu}_a$
consisting  of maps in $ \mathcal A^{\underline \nu}$
so that each map $g\in \mathcal P_{f,a}^{\underline \nu}$ has a unique parabolic periodic point 
of period $n$ near $a$ and of the same type as $a_0$, and if each $g\in   \mathcal A^{\underline \nu}_a$
with this property is in $\mathcal P_{f,a}^{\underline \nu}$.

If $a_0$ is of saddle-node or period-doubling type, then $\mathcal P_{f,a}^{\underline \nu}$ has 
codimension one, and if it  is of pitchfork type it has codimension two. 

If $f$ has $k_0$ parabolic periodic points of saddle-node or period-doubling type and 
$k_1$ periodic points of pitchfork type, then the analogous statement holds
and the codimension of $\mathcal P_{f,a}^{\underline \nu}$ is equal to $k_0+2k_1$. 
\end{thm} 
\begin{pf} For any vector field $v$ define $v^n=\frac{d}{dt} (f+tv)^n|_{t=0}$. 
Assume that  $a_0$ is a simple parabolic 
periodic point of period $n$.

If $Df^n(a_0)=1$ and $a_0$ is of saddle-node type  
then $f^n(x)=a_0+(x-a_0) + \alpha(x-a_0)^2+O(x-a_0)^3$ with $\alpha\ne 0$. 
In this case, choose a vector field $v$ so that $v^n(a_0)\ne 0$ (it is easy to 
see that this possible).  Let $B_\epsilon(a_0)$ be a small neighbourhood of $a_0$ and 
define $\Psi$ on a neighbourhood $\mathcal O\times (-\epsilon,\epsilon) \times B_\epsilon(a_0)$ of $(f,0,a_0)$ by
$$\Psi(g,t,x)=((g+tv)^n(x)-x, D(g+tv)^n(x)-1).$$ 
Here $\mathcal O$ is an open subset  of the real analytic Banach manifold $\mathcal A^{\underline \nu}_a$.
  Let us show that for each $g$ near $f$ there exist a unique $t$ and $x$
so that   $\Psi(g,t,x)=(0,0)$. Note that this implies that 
$x$ is a periodic point of $g+tv$ with multiplier $1$ which is 
of saddle-node type if $g$ is close to $f$.
The partial derivative matrix
of $\Psi$ w.r.t. $t$ and $x$ at $(g,t,x)=(f,0,a_0)$ is equal to 
$$\left( \begin{array}{cc} v^n(a)  & 0  \\ 0 & 2\alpha \end{array}\right).$$
Since $v^n(a_0)\ne 0$ and $\alpha\ne 0$, the Implicit Function Theorem gives for each $g\in \mathcal O$
the existence of  a unique solution $t(g),x(g)$ 
of the equation $\Psi(g,t(g),x(g))=(0,0)$. Here $t(g)$ and $x(g)$ depend 
smoothly on $g$. 
Thus we obtain that there exists a codimension-one manifold  with the desired properties, namely 
the space of maps of the form $g+t(g)v$. In fact, if we 
let $T\mathcal P_{f,a}$ be the linear space of vector fields $\tilde v$ so that $\tilde v^n(a_0)=0$
then $\mathcal P_{f,a}$ can be locally written 
as $\{f+\tilde v+\tilde t(\tilde v)v; \tilde v\in T\mathcal P_{f,a}\}$ where $\tilde t$ 
is an analytic function.  It is easy to see that $D\tilde t(0)=0$ and so
 the tangent space of $\mathcal P_{f,a}$  is equal to $T\mathcal P_{f,a}$. 

If  $Df^n(a_0)=-1$  and $a_0$ is of period-doubling type 
then $f^n(x)=a_0 - (x-a_0) + \alpha(x-a_0)^2+O(x-a_0)^3$.
In this case we choose a vector field $v$ so that  
$$2 Dv^n(a_0) +v^n(a_0) D^2f^n( a_0)\ne 0$$
and argue as in the saddle-node case. Here the above condition on $v$ follows
from the calculation in Remark~\ref{rem:perioddoubling}.

If $Df^n(a_0)=1$ and $f$ is attracting from both sides
then $f^n(x)=a_0+(x-a_0) + \alpha(x-a_0)^3+O(x-a_0)^4$ with $\alpha< 0$. 
In this case let $v$ and $w$ be independent vector fields so that 
$v^n(a_0)\ne 0$, $Dv^n(a_0)=0$, $w^n(a_0)=0$ and $Dw^n(a_0)\ne 0$. (It is easy to see that this is possible.)
Next define 
$$\Psi(g,t,s,x)=((g+sv+tw)^n(x)-x, D(g+sv+tw)^n(x)-1,D^2(g+sv+tw)).$$ 
The partial derivative matrix 
of $\Psi$ w.r.t. $s,t,x$ at $(g,s,t,x)=(f,0,0,a_0)$ is equal to 
$$\left( \begin{array}{ccc} v^n(a_0)  & 0 & 0  \\ 0 & Dw^n(a_0)  & 0 \\
0 & 0 & 6 \alpha  \end{array}\right).$$
By construction this matrix is invertible, and so there exists a unique $s(g),t(g),x(g)$ so that 
$\Psi(g,s(g),t(g),x(g))=(0,0,a_0)$. It follows that each map of the form $g+s(g)v+t(g)w$ 
has a pitchfork parabolic point at $x(g)$. Thus, using the 
Implicit Function Theorem we obtain a codimension-two manifold
of maps with a parabolic periodic point of pitchfork type. The tangent space of this manifold
is the space of vector fields $\tilde v$ so that $\tilde v^n(a_0)=D\tilde v^n(a_0)=0$. 
\end{pf} 

\begin{rem} A similar analysis can be found on page 510 \cite{ALM}. 
A more general result in this setting for arbitrary (non-degenerate) 
parabolic periodic points can be found in the Main Theorem  in \cite{LSvS2}.  
Note that if $a_0$ is of pitchfork type then it is  degenerate in the sense of \cite{LSvS2}.
\end{rem} 

\begin{rem} One has analogous results within the space of polynomials, rational maps, maps of finite type
or more general families of maps, see  \cite{Epstein,LSvS1,LSvS2,LSvS3}. However, here the situation is simpler
as it is easier to show that the map $\Psi$ from above is a submersion in the present case. This is because
here we can construct the appropriate transversal vector fields to the manifold by hand because the space 
$\mathcal A^{\underline \nu}_a$ is much larger than the space of rational maps. 
\end{rem}

\begin{cor} \label{cor:semihyperbolic}  The conclusion of Theorem~\ref{thm:manifold-specialmaps}  also holds if we replace 
 {\em hyperbolic periodic point} throughout by {\em hyperbolic or simple parabolic 
periodic point}. 
\end{cor} 
\begin{pf} 
This follows by combining the proofs of the previous theorem with that of Theorem~\ref{thm:manifold-specialmaps}.
\end{pf}

\section{Mating pruned polynomial-like mappings}\label{sec:mating}

%\textcolor{blue}{$F_\bullet:U_{F_\bullet}\to U_{F_\bullet}' $  is used so that we have a qc map
%on a neighbourhood of $K_F$. Is there an easier way to get this by pulling back the map near the pruning points???} 

The following theorem is the analogue of the mating theorems  \cite{DH} and \cite{Lyubich}
for polynomial-like mappings in our context of  pruned  polynomial-like mappings. 
 
\begin{thm}
  \label{thm:mating tpl}
Let $F:U_F\to U_F'$ and $G:U_G\to U'_G$ be pruned polynomial-like
mappings with $Q(F)=Q(G)$. %$F,G\in \mathcal{PL}_{Q}$.
%Moreover, assume that $F$ and $G$ are prunings of  \textcolor{red}{QQQ $F_e$}
%$F_\bullet:U_{F_\bullet}\to U_{F_\bullet}' $ and
%$G:U_{G}\to U_{G}',$
%respectively with $U_F' \Subset U_{F_\bullet}'$ and $U_G'\Subset U_{G}'$ with 
%$Q(F_\bullet)=Q(G)$.
%so that \textcolor{red}{$Q(F_\bullet)=Q(G)$.} %\in\mathcal{PL}_{Q^*}.$}
Then there exists a unique pruned polynomial-like mapping $\widetilde F:U_{\tilde F}\to U_{\tilde F}'$,
so that 
\begin{enumerate}
\item 
$\widetilde F$ is hybrid conjugate  to $F$. 
\item  $\widetilde F$ and $G$ have the same external mappings. 
\end{enumerate} 
\end{thm}

%Here we say $$\widetilde F\in \mathcal T_{F^*}$$ if and only if
% $\widetilde F$ and $F^*$ are topologically conjugate.  The map $\widetilde F$ is topologically conjugate to $F$ and has the same 
% external mapping as $G^*$, 
 We call   $\widetilde F$ a {\em mating} of $F$ and $G$.

%\medskip
%\noindent\textit{Remark.
%  We need to pass to a further pruning in this proposition to obtain
%a mapping $\widetilde F$ which is analytic on a neighbourhood of its pruned filled 
%Julia set via the Measurable Riemann Mapping Theorem. To obtain this
%conclusion, we require that $\widetilde F$ leaves invariant a Beltrami
%on a neighbourhood of its filled Julia set, and we only have such a
%neighbourhood after passing to a further pruning. It is also worth noting
%that in the presence of attracting or parabolic basins, a further
%pruning of this form may not exist. Further prunings always exist when all cycles are repelling.
%}
\begin{pf}
Since $Q(F)=Q(G)$ we obtain by Proposition~\ref{prop:qcconjpartial}
a $\R$-symmetric quasiconformal map $H\colon (U_F\cup U'_F)\setminus K_F \to (U_G\cup U'_G)\setminus K_G$
which conjugates  $F\colon U_F\to U'_F$
and $G\colon U_G\to U_G'$ wherever this is defined. 
Note that $K_F\cap \partial U_F$ 
consists of a finite number of non-real preimages of the pruning points. 
Therefore we can extend $H$ near each of these points $K_F\cap \partial U_F$  by taking  preimages under 
$F$ thus obtaining a qc conjugacy $H\colon U_F^*\setminus K_F\to U_G^*\setminus K_G$ between $F$ and $G$
(on these sets) where $U^*_F$, $U^*_G$ 
are neighbourhood of $K_F$ resp. $K_G$.  Here we use that  $F,G$ 
extend holomorphically to neighbourhoods of the closures of their domains. 
%\colon U_F\to U'_F$ extends holomorphically to a neighbourhood of $\overline U_F$.
 Since $\bar \C \setminus (\partial U_F\cup \partial U_F')$
and   $\bar \C \setminus (\partial U_G\cup \partial U_G')$ are quasi disks (and one choose $U^*_F,U_G^*$
conveniently),  $H$ has a qc extension  $H\colon \bar \C\setminus K_F \to \bar \C \setminus K_G$.

%
%~\ref{prop:qc-conjucy-external} 
%a $\partial \D$-symmetric quasiconformal map
% $H:(\V'_{\F},\V_{\F},B_{\F},B'_{\F})\rightarrow (\V'_{\G},\V_{\G},B_{\G},B'_{\G})$ which is 
% a conjugacy between $\F_X\colon \V_{\F}\cup B_{\F} \to  \V'_{\F}\cup B'_{\F}$
% and $\G_X\colon \V_{\G}\cup B_{\G} \to  \V'_{\G}\cup B'_{\G}$.
%So $H\circ \F_X = \G_X \circ H$ on the domain of $H$.
%Extend $H$ to a $\partial \D$-symmetric $H\colon \C\to \C$.
%\trevor{\textcolor{red}{Do we need to use $\F_X$ here?}}

%\textcolor{red}{Let us first assume that all periodic points of $F$ and $G$ are repelling.}
%Since $G\in\mathcal{PL}_Q,$
%we have that the external mapping
%$\F \colon V_{\hat F} \to  V'_{\hat F}$ and the ray structure 
%%$F_{\phi_F}\colon V_{\phi_F}\to  V'_{\phi_F}$ 
%$\Gamma_{\hat F}$ associated to $F\colon U_F\to U'_F$ is partially 
%conjugate to the external map $\G\colon V_{\hat G} \to V'_{\hat G}$
%together with  $\Gamma_{\hat G}$ associated to $G\colon U_G\to U_G'$. 
%By choosing a $\partial \D$-symmetric qc homeomorphism $H_0\colon (V_{\F},V'_{\F},\Gamma_{\F})\to
% (V_{\G},V'_{\G},\Gamma_{\G})$, and then pulling back this conjugacy
% we obtain a qc conjugacy $H\colon (V_{\F},V'_{\F},\Gamma_{\F})\to (V_{\G},V'_{\G},\Gamma_{\G})$
%between the external mappings 
%$\F:V_{\F}\to V'_{\F}$ and 
%$\G:V_{\G}\to V'_{\G}$. Extend $H$ to a $\partial \D$-symmetric $H\colon \C\to \C$. 
%%   Let $\hat A=\psi_F^{-1}(A)$ and $\hat A_0=\psi_F^{-1}(A_0).$

Let us consider the standard conformal structure $\sigma$ on $\bar \C$.
%Pulling it back by the Riemann mapping,
%$$\psi_{G}\colon \C\setminus \bar \D \to \C\setminus K(G),$$  we again obtain the standard conformal
%structure $\sigma$ on $\mathbb C\setminus \overline{\mathbb D}$.
%Reflect it about the unit circle, and push it forward by $H$ to obtain 
Then $\nu=H^*\sigma$ is a conformal structure on $\bar \C \setminus K_F$,
which is invariant under $F$ on its domain. %G_{\phi_F}$.
%Pushing $\nu$ forward by the Riemann mapping, 
%$$\psi_{F}\colon \C\setminus \bar \D \to \C\setminus K(F)$$ 
%we obtain a conformal structure 
%$(\psi_{F})_*\nu$ on
%$\mathbb C\setminus K(F)$. %\textcolor{red}{Or should this be $(\psi_{F})_*\nu$?} 
We extend $\nu$ to $\C$ by the standard conformal structure $\sigma$.
Thus we obtain a conformal structure $\mu$ on $\bar{\mathbb{C}},$
which is invariant under $F$ on $U_{F}^*\setminus K_F$. 
%(Note that by assumption $U'_{F}$ is compactly contained in $U_{F}$.)
We have that $\mu=H^*\sigma$ on $\C\setminus K_F$.

Applying the Measurable Riemann Mapping Theorem, see \cite{AB}, we
straighten $\mu$, and obtain a qc homeomorphism $H_\mu\colon \C\to \C$ 
which is $\R$-symmetric and sends the ellipse field defined by $\mu$
to the standard complex structure.  
%We can and will normalise $H_\mu$ so that $H_\mu(-1)= -1$, $H_\mu(1)=1$.
%By definition $H_\mu=H$. 
%which implies that $\widetilde F(\partial I)\subset \partial I$. 
%\trevor{\textcolor{red}{QQQQQQ CHECK}}
It follows that $\widetilde F:= H_\mu \circ F \circ H_\mu^{-1}$ is  
$\R$-symmetric conformal pruned polynomial-like mapping on $U_{\tilde F}=H_\mu U_{F}$.
%  which is qc conjugate to $F$
%on this set. 
We can and will normalise $H_\mu$ so that $H_\mu(-1)= -1$, $H_\mu(1)=1$
which implies that $\widetilde F(\partial I)\subset \partial I$.  Since $\mu=0$ on $K_F$, $\bar \partial H_\mu=0$ a.e. on $K_F$ and  it follows that $\tilde F$ is hybrid conjugate to $F$. This proves property (1).

Note that by construction $H\circ H_\mu^{-1}\colon \bar \C \setminus K_{\widetilde F} \to 
\bar \C \setminus K_G$ is conformal.  Denote by 
$$\phi_F \colon  \overline \C
\setminus K_{F,X_F} \to \overline \C \setminus \overline{ \mathbb D}\mbox{ and }
\phi_{G} \colon  \overline \C
\setminus K_{G,X_G} \to \overline \C \setminus \overline{ \mathbb D}$$  
the normalized Riemann mappings, then 
$\hat H = \phi_{G} \circ H\circ \phi_F^{-1}$ is a qc-conjugation between $\hat F_X$ and $\hat G_X$
near $\partial \D$, i.e.  $\hat  H\circ \F_X =\G_X\circ \hat H$
since $\hat F_X = \phi_F \circ F \circ \phi_F^{-1}$ and $\hat G_X = \phi_G \circ G \circ \phi_G^{-1}$.
Hence $\Psi\colon  \bar \C \setminus K_{\widetilde F}\to \bar \C\setminus \bar \D$
%=\hat H\circ \phi_{F}\circ H_\mu^{-1}
%\colon U_{\tilde F} \setminus K_{\widetilde F}\to \C\setminus \bar \D$
%is conformal, maps $-1$ to $-1$  (and $1$ to $1$)  (and so 
%agrees with $\phi_{\tilde F}$). This holds since 
defined by $\Psi:=\hat H\circ \phi_{F}\circ H_\mu^{-1}= \phi_{G} \circ H\circ \phi_F^{-1} \circ \phi_F \circ H_\mu^{-1}=
\phi_G \circ H \circ H_\mu^{-1}$ is conformal. 
Due  to the normalisations, $\Psi$ therefore agrees with $\phi_{\tilde F}$.
 Moreover
 $\Psi \circ \tilde F \circ \Psi^{-1}= (\hat H\circ  \phi_{F}\circ H_\mu^{-1}) \circ (H_\mu \circ F \circ H_\mu^{-1}) \circ (H_\mu \circ \phi_{F}^{-1} \circ \hat H^{-1}) = \hat H \circ (\phi_{F} \circ F \circ \phi_{F}^{-1}) \circ \hat H^{-1}=
\hat  H\circ \F_X \circ \hat H^{-1}=\G_X$, proving property (2).  

Let us  now assume that $\tilde F_1,\tilde F_2$ are two such maps. 
Then $\tilde F_1,\tilde F_2$ are qc-conjugate and by the remark
below this proof the conjugacy has the property that $\bar \partial H=0$
on $K(\tilde F_1)$. At the same the external maps of $\tilde F_1,\tilde F_2$
are the same, and so $\tilde F_1,\tilde F_2$ are conformally conjugate
outside their filled Julia sets. By construction these maps agree
on their filled Julia sets, and so $\tilde F_1,\tilde F_2$ are 
conformally conjugate. Since the conjugacy send $\pm 1$ to $\pm 1$
it follows that $\tilde F_1=\tilde F_2$.
\end{pf}

%\noindent 
%\textit{Remark.}  The reason we assume in  Theorem \ref{thm:mating tpl}
%that $F$ is the restriction of a pruned polynomial-like mapping 
%$F_\bullet$ is that this allows us to obtain a conformal structure that is invariant on a neighbourhood of a filled pruned Julia set of $F$. Without this condition, we would not know that the mapping $\tilde F$
%in the conclusion of Theorem~\ref{thm:mating tpl}, given by the Measurable
%Riemann Mapping Theorem is analytic.

\subsection{Hybrid classes are locally conformally equivalent} 

Using the previous mating result we obtain:

\begin{thm}
  \label{thm:tpl hybrid classes}
 % Let $f,g\in \mathcal A^{\underline \nu}_a$ and 
Suppose that $F_0:U_{F_0}\to U'_{F_0}$ and $G_0:U_{G_0}\to U_{G_0}'$ are pruned
polynomial-like mappings %extensions 
%of $f,g$, so that the domains of $F,G$ are 
%compactly contained in $\Omega_a$ and 
so that $Q(F_0)=Q(G_0)$.
Then 
\begin{itemize}
\item the hybrid classes of $F_0$ and $G_0$ 
are homeomorphic;
\item  if the hybrid class of $G_0$ has an analytic structure, then the one for $F_0$ also 
has an analytic structure.
\end{itemize} 
%
%Moreover, if the  real-hybrid class of $g$ is a real analytic manifold, then the real-hybrid class for $f$ is 
%an immersed analytic submanifold of $\mathcal A^{\underline \nu}_a$. 
\end{thm}

\begin{pf} 
Compare \cite[Lemma 4.3]{Lyubich}. We recall the following:

\medskip
\noindent\textbf{Fact.} If $F_\lambda = H_\lambda \circ F_0 \circ H_\lambda^{-1},$ $\lambda\in\Lambda\subset \mathbb C$ is a family of holomorphic mappings and $H_\lambda$ is a holomorphic family of qc mappings, then $F_\lambda$ is holomorphic in $\lambda$.

\medskip

By the Mating Theorem,  Theorem~\ref{thm:mating tpl},
we have an invertible mapping
$\psi:\mathcal H_{G_0}\to \mathcal H_{F_0}$. Let us 
prove that it is continuous (the proof that $\psi^{-1}$
is continuous is the same). 

%\trevor{\textcolor{red}{To prove continuity, it is enough to prove that $\psi$ restricted 
%to $\mathcal T_{G_0}\cap \mathcal  A^{\underline \nu}_a$ is continuous 
%where $a$ is so that the domain of pruned polynomial-like mapping
%$G_0$ is compactly contained in $\Omega_{a}$.}}  
Note that there exists a neighbourhood 
$\mathcal U$ of $G_{0}$   in the Caratheodory topology defined on 
the space of pruned polynomial-like mappings (see Definition~\ref{CaratheodoryTopology}), 
 % in $\mathcal A^{\underline \nu}_a$  for some 
so that  for each $G\in \mathcal U$ in the hybrid class of $G_0$  
%can be considered as pruned polynomial-like map 
%so that 
%there exists a pruned polynomial like mapping $G$ so that 
%$Q(G)=Q(G_{0})$ and so that 
%Moreover, for all  $g\in \mathcal U$ 
we
can select holomorphically moving domains
$U_{G}, U'_{G}$ so that
$G:U_{G}\to U'_{G}$ is a pruned  polynomial-like mapping and $Q(G)=Q(G_0)$. 
This holomorphic motion defines a family of
conformal structures on $U'_{G_{0}}\setminus U_{G_{0}}$.
More precisely, 
%From the analytic dependence in the Measurable Riemann Mapping Theorem,
we obtain  a holomorphic motion %s %have that
$H_{G}:(U_{G_{0}},U_{G_{0}}')\to (U_G,U'_{G})$ 
of  $\partial U_{G_{0}}\cup\Gamma_{G_{0}}$ over $\mathcal U$
conjugating 
$G_{0}:U_{G_{0}}\to U_{G_{0}}'$ with  $G:U_{G}\to U'_{G}$
on $\partial U_{G_{0}}\cup\Gamma_{G_{0}}$. 

Let us start with the standard conformal structure on
$\mathbb C \setminus U_{G}$.  Pulling it back by $H_{G}$, 
we obtain a conformal structure $H_{G}^*\sigma$ on
$\mathbb C\setminus U'_{G_0}$. 
Since the external mappings of $F_0$ and $G_{0}$
are qc conjugate, there exists a qc mapping $H:U_{F_0}\to U_{G_0}$ 
that conjugates $F_0$ with $G_{0}$ on $\partial U_{G_0}\cup\Gamma_{G_0}$.
Now, $H^*H_{G_{0}}^*\sigma$ is a conformal structure
defined on $U_{F_0}'\setminus U_{F_0}$.
Pulling it back by $F_0^n, n=0,1,2,\dots,$ 
we obtain an invariant
Beltrami differential defined on $U_{F_0}'\setminus K_{F_0},$
and extending this
to $\mathbb C$ by $\sigma$, we 
obtain a Beltrami differential $\nu_{G}$,
which is invariant on a neighbourhood of $K_{F_0}$,
and which depends holomorphically  on $G$.
By the Measurable Riemann Mapping Theorem and the fact above
we obtain
a family $F_{G}:U_{F_{G}}\to U'_{F_{G}}$ of
pruned polynomial-like mappings in $\mathcal H_{F_0},$
which depends continuously on $G$.
So if $G_\lambda$ is a holomorphic family of mappings of pruned polynomial-like mappings
then % $\mathcal H_{G_0}\cap \mathcal A^{\underline \nu}_a$,  then  
$F_{G_\lambda}$ also depends analytically on $\lambda$. 
% On the other hand, the filled Julia $K(F_0)$ varies continuously with respect to
% $\tilde F_0 \in\mathcal T_{F_0}$. Hence the normalized Riemann mapping
% $\psi_{\tilde F}:\mathbb C\setminus\overline{\mathbb D}\to \mathbb C\setminus K(\tilde F_0)$
% depends continuously on $\tilde F_0$ in the compact-open topology.
% The domain of $F_0^*$ is compactly contained in
% the domain of $F_0$. Thus we have that the mating depends continuously
% on $\tilde F_0$.
\end{pf}

%\textcolor{red}{
%\begin{cor}
%  \label{cor:hybrid classes converge}
% % Let $f,g\in \mathcal A^{\underline \nu}_a$ and 
%Suppose that $F_0:U_{F_0}\to U'_{F_0}$ is the extension
%of an interval map $f_0$ with only hyperbolic periodic points. 
%Let $G_{0,n}:U_{G_{0,n}}\to U_{G_{0,n}}'$ be a sequence of pruned
%polynomial-like mappings  so that $G_{0,n}\to F_0$ as $n\to \infty$  \trevor{in the PLL topology} 
%and 
%so that $Q(F_0)=Q(G_0,n)$ for all $n\ge 0$.
%Then there exists a neighbourhood $\mathcal U$ of $F_0$ such that  
%$\mathcal H_{G_{0,n}}\cap \mathcal U$ converges
%to  $\mathcal H_{F_0}\cap \mathcal U$.
%\end{cor}}

%\trevor{Check this and make sense of this} 

Later on we will show that $\mathcal T_f$ can be viewed as the 
product of $\mathcal H_f$ and the Teichm\"uller spaces of some 
punctured tori, see Theorem~\ref{thm:manifold-Tf}.  
%\textcolor{red}{Integrate the proof in the previous one? ???} 
%\begin{cor}
%  Suppose that $F:U_F\to U'_F$ together with a collection of rays $\Gamma_F$ is a
%  pruned polynomial-like mapping. Then the hybrid class of any further pruning
%  of $F$ can be endowed with a complex analytic structure.
%\end{cor}
%\begin{pf}
%\textcolor{red}{New Proof using density of hyperbolicity.}	
%Let $Q$ denote the set of points in $K_{F}$ that are landing points of rays in
%$\Gamma_F,$ so that $F\in\mathcal{PL}_{Q}$. Let $F^*:U_{F^*}\to
%U'_{F^*}\in\mathcal{PL}_{L_{Q^*}}$ 
%be a further pruning of $F$. 
%By density of hyperbolicity, 
%there exists a pruned polynomial like mapping $F_0,$ 
%satisfying the hypotheses of Proposition~\ref{prop:tpl hybrid classes}
%and additionally for every $c\in\Crit(F_0)$ there exists
%$j_c>0$ so that $F_0^{j_c}(c)\in\Crit(F_0)$. 
%By Proposition~\ref{prop:tpl hybrid classes} one may transfer the
%complex structure on the hybrid class of $F_0$ to the hybrid class of $F$.
%\end{pf}

\section{Hybrid classes form immersed analytic manifold}
\label{sec:manifolds} 

Our goal now is to show that the real-hybrid class of a real analytic mapping
$f$ has the structure of an analytic manifold. To do this, we
show that it can be exhausted by a union of ``compatible'' complex
analytic manifolds. The difficulty that we need to overcome
germs of real analytic maps of the interval are equivalence classes of
pruned polynomial-like maps, and we need to take prunings arbitrarily 
close to the real line.

\begin{lem}
  \label{lem:para} 
Suppose that $f$ is a real analytic mapping of the interval, 
so that each periodic point of $f$ is either hyperbolic or is a simple parabolic periodic point. 
Assume that $f\in \mathcal A^{\underline \nu}_{a}$.
Then there exists a sequence of real analytic interval map $g_n\in \mathcal A^{\underline \nu}_{a}$  and 
truncated polynomial polynomial-like maps $F,G_n$  which are complex extensions of $f,g_n$
so that 
\begin{enumerate}
\item $Q(F)=Q(G_n)$; 
\item  the domains of $F,G_n$
are compactly contained in $\Omega_{a}$;
and so that $g_n\to f$ on $\overline{\Omega_{a'}}$ for each $0<a'<a$.     
\item If $f$ has only hyperbolic periodic points, then each map $g_n$ is hyperbolic;
\item 
If $f$ has parabolic periodic points, then we can choose $g_n\in \mathcal P_f$ (where $\mathcal P_f$ is as in Theorem~\ref{thm:parabolicmanifold})
so that each critical point of $g_n$ is either in the basin of a parabolic periodic point 
or a hyperbolic periodic point. % and so that $f,g$ have extensions $F,G$ so that $Q(F)=Q(G)$.
\end{enumerate} 
\end{lem}
\begin{pf} This follows from density of hyperbolicity, see \cite{KSS-density}. Here, if $f$ has
simple parabolic points we restrict to the manifold of maps with corresponding 
parabolic periodic from Theorem~\ref{thm:parabolicmanifold}. 
\end{pf} 

\begin{rem} 
Note that $f\in \mathcal A^{\underline \nu}$ may have periodic attractors which  do have critical points 
in their basin, and in this case $g_n$ will have the same property for $n$ large.   
\end{rem} 

\begin{thm}[Manifold structure of $\mathcal H_f^{\R}$] 
Assume that $f_0$ is a  real analytic mapping  so that all its 
periodic points  are hyperbolic or simple parabolic points. 
  Then the real-hybrid class of the germ of $f_0$ is an immersed real analytic submanifold of 
  $\mathcal A^{\underline \nu}$. 
%  the space of
%germs on analytic mappings on the interval. Moreover, it can be
%expressed as an increasing limit of embedded real analytic Banach
%submanifolds $X_n$,
%and each point of the hybrid class is contained in some $X_n$.
\end{thm}
\begin{pf} 
Choose a sequence of finite pruning data $Q_n\subset \partial \D$, $n\ge 1$,  associated to $f_0$, so 
that $Q_{n+1}\supset Q_n$ and $Q_n$ is an admissible pruning set for $f_0$ for each $n\ge 1$. 
In particular,  $f_0$ has  pruned polynomial-like extensions $F_{0,n}\colon U_{0,n}\to U_{0,n}'$
for each $n$ (corresponding to pruning data $Q_n$). 
We can choose $Q_n$ so that each $f\in \mathcal H_{f_0}^{\R}$
has a pruned polynomial-like extension $F$ corresponding to some $Q_n$.  
%Assume $f_0\in \mathcal A^{\underline \nu}_a$ and let $F_0\colon U_0\to U_0'$
%be a pruned polynomial-like extension of $f_0$.  
 % and assume that  $f\in \mathcal A^{\underline \nu}_a$.
%Then there exist pruned polynomial-like extensions $F_0 \colon U_0 \to U_0'$ of $f_0$
%and $F\colon U\to U'$ of $f$ which are hybrid conjugate.  
From the previous lemma, if $f_0$ has no parabolic periodic points, 
there exists a sequence  $\{g_n\}$ of mappings converging to $f_0$ on $\overline{\Omega}_a$ 
with the properties that  $g_n$ is hyperbolic and so that 
$Q_n$ is an admissible pruning set for $g_n$. It follows that for 
each  $f\in \mathcal H_{f_0}^{\R}$ and all $n$ sufficiently large, 
$f$ and $g_n$ have pruned polynomial-like extensions $F_{n}\colon U_{n}\to U'_{n}$ 
and  $G_n\colon U_{g_n}\to U_{g_n}'$   so that $Q_n:=Q(G_n)=Q(F_{n})$
and so that the domains of $F_n,G_n$ are compactly contained
in the domain $\Omega_a$ of analyticity of $f$ and $g_n$. 

 If $f_0$ has parabolic periodic points, but all of them of simple type, 
then we can assume that $g_n\in \mathcal P_f$ and that all critical points
of $g_n$ are in the basins of hyperbolic or simple parabolic periodic points. 
%Here we use that   $Q_{n+1}\supset Q_n$.  
Here $n$ (and so $Q_n$) depends on $f$ and in particular on the domain of analyticity of the map $f$ 
and thus on the size of the pruning intervals $J_i\supset J_i^* \ni f(c_i)$. 
%Each $f$ in $\mathcal T_{f_0}$, has a pruned polynomial-like
%extension, $F:U\to U'$, so that for some $n$, the external mappings of
%$g_n$ and $F$ are combinatorially compatible. 

By  Theorem~\ref{thm:tpl hybrid classes}, we can therefore parameterize
the real-hybrid class of $F_{0,n}\colon U_{0,n}\to U'_{0,n}$ by the real-hybrid class of   $G_n\colon U_{G_n}\to U_{G_n}'$.
By the properties of $g_n$,   from Theorem~\ref{thm:manifold-specialmaps} (or Corollary~\ref{cor:semihyperbolic} if $f_0$ has parabolic periodic points) it follows that 
$\mathcal H^{\R}_{g_n}\cap \mathcal A_{a'}$ is a real analytic Banach manifold for each $a'\in (0,a]$ sufficiently small. 
Therefore,  Theorem~\ref{thm:tpl hybrid classes} implies that the real-hybrid class of $F_{0,n}\colon U_{0,n}\to U'_{0,n}$ is 
an immersed analytic submanifold  $X_n$ of $\mathcal A^{\underline \nu}$. 
Moreover, $f_0$ is  in this hybrid class for all $n$ sufficiently large. 
%Let $X_{n}\subset\mathcal H_{f_0}$ denote the real trace 
%of the image of $\mathcal H_{F_0}$ under the inclusion mapping  (i.e. the real analytic interval maps) 
%corresponding to those parametrised by $G_n$ and with pruning data $Q_n$. 
Since $Q_n\subset Q_{n+1}$ we have $X_n\subset X_{n+1}$.
Since any $f\in \mathcal{H}_{f_0}$ is included in some $X_n$ we have $\cup X_n=\mathcal{H}_{f_0}$. 
\end{pf}

%\textcolor{red}{
%\begin{thm}[Manifold structure of $\mathcal H_f^{\R}\cap \mathcal A^{\underline \nu}_a$] 
%Assume that $f_0\in \mathcal A^{\underline \nu}_a$ so that all its 
%periodic points  are hyperbolic or simple parabolic points. 
%  Then the real-hybrid class of the germ of $f_0$ is an immersed Banach submanifold of 
%  $\mathcal A^{\underline \nu}_a$. 
%%  the space of
%%germs on analytic mappings on the interval. Moreover, it can be
%%expressed as an increasing limit of embedded real analytic Banach
%%submanifolds $X_n$,
%%and each point of the hybrid class is contained in some $X_n$.
%\end{thm}}

\section{Topological conjugacy classes form immersed analytic manifold}
\label{sec:topmanifolds} 
To discuss the conjugacy class of $\mathcal T_{f_0}$ it will be useful
to denote by  $\mathrm T(T^2_{\mathrm g})$ 
%$\mathcal T^2_{\mathrm g}$  
the {\em Teichm\"uller space}
of a torus $T^2$ with $\mathrm g$ marked points. This means that  
$\mathrm T(T^2_{\mathrm g})$
is the space of all pairs $(X,\phi)$ where $X$ is a 2-torus with $\mathrm g$ marked points
and  $\phi\colon T^2 \to X$ (mapping marked points to marked points) is 
a quasiconformal map so that $(X_1,\phi_1)$ is considered equivalent to $(X_2,\phi_2)$
if $\phi_2\phi_1^{-1}\colon X_1\to X_2$ is isotopic to a holomorphic map. 
If $g>0$ this space has real dimension $2\mathrm g$. If $\mathrm g=0$ then 
$T^2_0$ is not a hyperbolic surface, and the Teichmuller space of 
$T^2_0$ is equal to the the upper half-space $\mathbb H$ and so again its real dimension is two. 
For more on this see \cite{IT}.  

In the theorem below, we only count critical points whose orbits is infinite, 
so not critical points which are eventually mapped into the periodic orbit.

\begin{thm}[Manifold structure of   $\mathcal T_f$] 
\label{thm:manifold-Tf} 
 Assume that $f\in \mathcal A^{\underline \nu}$ has only hyperbolic periodic points.
Then $\mathcal T_f $ is conformally equivalent to the product $\mathcal H_f$ and sets of the form   
$\mathrm T(T^2_{\mathrm g_1})$: 
\begin{equation} 
\mathcal T_f \approx  \mathcal H_f \times  \mathrm T(T^2_{\mathrm g_1})\times \dots \mathrm T(T^2_{\mathrm g_a})
\label{eq:TisomH} \end{equation} 
and is a  real analytic manifold.  Here $a$ is equal to the number of hyperbolic periodic attractors,
and $\mathrm g_i$ is equal to the number of critical points  with disjoint  infinite orbits in the basin 
of the $i$-th periodic attractor. 
%and the codimension of  $\mathcal T_f$ in $\mathcal A^{\underline \nu}$
%is $\mbox{codim}(\mathcal H_f) - (A+B)$. 
%
%Here $A$ is the maximal number of critical points with distinct orbits which are contained 
%in the punctured basins of periodic attractors and $B$ is the number of periodic attractors which do not contain critical points in their punctured basins. 
% \textcolor{red}{A and B} 
\end{thm} 

\begin{pf}  If $f_0$ has no periodic attractors or parabolic periodic points, 
then $\mathcal H_{f_0}=\mathcal T_{f_0}$ and so this statement follows
from the previous corollary. So assume that 
$f_0$ has precisely $a>0$ attracting periodic orbits $O_1,\dots,O_a$ of period $m_1,\dots,m_a$. 
Pick a periodic point $p_i\in O_i$  and let  ${\mathrm g}_i$ be the number of critical points with disjoint orbits that are contained in $B(O_i)\setminus O_i$.  Now take a fundamental annulus $A_i$ around $p_i$ so that $f^{m_i}$ is a diffeomorphism on the disk bounded by $A_i$
and so that  each  orbit of a critical point entering the basin of $p_i$ (and which is not mapped to 
$O_i$) enters $A_i$. The modulus of the annulus is equal to that of the annulus $\{z: |\lambda_1| <|z|<1\}$
where $\lambda_i=Df^{m_i}(p_i)\in \R$. 
Identifying the inner and outer boundary of $A_i$ via the identification $z\mapsto f^{m_i}(z)$
on the boundary,  we obtain a torus $T_i^2$. 

Given a map $f\in \mathcal H_{f_0}$,  
choose a fundamental annulus $A_{i,f}$ near the periodic attractor $p_{i,f}$ (corresponding to the periodic attractor  $p_i$ of $f_0$) which is the holomorphic image of the fundamental annulus $A_i$ for $f_0$
(mapping iterates of critical points of $f_0$ in $A_i$ to corresponding ones for $f$). 
By identifying the inner and outer boundary points of $A_i$ by the dynamics of $f_0$
we obtain a torus $T_i^2$ with $\mathrm g_i$ marked points
and similarly $A_{i,f}$ induces a marked torus $T_{i,f}$ 
and we have a conformal map $\phi_f\colon T_{i,f} \to T_i$. 
Now pick an element of $\mathrm T(T^2_{\mathrm g_i})$, i.e. $X$ and a quasiconformal map  $\phi\colon T_i\to X$.
Thus we obtain a map $\phi\circ \phi_f\colon T_{i,f} \to X$.  Let $\mu_i$ be the 
Beltrami coefficient of this map on $A_{i,f}$. Using dynamics
this induces an invariant Beltrami coefficient $\mu_i$
on the basin of $p_{i,f}$. Now do this for each periodic attractor $p_{i,f}$ of $f$. 
Now extend the Beltrami coefficient on the union of the basins to an invariant Beltrami 
coefficient $\mu=\mu_{\phi_1,\dots,\phi_a}$ on all of $\C$, so that  it is zero outside  the basins. 
Using the Riemann mapping theorem, we obtain a unique qc map $h_\mu$ normalised so that 
$h_\mu(\pm 1)=\pm 1$. Because $\mu$ is $f$-invariant, 
the map  $f_\mu=h_\mu \circ f \circ h^{-1}_\mu$ is again a holomorphic map. 
If we take $\mu$ real-symmetric then $f_\mu$ will again be an interval map. 
Thus we obtain 
a map  $\Psi$ defined by 
$$\mathcal H_{f_0} \times \mathrm T(T^2_{\mathrm g_1})\times \dots \mathrm T(T^2_{\mathrm g_a}) \ni (f,\phi_1,\dots,\phi_a)  \mapsto f_{\mu(\phi_1,\dots,\phi_a)}\in \mathcal T_{f_0}.$$
Since $h_\mu$ depends holomorphically on all choices, this map is holomorphic as a map 
into $\mathcal A^{\underline \nu}$.

To see that the above map is injective, assume that 
$(f,\phi_1,\dots,\phi_a), (\tilde f,\tilde \phi_1,\dots,\tilde \phi_a)\in \mathcal H_{f_0} \times \mathrm T(T^2_{\mathrm g_1})\times \dots \mathrm T(T^2_{\mathrm g_a})$ are so that  
$\Psi(f,\phi_1,\dots,\phi_a)=\Psi(\tilde f,\tilde \phi_1,\dots,\tilde \phi_a)$. 
Then $f,\tilde f$ are hybrid conjugate. In particular the multipliers are corresponding 
periodic attractors are the same, and the 
$\bar \partial$ derivative of  $h_{\mu(\tilde f,\tilde \phi_1,\dots,\tilde \phi_a)}^{-1}\circ 
h_{\mu(f, \phi_1,\dots,\phi_a)}$ vanishes a.e. on the basin of $f$. 
By the Weyl lemma it follows that $h_{\mu(\tilde f,\tilde \phi_1,\dots,\tilde \phi_a)}^{-1}\circ 
h_{\mu(f, \phi_1,\dots,\phi_a)}$ is conformal. Hence $\phi_i=\tilde \phi_i$. 
%\begin{equation} 
%\mathcal T(f_0)\approx \mathcal H_{f_0} \times \mathcal T^{g_1}\times \dots \mathcal T^{g_a} .\end{equation} 
To see that $\Psi$ is surjective, note that any $g\in \mathcal T_{f_0}$
can be obtained in this manner by quasiconformal surgery. Indeed, for 
each periodic attracting periodic orbit, choose a periodic point $p$, 
and  take a fundamental annulus  $A_{i,g}$ surrounding $p$ and containing
the $\mathrm g_i$ forward iterates of all critical point in the basin of $p$ (and  define $A_{i,f_0}$ similarly. 
Choose diffeomorphisms $\phi_i\colon A_{i,f_0}\to A_{i,g}$ which preserve the real line, 
and so that the $f_0$-iterates points of $f_0$ in $A_{i,f_0}$ are mapped to the corresponding
critical iterates of $g$ in $A_{i,f_0}$. Next let $\mu$ be the Beltrami coefficient
associated to  $\phi_i$ in $A_{i,f_0}$, and extend $\mu$ so that it is invariant
and zero outside the basins of attracting periodic points.  
Let $h_\mu$ be the the normalised qc map with Beltrami coefficient $\mu$, 
and define $g_0=h_\mu \circ f_0 \circ h_\mu^{-1}$. As usual 
$g$ is hybrid conjugate to $f_0$ and  $g\in \mathcal A^{\underline \nu}_a$. 
Moreover,  by construction 
$$(g,\phi_1,\dots,\phi_a) \mapsto  g_{\mu(\phi_1,\dots,\phi_a)}=f.$$ 
\end{pf} 

\begin{rem} In Theorem~\ref{thm:dimension-top-hybrid-class}  we will determine the codimension
of this manifold as a subset of $\mathcal A^{\underline \nu}$. 
\end{rem} 
\begin{rem} The space $\mathcal T_f$ consists of real maps,
and so the marked points in the tori $T_i$ are real.  
%Moreover,  if $\lambda_i>0$ then 
%$A_i\cap \R$
%consists of two intervals, and $T_i\cap \R$ consists of two disjoint circles. 
%If $\lambda_i<0$ then  $T_i\cap \R=A_i/\sim$ consists of a circle. Here $\sim$
%is the equivalence relationship  $x\sim f^{m_i}(x)$.
The Teichm\"uller space of $T_i$ with $\mathrm g_i$ real marked points can be considered
as  an ordered collection of $\mathrm g_i$ real points in these tori. 
%Thus we can see 
%that this real Teichm\"uller space is determined by $\lambda_i$ and 
%an ordered set of points $g_i$ points on the circles. 
Thus we see that
this  real Teichm\"uller space is a simplex of real dimension $\mathrm g_i$. 
If $\mathrm g_i=0$ then the modulus of the fundamental annulus, and therefore
the conformal structure on the torus,  is determined by the multiplier $\lambda_i$.
\end{rem} 

\begin{rem} A similar discussion can also be found in Section 6 of \cite{McSu}. 
\end{rem}

\begin{rem} If $f$ has  parabolic periodic points, then instead
of the annuli $A_i$ the fundamental domains are the crescent shapes $S$ from 
Figure~\ref{fig:parabolic}. Identifying boundaries we obtain an infinite cylinder, 
i.e.  $\C_*$ with $\mathrm g_i$ marked points. Denote its Teichm\"uller 
space by $\mathcal T(\C_{*,\mathrm g_i})$.
 If all parabolic periodic points of $f$ are simple, then  we obtain that 
\begin{equation} 
\mathcal T_f \approx  \mathcal H_f \times  \mathrm T(T^2_{\mathrm g_1})\times \dots \mathrm T(T^2_{\mathrm g_a})
\times \mathrm T(\C_{*,g'_1}) \times \dots  \times \mathrm T(\C_{*,g'_p})
\label{eq:TisomH2} \end{equation} 
where $p$ is the number of parabolic periodic points and $\mathrm g'_1,\dots,\mathrm g'_p$ 
is the number of infinite orbits in the basins of the $i$-th parabolic periodic point. 
\end{rem}

\section{Infinitesimal theory, horizontal, vertical and transversal vector fields}
\label{sec:infinitesimal} 
% Before proving the results of this section, we
% will need some preparation.

So far we have shown that $\mathcal H_f$ and $\mathcal T_f$
have manifold structures. In the next few sections we will show
that they have the correct codimension.

A continuous vector field $\alpha$ on an \label{qcvectorfield} 
open set $U\subset\overline{\mathbb C}$ is called 
$K$-{\em quasiconformal}, abbreviated $K$-{\em qc}, 
if it has locally integrable
distributional partial derivatives $\partial \alpha$ and
$\bar\partial \alpha$, and $\|\bar\partial \alpha\|_{\infty}\leq K$.
A vector field is {\em quasiconformal} if it is 
$\kappa$-qc for some $\kappa$. 
We say that a qc vector field is normalized if
it vanishes on $\{-2,2,\infty\}$  or on $\{-1,1,\infty\}$.
If $\alpha$ is a continuous vector field,
defined on a closed set $X$, we say that $\alpha$
is {\em quasiconformal} if it extends to a normalized
qc vector field on $\overline{\mathbb{C}}$. 
%   If $X\cap\{-2,2,\infty\}=\O,$
%     then any qc vector field on $X$ admits an extension
%       to a normalized qc vector field on $\mathbb C$.
We define
% \textcolor{red}{\sout{a norm on the space of qc vector fields
%$\alpha$ on $X$
%that admit an extension to a normalized qc vector field on
%$\overline{\mathbb{C}}$ by}}
$$\|\alpha\|_{\mathrm{qc}}=\inf\| \xi \|_{\infty},$$
where the infimum is taken over all normalized
qc-extensions $\xi$ of $\alpha$ to the Riemann sphere.

% Quasiconformal vector fields are tangent at the 
% identity to holomorphic motions:
% \begin{lem}\cite[Lemma 2.10]{ALM}
%     \label{lem:ALM210}
%       Let $h_{\lambda}:X\to \C,$ $\lambda\in\mathbb D$ be a holomorphic
%         motion with base point 0. Then 
%           $$\alpha=\frac{d}{d\lambda}h_{\lambda}\Bigg|_{\lambda=0}$$
%             is a qc vector field on $X$. Moreover, if $X$ is an open set, then
%               $$\bar\partial\alpha=\frac{d}{d\lambda}\mu_{h_{\lambda}}\Bigg|_{\lambda=0},$$
%                 where $\mu_{h_{\lambda}}$ is the Beltrami differential associated to 
%                   $h_{\lambda}$.
%                   \end{lem}

%\textcolor{red}{Suppose that 
%$f:U\to \mathbb C$ is a holomorphic mapping,
%and that $v$ is a holomorphic vector field on $U$.
%Let $X_0\subset U$ be forward invariant.
%We say that a vector field $\alpha$ on $X_0$
%is \textcolor{red}{ {\em equivariant} on $X_0\subset U$}
%{\em with respect to the pair $(f,v)$} if for any $z\in X_0$,
%$$v(z)=\alpha(f(z))-f'(z)\alpha(z).$$
%When $f$ is $\R$ symmetric then we can and will  impose (here and below) 
%that both $\alpha$ and $v$ are also $\R$ symmetric. }

%If $X_0\subset U$, not necessarily forward invariant,
%let $Y_0=f(X_0),$ and suppose that $\beta$ is a vector field
%on $Y_0$. Then $\alpha$ is called the 
%{\em lift} of $\beta$ to $X_0$ if for $z\in X_0$,
%$$v(z)=\beta(f(z))-f'(z)\alpha(z).$$
%%In Lemma~\ref{lem:first lift} we will give conditions
%%on when a vector field can be lifted.
%\label{def:liftvector} 

\begin{defn}[Hybrid horizontal vector fields $E^h_f$] 
\label{def:horvector}  Let $f\in\mathcal A^{\underline \nu} $.
Then $E^h_f$ is the space of all
holomorphic vector fields $v\in T_f\mathcal A^{\underline \nu} $
defined in a neighbourhood of the interval
such that there exists a pruned polynomial-like extension $F\colon U\to U'$  of $f$
and a qc vector field $\alpha$ on $U$ so that
\begin{equation} 
v(z)=\alpha\circ F(z)-DF(z)\alpha(z)\mbox{ for } z\in U \label{eq:equivariant} 
\end{equation} 
and so that $\bar \partial  \alpha=0$ on $K_F$.
We will call such a vector field $v$ {\em hybrid-horizontal}. 
A vector fields $v$ satisfying (\ref{eq:equivariant}) is called {\em horizontal} if  $\bar \partial  \alpha=0$ on $K_F$ does not 
necessarily hold. \end{defn}

\begin{rem}  Define $v^n=\frac{d}{dt} (f+tv)^n|_{t=0}$ and 
assume that each critical point of $F$ is in the basin 
of a periodic attractor and that $v$ satisfies 
\begin{equation} Dv^n(p)(Df^n(p) - 1) = v^n(p)D^2f^n(p). 
\label{eq:multipliercond}
\end{equation}
for each attracting periodic orbit.
Then a computation shows that infinitesimally the multipliers of these periodic orbit 
stays the same for maps of the form $f+tv$. However, (\ref{eq:multipliercond}) and (\ref{eq:equivariant})  are not enough to guarantee
that $v\in \mathcal H^f$ if there are several critical points
in the basin of $p$, see the discussion in Theorem~\ref{thm:manifold-Tf}.  
\end{rem}

\begin{rem} \label{def202}  Definition \ref{def:horvector} of hybrid tangent 
does not depend on the choice
of $F\colon U\to U'$ in the following sense:  assume that $F\colon U\to U'$ and $\tilde F\colon \tilde U\to \tilde U'$
are pruned polynomial-like extension of $f$ so that $\tilde U\subset U$. 
If  (\ref{eq:equivariant}) holds on $U$ then it also holds on $\tilde U$, since then 
$\tilde F|\tilde U = F|\tilde U$.  Because of this, whether or not $v\in E^h_f$ is a condition 
on $v$ restricted to $I$.  This also become apparent from the next proposition. 
\end{rem} 

%Suppose that $F:U\to U'$ is a pruned polynomial-like mapping
%associated to $f$. We let $E^h_F$ be the set of all 
%holomorphic vector fields $v$ defined in $U$,
%so that there exists a qc vector field 
%$\alpha$ defined on $U\cup U'$ so that
%$$v=\alpha\circ F-F'\alpha \mbox{ on }U.$$
%

%\medskip
%\noindent\textit{Remark.}
%In \cite{Lyubich} for polynomial-like mappings
%it was only required that $\alpha$ be defined in a neighbourhood of 
%the filled Julia set; indeed, any such $\alpha$ can be 
%extended to all of $U$. \marginpar{This extension claim is not quite right, but
%  it is possible to extend a restriction of $\alpha.$ \textcolor{red}{QQQ} }

%\medskip

The next proposition motivates the above definition and 
shows that hybrid-horizontal vector fields $v$ can be used
to parametrise $\mathcal H_f$ in the same way as the exponential map 
does in a Riemannian manifold. %, and that $T_f\mathcal H_f=E^h_f$.

\begin{prop}[$T_f\mathcal H_f=E^h_f$] \label{prop:41}
  Given each $v\in E^h_f$ %$v\in T_f\mathcal H_f$ 
  there exists a one-parameter family
of maps 
$f_{t,v}\in\mathcal H_f$
with $f_{0,v}=f$, depending analytically on $t$ and so that $\dfrac{d}{dt} f_{t,v}\big\vert_{t=0}=v$.
Vice versa, for each $g\in \mathcal H_f$ near $f$ there exists $v\in E^h_f$ so that 
$f_{1,v}=g$. 

Similarly, for each $g\in \mathcal T_f$ near $f$ there exists a horizontal vector field 
$v$ so that $f_{1,v}=g$, and each horizontal vector field $v$ corresponds to a
family $f_{t,v}\in \mathcal T_f$.  
\end{prop}
\begin{pf}
  Assume that $v\in E^h_f$. %$v\in T_f\mathcal H_f$. 
 Then $v$ is holomorphic and there exists a pruned polynomial-like extension $F\colon U\to U'$ and 
 qc vector field $\alpha$ so that $v(z)= \alpha(F(z))-F'(z) \alpha(z)$ for $z\in U$. 
Taking the $\bar \partial$ derivative, we see that $0=\bar \partial \alpha(F(z))\overline{F'(z)} -F'(z)\bar \partial \alpha(z)$ for all $z \in U$. 
i.e.  $F^*\mu=\mu$. 
Let $H_t$ be the normalised qc homeomorphism associated to $t\mu$. 
%\textcolor{red}{Let $f_t=h_t\circ f\circ h_t^{-1}$. It follows from linearising the equation $f_t\circ h_t = h_t\circ f$ that} $\alpha=\dfrac{d}{dt} h_t \big\vert_{t=0}$. 
%$\textcolor{purple}{QQQ add reference QQQ}
Since $\mu$ is an invariant linefield, we obtain that there exists a real analytic map $f_t$ with a pruned polynomial-like extension $F_t\colon U_t\to U_t'$ so that $F_t\circ H_t = H_t\circ F$. 
Taking the $\frac{d}{dt}$ derivative in  $F_t\circ H_t = H_t\circ F$,  
writing $\alpha=\frac{d}{dt} H_t \big\vert_{t=0}$, $v=  \frac{d}{dt} F_t \big\vert_{t=0}$
 we see that $v(z)+ DF(z) \alpha(z) = \alpha\circ F(z)$ 
as claimed. 

Vice versa, if $g\in \mathcal H_f$ then from Theorem~\ref{thm:qcconj} 
the maps $f,g$  have pruned polynomial-like extensions 
$F\colon U_F\to U_F'$
and $G\colon U_G\to U_G'$, which are quasiconformally conjugate
and so that the $\bar \partial $ of the conjugacy vanishes on $K_F$. 
Let $\mu$ be the Beltrami coefficient associated to this
qc conjugacy. Letting $h_t$ be the qc conjugacy associated to $t\mu$ and 
$f_t=h_t\circ f \circ h_t^{-1}$,  we obtain the vector field 
$v(z)=\frac{d}{dt}f_t(z)\vert_{t=0}$ which is in $E^h_f$ so that $f_{1,v}=g$. 
\end{pf}

\begin{rem} \label{rem:lemma19AB} 
Let $\mu_t$ be the Beltrami coefficient of a holomorphic motion $h_t$, and let $\alpha=\frac{d}{dt}h_t \vert_{t=0}$. Then 
$\mu_t$ depends holomorphically on $t$, $\alpha$ is a qc vector field 
 and  $\bar \partial \alpha = \frac{d\mu_t}{dt}\vert_{t=0}$, see 
%{|_{t=0}}$, see 
\cite[Lemma 2.10]{ALM} and  \cite[Lemma 19]{AB}. 
\end{rem}

\begin{defn}[Hybrid-vertical vector fields  $E^v_f$] 
\label{def:verticalvect}  A holomorphic  vector field $v\in T_f\mathcal A^{\underline \nu}$ 
is {\em hybrid-vertical} if there exists a
pruned
polynomial-like extension $F:U_F\to U'_F$ of $f$ and a 
holomorphic vector field $\alpha$ on
$\overline{\mathbb{C}}\setminus K_F$ vanishing at
$\infty$ such that
\begin{equation} 
v(z)=\alpha\circ F(z)-DF(z)\alpha(z) \mbox{ for } z\in U_F\setminus K_F. \label{eq:eqverticalv}  \end{equation} 
Let $E^v_F$ denote the space of hybrid-vertical vector fields at $f$ corresponding to $F$. 
\end{defn} 

%\textcolor{red}{This is equivalent to saying that $\alpha$ is a holomorphic 
%vector field on $\C\setminus K(F)$ with at most a simple pole at $\infty$}

\begin{rem} \label{rem:vertical-depF}  The definition of  $E^v_F$ seems to depend on $F$. 
Indeed, suppose that $F\colon U\to U'$ and $\tilde F\colon \tilde U\to \tilde U'$
are both pruned polynomial-like extensions of $f$ and that $\tilde U\subset U$. 
Then a vertical vector field $v$ associated to $F$ is not necessarily 
one for $\tilde F$. This is because $\tilde U\setminus K_{\tilde F}$ is not necessarily 
contained in $U\setminus K_F$ and so $\alpha$ satisfying (\ref{eq:eqverticalv}) does
not necessarily induces $\tilde \alpha$ on $\tilde U\setminus K_{\tilde F}$.  
\end{rem} 

\begin{defn}[Hybrid-transversal]  
We say that a vector field $v\in T_f\mathcal A_{a}^{\underline \nu}$ 
is {\em hybrid-transversal} if  any family $f_t \in \mathcal A_{a}^{\underline \nu}$ for which 
$\dfrac{d}{dt} f_t=v$ has the property that $f_t$ is not hybrid conjugate to $f$
for $|t|$ small and $t\ne 0$. 
\end{defn} 

\begin{rem} Assume that $f$ has no parabolic periodic points. 
Then, taking into account the Teichm\"uller spaces $\mathcal T^{\mathrm g_i}$ 
from Theorem~\ref{thm:manifold-Tf},  
it is also possible to define the notion of a topological-vertical vector field, namely a hybrid-vertical vector field
with the additional property that it infinitesimally preserves the positions of the  infinite critical orbits in 
the basins of periodic attractors (in terms of the Teichm\"uller space $\mathcal T^{\mathrm g_i}$) and the multipliers of these periodic attractors. 
We will not develop or need this description.  
\end{rem}

\begin{prop}
  \label{prop:splitting}
For $f\in\mathcal A^{\mathcal \nu}_a$,  then for each pruned polynomial-like
extension $F$ of $f$ we have  
$T_f\mathcal A^{\mathcal \nu}=E^h_f\oplus E^v_F.$
\end{prop}
\begin{pf}
  Suppose that $v\in T_f\mathcal A^{\underline \nu}$
is a holomorphic vector field defined in a neighbourhood
$\Omega_a$ of
the interval.
Select a pruned polynomial-like extension 
$F:U\to U'$
of $f$ so that 
$v$ is well-defined in a neighbourhood of $U\cup U'$.  Let $\Gamma$ be the 
set of curves associated to $F$. 
 
We claim that there exists a smooth vector field $w$
defined in $\overline{\mathbb{C}}\setminus K_F$
vanishing near $\infty$ and such that
\begin{equation} v(z)=w(F(z))-DF(z)w(z),\label{eq:vw} \end{equation} 
for all $z\in U$ where this is well-defined. 
%\quad\mbox{for}\ z\in U\cap \mathrm{Dom}(w). \label{eq:vw} \end{equation} 
To construct such a $w$ we do the following.
Let $\gamma$ denote one of the curves in $\Gamma$ landing at a periodic point $\alpha$
of $F$.
First we define $w$ in a neighbourhood, $O\cap \gamma$, of
the fundamental domain $\gamma\cap (U'\setminus U)$
so that (\ref{eq:vw}) holds on $O\cap \gamma$.
Then we extend $w$ to $\gamma\cap U$ using the relation:
$v(F^{n-1}z)=w(F^n(z))-DF^n(z)w(z),$
and then to 
$\cup_{n=0}^N F^{-n}(\gamma),$ where $N$ is chosen so that
the complex pullbacks of $F^{-N+1}(\gamma)$ intersect
$\partial U'$. This defines $w$ on
$\partial U'\cap F^{-(N+1)}(\gamma)$. Do this for all of such curves in $\Gamma$.
Next extend $w$ smoothly
to the rest of $\partial U'$, and define $w$ on $\partial U$ by
$v(z)=w(F(z))-DF(z)w(z)$. Notice that this definition
agrees with the previous one on $\partial U\cap\partial U'$.
Now extend $w$ to $\overline{\mathbb{C}}\setminus U$ smoothly and so
that it vanishes at infinity,
and to $\mathbb C\setminus K_F$ by 
$v(F^{n-1}z)=w(F^n(z))-DF^n(z)w(z).$
Thus we have constructed $w$ as in equation (\ref{eq:vw}). 

Let us consider the Beltrami differential
$\mu=\bar\partial w$ in $\overline{\mathbb{C}}\setminus K_F$,
extended by 0 to $K_F$.
Since $v$ is holomorphic on $U$,
$\mu$ is $F$ invariant on $U$:
$$\bar \partial v=0=\bar \partial (w\circ F)-\bar\partial (DF\cdot w),$$
$$0=( \bar\partial w\circ F)  DF- DF\, \bar\partial w,\quad
\mbox{so}\quad \bar\partial w\circ F = \bar\partial w.$$
Thus we have that $\mu$ has bounded $L^{\infty}$-norm
on the sphere. Now, we solve the 
$\bar\partial$-problem: $\bar\partial u=\mu$ where
$u$ is a qc vector field on $\overline{ \mathbb C}$ 
vanishing at $\infty$.

The vector field $v^h=u\circ F-DFu$ on $U$ is 
holomorphic, since $\mu$ is $F$-invariant.
Since $\bar \partial u=\mu =0$ on $K_F$, $v^h$ is hybrid-horizontal.
Let $\alpha=w-u$. Since $\bar\partial\alpha=0$ on $\overline{\mathbb{C}}\setminus K_F$, 
$\alpha$ is holomorphic on this set 
$\overline{\mathbb{C}}\setminus K_F$ and $\alpha$ vanishes at
$\infty$.
Moreover,
$v-v^h=\alpha\circ F-DF\alpha$ on $U$, 
so $v-v^h$ is a vertical vector field.

\medskip
\noindent{\em   Uniqueness of the splitting:}
Assume that there exists a vector field 
 $v \in E^h_f\cap E^v_f,$ $v\neq 0$.
Then there exists a pruned polynomial-like 
representative $F:U\to U'$ of $f$, a qc vector field $w$
defined in the neighbourhood  $U'$ of $K_F$ 
so that $v=w\circ F-DFw$  so that $\bar \partial w=0$ on $K_f$
and a holomorphic vector field $\alpha$
in $\overline{\mathbb{C}}\setminus K_F$, which vanishes at $\infty$
so that $v=\alpha\circ F-DF\alpha.$
Let us consider vector field 
 $u=w-\alpha$  defined on $U'\setminus K_F$.
Since $u\circ F^n = (F^n)' u$ it follows that 
%
%it is bounded with respect to 
%the hyperbolic norm on $U$. \textcolor{red}{SSS $u\circ F^n = DF^n u$. Why does this imply 
%$u$ is bounded??? } 
%Thus we have that
$|u(z)|\rightarrow 0$ as $z\rightarrow K_F$, $z\in U$.
It is important to observe that we have the same estimate 
at the points in $K_F\cap\partial U$. % (still by $F$-invariance).
Thus we have that $u$ admits a continuous extension to
$\overline U$ which vanishes on $K_F$.
Hence we have that $w$ and $\alpha$ agree on the 
pruned Julia set. 
Let $\beta$ be the vector field which is equal to 
$w$ on $K_F$ and $\alpha$ on 
$\overline{\mathbb{C}}\setminus K_F$. Then
$\beta$ has distributional
derivatives of class $L^2$ and $\bar\partial\beta=0$.
Thus by Weyl's Lemma, $\beta$ is a holomorphic vector field
defined on $\overline{\mathbb C}$, and since $\beta$ vanishes at
$\infty$, it is linear. 
Hence we can write $\beta(z)=[az+b]/dz$. Thus 
%$$v(z)=[a f(z)+b-f'(z)(az+b)]/\partial z,$$
$$v(z)=v_0(z) /\partial z \mbox{ where } v_0(z)=[a f(z)+b-f'(z)(az+b)].$$
From this it is easy to see that $v=0$. 
Indeed, note that $f_t$ is a family satisfying $f_t(\partial I)\subset \partial I$. 
Hence $\frac{df_t}{dt}=v$ implies $v(-1)=v(1)=0$.  Let us show that this implies $a=b=0$. 
To do this we need to consider several possibilities. 

 (i) If $f(-1)=-1$ and $f(1)=1$ then   $v_0(-1)=af(-1)+b - f'(-1)(-a+b) =(b-a)(1-f'(-1))=0$ 
and  $v_0(1)=af(1)+b - f'(1)(a+b)=(a+b)(1-f'(1))=0$.  Since we assume that 
$f$ has no parabolic points this implies $a=b=0$. 

(ii)  If $f(-1)=-1$ and $f(1)=-1$
then   $v_0(-1)=af(-1)+b - f'(-1)(-a+b)=(b-a)(1-f'(-1))=0$, which gives $a=b$  
and therefore $v_0(1)=af(1)+b - f'(1)(a+b)=(b-a)-f'(1)(a+b)]=0$ implies $a+b=0$ because we
assume $f$ has no critical points in $\partial I$.

 (iii) If $f(-1)=1$ and $f(1)=-1$ then we get $v_0(-1)=af(-1)+b - f'(-1)(-a+b)=(a+b)-f'(-1)(b-a)=0$ 
and $v_0(1)=af(1)+b - f'(1)(a+b)=(b-a)-f'(1)(a+b)=0$.  Combined this gives 
$(b-a)(1-f'(-1)f'(1))=0$ which again implies $b=a$ and therefore, using the previous expressions,  $a=b=0$.
\end{pf}

\section{Estimates for vertical vector fields} 

\begin{lem}
  \label{lem:44}
 Let $F\colon \E \to \E'$ be the part of the pruned polynomial-like mapping 
 $F\colon  \E\cup B\to \E'\cup B'$ coming from an 
 expanding Markov structure of the external map. 
Then there exist $\lambda>1$
and $N\in\mathbb N$ such that 
$|DF^N(z)|>\lambda$ for all $z\in F^{-N}(\partial \E)$ so that $z,\dots,F^{N-1}(z)\in \E$.
\end{lem}
\begin{rem}
The set $E$ also contains the basins of periodic attractors, 
%This lemma also holds if $f$ has periodic attractors,
 provided these all can be included
in the set $K_{X,O}$.  If $f$ has periodic attractors
with {\lq}large{\rq} basins then its global pruned polynomial-like extension 
is of the form $F\colon \E\cup B \to \E'\cup B'$, and then the above expansion
holds on $F^{-N}(\partial U)$, and not on $F^{-N}(\partial B)$. 
\end{rem} 
\begin{pf}  If $F^n(z)$ is in the interior of $\E'$ then 
take a balls $D\Subset D'\subset \E'$ containing $F^n(z)$.
By the Koebe Distortion Theorem and since puzzle pieces shrink in diameter to zero, 
it follows that there if $n$ is large enough then 
$|DF^n(z)|\ge 2$. If $F^n(z) \in \partial \E\cap \partial \E'$ then $F^n(z)\in \Gamma$.
If $F^n(z)$ is contained in a curve $\gamma\subset \Gamma$ 
which is part of the expanding Markov structure, then there exists 
$n_0$ which is independent of $z$ so that $z,\dots,F^{n-n_0}(z)\in \E$. 
Hence the proof goes as before in this case. If $\gamma$ is part of the attracting structure, then 
$\gamma$ is an invariant curve going through a repelling periodic point then it is possible
that $F^k(z),\dots,F^{n-1}(z)\in \gamma$ and $z,\dots,F^{k-1}(z)\in \E$
and we also get expansion because the multiplier at the periodic point is repelling
and because of the first part of the argument.
\end{pf}

Now we are in a position to prove the following proposition,
which is a modification of
\cite[Lemma 4.10]{Lyubich}. To 
deal with the fact that pruned polynomial-like mappings 
do not have moduli bounds,
we use the expansion along the boundaries of the puzzle pieces.
This lemma is one of the key technical tools in obtaining 
a lower bound on the codimension of $\mathcal H_f$.
In the statement, we use one more pullback than
in \cite[Lemma 4.10]{Lyubich} since if $F:U\to U'$ is a pruned polynomial-like mappings, we need not have that $U'\supset U$; however, we do have that $U\supset F^{-1}(U)$.

\begin{prop}[Control for vertical vector fields] 
  \label{prop:avbound}
Let $F\colon U=\E\cup B\to \E'\cup B'=U'$ be a global pruned polynomial-like mapping.
Let $W''=\E$, $W'=F^{-1}(W'')$ and $W=F^{-1}(W')$. 
Let $v$ and $\alpha$ %(z)/dz$ and $\alpha(z)/dz$ 
be holomorphic vector fields satisfying
$$v(z)=\alpha\circ F(z)-DF(z) \alpha(z) ,\quad z\in W',$$
where $v$ is holomorphic in $W''$ and $\alpha$ is 
holomorphic in $\overline{\mathbb C}\setminus K_F.$
Then there exists a constant $C$ depending on $\lambda,$
the constant from Lemma~\ref{lem:44} and the extremal widths of 
rectangles comprising $U'\setminus U$,
such that 
$$\|\alpha\|_{\overline{\mathbb C}\setminus W}
\leq C\|v\|_{W'}\;\mbox{and}\;
\|v\|_W\geq C^{-1}\|v\|_{W'}.$$
\end{prop}

\begin{pf}
  Let $\gamma=\partial W$.
By Lemma~\ref{lem:44},
there exists $N\in\mathbb N$
and $\lambda>1$ so that for $z\in F^{-N}\gamma$,
$|DF^N(z)|>\lambda.$

We have that
$$\alpha\circ F^N(z)-DF^N(z)\alpha(z)=DF^N(z)\sum_{k=0}^{N-1}
\frac{v(F^k(z))}{DF^{k+1}(z)},$$
so 
$$
\alpha(z)  = 
%& =\dfrac{\alpha\circ F^N(z) -DF^N(z)\sum_{k=0}^{N-1}
%  \frac{v(F^k(z))}{DF^{k+1}(z)}}{DF^N(z)} =\\ 
%  & \\ %
%& =   \dfrac{\alpha\circ F^N(z) -\sum_{k=0}^{N-1}
%v(F^k(z))DF^{N-k-1}(F^{k+1}z)}{DF^N(z)}  \\
%& =  
  \dfrac{\alpha\circ F^N(z)} {DF^N(z)}-
\sum_{k=0}^{N-1} \dfrac{v(F^k(z))}{DF^{k+1}(z)}
 %\end{array}  .
 $$
Since $|DF^{k+1}(z)|$ is bounded away from zero for all $z\in F^{-N}\gamma$
and all $0\le k\le N-1$, this expression implies that there exists   a constant $A>0$ (which depends on $N$ and on  $F\colon U\to U'$), so that
 \begin{equation}
 \|\alpha\|_{F^{-N}\gamma}\leq 
\frac{\|\alpha\|_{\gamma}+\lambda A\|v\|_{W}}{\lambda}\leq \frac{\|\alpha\|_{\overline{\mathbb C}\setminus W}
  +\lambda A\|v\|_W}{\lambda},\label{eq:boundalpha_n} \end{equation} 
where the last inequality follows from the Maximum Principle and since $\partial W=\gamma$.
By the Maximum Principle 
we also have 
$$\|\alpha\|_{\overline{\mathbb C}\setminus W}\leq\|\alpha\|_{F^{-N}\gamma}.$$
Thus
$$\|\alpha\|_{\overline{\mathbb C}\setminus W}\leq 
\frac{\|\alpha\|_{\overline{\mathbb C}\setminus W}
  +\lambda A\|v\|_W}{\lambda},$$
and so
$$\|\alpha\|_{\overline{\mathbb C}\setminus W}\leq \frac{\lambda A\|v\|_W}{\lambda-1}.$$
This proves the first inequality. The second inequality follows from 
it and the Maximum Principle:
$$\|v\|_{V'}\leq \|\alpha\|_{\overline{\mathbb C}\setminus V''}
+\|DF\|_{V'}\|\alpha\|_{\overline{\mathbb C}\setminus V'}$$
$$\leq \|\alpha\|_{\overline C\setminus W}(1+\|DF\|_{W'})
\leq \frac{\lambda A}{\lambda-1} (1+\|DF\|_{W'}) \|v\|_W .$$
\end{pf}

\begin{rem} 
The previous proposition is the main reason why we can deal with 
 vertical vector fields  more easily than in  \cite{ALM}
 where puzzle mappings are introduced which have infinitely many domains
 and which form a necklace neighbourhood of $I$. 
\end{rem}

\section{Estimates for horizontal vectors fields}
\label{sec:keylemma} 

In this section we will prove the following result, which is called the Key Estimate in \cite{ALM}. 
Here the proof follows easily from Part A of our paper.

%%%% KEY Lemma 
\begin{lem}[Control for horizontal vector fields] 
\label{lem:key lemma}
Suppose that $f\in\mathcal A_a^{\underline \nu}$ has only hyperbolic periodic orbits. 
Then there exists a neighbourhood 
$\mathcal W \subset\mathcal{A}_a^{\underline \nu}$ of $f$ and $C>0$, so that for any $g\in\mathcal W$ 
and any $v\in T_g\mathcal H_g,$ there exist a pruned polynomial-like map $G\colon U_g\to U_g'$
and a qc vector field $\alpha$  so that 
$$v(z)=\alpha\circ G(z)-DG(z)\alpha(z)\quad\mbox{for}\; z\in U_g,$$
so that 
$$\|\bar\partial\alpha\|_{\mathrm{qc}}\leq C\|v\|_{\infty}.$$
\end{lem}
\begin{pf} Take $v\in T_g\mathcal H_g$. Then by Proposition~\ref{prop:41} there exists a family 
$g_t\in \mathcal H_g$ with $\dfrac{d}{dt} g_{t}\big\vert_{t=0}=v$. Let $G_t$ and $G$ be pruned polynomial-like extensions 
of $g_t$ and $g$. By Theorem~\ref{thm:qcconj} there exists a family of qc maps $h_t$
so that $G_t=h_t\circ G \circ h_t^{-1}$. By Corollary~\ref{cor:keytypelemma} we have 
$$\varkappa(h_t) \le 1 + L ||g_t-g||_\infty , $$
where $\varkappa(h_t)$ is the qc dilatation of $h_t$ and as before $||\cdot ||_\infty$ stands for the supremum norm on $\overline \Omega_a$. 
So if we write $\alpha=\dfrac{dh}{dt}\vert_{t=0}$ we get $||\bar \partial \alpha||_{qc} \le L ||v||_\infty$.
As in the proof of Proposition~\ref{prop:41} we also 
 have $v(z)=\alpha\circ G(z)-G'(z)\alpha(z)$ for  $z\in U_g$. 
\end{pf}

\section{The codimension of  conjugacy classes}\label{sec:codimension} 

In Section~\ref{sec:manifolds} we showed that  the 
conjugacy class of a map $\mathcal{A}^{\underline \nu}$ is a real analytic manifold. 
 In this section, we show it  has the expected codimension in
$\mathcal{A}^{\underline \nu}$.

The lower bound for the dimension comes from Proposition~\ref{prop:avbound}.
To obtain an upper bound for the dimension we shall use the following lemma, 
which follows from Lemma~\ref{lem:key lemma} (which is our version of the   %Infinitesimal Pullback Argument and the 
``Key Estimate" of \cite{ALM}):  
%Theorem~\ref{thm:ipa} and 

\begin{lemma}[Continuity of tangent space] 
\label{lem: Eh closed}
If $f_n\in \mathcal A^{\mathcal \nu}_a$ converge to 
$f\in \mathcal A^{\mathcal \nu}_a$ where $f$ has only hyperbolic periodic points, and let $F_n\colon U_n\to U_n'$
be pruned polynomial-like extensions of $f_n$ so that $U_n'$ contains 
a $\delta$-neighbourhood of $I$ for each $n\ge 0$. Then if $v_n\in E_{f_n}^h$
is a sequence of horizontal vectors with $||v_n||_{U_n'}= 1$ 
then $v_n$ converges to a horizontal vector $v\in E_f^h$ with $||v||_{U'}\ge 1$ for some $U'\supset I$. 
\end{lemma} 
\begin{pf} Since $v_n\in  E_{f_n}^h$ there exists a 
qc vector field $\alpha_n$ on $U_n'$ so that 
$v_n(z)=\alpha_n (F_n(z)) - DF_n(z) \alpha_n(z)$ for all $z\in U_n$
and $v_n$ is holomorphic on $U_n$. 
By Lemma~\ref{lem:key lemma} (the Key Estimate) we have that $\alpha_n$ is a sequence of quasiconformal vector fields with bounded dilatation, thus by  the compactness lemma for qc vector fields, see for example \cite{ALM},
there exists qc vector field $\alpha$ so that $\alpha_n\to \alpha$ along some subsequence
and, as $||v_n||_{U_n'}=1$ there exists a holomorphic $v$ so that $v_n\to v$ on some
definite neighbourhood of $I$.  
Moreover, $v(z)=\alpha (F(z)) - DF(z) \alpha(z)$ for all $z$ in a definite neighbourhood of $I$  
and therefore $v$ is a horizontal vector field, $v\in E_f^h$. 
% \textcolor{red}{can you find bound on dilatation? yes}
\end{pf} 

We now show that the hybrid classes have the expected codimension,
{\it cf.} \cite[Theorem 4.11]{Lyubich} and also \cite[Theorem 10.4]{Sm-Fib}.

\begin{thm}\label{thm:dimension-top-hybrid-class} 
  Assume that $f\in \mathcal A^{\underline \nu}_a$ has only hyperbolic periodic points. Then 
  \begin{enumerate}
  \item $\mathcal H_f^{\R}$ is a real analytic manifold  whose 
 codimension in  $\mathcal A^{\mathcal \nu}_a$ is equal to $\nu_H=
\nu+\xi_{noness-att}$  where $\nu$ is the number of critical points of $f$ and $\xi_{noness-att}$
is the number of periodic attractors without (real) critical points in their basins - we call these non-essential attractors.
\item $\mathcal T_f$ is a real analytic
manifold whose  codimension  in $\mathcal A^{\mathcal \nu}_a$ is equal to 
$\nu_T=\nu-\zeta(f)$ where $\nu$ is the number of critical points of $f$
and $\zeta(f)$ is the maximal number of critical points in the  basins of periodic 
attractors with pairwise disjoint {\rm infinite} orbits. 
%
%$\nu'=
%\nu-\xi_{ess-att}$  where $\nu$ is the number of critical points of $f$
%and 
%$\xi_{ess-att}$ is the maximal integer $\xi$ so that there 
% exist $\xi$ critical points which are based in the basins of periodic attractors
% and  whose orbits are mutually disjoint. 
 \end{enumerate}
\end{thm}
 \begin{pf} Let us only consider the space $\mathcal H_f^{\R}$.  The dimension of the space $\mathcal T_f$
 follows from this  using Theorem~\ref{thm:manifold-Tf}.  

  \noindent\textit{Proof of the lower bound on the codimension.}
Recall that functions in $\mathcal A^{\underline \nu}_a$ are holomorphic on 
$\Omega_a=\{z\in\mathbb C:|z-I|<a\}.$
Suppose that $f_n\in \mathcal A^{\underline \nu}_a$ is a sequence of semi-hyperbolic  mappings
$f_n\rightarrow f$ on $\overline{\Omega_a}$.  That this is possible, follows from density 
of hyperbolicity, see \cite{KSS-density}. 
%Then   $f\in \mathcal A^{\underline \nu}_a$} and there exists $a>0,$ so that 
%$f_n$ converges to $f$ on compact subsets of 
%$\Omega_a$. Taking $a$ slightly smaller, we can assume that
%$f_n$ converges to $f$ on $\overline{\Omega}_a$.

To each $f_n$ we associate a pruned polynomial like 
mapping $F_n:U_n:=\E_n\cup B_n\to \E'_n\cup B'_n=:U_n'$ and to $f$ we associate a pruned polynomial-like
mapping $F:U:=\E\cup B \to \E'\cup B'=:U'$. We choose them so that 
$U'\subset \Omega_a$ and $U'_n\subset \Omega_a$ 
for all $n$ sufficiently large, and so that $F_n\colon U_n\to U_n'$
converges to $F\colon U\to U'$ (and so in particular 
$U_n$ contains each compact subset of $U$ 
for $n$ large (and similarly for $U_n'$ and $U'$).
By Theorem~\ref{thm:manifold-specialmaps}, for each $F_n$ the space of  
 vectors vertical to $\mathcal H_{F_n}$ has dimension $\nu_H$, 
 and let  $\{v_n^1,v_n^2,\dots, v_n^{\nu_H}\}$ be a basis. 
Moreover, we can assume that these vectors have unit length and
by a Theorem of Riesz, that they are
almost orthogonal in the sense that 
\begin{equation}
\dist(v^i_n,\mathrm{span}\{v_n^1,\dots,v_n^{i-1}\})>\frac{1}{2},
\mbox{ for }i=2,3,\dots, \nu_H,\label{eq:almostortho}\end{equation} 
where the distance is in the space of bounded 
holomorphic functions on $U_n$ which extend continuously to $\overline U_n$.

Let us show that the unit ball in the 
vertical direction is compact, {\em cf.} \cite[Corollary~10.1]{Sm-Fib}:
Suppose that $w_n$ is a sequence of
vertical vector fields at $F$ of unit length.
Then there exist $\alpha_n$
holomorphic on 
$\overline{\mathbb C}\setminus K_F$, vanishing at $\infty$
and satisfying
%$$w_n=\alpha_n\circ F-DF\alpha_n.$$
%\textcolor{red}{
$$w_n=\alpha_n\circ F_n-DF_n \alpha_n \mbox{ on }\E_n.$$
%}
Notice that because of the  Maximum Principle and equation (\ref{eq:boundalpha_n}), 
for each $i$ there exists $C_i$ so that 
%$$||\alpha_n||_{\overline{\mathbb C}\setminus F^{-i}(U')}\le C_i\mbox{ for all }n\ge 0.$$ 
%\trevor{\textcolor{red}{
$$||\alpha_n||_{\overline{\mathbb C}\setminus F_n^{-i}(\E'_n)}\le C_i\mbox{ for all }n\ge 0.$$
%}}
%we only need to check it on the preimages of $f^{-i}
%but there it follows immediately from the formula,
%$$\alpha_n(z)=\frac{\alpha_n\circ F^{k+1}(z)-w_n(F^k(z))}{DF^{k+1}(z)}.$$
%So there exists a sequence of constants
%$C_i>0$ so that for each $n$ sufficiently large,
%$$\|\alpha_n\|_{\overline{\mathbb C}\setminus F^{-i}(V)}\leq C_i.$$
Hence, because the basins $B_n$  are contained in the interior of $K(F_n)$,
there exists a subsequence $\alpha_{n_j}$  
which converges on compact subsets of $\mathbb C\setminus K_F$
to a holomorphic vector field $\alpha$, so that the associated
vector fields $w_{n_j}$ converge to a vector field
$w$ and so that  $w=\alpha\circ F-DF\alpha$ on $U\setminus K_F$. Thus the unit 
ball in the vertical direction is compact.

It follows that we can assume that
for each $i,$ $v_n^i\rightarrow v^i,$ uniformly on compact subsets
of $U$.
By Proposition~\ref{prop:avbound}, we have that
$v^i$ cannot be zero. 
Finally since the vectors $v_n^i$ are almost
orthogonal, see equation (\ref{eq:almostortho}), the collection $\{v^1,\dots,v^{\nu_H}\}$ is 
linearly independent.
\checkmark

\medskip
\noindent\textit{Proof of the upper bound on the codimension.}
Assume by contradiction that for $n$ sufficiently large we have that 
$\mathrm{codim}(\mathcal H_{F}) > \mathrm{codim}(\mathcal H_{F_n})$.
Because of Proposition~\ref{prop:splitting} it then follows that  
 there exists a non-zero vertical vector field $v$ in the tangent space to $F$, which can be approximated by horizontal vector fields $v_n$ tangent to $F_n$. This contradicts Lemma~\ref{lem: Eh closed}.
\end{pf}

\section{Hybrid conjugacies are embedded manifolds}\label{sec:embedded} 

\begin{thm} Take $f\in \mathcal A_a^{\underline \nu}$, assume 
that $f$ has only hyperbolic periodic points and let $v$ be a hybrid-vertical 
vector field. Then there exists $\epsilon>0$ so that the family 
$f_t=f+tv$ only intersects each hybrid class once. 
\end{thm} 
\begin{pf}  This follows immediately from the next lemma. 
\end{pf} 

\begin{lemma}  Assume that $f_n,g_n\in \mathcal A_a^{\underline \nu}$ are real-hybrid conjugate, 
converge to $f\in  \mathcal A_a^{\underline \nu}$ and that $f$ has only hyperbolic periodic orbits. 
Then any limit $v$ of $$\dfrac{f_n-g_n}{||f_n-g_n||_a}$$
is a hybrid-tangent vector to $f$, i.e.  $v\in T_f \mathcal H_f^{\R}$. 
\end{lemma} 
\begin{pf} See  \cite[Lemma 8.1]{ALM}.  
%\trevor{\textcolor{red}{Add detail??}}  
\end{pf}

\begin{rem} The above Theorem and lemma do not hold unless we assume that $f$ has only hyperbolic periodic points, 
see \cite[page 453 - footnote]{ALM}. In  \cite{CvS2} we will elaborate on this. 
 \end{rem} 

\section{Hybrid classes are Banach manifolds} \label{sec:hybridbanach} 
In Theorem~\ref{thm:dimension-top-hybrid-class}  we showed that  $\mathcal H_{f}$ is a real analytic manifold, in the germ sense. %  and also in the PPL sense.
Let us now improve this statement by showing:

\begin{thm}\label{thm:hybridbanach}   For each $f_0\in  \mathcal A_a^{\underline \nu}$, 
\begin{enumerate} 
\item $\mathcal H_{f_0}^\R \cap   \mathcal A_a^{\underline \nu}$ is real Banach submanifold of $\mathcal A_a^{\underline \nu}$ of codimension $\nu_H= \nu+\xi_{noness-att}$;
\item $\mathcal T_{f_0}^\R \cap   \mathcal A_a^{\underline \nu}$ is real Banach submanifold of $\mathcal A_a^{\underline \nu}$ of codimension $\nu_T=\nu-\zeta(f)$.
\end{enumerate} 
\end{thm} 
\begin{pf} In  Section~\ref{sec:manifoldstructure} it is shown that 
any $f\in \mathcal H_{f_0}^\R$ 
has a pruned polynomial-like extension $F\colon U\to U'$
so that its hybrid conjugacy class is conformally equivalent to 
the hybrid conjugacy class of $G\colon U_{G}\to U_{G}'$
where $G$ is the extension of a real analytic semi-hyperbolic  map 
$g\in  \mathcal A_a^{\underline \nu}$. 
This implies that %near $f_0$
 there exists $a'>0$ so that  $f\in \mathcal A_{a'}^{\underline \nu}$ and so that 
$\mathcal H_{f_0} \cap   \mathcal A_{a'}^{\underline \nu}$
is a Banach manifold of codimension-$\nu_H$  near $f$
where $\nu_H= \nu+\xi_{noness-att}$. 
From this we obtain that there exists a real analytic  function 
$$\Psi\colon \mathcal A_{a'}^{\underline \nu}\to \R^{\nu_H}$$
defined near $f$, 
which has maximal rank at $f$ and so that on some neighbourhood of $f$ in $\mathcal A_{a'}^{\underline \nu}$
one has 
$$\mathcal H_{f_0}^\R \cap   \mathcal A_{a'}^{\underline \nu}=\Psi^{-1}(0).$$  
Since $ \mathcal A_{a}^{\underline \nu}$ is dense in  $\mathcal A_{a'}^{\underline \nu}$
it follows that the restriction  $\tilde \Psi$ of $\Psi$ to 
$\mathcal A_{a}^{\underline \nu}$ also has maximal rank at each $\tilde f\in \mathcal  A_{a}^{\underline \nu}$ 
in a neighbourhood of $f$ in $\mathcal  A_{a'}^{\underline \nu}$. 
Moreover, whether $\tilde f$ is in $\mathcal H_{f_0}^\R$ only depends on 
$\tilde f$ restricted to $I$. It follows that 
$\mathcal H_{f_0}^\R \cap   \mathcal A_{a}^{\underline \nu}=\tilde \Psi^{-1}(0)$
and since $\tilde \Psi\colon  \mathcal A_{a'}^{\underline \nu}\to \R^{\nu_H}$ has maximal 
rank, it follows that $\mathcal H_{f_0}^\R \cap   \mathcal A_{a}^{\underline \nu}$ is a Banach manifold. 
The second assertion follows similarly. 
\end{pf}

% Take a sequence $g_n\in  \mathcal A_a^{\underline \nu}$ of semi-hyperbolic  maps
%converging to $f_0\in  \mathcal A_a^{\underline \nu}$.
%In Section~\ref{sec:manifoldstructure} it is shown that  $\mathcal H_{g_n} \cap   \mathcal A_a^{\underline \nu}$
%is a Banach manifold. Moreover,  by  Corollary~\ref{cor:hybrid classes converge}
%we have that there exists a neighbourhood $\mathcal U$ of $f_0$ in  $\mathcal A_a^{\underline \nu}$
%so that $\mathcal H_{g_n} \cap   \mathcal A_a^{\underline \nu}$  converges 
%to $\mathcal H_{f_0} \cap   \mathcal A_a^{\underline \nu}$.  Now vectors $v_1,\dots,v_{\nu_H}$
%spanning the linear space $E^v_{f_0}$.
%Since $\mathcal H_{g_n} \cap   \mathcal A_a^{\underline \nu}$  converges 
%to $\mathcal H_{f_0} \cap   \mathcal A_a^{\underline \nu}$, and by the Implicit Theorem 
%there exists  $n$ so that for each $f\in  \mathcal H_{f_0} \cap   \mathcal A_a^{\underline \nu}$
%there exist unique  $g\in \mathcal H_{g_n} \cap   \mathcal A_a^{\underline \nu}$ and $t_i\in \R$
%so that $f=g+t_1\nu_1+\dots+t_{\nu_H} \nu_{\nu_H}$. Therefore 
%$\mathcal H_{f_0} \cap   \mathcal A_a^{\underline \nu}\cap \mathcal U$ is the graph over a Banach manifold. 
%\end{pf} }

\section{Conjugacy classes of real analytic maps are path connected} 
\label{sec:pathconnected}

In this section we obtain Theorem~\ref{thm:connected} by proving the following: 

\begin{thm}\label{thm:pathconnected}
  Let $f,\tilde f\colon I\to I$ be two real analytic maps
in $\mathcal A^{\underline\nu}$ with all periodic points hyperbolic, which are
real-hybrid conjugate (via an order preserving conjugacy).
Then there exists a family of real maps $f_t\in\mathcal
A^{\underline\nu},$
$t\in [0,1]$ depending
real analytically on $t\in [0,1]$ so that $f_0=f$, $f_1=\tilde f$ 
and so that $f_t$ is real-hybrid  conjugate to $f$ and $\tilde f$ for each $t\in [0,1]$. 

Similarly, if $f,\tilde f\in \mathcal \A^{\underline \nu}$ are topologically conjugate and only 
have hyperbolic periodic orbits, 
 then there exists a real analytic path connecting $f,\tilde f$ in $\mathcal T(f)$. 
\end{thm}

%\begin{thm}\label{thm:pathconnected}
%  Let $f,\tilde f\colon I\to I$ be two real analytic maps
%in $\mathcal A^{\underline\nu}_a$ with all periodic points hyperbolic, which are
%real-hybrid conjugate (via an order preserving conjugacy).
%Then there exists a family of real maps $f_t\in\mathcal
%A^{\underline\nu}_a,$
%$t\in [0,1]$ depending
%real analytically on $t\in [0,1]$ so that $f_0=f$, $f_1=\tilde f$ 
%and so that $f_t$ is real-hybrid  conjugate to $f$ and $\tilde f$ for each $t\in [0,1]$. 
%
%Similarly, if $f,\tilde f\in \mathcal \A^{\underline \nu}_a$ are topologically conjugate and only 
%have hyperbolic periodic orbits, 
% then there exists a real analytic path connecting $f,\tilde f$ in $\mathcal T(f)$. 
%\end{thm}

\begin{pf} From \cite{CvS} there exists a quasisymmetric conjugacy between $f$ and $\tilde f$.
  Let $F\colon U_F\to U'_F$, $\tilde F\colon U_{\tilde F}\to U'_{\tilde F}$ be the
pruned polynomial-like
extensions of $f, \tilde f$ from Theorem~\ref{thm: pruned-pol-like-map}.
Choose a quasiconformal map
$h\colon U'_F \to U'_{\tilde F}$ so that $h$
maps $\partial U_F$ to $\partial U_{\tilde F}$ and $\Gamma_F$ to
$\Gamma_{\tilde F}$,
and so that $h$ is a conjugacy on these
sets and also so that $h$ 
agrees with the qs-conjugacy between $f$ and $\tilde f$ on the real line. 
We can also ensure that $h$ is a conformal conjugacy near hyperbolic 
periodic attractors. 
Now let $H_0=h$ and define $H_{n+1}$ by $\tilde F\circ H_{n+1}=H_n\circ
F$. Since $F, \tilde F$ are conformal, critical values of $F$ are mapped to
critical values of $\tilde F$,
and the conjugacy relation holds, 
$H_{n+1}$ is well-defined and has the same quasiconformal dilatation as $H_n$.  It follows that 
there exists a subsequence of $H_n$ converging to some quasiconformal map $H$. 

Let $Y_n=U'_F \setminus \cup_{i=0}^n (F^{-i}(U'_f\setminus U_F))$.  Then for each $n\ge 0$, $Y_{n+1}\subset Y_n$ 
and $H_{n+1}$ agrees with $H_n$ outside $Y_n$ and also on $\gamma_f$. 
Because of the last assertion of Proposition~\ref{prop:expandingcircle} it follows that
$\cap_{n\ge 0} Y_n$ has empty interior and so for each point $z$ outside this set there exists $n$ so that $H_{n+i}(z)=H_n(z)$
for all $i\ge 0$.  It follows that any convergent subsequence of $H_n$ converges to the same map $H\colon U_F\cup U'_F \to 
U_{\tilde F}\cup U'_{\tilde F}$ and therefore that $\tilde F\circ H=H\circ F$. 

Let $\mu$ be the Beltrami-coefficient associated to $H$  on $U_F\cup U'_F$ (describing the ellipse field obtained by the pullback of the 
standard circle field under $H$). Then $F$ preserves the ellipse field defined by $\mu$ 
since $\tilde F$ is conformal and since $F=H^{-1}\circ \tilde F \circ H$. Now extend $\mu$ to $\bar \C$
by setting $\mu=0$  on $\mathcal \C \setminus (U_F\cup U'_F)$
and let $H_{t\mu}$ be the corresponding
quasiconformal map corresponding to 
the Beltrami coefficient $t\mu$ normalised so that $H_{t\mu}(I)=I$ and $H_{t\mu} (\infty)=\infty$ (here we use the Measurable Riemann Mapping Theorem). 
Since $F$ preserves the ellipse field defined by $\mu$
and $F$ is conformal, $F$ also preserves the ellipse field defined by $t\mu$ for each $t\in [0,1]$. 
It follows that $G_t:=H_{t\mu} \circ F \circ H_{t\mu} ^{-1}$ is  a conformal
map on $H_{t\mu}(U_F)$. We can ensure that $\mu$ is $z\mapsto \bar z$ symmetric, 
and so we can assume that $H_{t\mu}$ is real for $t$ real. 
In particular $g_t=F|_I$ is a family of analytic maps of the
interval.
\end{pf}

We also have the following:  

\begin{prop}  Assume that $f_0\in  \mathcal A^{\underline\nu}_a$. 
Then $\mathcal H^{\R}_{f_0}\cap \mathcal A^{\underline\nu}_a$ is path connected. 
\label{prop:pathconnected} 
\end{prop}
\begin{pf} 
To see this, take $f,\tilde f\in  \mathcal H^{\R}_{f_0}\cap \mathcal A^{\underline\nu}_a$ and construct
a path $g_t\in  \mathcal H^{\R}_{f_0}$ connecting these maps as in the  proof of the previous theorem.  It could happen that $g_t$ is not analytic on $\Omega_a$. To rectify
this, we use the argument from the proof of \cite[Theorem 9.2]{ALM}, 
which we include for completeness. This argument shows that the curve $\{g_t\}$ can be approximated
by a curve in $\mathcal H_{f} \cap \mathcal A_a^{\underline \nu}$ connecting $f_0=f$ and $f_1=\tilde f$.
Indeed, let $a'\in (0,a)$ be so that for every $t\in[0,1]$, 
$\{g_t\}\subset\mathcal A_{a'}^{\underline\nu}$.
Let $\Pi_{a,a'}:\mathcal A^{\underline\nu}_a\to \mathcal
A^{\underline\nu}_{a'}$
be the inclusion mapping, given by the restriction of $f$ to $\Omega_{a'}$.
Next, consider a one parameter real analytic family of vector fields
$\{\hat v_t\}$ in $T\mathcal A^{\nu}_{a'}$ so that for each $t$, $\hat v_t
=(v_t^1,\dots, v_t^{\nu'})$ a basis for $E^v_{g_t}$, the vector space 
vertical  to $T_{g_t}(\mathcal H_{f_0}\cap \mathcal A_{a'}^{\underline \nu})$. 
Recall that $\nu' = \nu + \nu_{noness}$. (Here we use Theorem~\ref{thm:hybridbanach}.)

%\textcolor{red}{$\nu$ or $\nu'$ of these vectors?}

Define 
\begin{equation} P(t,\hat s)=g_t+\hat s \cdot \hat
v_t\in\mathcal{A}^{\underline\nu}_{a'},\;(t,\hat s)
\in[0,1]\times[-1,1]^{\nu'} \end{equation} 
and let $i_0:[0,1]\to[0,1]\times[-1,1]^{\nu'},$ where $i_0(t)=(t,0,\dots, 0)$.
Let $\Phi:[0,1]\times[-1,1]^{\nu'}\to\mathcal{A}^{\underline\nu}_a$
be a real analytic family so that
$\Phi(0,0,\dots,0)=P(0,0,\dots,0)=f_0,$ $\Phi(1,0,\dots,0)=P(1,0,\dots, 0) = f_1$ and
$\Pi_{a,a'}\circ\Phi$ is $C^1$ close to $P$.

To see that such a $\Phi$ exists we argue as follows: 
For each $a^*\in (0,a')$ and each $\epsilon>0$
there exists an integer $k>0$ so that the $k$-th order polynomial $J^kg_t$ which 
agrees with the first $k$ term of the Taylor expansion of $g_t$ at, say, the critical point $c_1$, 
has the property $||J^kg_t - g_t||_{\Omega_{a^*}}<\epsilon$. 
(One can see this by considering the $a'$-ball around each $x\in I$
and looking at the power series expansion.)

By the Implicit Function Theorem, there exists a real analytic curve 
$\zeta:[0,1]\to[0,1]\times[-1,1]^{\nu'}$, $C^1$ close to $i_0$, so that 
$\Pi_{a,a'}\circ \Phi\circ \zeta$ is contained in $\mathcal
H_{f}\cap\mathcal A^{\nu}_{a'}$ and so that the first coordinate of  $\zeta(t)$ is equal to $t$. 
Since $\Pi_{a,a'}\circ\Phi$ is $C^1$ close to $P$,
we have that the curve $\Pi_{a,a'}\circ\Phi(\{0\}\times[-1,1])$ only intersects the real-hybrid class of $f_0$ at $f_0$, and $\Pi_{a,a'}\circ\Phi(\{1\}\times[-1,1])$ only intersects the hybrid class of $f_0$ at $f_1$.
Thus $\Phi(\zeta(0)) = f_0$ and $\Phi(\zeta(1)) = f_1$.
Since $f_0, f_1\in\mathcal A^{\nu}_a$, it follows that 
$\Phi\circ\zeta:[0,1]\to\mathcal A^{\underline\nu}_a$ is a real-analytic path connecting $f_0$ to $f_1$ in $\mathcal H_f \cap\mathcal A^{\nu}_{a}$.
\end{pf} 
\begin{rem} In the previous theorem we cannot drop the assumption that $f,\tilde f$
only hyperbolic periodic orbits. This is because it is possible that 
$f, \tilde f$ are topologically conjugate on the interval $I$, 
so that $f$ has only repelling periodic orbits and $\tilde f$ 
has some periodic orbits which topologically repelling and parabolic
(up to coordinate change, an iterate takes the form $x\mapsto x+x^{2n+1}$). 
On the other hand, the proof goes through if we consider maps for which 
each parabolic periodic point is simple.
\end{rem} 

\section{Hybrid conjugacies form a partial lamination} 

In this section we will prove Theorem~\ref{thm:lamination}:

\begin{thm}[Partial lamination] 
	Every $f\in\mathcal{A}^{\underline{\nu}}$ without parabolic points has a neighbourhood which is laminated by real-hybrid conjugacy classes. More precisely, for each neighbourhood  $\mathcal V_2$ of $f$ 
there exists  
%	neighbourhood $\mathcal V$ of $f$ in $\mathcal{A}_a^{\underline{\nu}}$ so that for any $g\in \mathcal V$
%	the intersection of $\mathcal H_g$ with $\mathcal V$ is pathconnected. 
neighbourhood $\mathcal V_1\subset \mathcal V_2$ of $f$ 
so that for each 
	$g_0,g_1\in \mathcal V_1\cap \mathcal H_f$ then there exists a path 
	$g_t\in \mathcal{A}^{\underline{\nu}}$, $t\in [0,1]$ inside 
	$\mathcal V_2\cap \mathcal H_g$  connecting $g_0,g_1$.  	
	\end{thm}
\begin{pf} Consider the pruned polynomial-like extension $F\colon \E\cup B \to \E'\cup B'$ of $f$. 
Since $f$ has no parabolic periodic points, this pruned polynomial-like extension 
persists over a neighbourhood $\mathcal V_1$ of $f$  (of not necessarily 
real pruned polynomial-like maps)  by holomorphic motion. It follows that there exist pruned polynonial-like extensions
$G_i\colon U_{G_i}\cup B_{G_i} \to U_{G_i}\cup B'_{G_i}$ of $g_i$ which are obtainied
by holomorphic motion from the pruned polynomial-like extension $F$. By \cite[Theorem 2]{BR} 
it follows that for each $\epsilon>0$
there exists a neighbourhood $\mathcal V_1$ so that   the dilatation  $K$ 
of the {\lq}external{\rq} qc-conjugacy  from Proposition~\ref{prop:qcconjpartial} is at most 
$1+\epsilon$. As in the proof of Theorem~\ref{thm:pullback} and Corollary~\ref{cor:keytypelemma}
it follows that there exists  arc $g_t$ connecting $g_0$ and $g_1$ whose diameter is small  if $\epsilon>0$ is close to zero. 
\end{pf}

\section{Conjugacy classes of real analytic maps  are contractible} 
\label{sec:contractible}

Take $f_0\in \mathcal A^{\underline \nu}_a$. 
In order to show that the space $ \mathcal T_{f_0}$ is contractible, 
we will first need to show that the pruned polynomial-like structure 
persists on all of  $\mathcal T_{f_0}$
for some $a'\in (0,a)$. To do this on the entire infinite dimensional space
$\mathcal H_{f_0}$ we  cannot use holomorphic motions,
and therefore will use  the notion of 
quasiconformal motions, see \cite{ST}.

\def\T{\mathcal{S}}
\begin{defn}\label{def:qcmotion}  Let $\T$ be a connected topological Hausdorff space and $X\subset \C$. Then
the a map $H\colon \T\times X\to \C$ is called a {\em quasiconformal motion}  if, 
writing $H_t(z):=H(t,z)$, the following holds: 
\begin{enumerate}
\item  for some base point  $t_0\in \T$ we have $H_{t_0}=id$;  
\item for any $t\in \T$ and any $\epsilon>0$ there exists a neighbourhood $U$ of $t$
such that for all  $t',t''\in U$ and for all quadruples $x_1,x_2,x_3,x_4$ 
of distinct points in $X$ the cross-ratios of $H_{t'}(x_1),H_{t'}(x_2),H_{t'}(x_3),H_{t'}(x_4)$
and of $H_{t''}(x_1),H_{t''}(x_2),H_{t''}(x_3),H_{t''}(x_4)$
all lie within an $\epsilon$-ball in the Poincar\'e 
metric of $\C\setminus \{0,1\}$.
\end{enumerate} 
\end{defn} 

In our setting we will take $\T= \mathcal H_{f_0}\cap \mathcal A^{\underline \nu}_a$. %\mathcal A^{\underline \nu}_{a}\cap  \textcolor{red}{QQQ} 

\medskip
Since we cannot use holomorphic motions,  will use the following result of Douady-Earle
on extensions of qs maps on the boundary of a disk.  To state
it will be useful to associate to a qc map $H\colon \D \to \D$ its dilatation: 
$$\mu_{H}=\bar \partial H / \partial H.$$  
 
\begin{thm} [Douady-Earle extension]\label{DE-extension} 
 Let $h\colon \partial \D \to \partial \D$ and define 
$$G(z,w)= \dfrac{1}{2\pi} \int_{\partial \D} \dfrac{h(\zeta)-w}{1-\bar w h(\zeta)}\dfrac{1-|z|^2}{|z-\zeta|^2}\, 
|d\zeta|.$$
 \begin{enumerate}
\item 
Given $z\in \D$ there exists a unique $w$ so that $G(z,w)=0$. Define $H\colon \D\to \D$
by $H_h(z)=w$.  Then  $H_h$ is a real analytic diffeomorphism on $\D$. 
\item 
The map  $H_h \colon  \D \to  \D$ extends to a continuous map 
of $\overline{\D}$ to $\overline{\D}$ so that $H_h|\partial \D=h$.
\item The quasiconformal dilatation of $H_h\colon \D\to \D$
is bounded if $h\colon \partial \D\to \partial \D$ is quasisymmetric.
%\item Pick $z\in \D$ and define 
%$$\mu_h:=\bar \partial H_h / \partial H_h\in L^\infty(\D,\C).$$ Then 
%$h\mapsto \mu_h(z)\in \D$ is continuous
%when we endow the set of circle homeomorphisms $\partial \D \to \partial \D$ with the supremum norm. 
\item for any $\epsilon>0$ there exists $\delta>0$ so that 
if there exists a qc extension of $h$ whose dilation is $\le K$ then 
the dilatation of $H_h$ is at most $K^{3+\epsilon}$ provided $K\le 1+\delta$. 
\item  for any $\epsilon>0$ and any qs map $h_1\colon  \partial \D\to \partial \D$
 there exists $\delta>0$, so  that if $h_2\colon \partial \D\to \partial \D$
is a qs maps so that the map $h_2\circ h_1^{-1}$ has a qc extension $H$ to $\D$ with 
$||\mu_H||_\infty\le \delta$  then $||\mu_{H_{h_1}}-\mu_{H_{h_2}}||_{\infty}\le \epsilon$. 
\end{enumerate} 
\end{thm} 
\begin{pf} 
Items (1)-(3) are proved in \cite{DE}, see also \cite{Hub}. 
Item (4) is Corollary 2 on page 41 of \cite{DE}.  Item (5) generalises item (4) and 
follows  similarly. Indeed, let $M$ be the open unit ball in the Banach space $L^\infty(\D,\C)$.
 Given $\mu\in M$, by the Measurable Riemann Mapping Theorem, 
 there exists a unique qc map $$\phi^\mu\colon \D\to \D$$ fixing $\pm 1$ and $i$ 
 so that $\mu= \mu_{\phi^\mu}$ (i.e. $\bar \partial \phi^\mu = \mu \cdot \partial \phi^\mu$).  
 Set $h^\mu=\phi^\mu|\partial \D$. Now define 
$\sigma \colon M\to M$ by 
$$\sigma(\mu)=\mu_{H_{h^\mu}}.$$
So $\sigma$ assigns to the Beltrami coefficient of any extension 
of a qs map, the Beltrami coefficient of the Douady-Earle extension of this qs map. 
In particular, 
$$\sigma(\mu_{H_{h_1}}) = \mu_{H_{h_1}}.
$$
By assumption there exists a map $H\colon \D\to \D$ so that $H\circ h_1=h_2$ on $\partial \D$
and so that $||\mu_H||_\infty \le \epsilon$.
Hence $$\sigma (\mu_{H \circ H_{h_1}}) = \mu_{H_{h_2}}.$$  
By \cite[p 182. Theorem 5.5.6]{AIM} 
$$\mu_{\psi\circ \varphi^{-1}}(w)= \dfrac{\mu_\psi(z)-\mu_\varphi(z)}{1-\mu_\psi(z)\overline{\mu_\varphi(z)}} \cdot \left(  \dfrac{\varphi_z(z)}{|\varphi_z(z)|} \right)^2.$$ 
and taking $\psi=H \circ H_{h_1}$ and $\varphi = H_{h_1}$
we get from this formula that 
$$||\mu_{H_{h_1}}-\mu_{H\circ H_{h_1}}||_\infty \le ||\mu_H||_\infty \le  
\epsilon.$$
Since $\sigma$ is continuous, there exists $\delta$ so that 
$$||\mu_{H_{h_1}} - \mu_{H_{h_2}}||_\infty = ||\sigma (\mu_{H_{h_1}}) - \sigma (\mu_{H \circ H_{h_1}})||_\infty<\delta.$$
\end{pf}

\begin{prop}\label{prop:globalchoiceU} 
 Let $F_0\colon U_0\to U_0'$ be a pruned polynomial-like mapping with rays $\Gamma_0$. 
For each $\kappa>1$ 
one can redefine the domains of $F_0$ 
by possibly shortening the rays $\Gamma_0$ and lowering
the roofs,  and thus obtain an equivalent pruned polynomial-like map $F_{0}\colon U_{F_0}\to U_{F_0}'$ so that the following holds.  \begin{itemize}
\item 
 Let $\mathcal H_{F_0}(\kappa)$ be the set of maps
in $\mathcal H_{F_0}^{\R}$ which are $\kappa'$-qc conjugate to  $F_{0}\colon U_{F_0}\to U_{F_0}'$  for some $\kappa'<\kappa$. Then $\mathcal H_{F_0}^{\R}(\kappa)$ is an open subset of $\mathcal H_{F_0}^{\R}$, 
$\mathcal H_{F_0}^{\R}(\kappa)\subset \mathcal H_{F_0}^{\R}(\kappa)$ for $\kappa<\kappa^*$
and $\cup_{\kappa>0} \mathcal H_{F_0}^{\R}(\kappa)=\mathcal H_{F_0}^{\R}$. 
\item  For each  $G\in \mathcal H_{F_0}^{\R}(\kappa)$   
there exist a pruned polynomial-like mapping  $$G:U_{G}\to U'_{G}$$ 
with rays $\Gamma_{G}$  
and a quasiconformal motion
$$H\colon \mathcal H_{F_0}^{\R}(\kappa) \times X_n  \to \C$$
 of   $X_n=\overline{U_{F_0}\cup  U'_{F_0}\cup \Gamma_{F_0}}$   over $\mathcal H_{F_0}(\kappa)$
 so that \begin{enumerate}
\item
$\partial U_{G}=H_F(\partial U_{G_0}), \partial U'_{G}=H_G(\partial U_{F_0})$
 where we write $H_G(\cdot )=H(G,\cdot)$. 
%\item  For each $f\in T$ write $H_f(z)=H(f,z)$. Then 
%%There exists an extension of $H_F$  to a  quasiconformal map $H_F\colon \C\to \C$ so that 
%$H_f(U_{f_0})=U_f$, $H_f(U_{f_)}')=U'_f$, 
$$H_{G}\circ F_0(z)=G\circ H_G(z) \mbox{ for all }z\in \partial U_{F_0}\cup \Gamma_{F_0}$$
and so 
\item that the Beltrami coefficient of $H_G$ depends continuously on $G\in \mathcal H_{F_0}(\kappa)$, i.e. 
$G\mapsto \mu_{H_G}\in L^\infty(U_{F_0}\cup U_{F_0}',\D)$ is continuous. 
\end{enumerate} 
\end{itemize} 
\end{prop} 

\begin{pf} 
%That $\T$ is a Banach manifold follows from Theorem~\ref{thm:manifold-specialmaps}
%and that $\T_n$ is an open subset of $\T$ follows from 
%Corollary~\ref{cor:keytypelemma}. 
%By definition  the pruning data $Q_n$ can be used for 
%each $f\in \T_n$. In other words, one has a pruned Julia set $K_{f,n}$ which is associated to
%$f$ corresponding to the pruning data $Q_n$ and the set $X_n$. 
%This means that one can associate an external map $\hat f_{X_n}$ to $f$ and $X_n$. 
One can choose the smooth curves $\Gamma_{\hat f_X}$ from Lemma~\ref{lem:raysrepelling} 
so that they are orthogonal to $\partial \D$. Indeed, consider a linearisation 
at a repelling periodic point $p$ of $F_X|\partial \D$. The fact that $\hat f_X$ preserves $\partial \D$
implies that the multiplier at $p$ is real. Using the linearisation we immediately see that there is a unique smooth invariant curve through $p$ which is orthogonal to $\partial \D$. 
This ensures that the curves $\hat \Gamma_X$ are uniquely determined (by $f$ and the intervals
$J$)  apart from their length.   (In fact, a lower bound for the length of these smooth curves through $\partial \D$
is determined by the multiplier at the repelling periodic points and the upper bounds of $|D\hat f_X|$ 
and $|D^2\hat f_X|$.)  

Note that $\mathcal H_{F_0}^{\R}(\kappa)$ can be considered
as a subset of $\mathcal A^{\underline \nu}_a$ which is a metric space and therefore admits a partition of unity. 
%see Lemma~\ref{lem:partitionunity}.  
Hence, using a partition of unity argument one can choose 
$\Gamma_{\hat f_X}$ so that the arc length of each of these curves depends continuously on 
$f\in T_n\subset   \mathcal H_{f_0}^{\R}\cap \mathcal A^{\underline \nu}_{a}$. 
Similarly, one can choose the {\lq}roof{\rq} of the sets $\hat V_{\f_X}$ to depend continuously on 
 $f\in T_n\subset   \mathcal H_{f_0}^{\R}\cap \mathcal A^{\underline \nu}_{a}$. Here we will choose these roof curves
 to be circles near $\partial \D$. 
 
 Note that the normalised Riemann mapping from $\C\setminus K_X(f) \to \C\setminus \D$
 depends continuously on $K_{X}(f)$ and therefore on $f$, see Theorem~\ref{thm:KXcont}.  
 It follows that  the sets $\partial U_F$ and $\partial U'_F$, $\Gamma_F$  corresponding 
 tot $\partial V_{\f_X}$ and $\partial V_{\f_X}$ and the set $\Gamma_{\f_X}$ also move continuously with $F$.
In other words, the pruned polynomial-like maps 
 $F\colon U_F \to U'_F$   also move continuously with $f$ (in the Carath\'edory topology). 

Now parametrise the fundamental domains in $\Gamma_F$ by arc length (and on all of
$\Gamma_F$  dynamically).  Similarly parametrise  $\partial U_F$ and $\partial U'_F$ through arc length. 
Using this, one can   define maps $H_F$ on $\partial U_{F_0}\cup \partial U'_{F_0}$
to $\partial U_{F}\cup \partial U'_{F}$
so that  $H_F$ conjugates $F$ and $F_0$ on these sets.  
Similarly, if necessary,  we do this for the attracting structures $B_f$ associated to $f$. 
By construction, the sets $U_F$, $U'_F$,  $U'_F\setminus U_F$, $U_F\setminus U_F'$ (and similarly $B_F$ etc) 
form quasidiscs for each $f$. 

This means that we can use the Douady-Earle extension  to extend $H_F$ to a quasiconformal $H_F\colon U_{F_0}\cup U'_{F_0} \to U_{F}\cup U'_{F}$. Here we use that each the finitely many components of $\C\setminus (\partial U_F \cup \partial U'_F)$ is a quasidisc, and so we can use the  Douady-Earle extension on each of these
separately. Thus we obtain a quasiconformal map  $H_F\colon U_{F_0} \to U_F$  which depends continuously on $f$ in the required sense, see item 5 in Theorem~\ref{DE-extension}. 
\end{pf}

%We say that a holomorphic motion is {\em normalised} if it fixes the points $\pm 1$. \textcolor{red}{QQQ} 

\begin{rem} 
%\begin{enumerate}
%\item  %
The reason we cannot use a holomorphic motion here is because
it needs to be defined over all of the   infinite dimensional space  $\mathcal H_{F_0}^{\R}(\kappa)$. 
This forces us to use a partition of unity argument. Thus we obtain deformations which are not holomorphic in $t\in \mathcal H_{F_0}^{\R}(\kappa)$. 
%Of course another reason is that $\mathcal H_{F_0}^{\R}$ is not a complex manifold, but a real analytic manifold. 
\end{rem}

\begin{thm}\label{thm:contractible}   Assume that $F_0\colon U_0\to U_0'$ be a pruned polynomial-like mapping.
Then $\mathcal H_{F_0}^{\R}$ is contractible. 
\end{thm} 
\begin{pf}
Let $n>1$ be an integer. Let us first show that $\mathcal H_{F_0}^{\R}(n)$ is contractible. % is contractible in $S_{n+1}$. 
 From the previous proposition and the pullback argument 
from Theorem~\ref{thm:pullback} we obtain for each  $\mathcal H_{F_0}^{\R}(n)$
 a qc conjugacy $H_G$ between $G\colon U_{G}\to U'_{G}$ and $F_0\colon U_{0}\to U'_{0}$
so that the map $G\mapsto \mu_{H_G}  \in L^\infty(U_{F_0}\cup U_{F_0}',\D)$ 
is continuous.  Let $H_{t\mu_{H_G}}$
be the qc conjugacy associated to $t\mu_{H_G}$ normalised so that $H_{t\mu_{H_G}}(\pm 1)=\pm 1$
and $H_{t\mu_{H_G}}(\infty)=\infty$. Thus we get a  new map 
$$R_t(G)=H_{t\mu_{H_G} }\circ F_0 \circ H_{t\mu_{H_G}}^{-1}\in \mathcal H_{F_0}^{\R}$$   
depending analytically on $t$ and so that  
$R_0(G)=G$ and $R_1(G)=F_0$.  
Since $G\mapsto \mu_{H_G}\in L^\infty$ is continuous, the map $(t,G) \mapsto H_{t\mu_G}$
depends continuously on $t$ and $G$. Hence the retract $(t,G)\to R_t(G)$ is also continuous. 

Since $\mathcal H_{F_0}^{\R}(n)$ can be viewed as a subset of 
$\mathcal A^{\underline \nu}_a$ the result follows from the next theorem. 
%Note that $R_t(G)$ may not be in $ \mathcal A^{\underline \nu}_a$ for all $t\in [0,1]$.
%However, $R_0(F)$ and $R_1(F)$ are in  $ \mathcal A^{\underline \nu}_a$.
%Since $H_{t\mu_{H_f}}$ is at most a $n$-qc homeomorphism for each $t\in [0,1]$, 
%it follows that there exists $a'\in (0,a)$ so that  $R_t(F)\in  \mathcal H_{f_0} \cap \mathcal A^{\underline \nu}_{a'}$ for all $t\in [0,1]$. 
%In other words, $R_t(F)$ has an analytic extension to $\Omega_{a'}$. We do not claim that the domain of 
%the pruned polynomial-like map $R_t(F)$ is inside $\Omega_{a}$. 
\end{pf}

\begin{thm}[\cite{AE}]\label{thm:AE} 
If a normal space $\T$ is the union of a sequence of open subsets $\T_n$
such that $\overline \T_n\subset \T_{n+1}$ and $\T_n$ contracts to a point
in $\T_{n+1}$ for each $n\ge 1$, then $\T$ is contractible. 
\end{thm}

\begin{rem} By Lemma~\ref{lem:partitionunity}, for any open set $O$ 
in the space $\mathcal A^{\underline \nu}$ with the real analytic topology
and any $a>0$ there exists $g\in O\setminus \mathcal A^{\underline \nu}_a$. 
 That is why the above argument is insufficient to conclude that
 $\mathcal H^{\underline \nu}_f$ is contractible.  
\end{rem} 

We do however have the following:

\begin{thm} Let $f_0\in  \mathcal A^{\underline \nu}_a$. 
Then $\mathcal H_{f_0}\cap \mathcal A^{\underline \nu}_a$ is contractible. 
\end{thm} 
\begin{pf}  Choose a sequence  of pruning intervals $J_{f_0,n}$ for $f_0$ 
and let $X_n,Q_n$ be the corresponding pruning data. Choose
these pruning intervals so that $Q_{n+1}\supset Q_n$. 
Let $F_{0,n}\colon U_{F_0,n}\to U_{F_0,n}'$
be  a pruned   polynomial-like mapping which extends $f_0\colon I\to I$ corresponding
to pruning data $Q_n$, 
so that $F_{0,n+1}$ is a restriction of $F_{0,n}$. 

Let $\T_n$ be the set of maps $f\in \T := \mathcal H_{f_0}\cap \mathcal A^{\underline \nu}_a$
so that $f$ has a pruned polynomial-like mapping extension $F\colon U_{F,n}\to U_{F,n}'$ corresponding
to pruning data $Q_n$ and so that  $F_n$ and $F_{0,n}$ are qc conjugate with a
conjugacy which has dilatation $<n$. 
In other words, one has a pruned Julia set $K_{f,n}$ which is associated to
$f$ corresponding to the pruning data $Q_n$ and the set $X_n$. 
This means that one can associate an external map $\hat f_{X_n}$ to $f$ and $X_n$. 
 
Let us first show that $\T_n$ is contractible in $\T_{n+1}$. 
 From  Proposition~\ref{prop:globalchoiceU}  and the Pullback Argument  (see Theorem~\ref{thm:pullback}) 
 we obtain for each $f\in \T_n$ a qc conjugacy $H_F$ between $F\colon U_{F,n}\to U'_{F,n}$ and $F_0\colon U_{0,n}\to U'_{0,n}$
so that the map $F\mapsto \mu_{H_F}  \in L^\infty(U_{F_0}\cup U_{F_0}',\D)$ 
is continuous.  Let $H_{t\mu_{H_F}}$
be the qc conjugacy associated to $t\mu_{H_F}$ normalised so that $H_{t\mu_{H_F}}(\pm 1)=\pm 1$
and $H_{t\mu_{H_F}}(\infty)=\infty$. Thus we get a  new map 
$$R_t(F)=H_{t\mu_{H_F} }\circ F_0 \circ H_{t\mu_{H_F}}^{-1}\in \mathcal H_{F_0}$$   
depending analytically on $t$ and so that  
$R_0(F)=F_0$ and $R_1(F)=F$.  
Since $F\mapsto \mu_f\in L^\infty$ is continuous, the map $(t,F) \mapsto H_{t\mu_F}$
depends continuously on $t$ and $F$. Hence the retract $(t,F)\to R_t(F)$ is also continuous. 
Note that $R_t(F)$ may not be in $ \mathcal A^{\underline \nu}_a$ for all $t\in [0,1]$.
However, $R_0(F)$ and $R_1(F)$ are in  $ \mathcal A^{\underline \nu}_a$.
Since $H_{t\mu_{H_F}}$ is at most a $n$-qc homeomorphism for each $t\in [0,1]$, 
it follows that there exists $a'\in (0,a)$ so that  $R_t(F)\in  \mathcal H_{f_0} \cap \mathcal A^{\underline \nu}_{a'}$ for all $t\in [0,1]$. 
In other words, $R_t(F)$ has an analytic extension to $\Omega_{a'}$. We do not claim that the domain of 
the pruned polynomial-like map $R_t(F)$ is inside $\Omega_{a}$. 

To obtain a family $\tilde R_t(F)$ so that $R_t(F)\in  \mathcal H_{f_0} \cap \mathcal A^{\underline \nu}_{a}$, 
we argue as in Proposition~\ref{prop:pathconnected}. 
Indeed,  choose vector fields
$\{\hat v_f\}$ in $T_f \mathcal A^{\nu}_{a'}$ depending smoothly on $f\in \mathcal H_{f_0}\cap \mathcal A_{a'}^{\underline \nu}$
and so that  the vectors $\hat v_f =(v_f^1,\dots, v_f^{\nu'})$ together with  $T_{f}(\mathcal H_{f_0}\cap \mathcal A_{a'}^{\underline \nu})$ 
span the  tangent space  $T_{f}(\mathcal A_{a'}^{\underline \nu})$.  

Next define 
\begin{equation} P(f,t,\hat s)=R_t(f)+\hat s \cdot \hat
v_{R_t(f)} \in\mathcal{A}^{\underline\nu}_{a'},\;(f,t,\hat s)
\in \T_n \times [0,1]\times [-1,1]^{\nu'} \end{equation} 
and let $i_0: \T_n \times [0,1] \to \T_n\times [0,1]  \times[-1,1]^{\nu'},$ where $i_0(f,t)=(f,t,0,\dots, 0)$ for all $f\in \T_n$ 
and $t\in [0,1]$.
Let $\Phi: \T_n \times[-1,1]^{\nu'}\to \mathcal{A}^{\underline\nu}_a$
be a real analytic family so that
$\Phi(f,0,\dots,0)=P(f,0,\dots,0)=f$ for each $f$ and so that 
$\Pi_{a,a'}\circ\Phi$ is $C^1$ close to $P$. This can be done, as in  Proposition~\ref{prop:pathconnected},
by a polynomial approximation. 

By  Theorem~\ref{thm:hybridbanach} the space
 $\T_n$ is a Banach manifold.
Therefore,  using the Implicit Function Theorem (on Banach manifolds),  
we obtain $\zeta:\T_n\times [0,1]\to \T_n \times [0,1]\times[-1,1]^{\nu'}$, $C^1$ close to $i_0$
so that $\Pi_{a,a'}\circ \Phi\circ \zeta$ is contained in $\mathcal H_{f_0}\cap\mathcal A^{\nu}_{a'}$.
As before $\Phi\circ\zeta:[0,1]\to\mathcal A^{\underline\nu}_a$ is a  smooth 
 path connecting $f$ to $f_0$ in $\mathcal H_{f_0} \cap\mathcal A^{\nu}_{a}$.
Since $\Phi\circ \zeta$ is $C^1$ close to $P$ 
 fact lies in $\T_{n+1}$.  
 It follows that $\T_n$ is contractible inside $\T_{n+1}$.    Theorem \ref{thm:AE}  therefore implies
 that $\T$ is contractible. 
\end{pf}

\part*{Part C: Open questions} 

Let $Pol^{\R,d}$ be the space of real polynomials of degree $d$. 

\begin{question} Let $f\in \mathcal A^{\underline \nu}_a\cap Pol^{\R,d}$. 
Is $\mathcal T_{f,a}^{\underline \nu}\cap Pol^{\R,d}$ a real analytic manifold? 
\end{question} 
\medskip
The answer to this question is affirmative,  if real critical points  of $f$ have finite orbits, due
to transversality, see \cite{Epstein,LSvS3}.  On the other hand, it is obviously not enough that 
 $\mathcal T_{f,a}^{\underline \nu}$ is a real analytic manifold, see Theorem \ref{thm:manifold}, to conclude 
 that $\mathcal T_{f,a}^{\underline \nu}\cap Pol^{\R,d}$ is a real analytic manifold. 

In a similar vain, we have shown  in Theorem~\ref{thm:connected}
that two  real polynomials which are real-hybrid conjugate, 
can be connected by a one-parameter family of real analytic maps of the interval within the same real-hybrid class. 
It would be interesting to know whether one can find a one-parameter family of such polynomials. 

A partial answer to this question was given in \cite{ChvS} where the space $U^{2d}$ was considered
for $d=2$. Here $\mathcal U^{2d}$ is the space of real  polynomials of degree $2d$ 
with a unique critical point on the real line. So these maps are unimodal maps of the real line.

\begin{question} Assume that the real critical points of  $f_0,f_1\in  \mathcal U^4$ are periodic
and that $f_0,f_1$ are topologically conjugate on $\R$. Does this imply that there exists a continuous 
family $f_t$ so that  every map $f_t$ is also topologically conjugate to $f_0,f_1$ on $\R$. 
\end{question} 

In \cite{CvS} it was shown that for $d=2$ the answer if affirmative. This implies that 
the sets of maps in $\mathcal U^4$ with a given topological entropy form connected subsets, 
giving a partial generalisation of \cite{MTr, BvS} to polynomials with non-real critical points.

\part*{Appendices and References} 
\appendix

\section{Local connectivity of the pruned Julia set and
complex box mappings which include domains that do not intersect the real line} 
%Proof of Theorem~\ref{thm:pullbackcomplex}

The aim of this appendix is to prove Theorem~\ref{thm:KX}. 
% The main idea will be to use  complex box mappings which include domains that do not intersect the real line. 
In order to do this
we need the following notation. 
Given an interval $I\subset \R$, and $\theta\in(0,\pi)$, we denote by $D^{+}_{\theta}(I)$ (respectively
$D^{-}_{\theta}(I)$) the region in the upper (respectively
lower) half-plane bounded by $I$ together with the circle arc
subtending $I$ that meets the real axis with external angle $\theta$
at each boundary point of $I$. We let $D_{\theta}(I)=
D^{+}_{\theta}(I)\cup D^{-}_{\theta}(I)\cup I$ and call this set a {\em Poincar\'e disc}. 
If $\theta>\pi/2$  we shall call this a {\em lens domain}. 
\label{defn:poincare disk}
This set corresponds to the set of points with a fixed distance to $I$ in the Poincar\'e metric
in $\C_I$. Let $\Crit_{na}(f)$ be the set of critical points which are
not in the basin of periodic attractors.

\begin{thm}\label{thm:pullbackcomplex} 
 There exist complex neighbourhoods $W\subset W'$ of  $\Crit_{na}(f)$ 
 and $\delta>0$  so that 
if  $J_i$ are intervals containing $f(c_i)$ with $|J_i|<\delta$ for all $i$ and taking $J=\cup J_i$,   
the following holds.  Let $R_{K_1^*}$ be the first return map to $K_1^*:=\Comp_\Crit f^{-1}(J)$.
Then $R_{K_1^*}$ has an extension to a possibly 
multi-valued mapping $F\colon W \rightarrow W'$. 
More precisely,
\begin{enumerate}
\item $W'\supset  K_1^*$ and $W'$ is a union of pairwise disjoint Poincar\'e
  domains  each based on a real interval;
\item $W\subset W'$, i.e.  for each $x\in K_1^*$ and $c\in \Crit_{na}(f)$ and each $n> 0$ with  $f^n(x)\in W'_c$ is so that 
$W'_c$ is a component of $W'$ then  we have
 $$\Comp_x f^{-n}
 (W'_c)\subset W'.$$ 
\end{enumerate} 
\end{thm} 

Note that $K_1^*:=\Comp_\Crit f^{-1}(J)$ is not contained in $\R$, but that 
$f(K_1^*)\subset \R$. 

\begin{rem}  It follows that if $A\subset I$ is an attracting periodic orbit, so that 
the immediate basin of each $a\in A\cap W'$ is contained in $W'$
then for each $x\in K_1^*$ which is eventually mapped into $A$, then 
the component of the basin of $A$ containing $x$ is also contained in $W'$.  
\end{rem} 

\subsection{Terminology} 
We say that a puzzle piece is $\omega(c)$-\emph{critical} if it contains a critical point in $\omega(c)$.
Let $P$ be an $\omega(c)$-critical puzzle piece. An $\omega(c)$-critical puzzle piece $Q$ is a called a
\emph{child}\label{child} of $P$ if it is a unicritical pullback of
$P$; that is, there exists
a positive integer $n$ such that $Q$ is a component of $f^{-n}(P)$
containing a critical point in $\omega(c)$, and there exists a puzzle piece
$Q'\supset f(Q)$ such that the map $f^{n-1}\colon Q'\rightarrow P$ is a
diffeomorphism. 

A map $f$ is called \emph{persistently recurrent on $\omega(c)$}\label{persistently recurrent}
if $c$ is recurrent, non-periodic and each $\omega(c)$-critical puzzle
piece has only finitely many children. If $c$ is recurrent, but not persistently recurrent
then we say that $c$ is \emph{reluctantly recurrent}. 

It will be convenient to define $$\LL_xV:=\Comp_x f^{-n}(V)\mbox{ and } \hat \LL_xV:=\Comp_x f^{-n'}(V)
$$ where $n>0$ and $n'\ge 0$ 
are the smallest integers so that $f^n(x),f^{n'}(x)\in V$. 
Let $\rho>0$.
A puzzle piece $P$ around a persistently recurrent critical point 
is called $\rho$-\emph{nice}\label{rho-nice domain}
if for any $x\in P\cap\omega(c)$ one has
$\mod(P\setminus\mathcal{L}_x(P))\geq \rho$,
and $\rho$-\emph{free}\label{rho-free domain} if there are puzzle pieces $P^+\supset P\supset P^-$
such that $(P^+\setminus P^-)\cap\omega(c)=\emptyset$,
$\mod(P^+\setminus P)\geq \rho$ and $\mod(P\setminus P^-)\geq\rho$.
We refer to the annulus $P^+\setminus P^-$, which is disjoint from
$\omega(c)$, as \emph{free space}.
We say that a simply connected domain $U$ has 
$\rho$-\emph{bounded geometry with respect to $x\in U$}
if the Euclidian ball \mbox{$B(x,\rho\cdot\mathrm{diam}(U))\subset U$}.
A domain $U$ is said to have 
$\rho$-\emph{bounded geometry}\label{rho-bounded geometry} if there is an $x\in U$ such that $U$ has $\rho$-bounded geometry with respect to 
$x$. We say that $K'$ is a $\rho$-scaled neighbourhood of $K$ if
each component of $K'\setminus K$ has length $\ge \rho|K|$.

\subsection{Proof of Theorem~\ref{thm:pullbackcomplex}}

\begin{pf} 
This result uses the complex bounds from \cite{CvST,CvS} but is more
general in two ways: we will consider mappings  $F\colon U\to V$  which could include  
domains that do not intersect the real line and also $V$ contains
a neighbourhood of all critical points simultaneously. To construct $F$ decompose the set of critical points $\Crit_{na}(f)$ into the following  four sets: 
\label{Crr-etc}
\begin{enumerate} 
\item $\Crit_{nr}$: non-recurrent critical points;
\item $ \Crit_{rr}$: reluctantly recurrent critical points;  
\item $\Crit_{pr}$: persistently recurrent but not infinitely renormalizable critical points; 
\item $\Crit_{ren}$:  infinitely renormalizable critical points. 
\end{enumerate} 
Defining $\Crit'$ to be  $\Crit_{pr}$ or $\Crit_{nr}$, let 
 $\Crit'_2,\Crit'_1$ be the set of critical points $c\in \Crit'$ so that  $\LL_{c}I_{\Crit'}=I_{\Crit'}(c)$ resp. so that 
$\LL_{c}I_{\Crit'} \Subset I_{\Crit'}(c)$.
To prove Theorem~\ref{thm:pullbackcomplex} we need to combine the following complex bounds around each of  these sets in 
such a way that the appropriate non-real domains are included.

\begin{thm}[Complex bounds in the persistently recurrent case] \label{thm:cbpr} 
Suppose that $c\in \Crit_{pr}\cup \Crit_{ren}$. 
Then 
there exist $\rho_0>0$ and combinatorially defined intervals (puzzle pieces) $I\owns c$ of arbitrarily small diameter so that the following holds. Let
$$\hat{I}:=\bigcup_{c'\in\mathrm{Crit}(f)\cap\omega(c)}\hat{\mathcal{L}}_{c'}(I).$$
 \begin{enumerate}
\item Suppose that $f$ is non-renormalizable. Then the first return map to $\hat{I}$ extends to a complex box mapping 
$$F\colon U\rightarrow V\mbox{ so that } V\cap\mathbb{R}=\hat{I}\mbox{ and }$$
\begin{itemize}
\item for each component $U_i$ of $U$, $F|U_i$ has at most one
  critical point,
\item each component of $V$ is $\rho_0$-nice and $\rho_0$-free,
\item each component of $V$ has $\rho_0$-bounded geometry.
%\item on each component $V\in\mathcal{V}$, $F$ extends to the puzzle piece $P^{+}\supset V$ in the sense defined in Subsection \hyperref[extendible]{\ref{extendible}},
\end{itemize}
\item Suppose that $I$ is a terminating interval for $f$. Then the
  return map to $I^{\infty}$ extends to a
polynomial-like map $F\colon U\rightarrow V$ such that $\mathrm{mod}(V\setminus U)>\rho_0$.
\end{enumerate} 
\end{thm}
\begin{pf} This is Theorem 1.1 in \cite{CvST}. 
\end{pf} 

\begin{rem}\label{rem:A1}  In the first case there exists $\theta_0\in (0,\pi/2)$ 
so that 
$V_c \subset D_{\theta_0}(I_c^*)$ for each $c \in \Crit_{pr}\cup \Crit_{ren}$,
where $I_c^{*}$ is a neighbourhood of $I_c=V_c\cap \R$ for which 
$I_c^*\cap \omega(c)\subset I_c$.  Such a $\theta_0$ exists,  because 
$V_c$ is $\rho_0$-free w.r.t. $\omega(c)$ and since $V_c$ has bounded
geometry. Since $V_c$ has $\rho_0$-bounded geometry, 
we can assume that $I_c^*$ is small when $I$ is small. 
\end{rem} 

\begin{rem}\label{rem:A2}  In the 2nd case, by replacing $V$ by the disc $V'$ bounded by the 
core curve in $V\setminus U$ and replacing $U$ by the pullback $U'$
of $V'$,  we can assume that $V$ is $\rho$-nice, $\rho$-free and has
$\rho$-bounded geometry (all w.r.t. to $\omega(c)$).  
\end{rem}

\begin{thm}[Complex bounds in the reluctantly recurrent and non-recurrent case]  \label{thm:cbrr} 
Let $\Crit'$ be either $\Crit_{rr}$ or $\Crit_{nr}$. 
Then there exists $\theta_0\in (0,\pi/2)$ and arbitrarily small combinatorially defined real
neighbourhoods $I_{\Crit'}$ of $\Crit'$ such that 
\begin{itemize} 
\item the return mapping to $I_{\Crit'}$
extends to a complex box mapping
$F_{\Crit'} \colon U_{I_{\Crit'}}\rightarrow V_{I_{\Crit'}}$ where
$V_{I_{\Crit'}} \subset  \cup_{c'\in \Crit'} D_\theta(I_{\Crit'})$ 
where $\theta\in (\theta_0,\pi/2)$;
\item 
$V_{c}=D_{\pi/2}(I_{c})$ for each $c\in \Crit'_1$. 
\end{itemize} 
%a union of Poincar\'e lens domains, 
%$(\bm{V}_{I_{\Crit'}}\setminus\bm{U})_{I_{\Crit'}}\cap\mathbb H^+$ is a
%quasidisk and this complex box mapping satisfies the gap and extension
%properties.
\end{thm}
\begin{pf} In the setting that all critical points have even order this is Theorem 3 and the first line 
of the proof of  Proposition  1 in \cite{KSS-density}.
In \cite[Theorem 5.3]{CvS} this is extended to the general case. 
\end{pf}

We will also use that diffeomorphic pullback of Poincar\'e discs remain under control:

\begin{lem}
\label{lem:angle control along diffeos}
For any $\theta\in(0,\pi/2),$
there exist $\varepsilon>0$ and $\tilde\theta\in(0,\pi/2)$ such that
the following holds. Suppose that $|J_s|<\varepsilon$,  
and $f^s\colon J_0\rightarrow J_s$ is a diffeomorphism. Let
$\{J_{j}\}_{j=0}^s$ be the chain such that
$J_{j}=\Comp_{f^{j}(J_0)}f^{-(s-j)}(J_s)$.
Let $U_s=D_{\theta}(J_s)$, and set $U_{j}=\Comp_{J_j}(f^{-(s-j)}(U_s))$  for $j=0,\dots,s$. 
Then $U_j\subset D_{\tilde \theta}(J_j)$,  where 
can make the difference $\theta-\tilde\theta$ as small as we like by taking $\varepsilon>0$ sufficiently small.
\end{lem}
\begin{pf}  This builds on result  follows from  \cite{GSS} and
 \cite[Theorem B]{LS}. Under the additional 
assumption that  $J_0\cap \omega(c)\ne \emptyset$ this lemma is precisely Lemma 5.9 \cite{CvST}.
If  $J_0\cap \omega(c)\ne \emptyset$  then it follows  immediately from  Lemmas 4.4 and 4.5 in \cite{CvS}. 
\end{pf} 

The main issue in proving Theorem~\ref{thm:pullbackcomplex} is to 
pullback ranges around one type of critical point back to the range near another type of critical point. 
Note that in the above theorems the intervals $I$ are all nice intervals, whose boundaries are periodic 
or pre-periodic points.

By definition the forward orbit of $c\in \Crit_{pr}$ does not accumulate onto 
 $\Crit_{nr}\cup \Crit_{rr}\cup \Crit_{ren}$ and similarly the forward orbit of $c\in \Crit_{ren}$ does not 
 accumulate to $\Crit_{nr}\cup \Crit_{rr}\cup \Crit_{pr}$. So we can choose the neighbourhoods 
 $I_{\Crit_{nr}}$, 
 $I_{\Crit_{rr}}$,  $I_{\Crit_{pr}}$ and $I_{\Crit_{ren}}$ so small that forward iterates of 
 $c\in \Crit_{ren}$   avoid $I_{\Crit_{nr}}$, 
 $I_{\Crit_{rr}}$, $I_{\Crit_{pr}}$, $I_{c'}$ for any $c'\in \Crit$ with $c'\ne \omega(c)$ and 
 similarly so    that  forward iterates of 
 $c\in \Crit_{pr}$ avoid  $I_{\Crit_{nr}}$, 
 $I_{\Crit_{rr}}$, $I_{\Crit_{ren}}$ and  $I_{c'}$ for any $c'\in \Crit$ with $c'\ne \omega(c)$.
 
 Let $\theta_0>0$ be as in Remark~\ref{rem:A2}. 
Let $V_c$ be the corresponding ranges of the complex box mapping containing $c$
constructed in Theorems~\ref{thm:cbpr} and \ref{thm:cbrr}. 
Choose for each $c\in \Crit$ the interval $J(c)$ around $f(c)$ 
so small that $\Comp_c f^{-1}(J)\subset V_c$. The main issue now will be to deal with the case that 
a point $x$ enters a  component of $V$ which does not contain $x$.

We claim that for each $c\in \Crit_{pr}\cup \Crit_{ren}$ one can choose  $I_c$  so small 
so that % the corresponding set $\bm{V}_{I_c}$ is so small so that 
for each  $c'\in \Crit':=\Crit_{rr}\cup \Crit_{nr}$ and each 
$x\in \Comp_{c'}f^{-1}(J)$ % \subset V_{I_{\Crit'}}$ 
with $f^n(x)\in  V_{I_{c}}$ for some $n>0$, one has $\LL_{x} ( V_{I_c})\subset V_{I_{c'}}$. Similarly, we claim that 
if some forward iterate of $x\in  I_{c}$ for some $c\in \Crit_{pr}$ enters $V_{I_{c}}$ for some 
$c\in \Crit_{ren}$ then $\LL_{x} ( V_{I_c})\subset V_{I_{c'}}$. Let us prove
the first claim; the second claim goes in the same way. So take
 $x\in \Comp_{c'}f^{-1}(J) $ and assume that $f^n(x)\in V_c$ where $c\in  \Crit_{pr}\cup \Crit_{ren}$.
Then take $k<n$ maximal so that $f^k(x)\in V_{I_{\Crit'}}$ and let $k<k'\le n$ be minimal 
so that  whenever $f^i(x)$, $k'\le i \le n$ is contained in $V_{I_c}$ then it is contained
in a domain of  $F\colon U_{I_c}\rightarrow  V_{I_c}$  which intersects  $\omega(c)$.  If such an integer $k'$ does not exist then set $k'=n$ anyway. 
 By Theorem~\ref{thm:cbpr}, 
$U_{I_c}\subset V_{I_c}$ and so  we have that 
$\LL_{f^{k'}(x)} ( V_{I_c})\subset V_{I_{c}}$ if $k'<n$
and if $k'=n$ then we have by assumption $f^{k'}(x)\in V_{I_{c}}$.
By Remark~\ref{rem:A2} we have 
$$f^{k'}(x)\in V_{I_{c}}\subset D_{\theta_0}(I_c^*)$$
where $I_c^*$ is so that $(I_c^*\setminus I_c)\cap \omega(c)=\emptyset$. 

Because of the choices of the intervals
$I_c$, since  $(I_c^*\setminus I_c)\cap \omega(c)=\emptyset$
 and by the definition of $k,k'$ it follows that the pullback of $D_{\theta_0}(I_c^*)$ to $f^{k}(x)$ by $f^{-(k'-k)}$ does not 
meet any critical points from $\omega(c)$, $c\in  \Crit_{pr}\cup \Crit_{ren}$.
So the  map $f^{k'-k} \colon \LL_{f^k(x)} ( D_{\theta_0}(I_c^*)) \to D_{\theta_0}(I_c^*)$
is a diffeomorphism.   Hence, by Lemma~\ref{lem:angle control along diffeos} and by possibly 
choosing the interval $I_c$ smaller, we may assume that the set  $\LL_{f^k(x)} ( V_{I_c})\subset 
D_{\theta_0/2}( \LL_{f^k(x)}(I_c^*))$.

Let us show that for each $\rho_1>0$ there exists $\epsilon>0$ 
so that if  $|I_c|<\epsilon$ then  the $\rho_1$-scaled neighbourhood of  $\LL_{f^k(x)}(I_c^*)$ is contained in
the corresponding component of $I_{\Crit'}$. 
Indeed, if $I_c$ is small, then $I_c^*$ is small.
It follows that all preimages of $I_c^*$ are also small - this follows the Contraction Principle, see \cite[Section IV.5]{dMvS} (this principle  is related to the absence of wandering intervals).
So  $\LL_{f^k(x)}(I_c^*)$ is  small compared to $I_{c'}$ if $\epsilon>0$ is sufficiently 
small. If the assertion fails, then 
 $\LL_{f^k(x)}(I_c^*)$ is contained in a very small interval containing a boundary point of $I_{\Crit'}$. 
Note that boundary points of $I_{\Crit'} $ are pre-images of some periodic point  $p$. 
Hence some iterate $K$ of  $\LL_{f^k(x)}(I_c^*)$ is near the periodic point $p$ 
before it is mapped  onto $I_c^*$. Because of this, 
taking $\epsilon>0$ small (and therefore $I_c$ and $I_c^*$ small) will ensure
that $K$ occupies only a very small part of a fundamental domain
of the periodic point $p$. This in turn implies that  a $\rho_1$-scaled neighbourhood of 
 $\LL_{f^k(x)} (I_c^*)$ is contained in $I_{\Crit'}$, provided $\epsilon>0$ is small enough. 

From the choice of $k$ it follows that $f^k(x)$ must be contained in $I_{c'}$
with $c'\in \Crit'_1$ because the components $I_{c'}$ with $c'\in \Crit'_2$ 
are mapped onto a component $I_{c'}$ with $c'\in \Crit'_1$. Here we use that
points can only escape the box mapping around $\Crit'$ (to the persistently recurrent $c$) 
via the components $I_{c'}$ with $c'\in \Crit_1$.   Because a $\rho_1$-scaled neighbourhood
of $ \LL_{f^k(x)}(I_c^*)$ is contained in $I_{c'}$, this and 
the 2nd part of Theorem~\ref{thm:cbrr} implies that 
 $D_{\theta_0/2}( \LL_{f^k(x)}(I_c^*)) \subset  V_{I_{c'}}$
 (provided $\rho_1$ is sufficiently large).  It follows that $\LL_{f^k(x)}(V_c)\subset V_{I_{c'}}$. 
 Here,   if $k=0$, then  we use  that $J$ and therefore $J^{-1}$ can be taken arbitrarily small
 (this argument is needed and sufficient as $x$ may not be real but $f(x)\in J$ is). 
 If $k>0$ then we pullback using Theorem~\ref{thm:cbrr} till we reach $x$. 
  This shows how to pullback the ranges of the complex box mappings around the critical points in $\Crit_{pr}, \Crit_{ren}$ back
into range near $\Crit_{nr}, \Crit_{rr}$ (and also how to pullback range around $\Crit_{ren}$ back into the range
around near $\Crit_{pr}$). 

Now we show how to pullback the range near $\Crit_{nr}$ and $\Crit_{rr}$ back 
into the range near $c\in \Crit_{pr}$. This will involve choosing the neighbourhoods  $J(c)$ of $f(c)$ sufficiently small when   
$c\in \Crit_{pr}$. 

So take $x\in \Comp_c f^{-1}(J)$ and a (minimal) $n(x)$ so that $f^{n(x)}(x)\in I_{c'}$ for some $c'\in  \Crit_{rr}\cup \Crit_{nr}$.
By assumption %Since the set $\omega(c)$ is a minimal  Cantor set,  
the forward orbit of $c$ does not 
enter $I_{\Crit_{rr}}\cup I_{\Crit_{nr}}$. In particular, $f^{n(x)-1}\colon \LL_{f(x)}  V_{I_{c'}} \to   V_{I_{c'}}$ is a diffeomorphism and $n(x)$ is very large  if $J(c)$ is small. It follows that 
the real trace of  $\LL_{f(x)}  V_{I_{c'}}$ is  small when $J(c)$ is small, and so the components of $f^{-1} (\LL_{f(x)}  V_{I_{c'}})$ that are close to $c$ are all contained
in $ V_{I_{c}}$.   If $f^n(x)\in  I_{c'}$ for a non-minimal $n$ then we apply Theorem~\ref{thm:cbrr}   to $f^{n-n(x)}$
and again  $\Comp_x f^{-n}  (V_{I_{c'}})\subset  V_{I_{c}}$.  

This completes the proof of Theorem~\ref{thm:pullbackcomplex}. 
Note that although $\LL_{f(x)}  V_{I_{c'}}$ intersects the real line, 
 some components of the form $f^{-1} (\LL_{f(x)}  V_{I_{c'}})$  may be 
 disjoint from the real line. 
 \end{pf}

\subsection{Proof of Theorem~\ref{thm:KX}} 
Let $F\colon W\to W'$ be the map from 
Theorem~\ref{thm:pullbackcomplex}. This theorem shows that $K_X(f)\cap \partial W\subset \R$. 
Using Lemma~\ref{lem:angle control along diffeos} statements  (1) and the first part of (2) follow. 
In particular, for each $a'<a$ we can arrange it so that $K_X(f)\subset \Omega_{a'}$. 

{\bf $K_X(f)$ has no interior.} To see this, assume by contradiction that $K_X(f)$ contains
an open set $U$.  Then $U$ intersects $K_{N} $ for some $N$. Since $K_X(f)$
is forward invariant, it follows that there exists an open set $U'$ which intersects $I$
and which contains $K_X(f)$. However, this is impossible: consider a periodic point $p$
in $U'\cap I$ of period $k$ and construct a periodic {\lq}ray{\rq} $\gamma$ through $p$,  
i.e. so that  $f^{k}(\gamma)\supset \gamma$, and so that  $f^{k}(\gamma)\setminus \gamma$
is in the complement of $\Omega_{a'}$. In particular, 
each point  in $\gamma\setminus \{p\}$ is eventually mapped outside $\Omega_{a'}$. 
Since  $K_X(f)\subset \Omega_{a'}$ and since
$K_X(f)$ is forward invariant this implies that no point in $\gamma\setminus \{p\}$
can be contained in $K_X(f)$. This contradicts that $U'$ contains $K_X(f)$.

{\bf $K_X(f)$ is full.} Suppose by contradiction that some $\C\setminus K_X(f)$
contains a bounded component $U$. Then the $U$ is the interior of a topological disc $D$ bounded
by a set of the form
 $\overline{ \cup_{N\ge 0} \tau_N}$ where $\tau_N$ is a sequence of 
 subtrees  $\tau_{N}$  of the finite tree $K_N$  so that $\tau_{N}\subset \tau_{N+1}$ for all $N\ge 0$ and so that  
some endpoints $x_N,y_n$ of $\tau_N$ both converge to the same point $x$. 
Then, using the same argument as in the previous paragraph, some iterate 
of $D$ intersects the real line. In other words, the closure $D'$ of some component $U'$
of $\C\setminus K_X(f)$ intersects $I$ in an arc $J$. Since $K_X(f)$ is forward invariant, 
and since $f\colon I\to I$ has no wandering intervals, it follows that $U'$ must be periodic. 
Hence $D'\cap I$ contains an attracting periodic point. However, by the choice 
of the pruning intervals $J_i$, see the discussion in Example~\ref{example:K_X} of  Figure~\ref{fig:truncated2},
the set $K_N$ inside the basin of periodic attractors does not grow with $N$, 
and so the above situation cannot arise.

{\bf $K_X(f)$ is locally connected.} To show this consider the map $F\colon W\to W'$
and let $V$ be a component of $W$ and let $U$ be a pullback of $V$ by $f^n$.  
If $U$ contains a critical point $c\in \Crit_{pr}\cup \Crit_{ren}$
then $f^n$ is an iterate of the first return map of $F\colon W_{I_c} \to W'_{I_c}$. 
Since $K_X(f)\cap \partial W\subset \R$, it follows that $U\cap K_X(f)$ is connected. 
In general,  $f^n$ 
is a composition of a diffeomorphic iterate of $f$,  
an iterate of  the box mapping $F_{rr} \colon W_{I_{\Crit'}}\rightarrow W'_{I_{\Crit'}}$
from Theorem~\ref{thm:cbrr}  (possibly the identity), a diffeomorphic iterate and 
an iterate of $F\colon  W_{I_c} \to W'_{I_c}$
where $c\in \Crit_{pr}\cup \Crit_{ren}$ (possibly the identity).  
By Theorems~\ref{thm:cbpr} and \ref{thm:cbrr},  
puzzle pieces of $F_{rr} \colon W_{I_{\Crit'}}\rightarrow W'_{I_{\Crit'}}$ and of $F\colon W_{I_c} \to W'_{I_c}$ shrink in diameter as the depth
increases. Therefore, and because of Lemma~\ref{lem:angle control along diffeos}, 
it follows that the diameter of $U$ shrinks to zero as $n\to \infty$ unless there are critical points at which $f$ is infinitely renormalizable. In that case one can take a
shrinking sequence of renormalization intervals and also obtain that there exist arbitrarily small puzzle pieces around
each point. Since $\partial W$ only intersect $K_X(f)$ in the real line, it follows that each point $x\in K_X(f)$ 
is either contained in preimages of some component of $W$ of arbitrarily small diameter, 
or $x$ is real and does not accumulate to $\Crit$. In the latter case, by the (real) Ma\~n\'e theorem, see \cite{dMvS},
$x$ is contained in a hyperbolic set and again $x$ is contained in arbitrarily small neighbourhoods $W_n(x)$
so that $W_n(x)\cap K_X$ is connected. It follows that $K_X$ is locally connected. Hence there are rays landing at each point $z\in K_X(f)$. 

{\bf At most finitely many rays lands at any point of $K_X(f)$.}  Notice that 
$K_X(f)$ is the closure of a nested sequence of trees $K_{n+1}\supset K_n$. 
Let $z\in K_X(f)$ and assume that there are two rays $\gamma_1,\gamma_2$ landing on $z$. Then 
both components of $\C\setminus (\gamma_1\cup \{z\}\cup \gamma_2)$  have a non-empty 
intersection with $K_X$ and therefore with $K_n$ for some $n$.  As $K_n$ is a tree it follows 
that $z\in K_n$.  Moreover,  $K_{n+1}\setminus K_n$ does not contain arcs that are connected to endpoints of $K_n$. It follows that  if there are an infinite number of distinct  rays $\gamma_i$ landing on $z$, 
then when $i\ne j$ each component of $\C\setminus (\gamma_i\cup \{z\}\cup \gamma_j)$  
has a non-empty intersection with $K_n$ (where $n$ does not depend on $i,j$). 
Since $K_n$ is a finite tree with finite degree, this is impossible. This proves assertion (3).

{\bf Each non-real periodic point in $K_X(f)$ is repelling}. Since each non-real point $z\in K_X$ is contained in a sequence of puzzle pieces which lie nested, if $z$ is periodic then there exist topological discs $D'\supsetneq D\ni z$ with $D\cap \R=\emptyset$ 
and $n$ so that $f^n(z)=z$ and $f^n$ maps $D$ diffeomorphically to $D'$. Hence by the Schwarz Lemma
$z$ is a hyperbolic periodic point. \qed

\section{Continuity of the pruned Julia set} 

\begin{defn}
Let $(U_n,u_n)$ and $(U,u)$ be pointed open discs. We say that
$(U_n,u_n)\to (U,u)$ in the  {\em sense of  Carath\'eodory}  if 
(i) $u_n\to u$,  (ii) for each compact $K\subset U$, $K\subset U_n$ holds
for $n$ large and (iii) for any open connected set $N$ containing $u$, 
if $N\subset U_n$ for infinitely many $n$ then $N\subset U$. 
\end{defn}

Let $f\in \mathcal A^{\underline \nu}_a$ and let $f_n\in  \mathcal A^{\underline \nu}_a$
be so that  $f_n\to f$  on $\overline{\Omega_a}$.  Moreover, take 
the intervals $J_{n,i}\ni f_n(c_i)$  so small that Theorem~\ref{thm:KX} holds
and assume that $J_{n,i}\to J_i$ (i.e. the boundary points of $J_{n,i}$ converge 
to the boundary points of $J_i$).  Let $X,X_n$ correspond to $\partial J^{-1}$
and $\partial J^{-1}_n$ where $J_n=f_n^{-1}(J)\setminus \R$. 
Similarly choose $J_{n,i}^*\Subset J_{n,i}$ so that $J_{n,i}^*\ni f_n(c_i)$ and so that the boundary points
of $J_{n,i}^*$ are either (pre)periodic or contained in the basin of a periodic attractor. 
Pick a point $x\in \C\setminus K_X(f)$ and $x_n\in  \C\setminus K_{X_n}(f_n)$.

As before let $\psi_X\colon \overline \C \setminus \overline{ \mathbb D} \to \overline \C
\setminus K_X$ be the uniformising map that fixes $\infty$ and with real
derivative at $\infty,$ and write $\phi_X=\psi_X^{-1}\colon \overline \C \setminus
K_X\to \overline \C \setminus \overline{ \mathbb D}$. 

\begin{thm} \label{thm:KXcont} 
We have the following:
\begin{enumerate}
\item $K_{X_n}(f_n)\to K_X(f_n)$ in the Hausdorff topology;
\item $(\C\setminus K_{X_n}(f_n),x_n)$ converges in the Caratheodory topology to 
$(\C\setminus K_X(f),x)$; 
\item  $\phi_{X_n}\to \phi_X $; 
\item $\F_{X_n}\to \F_X$ on a neighbourhood of $\partial \D\setminus \hat J^*$ 
where $\hat J$ is the set corresponding to $J^*$. 
\end{enumerate}
\end{thm} 
\begin{pf} Statement (1) follows from the proof in the previous Appendix. Statement (2) 
follows from the definition of Carath\'eodory convergence. Statement (3) follows from 
Carath\'eodory's kernel theorem, see \cite{Pom}. Statement (4) follows from (3). 
\end{pf}

%\documentclass[11pt]{article}
%\documentclass[12pt, a4,oneside]{amsart}
%\let\latexput\put
%\usepackage{amsmath,amsfonts,amssymb,amsthm,bm,wasysym,url}%mathtext

%
%\begin{document}
%\input{macros.tex}
%\newtheorem{question}{Question}
%\newtheorem{defi}{Definition}
%\newtheorem{lemi}{Lemma}
%\newtheorem{exi}{Example}
%\newtheorem{remi}{Remark}
%\newtheorem{cori}{Corollary}
%\title{}
%\author{}
%\date{}
%\maketitle

\section{A topological and analytic structure on the space of real analytic functions}
\label{Append:toponA}

 Let $\mathcal A$ be the space of real  analytic functions on $I=[-1,1]$.
 Following \cite[Appendix 2]{Lyubich}, we will define in this section the real analytic topology  and an analytic 
 structure on $\mathcal A$.
This  space can be viewed as a space of germs in the following way:
  
 \begin{defn}[$\mathcal A$ as a space of germs]
 Let $\mathcal U$ be a collection of open sets in $\C$ containing $[-1,1]$
so that  for each $a>0$ there exists an open set $U\in \mathcal U$ so that $I\subset U\subset \Omega_a:=\{z\in \mathbb C: \dist(z,I)<a\}$. 
% and assume that for each $V$ neighbourhood of $I$ there exists $U\in \mathcal V$
% so that $U\subset V$. (For example $\mathcal U$ could be the collection 
% of open sets of the 
%orm  %(with $\mbox{dist}(x,y)$ say the Euclidean distance between $x) 
 For $U\in \mathcal U$, let $\mathcal A_U$ be   the Banach space of holomorphic maps on $U$ which extend continuously 
 to $\overline U$. If $U\subset V$ then define 
 $i_{U,V}\colon \mathcal A_V\to \mathcal A_U$ as the restriction function (for $f\in \mathcal A_V$ let $i_{U,V}f=f|U$). 
 We say that $f_1\in \mathcal A_U$ and $f_2\in \mathcal A_V$  are {\em equivalent} if there exists $W\in \mathcal U$
 with $W\subset U\cap V$
 so that $i_{W,U}f_1=i_{W,V}f_2$. Equivalence classes are called {\em germs} and the space
 of such germs is called the inductive limit of the Banach spaces $\mathcal A_U$. 
 We say that $f\in \mathcal A$ is in $\mathcal A_U$ if some representative of $f$
 is in  $\mathcal A_U$. 
%, and so in other words if $f$ extends to a map ini $\mathcal A_a$. 
 \end{defn} 

\begin{rem}\label{rem:conti}  $i_{U,V}\colon \mathcal A_V \to \mathcal A_U$ is continuous 
(when $U\subset V$). 
\end{rem} 
 
  \begin{defn}[The inverse limit topology on $\mathcal A$] 
In the  {\em $C^\omega$ topology}  (also called the real-analytic topology, or the inverse limit topology)  on $\mathcal A$ 
a set $\mathcal O$ is defined to be open if and only if $\mathcal O\cap \mathcal A_U$ is open for any $U\in \mathcal U$
 (in the Banach space topology 
on  $\mathcal A_U$). 
\end{defn}

Because of the assumption on $\mathcal U$ and since $i_{U,V}$ is continuous, 
the $C^\omega$ topology on $\mathcal A$ does {\em not depend} on the choice of $\mathcal U$.
It is often convenient to consider the collection $\mathcal U$ of open sets of the form $\Omega_a$
and to define $\mathcal A_a :=\mathcal A_{\Omega_a}$. 
%% provided 
%%for each $a>0$ there exists an open set $U\in \mathcal U$ so that $I\subset U\subset \Omega_a=\{z\in \mathbb C: \dist(z,I)<a\}$. 
%\end{lemi} 
%\begin{pf} This follows from Remark~\ref{rem:conti}. 
%\end{pf} 
%%If we would instead  of taken the collection $\mathcal U$ of all open sets  $U\supset I$
%It follows from Remark~\ref{rem:conti} that we would end up with the same topology on $\mathcal A$ if, 
%in the previous definition, we would take  an arbitrary collection $\mathcal U$ of open sets $U\subset I$ so that for any $a>0$ 
%there exists $U\in \mathcal U$ with $U\subset \Omega_a=\{z\in \mathbb C: \dist(z,I)<a\}$. 
% We say that $f\in \mathcal A$ is in $\mathcal A_U$ if some representative of $f$
% is in  $\mathcal A_U$. 
%, and so in other words if $f$ extends to a map ini $\mathcal A_a$. 
  
  For later use,  let $\mathcal A_U(f,R)$ be the $R$-ball around $f$ in $\mathcal A_U$, i.e., the set of 
 functions  $g\in \mathcal A_U$ so that $||g-f||_U:=\sup_{U}|g(z)-f(z)|\le R$.

The above topology is Hausdorff. The next lemma is  
part of \cite[Lemma 11.4]{Lyubich} (and is a special case of general results
for the inverse limit topology), and gives some desirable properties of this topology: 

\begin{lem}  \label{lem:compact} \leavevmode %The following holds: 
\begin{itemize}
\item 
For  $f_n,f\in \mathcal A$ we have 
$f_n\to f$ in the $C^\omega$ topology if and only if there exists $U\in \mathcal U$ so that $f_n,f\in \mathcal A_U$ and $f_n\to f$ on $\overline U$.
\item  If $\mathcal K\subset \mathcal A$ is compact then 
there exists $U\in \mathcal U$ so that $\mathcal K \subset \mathcal A_U$.
 \end{itemize} 
\end{lem} 
\begin{pf}  Let us prove this for the collection $\mathcal U$  of open sets $\Omega_a$, $a>0$. 
Assume $f_n\to f$ and assume (by contradiction) there exist $a_n\downarrow 0$
so that  $f_n\in \mathcal A_{a_n}\setminus \mathcal A_{a_{n-1}}$.  We may assume that $f_n\ne f$ for all $n$.
Then $\mathcal O =   \mathcal A \setminus \{f_1,f_2,\dots\}$ is open neighbourhood of $f$, which gives a contradiction.  
So let us prove the 2nd assertion: if $\mathcal K$ is compact, 
each infinite sequence $f_n\in \mathcal K$ must  have a cluster point. However, if the claimed assertion 
does not hold, then there exists $f_n$ and $a_n\downarrow 0$ so that $f_n\in \mathcal A\setminus \mathcal A_{a_n}$. 
By the first part of the  lemma such a sequence can't have a cluster point. 
\end{pf} 

On the other hand, 
\begin{lem}\leavevmode
\begin{itemize}
\item Let $U\subsetneq V$. Then $i_{U,V}\mathcal A_V$ is not an open subset of $ \mathcal A_U$.
\item For any $f\in \mathcal A_V$,  any $U\subsetneq V$ and any open set $\mathcal O\ni f$ 
in the real analytic topology, there exists  $g\in  \mathcal O \setminus \mathcal A_U$.
\end{itemize} 
\end{lem} 
\begin{pf} If $f\in \mathcal A_V$ then $\tilde f=i_{U,V}f \in \mathcal A_{U}$. 
Each open neighbourhood of $\tilde f$ in  $\mathcal A_{U}$  contains maps which are not
analytic on $V$.  
Hence  $i_{U,V}\mathcal A_V$ is not open in $\mathcal A_U$ (and so not open in $\mathcal A$).

The 2nd statement follows from this: since $f\in  \mathcal A_U$ implies  $f\in  \mathcal A_W$ for 
any $W\subsetneq U \subsetneq V$, 
the previous assertion implies that there exists a map $g$ in $\mathcal O$ which is not contained
in  $\mathcal A_U$. 
\end{pf} 
 
%Moreover, for  $f_n,f\in \mathcal A^{\underline \nu}$
%we say $f_n\to f$ if there exists $a>0$ so that $f_n,f\in \mathcal A^{\underline \nu}_a$ and $f_n\to f$ on $\overline \Omega_a$. 
%We say that $\mathcal O\subset  \mathcal A^{\underline \nu}$ is open 
%if for each $f\in \mathcal O$ and each $f_n\in \mathcal A^{\underline \nu}$ with $f_n\to f$
%we have that $f_n\in \mathcal O$ for all but finitely many $n$. 

Obviously when $U\subsetneq V$ are open subsets, then the following properties hold, 
cf. Properties C1-C2, P1-P3 in \cite[p412, p416 and p342]{Lyubich}: 
\begin{itemize}
\item[C1.] The image of $i_{U,V}\mathcal A_V$ is dense in $\mathcal A_U$. 
\item[C2.] the map $i_{U,V}$ is compact, i.e.  $i_{U,V} \mathcal A_V(f,R)$
pre-compact in $\mathcal A_U$ for any $R>0$. 
%Each open neighbourhood  $\mathcal O$ of $f\in \mathcal A$ contains, for each $a'>0$, 
%maps which are not analytic in $\Omega_{a'}$. % and so are not contained  in $\mathcal A_{a'}$. 
 \end{itemize}

The following result is implicitly contained in more general results
on inverse limit topology, see  \cite{Engel}. 

\begin{lem}\label{lem:partitionunity}  $\mathcal A$ is not locally compact,  not metrizable and not Baire.
However it is regular,  paracompact,  Lindel\"of and admits a partition of unity. 
 \end{lem}
We recall that a topological space is called {\em locally compact} if each point $x$ has a compact neighbourhood, i.e. there
 exists an open set $U$ and a compact set $K$ so that $x\in U \subset K$. It is {\em Baire} if
 the countable intersection of open and dense sets is dense. We say that a topological space is 
{\em regular} (or $T_3$) if any point and  closed set can be separated by open sets. 
It is {\em paracompact} if each cover contains a subcover which is locally finite. 
It is {\em  Lindel\"of} if each cover contains a countable subcover. 
%It is {\em normal } if two closed sets can be separated by open sets. 
  \begin{pf} According to the previous lemma, for each open set $\mathcal O\subset \A$ containing $f$ and each $a>0$ there exists 
$f_n\in \mathcal O\setminus \mathcal A_a$. By Lemma~\ref{lem:compact} this implies that no set $K\supset U$ can be compact. 
To see that $\mathcal A$ is not metrizable, assume that $d$ is a metric on $\mathcal A$. Take $f\in \mathcal A_a$
and a sequence of open balls $\mathcal O_n\ni f$ with diameter going to zero. 
Then, by the previous lemma,  $\mathcal O_n$ contains  some $f_n\in \mathcal A\setminus \mathcal A_{1/n}$.  So $d(f_n,f)\to 0$ but  there exists no
$\epsilon>0$ so that every $f_n$ is contained in $\mathcal A_\epsilon$, contradicting  Lemma~\ref{lem:compact}. 
 To see that  $\mathcal A$ is not Baire, choose  $a_n\downarrow 0$ and $f_n\in \mathcal A_{a_n}\setminus \mathcal A_{a_{n-1}}$ with $f_n\ne 0$. 
Then the  set $\mathcal O_n=\{f\in \mathcal A ; f\ne f_n\}$ is open and dense in  the real analytic topology. 
On the other hand,  $\cap \mathcal O_n$ is not dense, because  $\mathcal O =   \mathcal A \setminus \{f_1,f_2,\dots\}$ 
is an open set (containing $f\equiv 0$) which is disjoint from each $\mathcal O_n$. 

To see that $\mathcal A$ is regular, take $f\in \mathcal A_{a_0}$  and let $K_1=\{f\}$ 
and some closed set $K_2$ not containing $f$. Since $f\notin K_2$, $K_2$ is contained in a closed set of the form 
$K_2'=\cup_{0<a<a_0} \{g\in \mathcal A_a ;  \sup_{z\in \Omega_a}|g(z)-f(x)|\ge t(a)\}$ where $t(a)>0$. 
Then $U_1=\cup_{0<a<a_0} \{g\in \mathcal A_a ; \sup|g(z)-f(z)|<(1/4)t(a)\}$  and 
$U_2=\cup_{0<a<a_0} \{g\in \mathcal A_a ; \sup|g(z)-f(z)|\ge (3/4)t(a)\}$
are open disjoint sets containing $K_1$ and $K_2$, proving that $\mathcal A$ is regular. 

Let us next prove that $\mathcal A$ is the countable union of compact subsets. 
Indeed, for each $f\in \mathcal A$ there exists $n,k\in \mathbb N$ so that $f\in \mathcal A_{1/n}(0,k)$
and so in particular  $f\in i_{\Omega_{1/2n},\Omega_{1/n}}\mathcal A_{1/n}(0,k)$. It follows that 
$\mathcal A = \bigcup_{k,n} i_{\Omega_{1/2n},\Omega_{1/n}}\mathcal A_{1/n}(0,k)$. 
By Property C1 the set $i_{\Omega_{1/2n},\Omega_{1/n}}\mathcal A_{1/n}(0,k)$ is a precompact subset of $\mathcal A_{1/2n}$. 
So the closure $C_{n,k}$  of $i_{\Omega_{1/2n},\Omega_{1/n}}\mathcal A_{1/n}(0,k)$ in $\mathcal A_{1/2n}$
is a compact subset of  $\mathcal A_{1/2n}$. Let $\mathcal W$ be a union of open subsets of $\mathcal A$
covering $C_{n,k}$. Then by definition of the topology on $\mathcal A$, 
the sets  $W\cap \mathcal A_{1/2n}$, $W \in \mathcal W$ are open in $\mathcal A_{1/2n}$ and
cover $C_{n,k}$. It follows that a finite collection in $\mathcal W$ already covers $C_{n,k}\subset \mathcal A$, 
and so $C_{n,k}$ is also compact as a subset of $\mathcal A$. It follows that 
$\mathcal A$ is the countable union of compact subsets. 

Because of this and since $\mathcal A$ is regular,  it is a  Lindel\"of space, \cite[page 192]{Engel}. 
Hence it is  paracompact, see \cite[Theorem 5.1.1 on page 300]{Engel}.  
Therefore any open cover of  $\mathcal A$ has a partition of unity subordinate to it, see \cite[Theorem 5.1.9 on page 301]{Engel}. 
\end{pf} 

%It follows that 
%$\mathcal A$ is not second countable (i.e. does not have a local countable base of neighbourhoods), because if it was it would be metrizable by Urysohn's theorem (since $\mathcal A$ 
%is Hausdorff and regular). 
%
%Is $\mathcal A_a$ second countable? For a metric space the properties of 
%second-countable, separable, and Lindelöf are all equivalent.
%
%It is known that $H^\infty$, the space of functions on $\mathbb D$ with bounded supremum 
%norm is non-separable, see \url{https://math.stackexchange.com/questions/1689215/h-infty-is-not-separable}. 
%On the other hand, if one restricts to a compact subset one has separable: 
%\url{https://mathoverflow.net/questions/366683/supremum-over-which-sets-makes-h-infty-non-separable}. 

%\begin{exi} 
\begin{lem} The $C^\omega$ topology on $\mathcal A$ is finer than the $C^\infty$ topology on $\mathcal A$. 
\end{lem} 
\begin{pf} It is sufficient to show that there exists an open set $\mathcal O$
in the real analytic topology which is contained in the set $\{g; \sup_{x\in I} |D^kg(x)|<\epsilon\}$. 
So choose $\epsilon(a)>0$  so small that if $g\in \mathcal A_a$ is so that $|g(x)|<\epsilon(a)$ on $\overline \Omega_a$
then $\{g; \sup_{x\in I} |D^kg(x)|<\epsilon\}$. It follows that the open set 
$\mathcal O=\cup_a \{g\in \mathcal A_a; \sup_{x\in \Omega_a} |g(x)|<\epsilon(a)\}$ in the real analytic topology
has the desired properties,
\end{pf}

%\end{exi} 
%and that  $f\in \mathcal A^{\underline \nu}_a$ can be approximated
%(in the $C^k$ topology on $[-1,1]$) by real analytic maps
%$f_i\in \mathcal A^{\underline \nu}_{a_i}$ which are not 
%contained in $f_i\in \mathcal A^{\underline \nu}_{a'_i}$ where $0<a_i<a_i'$
%is so that $a_i'\to 0$. 
%The space $\mathcal A^{\underline \nu}$ can be viewed 
%as {\em germs} of interval maps, in the usual sense,  and $\mathcal A^{\underline \nu}=\lim \mathcal A^{\underline \nu}_a$. 

%\begin{rem} 
%\begin{itemize}
%\item regular Lindelöf space is paracompact
%\item regular means: 
%\item According to a theorem by Stone, any metric space is paracompact. 
%\item A locally compact Hausdorff second-countable space is paracompact.
%\item Any paracompact Hausdorff space admits a partition of unity. 
%\end{itemize} 
%\end{rem} 
%
%Question: Is $\mathcal A$ normal? (two disjoint closed sets
%have disjoint open neighbourhoods). 
%
%
%\bigskip 
%
%According to Lyubich: $\mathcal A$ does not have a local countable base of neighbourhoods, 
%$\mathcal A$ is not metrizable, is not a Frechet space. 
%NO local countable base of neighbourhoods 
%
%Frechet space: a locally convex metrizable topological vector space (TVS) that is complete as a TVS,
%Generalizes Banach spaces. 
%
%second countable: there exists a countable base  of neighbourhoods. 
%
%Urysohn's metrization theorem states that every second-countable, Hausdorff regular space is metrizable.
%
%
%

\begin{defn} 
$\mathcal M\subset \mathcal A$ is a  {\em real analytic manifold  modelled on a family of Banach spaces}, 
or simply {\em real analytic manifold}, 
if $\mathcal M$ is the union of the image of a family of injections $j_V\colon \mathcal O_U\to \mathcal M$
where  $U\in \mathcal U$ and where $\mathcal O_U$ is an open subset of  the Banach space $\mathcal A_U$. 
The set  $j_V( \mathcal O_U)$ is called a {\em Banach slice} of $\mathcal M$. 
\end{defn} 

In this paper we will choose in this definition for $\mathcal U$ one of the following collections of open sets: 
\begin{enumerate}
%\item all open sets $U\supset I$, or  
\item all sets of the form $\Omega_a$, or  
\item domains $U$ of pruned polynomial-like extensions of some real analytic $f\colon I\to \R$. 
\end{enumerate} 

\begin{defn}  A set $\mathcal X\subset \mathcal A$ is called {\em an immersed submanifold} if there exists a real  analytic manifold
$\mathcal M$ and a  map $i\colon \mathcal M\to  \mathcal A$  so that $Di(m)$
is a linear homeomorphism onto its range, and if $\mathcal X=i(\mathcal M)$. We say that 
$\mathcal X$ is a {\em embedded manifold} if $i\colon \mathcal M \to \mathcal X$ is a homeomorphism with the topology on $\mathcal X$ coming from the one on $\mathcal A$.
\end{defn} 

\begin{rem} The inclusion map $i\colon \mathcal M\to \mathcal A$ has the properties 
P1-P3 from \cite[Appendix 2]{Lyubich}. This implies that $\mathcal X$ can be given an intrinsic  
manifold structure and that the tangent space and codimension on $\mathcal X$ is well-defined, 
see~\cite{Lyubich}.
\end{rem} 

\begin{rem} Lemma~\ref{lem:compact} implies that $f_\lambda$, $\lambda \in \R^k$ is a family of
maps in $\mathcal A$ which depending continuously on $\lambda$ which is contained  an immersed submanifold 
$\mathcal X$ then for $\lambda\approx \lambda_0$ we have that $f_\lambda$ is contained in a real analytic Banach manifold. 
\end{rem} 

%\begin{question} 
%I would like to better understand the difference between this real analytic topology on $\mathcal A$ and the 
%$C^\infty$ topology on $\mathcal A$. 
%\end{question} 
%
%
%
%\begin{question} 
%In Section 42 of  the book {\em Convenient setting of global infinite dimensional analysis}
%by Kriegl and Michor they also define an analytic manifold structure on $\mathcal A$ which seems
%to come from the $C^\infty$ topology. How to understand the difference between 
%these manifold structures on $\mathcal A$? 
%\end{question} 

%\begin{thebibliography}{CCCCC}
%\bibitem[Eng]{Engel}
%R. Engelking, {\em General topology}, Translated from the Polish by the author. Second edition. Sigma Series in Pure Mathematics, 6. 
%Heldermann Verlag, Berlin, 1989. viii+529 pp.
%\bibitem[Lyu]{Lyubich} M. Lyubich, {\em Feigenbaum-Coullet-Tresser universality and Milnor's hairiness conjecture,} Ann, of Math. 
%  \textbf{149} (1999), 319-420. 
%\end{thebibliography} 

   %

\newpage

\section{Summary of notation and definitions} 

In general,  sets in $\C$ are denoted by symbols such as $U,U',V,V,O$
whereas functions spaces are denoted by symbols such as $\mathcal A, \mathcal B, \mathcal T, \mathcal H, 
\mathcal U$ etc. 
\bigskip 
\bigskip 

\noindent

\begin{tabular}{cc}
\begin{tabular}{|l|l|}
\hline
$\mathcal A^{\underline \nu}$, $\mathcal A^{\underline \nu}_a$,  $\Omega_a$  &  p. \pageref{page1}  \\
 \hline 
topological conjugacy class $\mathcal T_{f}$ & p. \pageref{Tf}  \\
%\hline 
% &  p. \pageref{notation} \\ 
%\hline 
% &  p.  \pageref{notation} \\
%\hline 
%%$\Crit(f)$  &  p. \pageref{notation} \\
%%\hline 
%Lens domain, Poincar\'e disc  & p. \pageref{notation} \\
\hline 
pruned polynomial-like& \\
 mapping $F\colon U\to U'$, $\Gamma_f$ & p. \pageref{prunedpollike}    \\
%\hline 
% &p. \pageref{prunedpollike}  \\
\hline 
pruned filled Julia set $K_f$ &p.~\pageref{prunedpollike}  \\ 
\hline 
notation $f\colon I\to I$, $F\colon U\to U'$ & \\
vs. $\hat f_X\colon \partial \D \to 
\partial \D$, $\hat F_X\colon \E \to  \E'$ & 
 p.  \pageref{notationfandfX}  \\ 
\hline 
$X_0,X$ and pruning intervals $J_i$ & p.~\pageref{XJi}  \\
\hline 
pruned Julia set 
$K_X(f)$ & p.~\pageref{KXf}  \\
%\hline 
%$\Crit_{nr},\Crit_{rr},\Crit_{pr},\Crit_{ren}$ & p.~\pageref{Crr-etc} \\
\hline 
external map $\hat f_X\colon \D\to \D$, $\hat X$ & p.~\pageref{subsec:ppl} \\
 \hline 
sign $\epsilon(f)$ & p.~\pageref{def:signf} \\ 
\hline 
$\hat X^*$, $Y$, $\hat Y$ and $J_i^*$ & p.~\pageref{subsec:furtherpruning}  \\
\hline 
$\Lambda_N$, $\Lambda_\infty'$ & p.~\pageref{defLambdaN},\pageref{def:Lambda'} \\
% \hline 
%the landing set $L_f$  & p.~\pageref{Lf} \\
\hline 
expanding Markov structure &  \\
$\hat F_X \colon {\V}\to {\V'}$ & p.~\pageref{def:expmarkov}  \\
\hline 
rays \& roofs/equipotential & p.~\pageref{def:ray-equipotential}  \\
\hline 
admissible pruning set $Q$ & p.~\pageref{def:admispruning}  \\
\hline 
pruning data and equivalence of & \\
pruned polynomial-like mappings &p.~\pageref{pruning-equivalence}   \\
\hline 
\end{tabular} &
\begin{tabular}{|l|l|}
%\hline 
%$\Lambda$ compact forward invariant set & p.~\pageref{lambda} 
% \\
%\hline 
%Markov structure of $\hat f_X\colon V\to \tilde V$ & p.~\pageref{lambda} 
%\\
%
%\hline 
%admissible pruning set $Q$ & p.~\pageref{def:admispruning}  \\
\hline 
attracting structure $F\colon B\to B'$ & p.~\pageref{polstructure-attractors} \\ 
\hline
global pruned polynomial-like maps & \\
with attractors & p.~\pageref{def:globalpruned}   \\
$F\colon \E\cup B\to \E'\cup B'$ & \\ 
\hline
${\mathcal B}_a^{\underline \nu} $  &  p.  \pageref{defBnua} \\
\hline 
simple parabolic point & p.~\pageref{def:simplepara} \\ 
\hline hybrid conjugacy class $\mathcal H^{\R}_f$ and $\mathcal H_{F}$ & p. \pageref{Hf},\pageref{HF}  \\
\hline
Real analytic topology  & p.~\pageref{realanalytictopology}  \\
\hline
Carath\'eodory topology  & p.~\pageref{CaratheodoryTopology}  \\
\hline 
manifold structure on $\mathcal A^{\underline \nu}$ & p. \pageref{manifold-structure},\pageref{Append:toponA}   \\
\hline 
hyperbolic, semi-hyperbolic map & p.~\pageref{def:semihyperbolic} \\
\hline
transversal vector field to & \\ topological conjugacy class &   p. \pageref{subsec:condhybrid} \\
of a hyperbolic map &  \\
\hline 
quasiconformal vector field &  p. \pageref{qcvectorfield}  \\
%\hline  equivariant, lift vector field  & p. \pageref{def:liftvector}   \\
\hline  horizontal vector field   & p. \pageref{def:horvector} \\
\hline 
vertical and transversal vector field &  p. \pageref{def:verticalvect}  \\ 
 \hline 
quasiconformal motion & p. \pageref{def:qcmotion}  \\
\hline 
Real-analytic manifolds & p. \pageref{subsec:real-analytic} \\
\hline 
%neighbourhood $Q$ in the Fatou coordinates & page 35, Figure 5
%
%$\tau_f$ p33
%
%$\Theta_t$ transition map p32
%
%$S_{\pm, f}$ crescent p32
%
%$\mathcal H_n$ space of maps with transition time $n$ p47
\end{tabular}
\end{tabular}


\newpage 

\begin{thebibliography}{CCCCC}

%\bibitem[AH]{AH} J.M. Anderson and A. Hinkkanen, 
%  {\em Quasiconformal self-mappings with smooth boundary values}, 
%Bull. London Math. Soc. {\bf 26} (1994), no. 6, 549-556. 

%\bibitem[Dyn]{Dyn} E. Dyn'kin, {\em Estimates for asymptotically conformal mappings}, 
%  Ann. Acad. Sci. Fenn. Math. {\bf 22} (1997), no. 2, 275-304.

   %      Let ? be a quasiconformal automorphism of the unit disk ? with boundary values ? on ?? and complex dilatation ?(z). Suppose esssup{|?(z)|: 1?t?|z|<1}=O(K(t)) as t?0, where K(t)?0 as t?0. The authors obtain various smoothness properties of ? while imposing upon K(t) some additional restrictions. Thus K(t)=t? for some ?>0 implies ??Lip? for all ?<?+1 (Theorem 1), whereas limt?0t?nK(t)=0 for all n?N results in ??C? (Theorem 2).

\bibitem[AB]{AB} L. Ahlfors and L. Bers, {\em Riemann mapping theorem for variable metrics}, Ann. Math.
{\bf 72} (1960), 385--404.

\bibitem[AE]{AE}  F.D. Ancel and R.D.  Edwards, {\em Is a monotone union of contractible open sets contractible?} Topology Appl. 214 (2016), 89--93. 

 \bibitem[AIM]{AIM} K. Astala, T. Iwaniec and G, Martin,  {\em Elliptic partial differential equations and quasiconformal mappings in the plane},  Princeton Mathematical Series, 48. Princeton University Press, Princeton, NJ, 2009. xviii+677 pp.


\bibitem[AL]{AL} A. Avila and M. Lyubich,  {\em The full renormalization horseshoe for unimodal maps of higher degree: exponential contraction along hybrid classes}, Publ. Math. Inst. Hautes \'Etudes Sci. No. 114 (2011), 171-223.

\bibitem[ALM]{ALM} A. Avila, M. Lyubich and W. de Melo, {\em Regular or stochastic dynamics in real analytic families of unimodal maps}, Invent. Math. \textbf{154} (2003), 451-550.


%\bibitem[Av]{Av} A. Avila, Private communication. 

\bibitem[BM]{BM} E. Bierstone and P.D. Milman, {\em Semianalytic and subanalytic sets}, 
Publ. Math. Inst. Hautes \'Etudes Sci. {\bf 67} (1988),5--42. 

\bibitem[BvS]{BvS} H. Bruin and S. van Strien, {\em Monotonicity of entropy for real multimodal maps},
 J. Amer. Math. Soc. 28 (2015), no. 1, 1--61.

\bibitem[BR]{BR} 
  L. Bers and H.L. Royden,
{\em Holomorphic families of injections,}
Acta Math. {\bf 157} (1986), 259--286

\bibitem[ChvS]{CvS} D. Cheraghi and S. van Strien, 
{\em Monotonicity of entropy in a three dimensional family of unimodal maps}, 
manuscript. 


%\bibitem[Bruin]{Bruin}
%H.  Bruin, {\em  Non-monotonicity of entropy of interval maps}. Phys. Lett. A 202 (1995), no. 5-6, 359--362.

\bibitem[C]{C} T. Clark,
  {\em Regular or stochastic dynamics in families of higher degree 
  unimodal maps,} 
Ergodic Theory Dynam. Systems, 
{\bf 34}(5) (2014), 1538--1556. 

\bibitem[CDKS]{CDKS} T. Clark, K. Drach, O. Kozlovski and S. van Strien,
  {\em The dynamics of complex box mappings}, To appear in Arnol'd Math. J.  {\bf 8} (2022), no. 2, 319--410.
  
\bibitem[CvS]{CvS} T. Clark and S. van Strien,
  {\em Quasi-symmetric rigidity},  \url{https://arxiv.org/abs/1805.09284}, Submitted for publication. 

\bibitem[CvS2]{CvS2} T. Clark and S. van Strien, 
{\em  Conjugacy classes laminate  $S$-unimodal  maps and applications}, In preparation. 
 
\bibitem[CvST]{CvST} T. Clark, S. van Strien and S. Trejo, 
  {\em Complex bounds for real maps}, Commun. Math. Phys. {\bf 355} (2017), 1001-1119.
  
  \bibitem[CDKvS]{CDKvS} 
T. Clark, K. Drach, O. Kozlovski and S. van Strien, {\em The dynamics of complex box mappings}, Arnold Mathematical Journal (2022) 8:319--410.

\bibitem[DE]{DE} A. Douady and C.J. Earle, {\em  Conformally natural extension of homeomorphisms of the circle}, Acta Math. {\bf  157} (1986) 23-48. 
  
 \bibitem[DH]{DH} A. Douady and J.H.  Hubbard {\em On the dynamics of polynomial-like mappings},
  Ann. Sci. \'Ecole Norm. Sup. (4) {\bf 18} (1985)  287--343.


%\bibitem[Kol]{Kol} S.  Koljada, {\em One-parameter families of mappings of an interval with a negative Schwartzian derivative which violate the monotonicity of bifurcations.}  Ukr. Mat. Z. {\bf 41} (1989) 258. 
%


\bibitem[Eng]{Engel}
R. Engelking, {\em General topology}, Translated from the Polish by the author. Second edition. Sigma Series in Pure Mathematics, 6. 
Heldermann Verlag, Berlin, 1989. viii+529 pp.

\bibitem[Eps]{Epstein}  A. Epstein, {\em Transversality in holomorphic dynamics}, \url{http://homepages.warwick.ac.uk/~mases/Transversality.pdf} and slides of talks available in
\url{https://icerm.brown.edu/materials/Slides/sp-s12-w1/Transversality_Principles_in_Holomorphic_Dynamics_\%5D_Adam_Epstein,_University_of_Warwick.pdf}


\bibitem[dFdM]{dFdM} E. de Faria and W. de Melo, {\em Rigidity of critical   circle mappings I},
J. Eur. Math. {\bf 13}, 2,  (2000),  343--370.

\bibitem[dFdM2]{dFdM2} E. de Faria and W. de Melo, {\em Rigidity of critical   circle mappings II},
J. Amer. Math. Soc. {\bf 1}, 4,  (1999), 339--392. 

%\bibitem[Dyn]{Dyn} E. Dyn'kin, {\em Estimates for asymptotically conformal mappings}, 
%  Ann. Acad. Sci. Fenn. Math. {\bf 22} (1997), no. 2, 275-304.

\bibitem[GS2]{GS2} J. Graczyk and G. \'Swi\c{a}tek, {\em Generic hyperbolicity in the logistic family}, Ann. of Math. \textbf{146} (1997), 1-52.

\bibitem[GSS]{GSS} J. Graczyk, D. Sands and G. \'Swi\c{a}tek, 
  {\em Decay of geometry for unimodal maps: negative Schwarzian case},
Ann. of Math. \textbf{161} (2005), 613-677.

%\bibitem[GSm1]{GSm1} C. Grotta-Ragazzo, D. Smania, 
%{\em Birkhoff sums as distributions I: Regularity}, 
%\url{https://arxiv.org/abs/2104.04806}

\bibitem[GSm]{GSm} C. Grotta-Ragazzo, D. Smania, 
{\em Birkhoff sums as distributions II: Applications to deformations of dynamical systems}, 
\url{https://arxiv.org/abs/2104.04820}. 

\bibitem[Hub]{Hub} J.H.  Hubbard {\em Teichm\"uller theory and applications to geometry, topology, and dynamics}, Vol. 1. Matrix Editions, Ithaca, NY, 2006. xx+459 pp. 

\bibitem[IT]{IT}  Y. Imayoshi and M. Taniguchi, {\em Introduction ot Teichm\"uller spaces}, Springer Verlag 1992.  

 \bibitem[Koz]{Koz} O.S. Kozlovski, {\em Axiom A maps are dense in the space of unimodal maps in the $C^k$ topology}, Ann. of Math. (2) 157 (2003), no. 1, 1--43.


\bibitem[KvS]{KvS} O. Kozlovski and S. van Strien, \emph{Local connectivity and quasi-conformal rigidity of non-renormalizable polynomials}, Proc. London Math. Soc. (3) \textbf{99} (2009), 275--296.
 
\bibitem[KSvS1]{KSS-rigidity} O. Kozlovski, W. Shen, S. van Strien, {\em Rigidity for real polynomials}, Ann. of Math. \textbf{165} (2007), 749-841.

\bibitem[KSvS2]{KSS-density} O. Kozlovski, W. Shen, S. van Strien, {\em Density of hyperbolicity in dimension one},  Ann. of Math. \textbf{166} (2007), 145-182.

\bibitem[KP]{KP} S. G. Krantz and H. R. Parks, {\em A Primer of Real Analytic Functions}, 2nd ed.,
Birkh\"auser, Boston, 2002. 


\bibitem[LSvS1]{LSvS1}  G. Levin, W. Shen and S van Strien, {\em  Positive transversality via transfer operators and holomorphic motions with applications to monotonicity for interval maps,}  Nonlinearity 33 (2020), no. 8, 3970--4012.

\bibitem[LSvS2]{LSvS2} G. Levin, W. Shen and S van Strien,
{\em Transversality in the setting of hyperbolic and parabolic maps,} 
J. Anal. Math. 141 (2020), no. 1, 247--284. 

\bibitem[LSvS3]{LSvS3} G. Levin, W. Shen and S van Strien, {\em Transversality for critical relations of families of rational maps: an elementary proof.} In {\lq}New trends in one-dimensional dynamics{\rq}, 
201--220, Springer Proc. Math. Stat., 285, Springer, Cham, (2019). 

\bibitem[LS]{LS} S. Li and W. Shen, {\em  An improved real $C^k$ Koebe principle}, Ergod. Th. \& Dynam. Sys. {\bf 30} (2010), 1485--1494. 

\bibitem[Ly1]{Lyubich} M. Lyubich, {\em Feigenbaum-Coullet-Tresser universality and Milnor's hairiness conjecture,} Ann, of Math. 
  \textbf{149} (1999), 319-420. 

\bibitem[Ly2]{Lyubich-quadratic dynamics} M. Lyubich, {\em Dynamics of quadratic polynomials I-II,} Acta Math. \textbf{178} (1997), 185-297.

\bibitem[MMS]{MMS} R. Ma\~n\'e, P. Sad and D. Sullivan,
{\em On the dynamics of rational maps},  Ann. Sci. \'Ecole Norm. Sup. (4) 16 (1983), no. 2, 193--217. 

\bibitem[MMvS]{MMvS}  M. Martens, W. de Melo and S. van Strien. {\em  Julia-Fatou-Sullivan theory for real one-dimensional dynamics},  Acta Math. 168 (1992), no. 3--4, 273--318. 

\bibitem[McM]{McM1} C.T. McMullen. {\em Complex dynamics and
  renormalization}, Annals of Mathematics Studies, 135. Princeton
University Press, Princeton. 1994.

\bibitem[McSu]{McSu} C.T. McMullen and D. Sullivan, {\em Quasiconformal homeomorphisms and dynamics. III. The Teichmüller space of a holomorphic dynamical system},  
Adv. Math. 135 (1998), no. 2, 351--395. 

\bibitem[dMvS]{dMvS} W. de Melo and S. van Strien. {\em One-Dimensional Dynamics}, Springer-Verlag. Berlin (1993).

\bibitem[Mil]{Mil} J. Milnor. {\em Dynamics in one complex variable}, Annals of Math Studies {\bf 160} 2006.

\bibitem[MTr]{MTr} J. Milnor and C. Tresser, 
{\em On entropy and monotonicity for real cubic maps}, 
With an appendix by Adrien Douady and Pierrette Sentenac. 
Comm. Math. Phys. 209 (2000), no. 1, 123--178. 

%\bibitem[MP]{MisPer} M. Misiurewicz and R. P\'erez,
%  {\em Real saddle-node bifurcation from a complex viewpoint,}
%Conform. Geom. Dyn. {\bf 12} (2008), 97-108.

%\bibitem[NY]{NY} H.E. Nusse and J.A. Yorke, {\em Period halving for $x_{n+1}=MF(x_n)$ where $F$ has negative Schwarzian derivative.}

\bibitem[Pom]{Pom} Ch. Pommerenke, {\em Boundary behaviour of conformal maps}, 
Grundlehren {\bf 299}, Springer Verlag. 

\bibitem[R]{R}  L. Rempe, {\em Rigidity of escaping dynamics for transcendental entire functions}, Acta Math. 203 (2009), no. 2, 235--267.

\bibitem[RvS]{Rempe-Strien} L.   Rempe-Gillen and S.  van Strien, {\em Density of hyperbolicity for classes of real transcendental entire functions and circle maps,} Duke Math. J. {\bf 164} (2015), no. 6, 1079--1137. 

\bibitem[She1]{Shen}
W. Shen, \emph{On the measurable dynamics of real rational functions}, Ergod. Theory Dyn. Syst. \textbf{23} (2003), 957--983.


%\bibitem[Shi]{Shi} M. Shishikura, {\em Bifurcation of parabolic fixed points}.  In: The Mandelbrot set, theme and variations, 325–363, London Math. Soc. Lecture Note Ser., 274, Cambridge Univ. Press, Cambridge, 2000. 
\bibitem[Sh]{Sh} M.  Shub, {\em Endomorphisms of compact differentiable manifolds}, 
Amer. J. Math. {\bf 91}, (1916), 175--199.  

\bibitem[Sm]{Sm-Fib} D. Smania,
  {\em Puzzle geometry and rigidity: the Fibonacci cycle is hyperbolic},
Jour. Amer. Math. Soc. {\bf 20}(3) (2007), 629--673.

\bibitem[Sm2]{Sm-shy} D. Smania,
{\em Solenoidal attractors with bounded combinatorics are shy},  Ann. of Math. (2) {\bf 191} (2020), no. 1, 1--79. 


\bibitem[vS1]{strien-robustchaos}
S. van Strien, \emph{Density of hyperbolicity and robust chaos within
  one-parameter families of smooth interval maps}, Proc. Amer. Math. Soc.
  \textbf{138} (2010), no.~12, 4443--4446.
  
 \bibitem[vS2]{strien-density-bis}
S. van Strien, \emph{Density of hyperbolicity holds within families of real analytic one-dimensional maps}, in preparation. 
  
\bibitem[vSV]{vSV} S. van Strien and E. Vargas,
  {\em Real bounds, ergodicity, and negative Schwarzian for multimodal
  maps,}
Jour. Amer. Math. Soc. {\bf 17}(4) (2004), 749--782.

\bibitem[Su1]{Su1} D. Sullivan, {\em Quasiconformal homeomorphisms and dynamics I: a solution of Fatou-Julia problem on wandering domains}, 
  Ann. of Math. {\bf 122} (1985), 401-418.
  
\bibitem[ST]{ST}  D. Sullivan and W. Thurston, {\em Extending holomorphic motions}, Acta Math. {\bf 157} (1986), 243--257.

%\bibitem[Zd]{Zd} A. Zdunik, Fund. Math. {\em Entropy of transformations of the unit interval}, 124 (1984) 235--241.

\end{thebibliography}
\end{document}